\documentclass[11pt,a4paper]{article}
\usepackage{makeidx}
\makeindex
\newtheorem{Theorem}{Theorem}[subsection]

\newtheorem{Proposition}[Theorem]{Proposition}
\newtheorem{prop}[Theorem]{Proposition}
\newtheorem{Corollary}[Theorem]{Corollary}

\newtheorem{Lemma}[Theorem]{Lemma}
\newtheorem{Axiom}{Axiom}
\newtheorem{Example}{Example}[subsection]
\newtheorem{Definition}{Definition}[subsection]
\newtheorem{Remark}{Remark}[subsection]
\newtheorem{remark}[Remark]{Remark}
\newenvironment{eproof}{\par\vspace{1mm}%
\par\bgroup\small{\bf Proof.}}{\qed\egroup\par
\vspace{1mm}%
}
\newenvironment{eproofi}[1]{\par\vspace{1mm}%
\par\bgroup\small{\bf Proof#1.}}{\qed\egroup\par
\vspace{1mm}%
}
\newdimen\mywidth
\mywidth=\textwidth
\advance\mywidth by -\parindent

\newcommand{\seperateline}{\par\vskip 2mm
\vbox{\smallskip\hbox to \hsize{}\hrule\smallskip}%
\vskip 2mm
\par}%
\newcommand{\fullline}{\vbox{\smallskip\hbox to \hsize{}\hrule\smallskip}}%
\newcommand{\halfseperateline}{\par
\hbox to \hsize{\hfil\cleaders\hrule\hskip 7cm\hfil}
\vspace{2mm}\par
}%
\newcommand{\slinewithtitle}[1]{\par
\hbox to \hsize{\hfil\cleaders\hrule\hskip 3cm#1\cleaders\hrule\hskip 3cm\hfil}
\vspace{2mm}\par
}%
\newcommand{\ba}{\both{{\bf a}}}
\newcommand{\bfzero}{\both{{\bf 0}}}
\newcommand{\be}{\both{{\bf e}}}
\newcommand{\bu}{\both{{\bf u}}}

\newcommand{\bx}{\both{{\mathbf x}}}

\newcommand{\both}[1]{{\ifmmode{#1}\else{$#1$}\fi}}

\newcommand{\integer}{\both{\mathbb Z}}

\newcommand{\itwo}{\both{\;(i=1,\;2)\;}}

\newcommand{\myBox}{\rule{4pt}{10pt}}
\newcommand{\nin}{\notin}
\newcommand{\naturalnumber}{\both{\mathbb N}} 
\newcommand{\nat}{\both{\naturalnumber}}

\newcommand{\pair}[1]{\both{{\langle#1\rangle}}}

\newcommand{\pow}{\both{{\bf pow}\;}}

\newcommand{\qed}{\hfill \myBox}

\newcommand{\ra}{\both{\rightarrow}}

\newcommand{\rational}{\both{\mathbb{Q}}}
\newcommand{\real}{\both{\mathbb{R}}}

\newcommand{\setii}[2]{\both{{ \left\{\;\left. {\strut#1} \;\right|\; {\strut#2} \;\right\}}}}
\newcommand{\seti}[1]{\both{{\left\{\; #1 \;\right\}}}}

\newcommand{\seq}[2]{\both{{#2_1,\cdots,#2_{#1}}}}
\newcommand{\seqinf}[1]{\both{{#1_1,#1_2,\cdots}}}

\newif\ifwithkotae
\withkotaetrue

\newcommand{\verylongarrow}[1]{%
\smash{\mathop{\hbox to 1cm{\rightarrowfill}}\limits^{#1}}}

\newcommand{\Frac}[2]{\displaystyle{\frac{\strut #1}{\strut #2}}}

\def\twov(#1,#2){\mtrxall{c}{#1\\#2}}
\def\twovv(#1,#2){(\mtrxall{c}{#1\\#2})}

\def\rrr{\mathrel|\kern-.5pt\joinrel\mathrel\leadsto}

\newsavebox{\circlebox}
\savebox{\circlebox}{\fontencoding{OMS}\selectfont\Large\char13}
\newlength{\circleboxwdht}

\newcommand{\mynext}[1]{#1^{+}}
\newcommand{\mylabel}[1]{\label{#1}}
\newcommand{\mylabelitem}[1]{\label{#1}}
\newcommand{\mylabeleq}[1]{\label{#1}}

\newcommand{\sabun}[1]{\both{\frac{\Delta{}#1}{\Delta x}}}

\newcommand{\hendoukansuu}[2]{%
\Frac{\partial#1}{\partial x_{#2}}%
}
 
\newcommand{\koukaisabuni}[2]{\both{\frac{\Delta^{#2}#1}{\Delta x^{#2}}}} 
\newcommand{\koukaibibun}[2]{\both{\Frac{\Delta^{#2}#1}{\Delta x^{#2}}}}

\newcommand{\hensabun}[2]{\both{\mathcalD_{#2}#1}}

\newcommand{\myif}{\quad \mbox{if  }}
\newcommand{\virtualline}{\both{(-\infty,\infty)}}

\newcommand{\kukanii}[2]{\both{[-#1,#1]_{\frac1{#2}}}}
\newcommand{\kukaniie}[2]{\both{[-#1,#1]_{#2}}}
\newcommand{\kukan}[1]{[0,1]_{\frac1{#1}}}
\newcommand{\kukane}[1]{[0,1]_{{#1}}}

\newcommand{\accessible}{accessible}

\newcommand{\inaccessible}{inaccessible}
\newcommand{\accessibility}{accessiblility}
\makeatletter
\@twosidetrue
\def\ps@lct{%
\let\@mkboth=\markboth
\def\@oddhead{%
\vbox{\hbox to \hsize{\oddpagehead}\medskip\hrule}}%
\def\@evenhead{%
\vbox{\hbox to \hsize{\evenpagehead}\medskip\hrule}}%
\def\@oddfoot{\oddpagefoot}
\def\@evenfoot{\evenpagefoot}
}
\def\@oddfoot{%
\vbox{
\hbox to \textwidth{\oddpagefoot}}
}
\def\@evenfoot{
{\vbox{
\hbox to \textwidth{\evenpagefoot}}}
}
\makeatother

\def\runningtitle{}
\def\oddpagehead{{\tiny Ver. \ver}\hfill \runningtitle\hspace{1cm}\thepage}
\def\evenpagehead{\thepage\hspace{1cm}\hfill}
\def\oddpagefoot{}
\def\evenpagefoot{}
\newcommand{\kosuu}[1]{{}^{\#}(#1)}
\newcommand{\support}[1]{|#1|}
\newcommand{\keywords}[1]{}
\newcommand{\subclass}[1]{}
\newcommand{\mathcalD}{{\mathcal{D}}}

\usepackage[all,cmtip]{xy}

\usepackage[pdftex]{graphicx}
\usepackage{hyperref}
\usepackage{amsmath}
\usepackage{amssymb}
\newtheorem{example}[subsubsection]{Example}
\renewenvironment{Example}{\begin{example}}{\qed\end{example}}
\newtheorem{definition}[subsubsection]{Definition}
\renewenvironment{Definition}{\begin{definition}}{\qed\end{definition}}
\renewenvironment{Remark}{\begin{remark}}{\qed\end{remark}}
\def\oddpagehead{\hfill \thepage}
\def\oddpagefoot{}
\def\evenpagehead{\oddpagehead}
\def\evenpagefoot{\oddpagefoot}
\title{Alternative Mathematics without Actual Infinity
\thanks{Thanks to Ritsumeikan University for the sabbathical leave 
which allowed the author to concentrate on doing research on this theme.}
}
\author{Toru Tsujishita
}
\date{2012.6.12}
\begin{document}
\maketitle
\addcontentsline{toc}{section}{Abstract}
\begin{abstract}
An alternative mathematics based on qualitative plurality of
finiteness is developed to make non-standard mathematics independent
of infinite set theory.  The vague concept ``accessibility'' is used
coherently within finite set theory whose separation axiom is
restricted to definite objective conditions. The weak equivalence
relations are defined as binary relations with sorites phenomena. Continua are collection with weak equivalence
relations called indistinguishability. The points of continua are the
proper classes of mutually indistinguishable elements and have
identities with sorites paradox. Four continua formed by huge binary
words are examined as a new type of continua. Ascoli-Arzela type
theorem is given as an example indicating the feasibility of
treating function spaces.

The real numbers are defined to be points on linear continuum and have
indefiniteness. Exponentiation is introduced by the Euler style and
basic properties are established. Basic calculus is developed and the
differentiability is captured by the behavior on a point. Main tools
of Lebesgue measure theory is obtained in a similar way as Loeb measure.

Differences from the current mathematics are examined,
such as the indefiniteness of natural numbers, qualitative plurality of finiteness, mathematical usage of vague concepts, the continuum as a primary inexhaustible entity
and the hitherto disregarded aspect of ``internal measurement'' in mathematics.

\subclass{26E35  
\and 03H15 
\and 03B52 
\and 03E70 
\and 54J05 
\and 05C99 
}
\end{abstract}
\clearpage
\addcontentsline{toc}{section}{Contents}
\tableofcontents
\newpage
\def\oddpagehead{\S\thesection\hspace{2mm}\runningtitle \hfill \thepage}
\def\oddpagefoot{}
\setcounter{section}{-1}
\section{Introdution}
\def\runningtitle{Introduction}
\subsection{Nonstandard Approach as a Genuine Alternative}
Mathematics has evolved by integrating paradoxes. The Galilei paradox
and the Sorites paradox represent two logical phenomena concerning
infinity. The former, the inevitable paradox of actual infinity that a
part is the same size as the total is incorporated in mathematics as
the very definition of infinite sets. The latter, the inevitable
paradox resulting from two incommensurable view points such as micro
vs macro is in harmony with the infinite phenomena daily experienced by
us and is taken in mathematics implicitly by nonstandard mathematics.

As a result mathematics has currently two methods of treating infinity, 
Cantorian set theory and nonstandard mathematics.
In contrast to the former which handles infinity as a definite concept, 
the latter handles infinity as being
``incarnated'' in finiteness thus removing the inconvenient dichotomies such 
as infinite vs finite and continuous vs discrete. 
Nonstandard mathematics shows new ways 
of making mathematical discourse more intuitive without losing logical
rigor and giving more flexible ways of constructing mathematical objects.
We may say that by discriminating between ``actual finiteness'' 
and ``ideal finiteness'', we obtain a better system of handling infinity
than the ``actual infinity'' offers.

Surely the nonstandard mathematics was born and has been bred in the
realm of Cantorian set theory. Various axiomatic systems for
nonstandard mathematics such as IST (Internal Set Theory) 
\cite{nelson77:_inter_set_theor} of
E.Nelson, RST (Relative Set Theory ) \cite{peraire1992théorie} 
of P\'eraire and EST (Enlargement Set Theory ) \cite{ballard1994foundational} 
of D. Ballard  are conservative
extensions of ZFC, so that the statements of usual mathematics proved
in the new axiomatic systems can be proved without them. It is natural
that many researchers considered the conservativeness the crucial
point since the significance of nonstandard mathematics at first was
able to be claimed only through its relation to current mathematics.
Besides one could believe the consistency of axiomatic systems of
nonstandard mathematics only through reducing it to that of the
standard systems.

However as long as it remains grafted to Cantorian set theory, the
nonstandard mathematics will not unveil its seminal significance as a
genuine alternative to modern mathematics and its potentiality will
not be fully brought to fruition. It seems high time to break the fetters 
and to make nonstandard mathematics independent of Cantorian set
theory. 
In fact, after 50 years after its birth, there seems to be 
widespread conviction that most of modern mathematics can be rebuild
more efficiently by nonstandard mathematics and that new wine must be
put in new bottle, namely the foundation of nonstandard mathematics
itself should be rebuild without recourse to infinite set theory.

In fact, already in 1991, P. Vop\v{e}nka \cite{vopenka1991philosophical}
clearly stated such a view as follows\footnote{For more quotations of similar views, see 
\S~\ref{sec:quotation}.}.
\begin{quote}
As long as this master-vassal relationship lasts, Non-standard Analysis cannot use all its potential, which lies mainly in new formalizations of various situations and not in new proofs of classical theorems. $\dots$ 
It is necessary to approach the study of natural infinity directly and not through its pale reflection as found within Cantor's Set Theory. Such a direct approach is what Alternative Set Theory attempts.
\end{quote}
E. Nelson \cite{nelson-2007-simplicity} points out the importance of thinking of nonstandard analysis as a genuine alternative to modern mathematics.
\begin{quote}
Heretofore nonstandard analysis has been used primarily to simplify
proofs of theorems. But it can also be used to simplify theories. There are
several reasons for doing this. First and foremost is the aesthetic impulse, to
create beauty. Second and very important is our obligation to the larger scientific
community, to make our theories more accessible to those who need to use
them. To simplify theories we need to have the courage to leave results in simple,
external form ------ fully to embrace nonstandard analysis as a new paradigm
for mathematics.
Much can be done with what may be called minimal nonstandard analysis.
\end{quote}

\subsection{Multiple Levels of Finiteness}
The crucial point of the nonstandard mathematics is to 
afford qualitative multitudeness of finiteness. 
Unfortunately one must take currently a long detour to actualize 
the qualitative multitudeness of finiteness in modern mathematics, 
because of its deep belief in the qualitative uniqueness of finiteness 
symbolized as ``the infinite set \nat''. 
However the presence of qualitatively different levels of finiteness 
is an undeniable state of affairs in real life 
and may be assumed as a fundamental principle 
much more secure than the belief in the $\nat$-dogma.

The importance of considering seriously the qualitatively different
kind of natural numbers has been stressed repeatedly by many mathematicians from
the middle of the last century. 
In 1952, E. Borel \cite{borel-1952} considered ``inaccessible numbers'' would be
important. Around 1960, E. Volpins \cite{volpin} claimed the
multitude of natural number sequences and in 1971 R. Parikh \cite{r.71:_exist_and_feasib_in_arith} pointed out various paradoxical phenomena resulting 
from the uniqueness of the natural number concept, e.g., construction of certain formulas which are shown to be provable but the proof is too long to be actually 
carried out. See \S\ref{sec:plurality-of-numbers} for more comments on these aspects.

Now there have been many trials to lay foundation of mathematics based on the
multitude of finiteness. Vop\v{e}nka's Alternative Set Theory is one of
the most elaborated approaches admitting only finite sets some of
which are huge containing actually all the ``concrete numbers''.
Similar systems are elaborated in Hyperfinite Set theory \cite{andreev2006theory} 
of Gordon et al..  

The outstanding variance of nonstandard mathematics from the
conventional mathematics is the acceptance of so called ``vague concept''
in mathematics\footnote{See \S\ref{sec:finitism} for a criticism to the Dammetts' arguments on
the incoherence of vague concepts.
}. The totality of accessible objects is indefinite
since the accessibility depends on the methods of access and even if a
method is fixed it is not clear how far we can access.  Hence one
cannot consider the collection of all the accessible numbers as 
a set and must treat it as a proper class like the totality of sets. 
However this collection is contained in the finite set of numbers 
less than an inaccessible number,
whence the notion of semisets of Vop\v{e}nka will play vital roles
in this new mathematics. See \S\ref{sec:vagueconcept} for more
points on concepts without extension.

In educational studies of mathematics, it has been
pointed out that the concept of ``measuring infinity''\cite{tall1980notion}\footnote{
This is another aspect of the natural infinity in the sense of P. Vop\v{e}nka. See \S~\ref{sec:discussions}.}
such as the hyperfinite numbers in nonstandard mathematics is more
intuitive than that of ``cardinal infinity'' of Cantorian set
theory \cite{tall2001infinity,brown2010step}.
Regrettably the usage of nonstandard mathematics in elementary
levels of university education is not workable at present
because of various artifacts in its usual framework 
resulting from the detour through infinite set theory\footnote{
There seems recently to be new trials \cite{hrbacek2010analysis}
with good results based on relative set theory}. 

However E. Nelson \cite{nelson-nonstandard} clearly showed that
``minimal nonstandard analysis'' captures directly the essence of a
deep mathematical theory in an elementary way without artificial arguments
when freed from the burden of the infinite sets theory. The preface
states clearly his intention as follows.
\begin{quote}
This work is an attempt to lay new foundations for probability theory,
using a tiny bit of nonstandard analysis. The mathematical background
required is little more than that which is taught in high school, and it is
my hope that it will make deep results from the modern theory of stochastic
processes readily available to anyone who can add, multiply, and reason.
What makes this possible is the decision to leave the results in nonstandard
form. Nonstandard analysts have a new way of thinking about
mathematics, and if it is not translated back into conventional terms then
it is seen to be remarkably elementary.

Mathematicians are quite rightly conservative and suspicious of new
ideas. They will ask whether the results developed here are as powerful as
the conventional results, and whether it is worth their while to learn nonstandard
methods. These questions are addressed in an appendix, which
assumes a much greater level of mathematical knowledge than does the
main text. But I want to emphasize that the main text stands on its own.
\end{quote}

Just as it took only a few decades for mathematicians to get
comfortable with the cardinal infinity, it may not take long that
discourse using the measuring infinity become common practice as tools
more fundamental and more versatile than the cardinal infinity. But it
will surely take at least a few decades and most mathematicians might
hesitate to take the risk of get involved in such a long range
uncertain project.  But various trials to develop such a genuine
alternative to modern mathematics are indispensable for healthy
evolution of future mathematics in view of the strong evidence of the
radical superiority of the alternative over the current
mathematics. Besides already mentioned contributions \cite{vopenka},
\cite{nelson-nonstandard}
there are many proposals and trials of alternative mathematics
based on similar intention such as 
\cite{sochor-differential-calculus-AST}, 
\cite{mycielski81:_analy_without_actua_infin}, 
\cite{harthong1983éléments}, 
\cite{beck80:_relat_between_algeb_and_analy,beck79:_simpl_sets_and_found_of_analy}
\cite{Lutz1987reveries},\cite{lutz1992force},
\cite{Laugwitz-omega-calculus},
\cite{diener1992application} to mention a few.
I hope this another trial would play some role, however small it may be,
to strengthen and quicken the movement to free nonstandard analysis from
current mathematics.

\subsection{Points of Conflicts with Modern Mathematics}
The followings are some of the features of our approach
radically different from the usual mathematics.
\begin{description}
\item[Sets are finite.] The usual ``infinite sets'' such as \nat{} and \rational{}  
are considered as proper classes so that the totality is not considered as 
a definite object. 

\item[Sorites Axiom.] A number $x$ is called accessible if there is a
  certain concrete method of obtaining it\footnote{For example there
    is a concrete Peano formula $P(u)$ such that $x$ is the minimal
    number satisfying $P(x)$.}. We postulate the existence of
  inaccessible numbers as the most basic axiom of our framework. 
  The accessible numbers form an nonending number series 
  which is closed under the operation 
  $x\mapsto x+1$ but differes from the total number series. 
  Accordingly, fundamental notions such as
  transitivity, equivalence relation, provability, compatibility, etc.
  become relative to the number series chosen. 

\item[The overspill axiom.] If an objective condition holds for all
  accessible numbers, then it holds also for an inaccessible number.
  Here a condition is called objective if it can be specified without
  the notion of accessibility.

\item[Vague conditions.] The vaguness of the accessibility
  prohibits us to regard the collection of accessible numbers as a set.
  It is a proper class contained in a finite set, called semiset in
  Alternative Set Theory of Vop\v{e}nka\cite{vopenka}.

\item[Continua are not infinite sets.] The real line is considered
  as the ``quotient'' of the proper class \rational{} by the
  indistinguishability relation defined by $r\approx r'$ if and only
  if $k|r-r'|<1$ for every accessible number $k$. Although this
  ``quotient'' is used only as a way of speech, we can represent 
  for example the ``unit interval'' $\setii{r\in \rational}{0\leq r\leq 1}/\approx$ by the
  quotient of the finite set $\setii{\frac{i}{\Omega}}{0\leq i\leq
    \Omega}/\approx$ with an inaccessible number $\Omega$. 
See \S\ref{sec:continuum} for more discussions on continuum.

\item[Functions not as arbitrary mappings.] A function on a proper
  class must be given by an explicit objective specification. However functions on sets are precisely the usual
  arbitary mappings since every map has an explicit specification as a
  finite table. A function on a semiset can be extended to a mapping
  defined on a set including $D$.  For example a sequence defined on
  the accessible numbers is uniquely extended up to a certain
  inaccessible number.
\end{description}

\subsection{Background}
We augment the above position by examining key differences 
between Cantorian infinity and ``Robinsonian infinity''.

\renewcommand{\runningtitle}{Background} 
\label{intro}
\label{sec:discussions}
\subsubsection{Qualitative Plurality of Numbers}
\label{sec:plurality-of-numbers}
\label{sec:quotation}
``The infinite set \nat'' has brought phenomenal evolution of
mathematics by its boundless productivity. However it still remains a 
pure dogma, without any supporting mathematical phenomena. On the contrary,
there have been found many mathematical observations against it such as
Skolem theorem and  G\"{o}del's incompleteness theorem
signifying respectively ontological and epistemological indefiniteness
of the collection of natural numbers. 
As a result various disbelief in ``the infinite set \nat'' has never 
vanished and quite a few mathematicians have stated strong views against it.

Perhaps one of the earliest positive criticism against it is
stated by E. Borel in \cite{borel-1952} where he pointed out 
the potential productivity of taking accessibility into account
as follows.
\begin{quote}
Il me semble que les math\'ematiciens, tout en conservant
le droit d'\'elaborer des th\'eories abstraites d\'eduites d'axiomes
arbitraires non contradictoires, ont int\'er\^et, eux aussi,
\`a distinguer, parmi les \^etres de raison qui sont la substance
de leur science, ceux qui sont v\'eritablement accessibles,
c'est-\`a-dire ont une individualit\'e, une personnalit\'e qui les
distingue sans \'equivoque; on est ainsi conduit \`a d\'efinir
avec pr\'ecision une science de l'accessible et du r\'eel, au del\`a
de laquelle il reste possible de d\'evelopper une science de
l'imaginaire et de l'imagin\'e, ces deux sciences pouvant, dans
certains cas, se pr\^eter un appui mutuel.\footnote{
``It seems mathematicians have interests in distinguishing really
accessible objects, namely, those which have individuality with
personality distinguishing them clearly from others, among the
intellectual objects which constitute the substance of their
discipline, keeping of course the right to elaborate the abstract
theories deduced from arbitrary consistent axioms. Thus one can
precisely define a science of accessibility and reality from
which it is possible to develop a science of imagination and imagined
objects, and in certain cases these two sciences can support each other.''
}.
\end{quote}

Around 1960, E.Volpin \cite{volpin} stated the radical view of the
multitudeness of natural number series which has given various impetus
to explore alternative mathematics freed from the dogma of ``the
infinite set \nat''.  An example is the seminal paper of
R. Parikh \cite{r.71:_exist_and_feasib_in_arith} which showed several
paradoxical consequences of the $\nat$-dogma and suggested the importance 
of taking the notion of ``feasibility'' into account in mathematics.
\begin{quote}
  Does the Bernays' number $67^{257^{729}}$ actually belong to every
  set which contains $0$ and is closed under the successor function?
  The conventional answer is yes but we have seen that there is a very
  large element of fantasy in conventional mathematics which one may
  accept if one finds it pleasant, but which one could equally
  sensibly (perhaps more sensibly) reject.
\end{quote}

Another example is an outline \cite{rashevskii1973dogma} by P.K.Rashevsky 
of radically different type of mathematical theory on numbers as follows.
\begin{quote}
What would correspond more to the spirit of physics would be a mathematical
  theory of the integers in which numbers, when they became very
  large, would acquire, in some sense, a ``blurred'' form and would
  not be strictly defined members of the sequence of natural numbers
  as we consider it. The existing theory is, so to speak,
  over-accurate: adding unity changes the number, but what does the
  addition of one molecule to the gas in a container change for the
  physicist? If we agree to accept these considerations even as a
  remote hint of the possibility of a new type of mathematical theory,
  then first and foremost, in this theory one would have to give up
  the idea that any term of the sequence of natural numbers is
  obtained by the successive addition of unity - an idea which is not,
  of course, formulated literally in the existing theory, but which is
  provoked indirectly by the principle of mathematical induction. It
  is probable that for ``very large'' numbers, the addition of unity
  should not, in general, change them (the objection that by
  successively adding unity it is possible to add on any number is not
  quoted, by force of what has been said above).
\end{quote}

See \cite{isles},\cite{mayberry},\cite{sazonov-feasible} for similar
views.\footnote{S.Yatabe observes in \cite{MR2500986} that sorites
phenomena is unavoidable for models of natural numbers 
in set theories in a non-classical logic.} 

Around the same period, although not directly connected with the above tide, 
A. Robinson\cite{robinson1966} created nonstandard analysis,
which took advantage of a mathematical phenomenon conflicting with
the $\nat$-dogma . As is often quoted, he comments on the last page of his
book\cite{robinson1966}
\begin{quote}
   Returning now to the theory of this book, we observe that it is
   presented, naturally, within the framework of contemporary
   Mathematics, and thus appears to affirm the existence of all sorts
   of infinitary entities. However, from a formalist point of view we
   may look at our theory syntactically and may consider that what we
   have done is to introduce new deductive procedures rather than new
   mathematical entities. Whatever our outlook and in spite of Leibniz'
   position, it appears to us today that the infinitely small and
   infinitely large numbers of a non-standard model of Analysis are
   neither more nor less real than, for example, the standard
   irrational numbers.
\end{quote}

Our main purpose is to give another support to the position that ``the
existence of all sorts of infinitary entities'' is not indispensable for
nonstandard mathematics. We try to show this by the strategy of 
developing core mathematics without infinite set theory 
taking the multitudeness of finiteness as the very basic axiom 
considered as more reliable than that of its uniqueness.

\paragraph{Quotations}
The followings are quotations from authors who take the position that 
nonstandard mathematics is a genuine alternative way of handling 
infinity and infinitesimals.

P. Vop\v{e}nka \cite{vopenka} wrote in 1976 
\begin{quote}
Cantor set theory is
  responsible for this detrimental growth of mathematics; on the other
  hand, it imposed limits for mathematics that cannot be surpassed
  easily. All structures studied by mathematics are a priori completed
  and rigid, and the mathematician's role is merely that of an
  observer describing them.  This is why mathematicians are so
  helpless in grasping essentially inexact things such as
  realizability, the relation of continuous and discrete, and so on.
\end{quote}
In 1991 \cite{vopenka1991philosophical}, he analyzed philosophically
the Cantorian set theory and called its infinity as ``classical'' 
and introduced the concept of ``natural infinity'' to capture the aspect 
of infinity present already in huge finite sets 
emerging from the ``horizon'' which bounds our ``view'', and write
\begin{quote}
  Even classical mathematics then studies natural infinity; however,
  it does so inappropriately. Classical mathematics is restricted by
  the accepted limitations, mainly by those inflicted on the horizon.
  The acceptance of the hypothesis that the sharpening process can lead
  to a complete sharpening does not extend the field of our study but
  rather to the contrary, restricts it. The study of situations where
  the sharpening process itself is essential is thus completely
  blocked.  To put it briefly, the laws that govern classical infinity
  are nothing more than a drastic restriction of the laws that govern
  natural infinity.
\end{quote} 
Incidentally the following remark in \cite{vopenka1991philosophical} 
on the nature of the ``horizon'' seems helpful to understand 
the main idea behind the concept of semisets.
\begin{quote}
The following three characteristics of the horizon are now
important for our theme. 
Firstly, we do not understand the horizon as the boundary of the world, 
but as a boundary of our view. So the world continues even beyond the horizon. 
Secondly, the horizon is not some line drawn and fixed in the world but 
it moves depending on the view in question, specifically on the degree 
of its sharpness. The further we manage to push the horizon, the sharper the view. 
Thirdly, for a phenomenon situated in front of the horizon, the 
closer it is to the horizon, the less definite it is.
\end{quote}

\noindent{}G.Reeb\cite{reeb1981math} wrote in 1981
\begin{quote}
Donc, contrairement \`a une l\'egende, il ne s'agit pas du tout de compl\'eter
\nat{}, par l'adjonction d'objets nouveaux, en un ensemble plus large $\nat{}^{*}$;
mais il s'agit de reconnaître que seulement quelques objets privil\'egi\'es de
\nat{}, en particulier $0, 1,2,3,4$ etc., m\'eritent le label standard\footnote{
``Therefore, contrary to the legend, it is not the question of augmenting 
$\nat$ by adding new objects to a larger set $\nat^{*}$ but it is only the matter
of recognizing that some priviledged elements of $\nat$, in particular $0,1,2,3,4$ etc,
are entitled to be labeled standard.''
}. 
\end{quote}
In 1983, J.Harthong \cite{harthong1983éléments} wrote
\begin{quote}
Je voudrais montrer
  dans cette communication que ....
si on admet que les entiere na\"is ne remplissent pas \nat,
la seule th\'eorie des ensembles finis suffit \^^ a rendre compte de toutes les
proproi\'et\'e du continu, et il est inutile de recourir \^^ a des ensembles non
d\'emombrables\footnote{``I would like to show in this communication that if the naive
integers do not fill \nat{} then only the finite set theory suffices to
treat all the properties of continuum and it is not
to necessary to have recourse to uncountable sets.''}.
\end{quote}
In 1985, A.G.Dragalin \cite{dragalin1985correctness} points out 
the inconsistency of feasibility can be tamed by taking into account
the qualitative difference of length of proofs.
\begin{quote}
We investigate theories with notions ``infinitely large'' and respectively 
``feasible'' numbers of various orders. These notions are 
inconsistent in a certain sense, so our theories turn out to be inconsistent in
an exact sense. Nevertheless, we show that by the short proofs in these theories we get true formulas.
\end{quote}
In 1996, R. Chuaqui and P.Suppes \cite{chuaqui1995free} consider it 
important to ignore the standard part operation.
\begin{quote}
To reflect the features mentioned above that are characteristic of works in
theoretical physics, the foundational approach we develop here has the following
properties:\\
(i) The formulation of the axioms is essentially a free-variable one with no use
of quantifiers.\\
(ii) We use infinitesimals in an elementary way drawn from nonstandard analysis, 
but the account here is axiomatically self-contained and deliberately elementary
in spirit.\\
(iii) Theorems are left only in approximate form; that is, strict equalities and
inequalities are replaced by approximate equalities and inequalities. In particular,
we use neither the notion of standard function nor the standard part function.
Such approximations are not explicit in physics, but can be viewed as implicit in
the way infinitesimals are used.
\end{quote}
In 2005, Y.P\'eraire \cite{peraire2005replacement} pointed out that
nonstandard analysis made it possible to express indefiniteness in 
mathematics.
\begin{quote}
The recent history of nonstandard mathematics is displayed so
as to reveal a modification in the used language as well that in the way the
referentiation of the statements is done. These changes could lead to bring the
mathematical language closer to a language of communication. The profusion
of constructions of sets can be limited thanks to a little richer vocabulary
making it possible “to express the indetermination”, indiscernibility, 
inaccessibility \dots  when it is necessary, and permit also to explore 
more precisely with the mathematical language, using a sort of translation of the ordinary language, some concepts about which the language of conventional mathematics
is almost dumb such as concepts of point, infinity or infinitesimality.
\end{quote}
In 2006, Hrba\v{c}eck et al. \cite{Hrbaceck-ultrasmall-2010} also recognizes
the key point of nonstandard mathematics is to incorporate vague concepts
with ``soritic properties'' into mathematics and write as follows.
\begin{quote}
There are many
  examples of ``soritic properties'' for which mathematical induction
  does not hold (``number of grains in a heap'', ``number that can be
  written down with pencil and paper in decimal notation'',
  ``macroscopic number'', ... ), but mathematicians traditionally take
  no account of them in their theories, with the excuse that such
  properties are vague.  We present here a mathematically rigorous
  theory in which a soritic property is put to constructive use. 
\end{quote}

\subsubsection{Properties without Extension}
\label{sec:vagueconcept}
The above quotations may be said to point out the essence of nonstandard mathematics
consists in the positive usage of indefiniteness in mathematics, 
which means the rejection of the monism of sets in modern mathematics. 
How is it possible to treat conditions without definite extensions?

Surely modern mathematics do not exclude conditions because it is
without extension. For example the condition $x\nin x$ is not
considered as nonsense even though we cannot consider its extension.
In fact, from purely formalistic points of view, a ``vague'' concept
has no difference from the usual ones provided the rule of its usage
is precisely given. In fact in the axiomatic formulation of
nonstandard mathematics such as the internal set theory
\cite{nelson04:_bookr}, the rules of the usage of the word
``standard'' is precisely given among which is the prohibition to
consider its extension.  It might be said that we have
already enough experience about reasoning coherently with conditions
without extensions at least formally.

However in order to ``do mathematics'' actually, purely formalistic
position is not helpful and it is beneficial if even vague concepts
have certain kind of extentions so that they have ``set theoretical''
meaning. It is P. Vop\v{e}nka \cite{vopenka} who found the notion of 
semisets which disclosed essential difference of nature between sets 
and ``external sets'' often used informally.

We can not only coherently and naively develop an alternative mathematics 
admitting properties without extension but also
enjoy its advantage over usual mathematics 
since we can treat infinitary concepts and continuum 
more naturally by keeping their indefiniteness.
In \cite{vopenka,vopenka1991philosophical}, P. Vop\v{e}nka
points out that infinite sets are not necessary to treat infinitary
phenomena\footnote{ ``We shall deal with the phenomenon of infinity in
  accordance with our experience, i.e., as a phenomenon involved in
  the observation of large, incomprehensible sets. We shall be no
  means use any ideas of actually infinite sets. Let us note that by
  eliminating actually infinite sets we do not deprive mathematics of
  the possibility of describing actually infinite sets sufficiently
  well in the case that they would prove to be useful.''  }. He also
points  out the merit of his alternative set theory which allows new
kind of natural concepts which are not available in usual
mathematics\footnote{ ``Our theory makes possible a natural
  mathematical treatment of notions that either have not yet been
  defined mathematically or that have been defined in n unsatisfactory
  way. As an example we have here the chapter dealing with motion.''
}.

\subsubsection{Coherence of Vague Concepts}
\label{sec:finitism}
We do not take the ultrafinitistic standpoint and admit the existence
of inaccessible numbers\footnote{
According to R. Tragesser \cite{tragesser98:_part_i}, ultrafinitistic aim
is not to restrict mathematics to concrete objects but to reconstruct
the idealization of mathematics properly. In this sense, our program
might be called ultrafinitistic.}. Just as infinite sets, huge numbers are ideal objects, but,
in contrast to Cantorian infinity, huge finiteness is philosophically
less problematic and intuitively more in harmony with naive concepts
of infinite quantities\footnote{ This view is supported by educational 
studies of university mathematics. For example J.Monaghan \cite{monaghan2001young}
says as follows.
\begin{quote}
  Cantor's transfinite universe became the infinite paradigm during
  the 20th Century. This affected educational studies, which tended to
  view children's responses against Cantorian ideas. Robinson's
  non-standard universe (Robinson, 1966) is equally authoritative
  (though not as well known) and it is a different paradigm. It offers
  researchers a release from a single paradigm and allows them to
  interpret children's ideas with reference to children's ideas
  instead of with reference to Cantorian ideas.
\end{quote}
}.

However the concepts of accessibility and hence that of hugeness interpreted
as inaccessibility are vague. Since
Frege, vague concepts have been considered as useless in mathematics
because of various incoherences associated to them.  For example,
M. Dummett \cite{dummett1975} argues for this Frege's position that
use of vague expressions is fundamentally incoherent and concludes as
follows.
\begin{quote}
Let us review the conclusions we have established so far.
\def\labelenumi{(\theenumi)}
\begin{enumerate}
\item\label{item:wang1} Where non-distinguishable difference is
  non-transitive, observational predicates are necessarily vague.

\item\label{item:wang2} Moreover, in this case, the use of such predicates is
  intrinsically inconsistent.

\item\label{item:wang3} Wang's paradox merely reflects this
  inconsistency.  What is in error is not the principles of reasoning
  involved, nor, as on our earlier diagnosis, the induction step. The
  induction step is correct, according to the rules of use governing
  vague predicates such as 'small': but these rules are themselves
  inconsistent, and hence the paradox. Our earlier model for the logic
  of vague expressions thus becomes useless: there can be no coherent
  such logic.

\item\label{item:wang4} The weakly infinite totalities which must underlie any strict
  finitist reconstruction of mathematics must be taken as seriously as
  the vague predicates of which they are defined to be the
  extensions. If conclusion (\ref{item:wang2}), that vague predicates of this kind
  are fundamentally incoherent, is rejected, then the conception of a
  weakly infinite but weakly finite totality must be accepted as
  legitimate. However, on the strength of conclusion (\ref{item:wang2}), weakly
  infinite totalities may likewise be rejected as spurious: this of
  course entails the repudiation of strict finitism as a viable
  philosophy of mathematics.
\end{enumerate}
\end{quote}

He identifies the condition of transitivity
$$a\approx b\approx c \mbox{ implies }a\approx c$$
with the multiple transitivity
\begin{equation}
\label{eq:transitivity}
a_1\approx a_2\approx a_3\approx \cdots\approx a_n
\mbox{ implies } a_{1}\approx a_n.
\end{equation}
This identification is based on
the tacit assumption that the notion of natural number is uniquely
determined, which is precisely the ultrafinitistic position
doubts. When there are two kinds of natural numbers, for example
feasible and unfeasible ones, it is possible to define the weak
transitive relations for which the multiple transitivity \eqref{eq:transitivity}
holds only for feasible $n$. Hence the conclusions (\ref{item:wang1}) and
(\ref{item:wang2}) are untenable if the assumption of the qualitative
uniqueness of finiteness is abandoned, which opens the possibility to
use weakly transitive relations consistently. Namely, weak
transitivity of non-distinguishable difference turns out to be one of
the corner stone of the new approach to continuum developed here.

As for the conclusion (\ref{item:wang3}), the key arguments against
the skepticism about induction is as follows.  Assume the
ultrafinitistic position that a proof is legitimate only when the
totality of the inferences is survayable.  A number $n$ is called
\textit{apodictic} if a proof, without induction principle, of length
less than or equal to $n$ is legitimate as a proof from
ultrafinitistic standpoint.  Then the condition of being apodictic is
inductive in the sense that $0$ is apodictic and if $n$ is apodictic
then $n+1$ is apodictic.  Moreover a number less that an apodictic
number is also apodictic. If a condition $F$ is inductive then $F(n)$
is true whenever $n$ is apodictic since there is the obvious proof of
$F(n)$ consisting of $n$ lines of modus ponen. Now choose an apodictic
number $k$ and define the condition $S(n)$ to be $n+k$ is
apodictic. Then $S$ is obviously inductive. Suppose there are an
apodictic number $n$ such that $n+k$ is not apodictic. Then $S(n)$ is
false by definition but since $S$ is inductive $S(n)$ is true, a
contradiction. Hence he concludes that the arguments against the
induction principle is not tenable and also implicitly that the notion
of apodictic is incoherent and hence the ultrafinitistic standpoint is
incoherent.

However the contradiction comes from the assumption that there are two
apodictic numbers $k,n$ such that $n+k$ is not apodictic. However this
is based on the tacit assumption that there are no nontrivial
inductive properties of numbers closed under the addition which
ultrafinitistic position doubts. Since not only the induction remains
problematic but also there is coherent usage of ``non-transitive
non-distinguishable difference'' the conclusion (\ref{item:wang3}) is
untenable.

Since conclusion (\ref{item:wang2}) is misleading, so is the
conclusion (\ref{item:wang4}). See \cite{magidor2007strict} for 
similar criticism against Dummett's arguments.

Thus Dummett's arguments against not only to ultrafinitism but also
to any alternative mathematics which use vague concepts is essentially
grounded on the basic assumption of modern mathematics that there is
unique concept of natural numbers, which is exactly the alternative 
approach in this paper negates.

In fact, the secret of
effectiveness of nonstandard analysis might be pin downed to the vague
concept ``standard'' which forbids formation of the set of standard
elements.

\index{$\nat$, class of natural numbers}

\subsubsection{Continuum}
\label{sec:continuum}
The infinite sets are considered indispensable to modern mathematics
since the continua are infinite sets. For example the interval, the
simplest continuum, is identified with ``a set of real numbers between
$0$ and $1$'' which have more elements than ``the set of natural
numbers''. However historically this atomic view regarding continuum 
as a mere aggregation of its points has been
criticized repeatedly from various points of view since ancient
times to today. 

H. Weyl gives in 1921 an overview\footnote{
``An atomistic view, taking the continuum to consist of individual
points, and a view that takes it to be impossible to understand the
continuum flux in this manner, have been opposing each other from time
immemorial. The atomistic one has a system of existing elements that
can be conceptually grasped, but it is incapable of explaining motion
and action. In it, all change must degenerate into appearance. The
second conception was not capable, in antiquity, and up to the time of
Galilei, to lift itself from the sphere of vague intuition to the one
of abstract concepts that would be suitable for a rational analysis of
reality. The solution that was finally achieved is the one whose
mathematical systematic form is given in the differential and integral
calculus.  Modern criticism of analysis is destroying this solution
from within, however, without being particularly conscious of the old
philosophical problems, and it lead to chaos and nonsense.  The two
rescue attempts discussed here revive the old antithesis in a sharper
and more clarified form. The previously described theory is radically
atomistic([I am saying this] in full awareness of the fact that, as it
is, this theory does not fully capture the intuitive continuum, the
idea being that the concepts are capable of grasping only rigid
existence.) Brouwer's theory, on the other hand, undertakes to do
justice to Becoming in a valid and tenable manner. \cite{bo:Mancosu98}''
} of the two opposing approach to continuum, culminating respectively
to Cantorian set theoretic approach and Brouwer's intuitionistic
approach.  He did not satisfied with the atomic approach to continuum
of his book \cite{Weyl-continuum} published in 1917 and recognized the
need to reconstruct it radically according to his philosophy, but he
regrets in the ``preface to the 1932 Reprint'' \footnote{
``It seems not to be out of the question that the limitation prescribed
in the present treaties-- i.e., unrestricted application of the concepts
"existence" and "universality" to the natural numbers, but not to 
sequences of natural numbers-- can once again be of fundamental significance.
It would not be possible, without radical rebuilding, to bring the content
of this monograph into harmony with my current beliefs -- and such a 
project would keep me from satisfying other demands on my time.''
} that he has no time to undertake it \cite{Weyl-continuum}.  

Now that topology has become one of the major disciplines of
mathematics, there seems to be quite a few mathematicians 
who, independently of the antagonism between classical logic vs intuitionistic one,  
consider continua as primitive objects. For example R. Thom amplifies the claim that contiuum
ontologically precedes discrete objects in \cite{MR1413523}.\footnote{It might be said that such viewpoints is reflected for example in the 
computational approach to topology such as \cite{DBLP:journals/corr/abs-1005-5685}. 
}
\begin{quote}
  Ici, je voudrais m'attaquer \`a un mythe profond\'ement ancr\'e dans la
  math\'ematique contemporaine, \`a savoir que le continu s'engendre
  (voire se d\'efinit) \`a partir de la g\'en\'erativit\'e de l'arithm\'etique,
  celle de la suite des entiers naturels. Je fais bien entendu
  allusion \`a la construction de Dedekind o\`u $\real$ se d\'efinit par
  compl\'etion des coupures d\'efinies sur les rationnels. J'estime, au
  contraire, que le continu arch\'etypique est un espace ayant la
  propri\'et\'e d'une homog\'en\'eit\'e qualitative parfaite.
\footnote{``Here I would like to attack a myth deeply anchored in modern mathematics
    which says that continuum is obtained from the generative
    feature of the arithmetics and the series of natural numbers.
    Of course I am referring to Dedekind construction which defines
    $\real$ by completion using the cuts on rationals.
    I consider on the contrary that archetype of continuum is a
    space with qualitatively complete homogeneity.''} 
\end{quote}

Our intention is to develop mathematics which takes as primitive both
the discrete and the continuous based on the plurality of finiteness.
This might be said to conform with the viewpoint of Brouwer who
stated as follows in his dissertation 1907 according to \cite{MR1075018}.
\begin{quote}
Since in the Primordial Intuition the continuous and the discrete
appear as inseparable complements, each with equal rights and
equally clear, it is impossible to avoid one as a primitive entity
and construct it from the other, posited as the independent
primitive.
\end{quote}

However, H. Weyl was dissapointed with the Brouwer's approach 
in which it is awkward to carry out usual mathematics as is remarked\footnote{
``Mathematics with Brouwer gains its highest intuitive clarity. He
succeeds in developing the beginnings of analysis in a natural
manner, all the time preserving the contact with intuition
much more closely than had been done before. It cannot be denied,
however, that in advancing to higher and more general theories
the inapplicability of the simple laws of classical logic
eventually results in an almost unbearable awkwardness. And
the mathematician watches with pain the larger
part of his towering edifice which he believed to
be built of concrete blocks dissolve into mist before his eyes.''
}in his book \cite{weyl-philosophy} published in 1949.

The success of nonstandard mathematics suggests high feasibility of
this approach. 

We regard intervals as primitive objects and the basic operation is
to divide them to subintervals which are similar in nature to the
total interval. This fractal nature is the essential feature of
continua. Dividing the subintervals again and again, we get many small
intervals with many points which bounds them. Although we can divide
only concrete number of times, the division process can be continued
to huge number of times in principle.
Thus we can imagine a set of huge finite number of infinitesimal
intervals each of which we cannot discriminate from the neighboring
ones. Moreover since the intervals obtained are so small that each
interval determines uniquely a position up to indistinguishability on the
interval. The subintervals are infinitesimal but are not identified
with the positions they determine. They are themselves continua which
have the same property as the initial interval.

So we arrive at an imaginary picture that intervals are composed of
huge number of infinitesimal intervals, which themselves can be
divided indefinitely into smaller intervals.
The result of a huge number of division can be described by the huge
finite set 
of the rationals coding the positions of the boundaries of
the resulting infinitesimal intervals. This huge finite set, called a representation
of the interval, has an indistinguishability binary relation satisfying 
the usual axiom of equivalence relation but the huge chain of
indistinguishable elements connects distinguishable elements and 
there arises the sorites paradox. 
So we define a \textit{rigid mesh continuum} 
as a huge finite set equipped with such a paradoxical equivalence relation,
called sorites relation, which exists by virtue of
Axiom~\ref{axiom:sorites}.

Thus continuum as a primitive entity should be represented as an
``equivalence class'' of rigid mesh continua, but we use more handy
formulation, for example, of the linear continuum $\real$ 
as the proper class $\rational$ of rational numbers equipped with the indistinguishability
relation.\footnote{A rational $r$ is infinitesimal and written $r\approx 0$ if
$|r|<\frac1k$ for every accessible number $k$ and two rational numbers are considered
indistinguishabile if their difference is infinitesimal.}
Many continua such as the real line and various types of
intervals are given as subcontinua defined by possibly vague
conditions.\footnote{ For the real line, the condition is the
  finiteness, namely the absolute value is less than an \accessible{}
  number. For the open intervals such as $(0,1)$, the condition is
  $0\prec x\prec 1$ where $x\prec y$ means that 
  that $y$ is visibly greater than $x$.}

Besides the Euclidean continua, a large class of continua is provided
by metric spaces with rational distance functions by defining the
indistinguishability $x\approx y$ as $d(x,y)\approx 0$. Symmetric graphs with
infinitesimal positive distances given to edges form a rich subclass
of metric continua.  This construction given for the first time by
L. van den Dries and A. J. Wilkie \cite{MR751150} plays vital roles in
their proof of the Gromov's theorem on groups of polynomial growth
using nonstandard method.  We note that Urysohn space
\cite{Vershik-Urysohn-1998-MR1691182} can be regarded as a ``universal
continuum'' which includes all metric contiua as subcontinua.

\subsection{Outline of Contents}
In Section \ref{sec:axioms}, we explain fundamentals in a naive way to
emphasize the approach is more easily assimilated than that of
infinite set theory. Then we treat directly ``continua'' represented
as a ``quotient'' of finite sets by weak equivalent relation of
indistinguishability in Section \ref{sec:basic-continuum}.  Usual
topological concepts are reformulated by the indistinguishability
relation in Section \ref{sec:topology-continua}.  

The following two sections discuss concrete examples of continua.
Section \ref{sec:binary-words} treats the continua arising from the
finite sets of the $01$-words of an inaccessible fixed length endowed
with various distances, which demonstrates the drastic increase of
freedom of construction in the new approach. Section
\ref{sec:continuous-maps} investigates the continuum of morphisms and
show the Ascoli-Arzela Theorem, with the purpose to demonstrate how
our framework can treat function spaces.

As a special case of continua, we treat ``real numbers'' as rational numbers under the
 weak equivalence relation of indistinguishability in Section \ref{sec:real-numbers}.
 Section \ref{sec:real-valued-functions} treats real valued functions 
and proves the mean value theorem and the maximum principle.
 The exponential functions is treated just like in the Euler's way. 

The calculus of one
 variable and multiple variables are treated respectively in
 Sections \ref{sec:calculus-single} and \ref{sec:calculus-multiple}.  A new
 feature is that the differentiability of a function controls its
 behavior only on large infinitesimals and the behavior on tiny
 infinitesimal neighborhood can be taken rather arbitrarily, which seems to
 open new freedom to represent functions. The integration is treated
 in a way similar to Loeb measure in Section \ref{sec:measure}. 

\newpage
\section{Fundamentals} 
\renewcommand{\runningtitle}{Fundamentals} 
\label{sec:axioms}
\subsection{Numbers}
We assume usual elementary arithmetic taught up to high school.  For
examples, we have natural numbers $0,1,2,3,\cdots$ sometimes simply 
called numbers, and the integers $0,\pm1,\pm2,\pm3,\cdots$ with the addition and multiplication
satisfying the axiom of rings. We have also the rational numbers
$\pm\frac{p}{q}$ with natural numbers $p,q\neq 0$ with the addition
and multiplication satisfying the axiom of field.

\subsubsection{Accessibility}
\label{subsub:accessibility}
A number is called \textit{\accessible{}} if it can be actually accessed
somehow.  For example, numbers which can be written by some notation
is \accessible.  Since such a naive concept of \accessibility{} has
inevitable vagueness, we use it as an undefined terminology obeying
rigorously the axioms which reflects naive meaning of
\accessibility\footnote{We remark that the concept of ``\accessible{}
  numbers'' is semantically vague but just as vague as that of
  ``numbers'' and less vague than that of ``infinite sets''.  }.

\begin{Axiom}[Accessible{} Numbers]
\label{axiom:axioms-concr-numb}
 \begin{enumerate}
 \item The numbers $0$ and $1$ are \accessible{}{}.
 \item\label{item:4-821} The sum and product of two \accessible{}{} numbers are \accessible{}{}.
 \item\label{item:3-821} Every number less than an \accessible{}{} number is \accessible{}{}\footnote{
Hence numbers less than a big numbers such as
$10^{10^{10}}$ are considered to be \accessible{}{}  although most of them
cannot be written explicitly.
}.
 \end{enumerate}
\end{Axiom}
\index{Axiom, \accessible{}{} numbers}

The following reflects the intuition that there are \inaccessible{} numbers. 
\index{Axiom, sorites}
\begin{Axiom}[Sorites Axiom]
\label{axiom:sorites}
There are numbers which are not \accessible{}{}.
\end{Axiom}

By Axiom \ref{axiom:axioms-concr-numb}, the number $0$ is \accessible{}{} and if
$n$ is \accessible{}{} then $n+1$ is \accessible{}{}, whence every numbers
are \accessible{} if the unrestricted induction principle is applied, 
contradicting to Axiom
\ref{axiom:sorites}. This is one version of the sorites paradox. We
weaken in \S\ref{subsec:induction} the induction principle in order to use the concept
\accessibility{} coherently. 

If a natural number $n$ is \accessible{}{}, we say $n$ is \textit{finite} and
write $n<\infty$. If a natural number $n$ is not \accessible{}{}, 
we say  $n$ is \textit{huge} and write $n\gg 1$. 
An integer is called \accessible{}{} if
its absolute value is \accessible{}.
\index{huge number}%
\index{rational number, huge}%
\index{rational number, finite}%
\index{rational number, accessible}%
We call a rational number $r$ is \textit{bounded from above} and write $r<\infty$
if there is an \accessible{}{} number $n$ with $r<n$ and $r$ is \textit{bounded from below}
and write $-\infty<r$ if $-r<\infty$.  A rational number is called \textit{finite}
if its absolute value is bounded from above.

We say a rational
number is \textit{\accessible} if it is written as $\pm\frac{p}{q}$ with 
accessible $p,q$. An \accessible{}{} rational number is finite but 
the  converse is not true. For example the rational number $\frac1{N}$ with
$N\gg 1$ is finite but is not \accessible{}.

\subsubsection{Rational Numbers}
\label{subsub:infinitesimal-rational}
We call a rational number $r$ \textit{infinitesimal} and write
$r\approx0$ if $|r|<\frac1k$ for all \accessible{}{} number $k$.
We say two rational numbers $r,s$ are \textit{indistinguishable} and write $r\approx s$
if $r-s$ is infinitesimal. 

We remark that the assertion ``for all \accessible{}{} number $k$ the condition $P(k)$ is true'' means that there is a proof of the assertion $P(k)$ with parameter $k$ 
which do not use peculiarity of $k$. See \S~\ref{sss:subclasses} for more elucidation
about this.
\index{infinitesimal}
\index{$r\approx s$}

Axiom \ref{axiom:sorites} implies
\begin{Proposition}
There are nonzero infinitesimal rational numbers.
\end{Proposition}
\begin{eproof}
Let $r=\frac1N$ with $N\gg1$. Let $k$ be an \accessible{}{} number. Then $k<N$ 
whence $kr<1$.  Hence $r$ is an infinitesimal but nonzero rational number.
\end{eproof}

For rational numbers $r,s$, 
we write $r\prec s$ and say that $r$ is \textit{visibly smaller than} $s$
if there is an \accessible{} number $k$ satisfying $r+\frac1k<s$.
\index{$r\prec s$}
\index{$r\preceq s$}
\index{visibly}
We write $r\preceq s$ if $r\prec s$ or $r\approx s$.

Note that $r\leq s$ implies $r\preceq s$ but the
other implication is generally false and 
$r\prec s$ implies $r<s$ but the other implication is generally false.
In fact if $\varepsilon$ is a positive infinitesimal, we have
$r\preceq r-\varepsilon$ and $r\not\prec r+\varepsilon$. 

Obviously we have
\begin{prop}
\label{prop:daishou}
\begin{enumerate}
\def\labelenumi{(\theenumi)}
  \item $r\prec s$ satisfies the transitivity.
  \item The conditions $r\prec s$, $r\approx s$,$s\prec r$ are mutually exclusive and
    just one of them is valid.
  \item The relation $\prec$ is $\approx$-congruent, 
namely, if $r\prec s$,$r\approx r'$ and $s\approx s'$ then $r'\prec s'$.
  \item The relation $\preceq$ is $\approx$-congruent, 
namely, if $r\preceq s$,$r\approx r'$ and $s\approx s'$ then $r'\preceq s'$.
  \item If $0\prec r,s$ then $0\prec r+s,rs$.
  \item If $0\prec r,s_1\prec s_2$ then $rs_1\prec rs_2$.
  \end{enumerate}
\end{prop}

\subsection{Sets and Classes}
\subsubsection{Basic Concepts}
A collection of objects is called a \textit{class} if its elements
have distinctiveness, namely, given two objects $x,y$ qualified
as its elements, it is possible to determine $x=y$ or
$x\neq y$.  If an object $x$ belongs to a class $X$, we  write $x\in X$.

A \textit{set} is a class $X$ with enumeration
$X=\seti{x_1,x_2,\cdots,x_n}$ for some number $n$.
An enumeration of a finite set
without repetition is called a \textit{tight enumeration}. 
\index{enumeration, tight}
The number of elements of a set $A$ is denoted by $\kosuu{A}$.  
\index{$\kosuu{A}$, number of elements}

We take usual naive set theory for granted 
with the exception of those concepts and propositions 
referring to infinite sets.

A class which is not a set is called a \textit{proper class}.
For example, the collections \nat, \integer{} and \rational{}
respectively of natural numbers, integers and rational numbers are
proper classes.

Let $X$ and $Y$ be classes. We say they are equal and write $X=Y$ if
and only if we can prove that every object belongs to $X$ if and only
if it belongs to $Y$.  We say $X$ is different from $Y$ and write 
$X\neq Y$ if and only if we can find an object $x$ either satisfying 
$x\in X$ and $x\nin Y$ or satisfying $x\nin X$ and $x\in Y$. Hence
it is not logically evident that either $X=Y$ or $X\neq Y$ holds.
Hence classes have no distinctiveness so that the collection of 
classes do not form a class.

\subsubsection{Subclasses and Subsets}
\label{sss:subclasses}
A class $Y$ is a \textit{subclass} of a class $X$ 
written as $Y\subset X$ 
if every element of $Y$ is also an element of $X$.

A subset of a class $X$ is a set with elements in $X$. 

\paragraph{Bounded Conditions}
We say a quantification is \textit{bounded} if it is either $\forall
x\in a$ or $\exists x\in a$ with $a$ being a set. A condition is
called \textit{definite}\footnote{Usually called
  $\Delta_0$-conditions.} if it has only bounded quantification. Since
sets can be exhausted, a definite condition $P$ has semantically
definite truth value and either $P$ is true or $P$ is false.
\index{quantification, bounded}
\index{condition, definite}
A condition on the class $\nat$ is bounded precisely when the
quantifications are of the form $\forall x\leq n$ or $\exists x\leq n$.

If $X$ is a proper class, the truth value of an unbounded condition such as
$\forall x\in X.P(x)$ or $\exists x\in X.P(x)$ cannot be determined
semantically, namely, by evaluating the truth value of $P(x)$ for each
$x\in X$ since a proper class cannot be exhausted by any procedures.
So we adopt the proof-theoretic interpretation that ``$\forall x\in X.P(x)$
is true'' means that the assertion $P(a)$ with the parameter $a$ has
a proof which is independent of the parameter $a$, and ``it is false''
means that the assumption that every $x\in X$ satisfies $P(x)$ implies
a contradiction.  For example if we have found an $a$ for which $P(a)$
is false, it is false.  Similarly ``$\exists x\in X.P(x)$ is true''
means that we have constructed an object $a$ satisfying $P(a)$ and
``it is false'' means that the existence of an object $a$ such that
$P(a)$ implies a contradiction.

Note that the condition $Y\subset X$ is not definite if $Y$ is a
proper class. Moreover for two proper subclasses $Y_i\subset X$ \itwo,
the condition $Y_1=Y_2$ is not definite. Hence the collection of
proper subclasses of $X$ is not a class. It will turn out that the
collection of subsets of $X$ is a class when $X$ is $\sigma$-finite 
in the sense defined in \S~\ref{sec:sigma-finite}.

\paragraph{Power Set}
The collection of subsets of a set forms a set as follows. 
Let $A$ be a set with a tight enumeration $\seti{\seq{n}{a}}$. 
An integer $k\in [1..2^{n}]$ defines a set $S_k\subset A$ 
by 
$$
S_k:=\setii{a_i}{\mbox{the binary expansion of $k-1$ has $1$ on the $i-1$-th position}}.
$$
Conversely, for each $B\subset A$, define
$k=\sum_{a_i\in A}2^{i-1}+1$. Then $S_k=B$. 
 
Hence the subsets of a set $A$ defines the power set $\pow(A)$ with 
the explicit enumeration $\setii{S_k}{k\in [1..2^{n}]}$.
We show in \S~\ref{ss:class-constructions}, the subsets of a $\sigma$-finite
class form a $\sigma$-finite class.

\subsubsection{Objective Conditions and Semisets}
A condition is called \textit{objective} if it is specified independently of
the concept of \accessibility.  An \textit{objective subclass} is a subclass
defined by an objective definite condition.

\begin{Remark}
\label{remark:definiteness-of-equality-of-subclasses}
If proper subclasses $A$ and $B$ are not objective, then 
the equality condition $A=B$ is not definite and the
collection of subclasses of a class is not a class generally. 
However the collection of objective subclasses form a class
since the equality condition of objective subclasses is definite.
\end{Remark}

A subclass of a set is called a \textit{semiset}.  A semiset which is
not a set is called a \textit{proper semiset}.  We write $A\sqsubset
x$ if $A$ is a semiset included in a set $x$.
A set including a proper semiset is called \textit{an environment set of} it.
Note that the intersection of two environment sets is also an environment set.
\index{condition, objective}
\index{semiset, environment set}
\index{semiset, proper}
\index{subclass, objective}

\begin{Axiom}[Objective seperation]
\label{axiom:objective-seperation}
An objective semiset is a set.
\end{Axiom}

The class $\nat_{acc}$ of \accessible{}{} numbers is a proper class
and hence is a proper semiset. Generally a proper semiset present
itself only when the defining condition depends on \accessibility{}
explicitly or implicitly. Thus proper semisets play vital roles in
the mathematical treatment of vague concepts such as \accessibility{}.

The proper semisets plays in our theory the similar role as is played
by the infinite sets in usual mathematics. The following is the key
tool in the arguments of the proper semisets.
\begin{Theorem}[General Overspill Principle]
\label{theorem:overspill}
Let $A$ be a proper semiset of a set $X$. Suppose every element 
of $A$ satisfies a definite objective condition $P$ on $X$. 
Then there is an $x\in X\setminus A$ satisfying $P$.
\end{Theorem}
\begin{eproof}
Axiom~\ref{axiom:objective-seperation} implies that 
$B=\setii{x\in X}{P(x)}$ is a subset which includes $A$. 
Since $A$ is not a set, the class $B\setminus A$ must be nonempty.
\end{eproof}

\subsubsection{Class Constructions}
\label{ss:class-constructions}
If $A,B$ are subclasses of $X$ defined respectively 
by definite conditions $P_A,P_{B}$, then usual Boolean operations 
$$A\bigcap B,A\bigcup B,A\setminus B$$ 
are defined respectively by the definite conditions ``$P_A$ and $P_B$'',
``$P_A$ or $P_B$''  and ``$P_A$ but not $P_B$'' and  
obeys usual algebraic laws of Boolean operations.

If $A_i$ ($i\in [1..n]$) are subclasses of a class $X$, then 
subclasses 
$\bigcup_{1\leq i\leq n}A_i$ and $\bigcap_{1\leq i\leq n}A_i$ of $X$ 
are defined respectively by the definite conditions $\forall i\leq n. x\in A_i$ and 
$\exists i\leq n. x\in A_i$.

If $X_i$ \itwo are classes, the product class $X_1\times X_2$ is 
defined as the collection of the ordered pair $\pair{x_1,x_2}$ of
$x_i\in X_i$ \itwo. The coproduct class $X_1\coprod X_2$ is
defined as the collection of $\pair{i,x_i}$ with $x_i\in X_i$ \itwo 
with the canonical inclusions $\imath_i:X_i\rightarrow X_1\coprod X_2$ \itwo
defined by $\imath_i(x_i)=\pair{i,x_i}$ \itwo.

If there is a rule to define a class $A_n$ for each $n\in \nat$ such
that $A_n\subset A_{n+1}$ for all $n$, we say $\setii{A_n}{n\in
  \nat}$ is an \textit{increasing family of classes}. Then the union class
$\bigcup_{n\in N}A_n$ is defined as the collection of the elements of
some $A_n$. There is a function $rank:\bigcup_{n\in N}A_n\rightarrow \nat$
defined by $rank(x)=\min\setii{k\leq m}{x\in A_k}$ for $x\in A_m$,
which satisfies $x\in A_{rank(x)}$. 

\subsubsection{Functions}
\label{sec:functions}
Let $X,Y$ be classes and suppose $f$ is a correspondence which assigns
$x\in X$ to $f(x)\in Y$ by an objective definite rule $R_f$. Here a rule
is called objective if it is specified without recourse to the concept of \accessibility{} and is called definite if the specification does not involve unbounded quantification. We say then that
$f$ is a \textit{function} from $X$ to $Y$ and write $f:X\rightarrow Y$.

Two functions $f,g:A\rightarrow B$ are called equal and written $f=g$
if 
$$\forall x\in A. f(x)=g(x)$$
is proved.
We say $f\neq g$ if we have found an $a\in A$ such that $f(a)\neq g(a)$.  
Since it is not logically obvious that either $f=g$ or $f\neq g$ is true,
functions defined on proper classes have generally 
no distinctiveness and do not form a class.
However Proposition~\ref{prop:continuum-functions} below shows that the collection of 
functions defined on a semiset forms a class.

\begin{Axiom}[Extenstion Axiom]
\label{axiom:extension}
If $f:A\rightarrow Y$ is a function from a proper semiset $A$
to a set $Y$, then there is an environment set $b$ and 
a function $g:b\rightarrow Y$ which coincides with $f$ on A.
\end{Axiom}
A rationale of this axiom is as follows. 
Suppose $x$ is an environment set of $A$. 
The condition on the elements of the set $x$ 
that the defining condition of $f$ has meaning 
and determines an element of $Y$ 
is objective and includes $A$, 
hence defines a set $b\subset x$ 
such that $A\sqsubset b\subset x$.

We call a class $Y$ \textit{set-like} if every function from a semiset to
$Y$ can be extended to an environment set of $Y$. 
\index{class@class, set-like}

Two extensions coincides on an appropriate environment set.
\begin{Proposition}
Let $f:A\rightarrow Y$ be a function from a proper semiset $A$ to a set-like class $Y$
and $f_i:a_i\rightarrow Y$ \itwo be its extensions. Then 
$f_1=f_2$ on an environment set included in $a_1\bigcap a_2$.
\end{Proposition}
\begin{eproof}
The objective condition $f_1(x)=f_2(x)$ on the elements $x\in a_1\bigcap a_2$
is satisfied on $A$, whence defines an environment set of $A$.
\end{eproof}
 
Moreover we can choose extensions of family of functions so that
their domains coincides.
\begin{Proposition}
\label{prop:common-domain-of-extension}
Let $A$ be a proper semiset and   
$$f_i:A\rightarrow Y\quad  i\in [1..n]$$
a family of functions to a set-like class $Y$.
Then there is a set $b$ and extensions $\tilde{f}_i$ ($i\in [1..n]$)
with domains $b$.
\end{Proposition}
\begin{eproof}
Just take any extensions of $f_i$'s and then restrict them to the intersection of
their domains.
\end{eproof}

\begin{Proposition}
\label{prop:continuum-functions}
If $X_1$ is a semiset and $X_2$ is a set-like class. 
Then the collection of functions from $X_1$ to $X_2$ forms a
class $Fun(X_1,X_2)$. 
\end{Proposition}
\begin{eproof}
If $X_1$ is a set, the equality of two functions is obviously definite.

Suppose $X_1$ is a proper semiset.
Let $f_i:X_1\rightarrow X_2$ \itwo be functions.
By Proposition~\ref{prop:common-domain-of-extension}
we can choose their extensions $\tilde{f}_i:b\rightarrow X_2$ \itwo  
with a common domain set $b$. Then the equality condition 
$$\forall x\in X_1 . f_1(x)=f_2(x),$$
which is unbounded since $X_1$ is proper, is equivalent 
to the bounded condition
\begin{equation}
\label{eq:20120228}
\exists c\subset b. \forall x\in c . \tilde{f}_1(x)=\tilde{f}_2(x),
\end{equation}
whence the indistinguishability is definite. 
Note that the validity of (\ref{eq:20120228}) is
independent of the choice of the extensions.

Hence the functions from $X_1$ to $X_2$  forms a class 
$Fun(X_1,X_2)$.

\end{eproof}

Suppose $X,Y$ are sets with tight enumerations $X=\seti{\seq{n}{x}}$
and $Y=\seti{\seq{m}{y}}$. Each $k\in [1..m^{n}]$
defines a function $f_k:X\rightarrow Y$ by the rule
$ f_k(x_i)=y_{j} $ if and only if $j-1$ is the number in the $i-1$-th
position of the $m$-arry expansion of $k-1$. Obviously any function
from $X$ to $Y$ is given as $f_k$ for some $k$, whence the 
collection of functions from $X$ to $Y$ forms a set $Y^{X}$ with
explicit enumeration $\seti{\seq{\kosuu{Y}^{\kosuu{X}}}{f}}$.

Thus we can define a function from a set $X$ to a set $Y$
by choosing an arbitrary element of $Y$ for each $x\in X$.
In particular we have the following choice principle.

We remark that if $f:X\rightarrow Y$ is a function between sets, then
it induces functions on the power sets. Namely, 
if $A\subset X$, $B\subset Y$ are subsets, then
$$
f(A):=\setii{f(a)}{a\in A}\subset Y,$$
$$f^{-1}(B):=\setii{x\in X}{f(x)\in B}\subset X
$$
are subsets. In particular, for $y\in Y$,
$$f^{-1}(y):=f^{-1}\seti{y}$$
is a subset.

\subsubsection{$\sigma$-finite Classes}
\label{sec:sigma-finite}
The union of an increasing family of sets is called \textit{$\sigma$-finite}. 
If the sequence $\seti{A_n}$ is strictly
increasing in the sense that $A_n\neq A_{n+1}$ for all $n$, then
$\bigcup_{n\in \nat}A_n$ is a proper class.

Two increasing sequences $\setii{A_n}{n\in \nat}$  and $\setii{B_n}{n\in \nat}$ are
called \textit{equivalent} if there are functions $f,g:\nat\ra \nat$
such that $A_n\subset B_{f(n)}$ and $B_n\subset A_{g(n)}$ hold
for all $n$. 
Obviously we have the following.
\begin{Proposition}
Let $\setii{A_n}{n\in \nat}$  and $\setii{B_n}{n\in \nat}$ by 
equivalent increasing sequences of sets. Then
$$\bigcup_{n\in \nat}A_n=\bigcup_{n\in \nat}B_n.$$
\end{Proposition}

For a $\sigma$-finite class $X$, an increasing sequence of
sets $\setii{A_n}{n\in \nat}$ with a function $\rho:X\rightarrow \nat$
called \textit{ranking} satisfying $x\in A_{\rho(x)}$ for all $x\in X$ is called 
its \textit{representation}. We write then $X=\bigcup_{n\in \nat}A_n$
with the implicit agreement of the existence of a ranking function. 
\index{class@$\sigma$-finite class, ranking}
\index{class@$\sigma$-finite class}
\index{class@$\sigma$-finite class, representation}

For example, $\nat,\integer,\rational$ are $\sigma$-finite classes with
the following representations.
\begin{eqnarray*}
\nat&=&\bigcup_{n\in \nat}[0..n],\\
\integer&=&\bigcup_{n\in \nat}[-n..n],\\
\rational&=&\bigcup_{n\in \nat}\rational_n.\\
\end{eqnarray*}
where
$$\rational_n:=\setii{\frac{p}{q}}{p,q\in [-n..n],q\neq 0}\subset \rational. $$

If $X=\bigcup_{n\in \nat}X_n$ is a representation of a $\sigma$-finite class 
then every subset is contained in some $X_n$. 
In fact if $a\subset X$, then $a\subset X_n$ with 
$$n:=\max\setii{\mbox{the rank of }u}{u\in a}.$$
Hence, the subsets of $X$ forms a $\sigma$-finite class $\pow(X)$ with
a representation
$$\pow(X)=\bigcup_{n\in \nat}\pow(X_n).$$

If $X$ is a set and $Y=\bigcup_nY_n$ is a $\sigma$-finite class, then
a function from $X$ to $Y$ 
is given by a map from $X$ to $Y_n$ for some $n$ and hence 
the collection of functions from $X$ to $Y$
forms the $\sigma$-finite class $Y^X$ with the representation $Y^X=\bigcup_nY_n^X$.

If a class $X$ is $\sigma$-finite, the elements of $X^n:=X^{[1..n]}$ is called
a sequence of length $n$ in $X$ and is written as $(\seq{n}{x})$. 
Hence, we have $\sigma$-finite classes $\nat^n,\integer^n,\rational^n$ of 
sequences of length $n$ for each $n$. 

Note that an objective subclass of a $\sigma$-finite class is $\sigma$-finite
by Axiom \ref{axiom:objective-seperation} namely the objective separation axiom.

If $X_i$ ($i\in [1..n]$) are $\sigma$-finite classes, 
their union class $\bigcup_{i\in [1..n]}X_i$ is obviously $\sigma$-finite.

The \textit{product class} $$\prod_{1\leq i\leq n}X_i$$ 
is defined as the collection of functions 
$$f:[1..n]\rightarrow \bigcup_{1\leq i\leq n}X_i$$
satisfying $f(i)\in X_i$ for all $i\in [1..n]$. 

The \textit{coproduct class} $$\coprod_{1\leq i\leq n}X_i$$ 
is defined as the collection of the pairs 
$$(i,x)\in [1..n]\times \bigcup_{1\leq j\leq n}X_j$$
such that $x\in X_i$. 

Obviously the product and the coproduct of $\sigma$-finite classes
are $\sigma$-finite.

\begin{Lemma}
If $X_i$'s are semisets, then the product $\prod_{1\leq i\leq n}X_i$ 
and the coproduct $\coprod_{1\leq i\leq n}X_i$ are semisets.
\end{Lemma}
\begin{eproof}
Suppose $X_i$ is a subclass of a set $a_i$ for $i\in [1..n]$. Then
the product class is a subclass of the set $\prod_{1\leq i\leq n}a_i$
and the coproduct class is a subclass of the set $\coprod_{1\leq i\leq n}a_i$.
\end{eproof}

\subsection{Induction Axioms}
\label{subsec:induction}

A condition $P(n)$ on numbers is called \textit{inductive} if 
$P(0)$ is true and for all $n$ if $P(n)$ is true then $P(n+1)$ is true.
Even when $P(x)$ is an inductive condition, we do not think it is
obvious that $P(a)$ is true for all $a\in \nat$ but rather we take it
as an evidence that we can refine the class $\nat$ so that the
condition $P(x)$ turns out to be true on it. In this way, we consider
the class $\nat$ open to continual refinements whenever new inductive
conditions are found. With the tacit understanding that such
refinements being done background automatically, we are convinced of
the validity of following axiom.
\begin{Axiom}[Strong Induction Axiom]
\label{axiom:definite induction}
\index{induction for definite conditions} If $P$ is
an objective definite inductive condition on the class \nat, 
$P(x)$ holds for every $x\in \nat$.
\end{Axiom}

A condition $P(n)$ on numbers, which may not be necessarily objective,
is called \textit{weakly inductive} if it satisfies that $P(0)$ is
true and for every \accessible{}{} number $n$, if $P(n)$ is true then
$P(n+1)$ is true.
\index{induction for vague conditions}
\begin{Axiom}[Weak Induction Axiom]
\label{axim:vague induction}
  If $P$ is a weakly inductive definite condition, 
  then $P(i)$ holds for every \accessible{}{} $i$. 
\end{Axiom}

Note that the condition that a primitive recursive function 
is totally defined is not bounded for most primitive 
recursive functions and hence we cannot show that 
they are total functions using Strong Induction Axiom. 

However the condition that a premitive recursive function is defined
for accessible numbers and has accessible values can by expressed by
bounded formula by virtue of the overspill principle, whence the weak
induction axiom shows that every primitive recursive function is
defined at least on accessible numbers with accessible values. 
We omit the detail.

The following is used frequently.
\begin{Lemma}
\label{lemma:prim-recurs-funct}
Let $f(x)$ be a primitive recursive function. Then for each huge $K$,
there is a huge $M$ such that $f(M)$ is defined and satisfies $f(M)<K$.
\end{Lemma}
\begin{eproof}
Since the set $\setii{n}{f(n)<K}$ contains all \accessible{} numbers,
it contains also a huge number $M$. 
\end{eproof}
For example, for every huge $K$, there is a huge $L$ with $L^{L}<K$.

\subsection{Overspill Principles}
The following is a special case of Theorem \ref{theorem:overspill}.
\begin{Theorem}[Overspill Principle]
\label{overspillaxiom}
Let $P$ be a definite objective condition on natural numbers. If all the
\accessible{}{} numbers satisfy the condition $P$, then there is a huge 
  number satisfying $P$.
\end{Theorem}

The contraposition of the theorem for the negation of $P$ gives the following
\begin{Corollary}
\label{cor:overspill} If all the \inaccessible{} numbers satisfy a definite
  objective condition $P$, then an \accessible{}{} number satisfy $P$.
\end{Corollary}

The following over spill principle is often used.
\begin{Theorem}
\label{th:overspill-equivalence}
Let $P$ be a definite objective condition.  
Then the following two conditions are equivalent:
\begin{description}
\item[(a)] There is an \accessible{}{} number above which every \accessible{}{} number satisfy $P$.
\item[(b)] There is a huge number under which every huge number satisfy $P$.
\end{description}
\end{Theorem}

\begin{eproof} Suppose (a) and every \accessible{}{} number greater than $k$
  satisfies $P$. Denote by $Q(n)$ the condition that
  every \accessible{}{} $x\in [k+1..n]$ satisfies $P$. 
  Since $Q$ is definite and 
  every \accessible{}{} number satisfies $Q$, the Theorem~\ref{overspillaxiom}
  implies $Q(K)$ for a huge number $K$, hence every
  huge number $I\leq K$ satisfies $P(I)$ since $I>k$, whence $(b)$.

Conversely suppose (b) and 
there is an $M\gg 1$ such that 
every huge numbers $I\leq M$ satisfies $P(I)$.
Denote by $R(x)$ the condition that every number $n\in [x..M]$ 
satisfies $P(n)$. Since $R$ is definite  and every huge number satisfies it,
the Corollary ~\ref{cor:overspill} implies that there is an  \accessible{}{} $i$ satisfying
$R(i)$, whence (a).
\end{eproof}

Taking its contraposition for the negation of $P$, we have
\begin{Theorem}
Suppose $P$ is a definite objective condition. Then the following conditions are equivalent.
\begin{description}
\item[(a)] For every \accessible{}{} number $n$, there is an \accessible{}{} number $x>n$ satisfying $P(x)$.
\item[(b)] For every huge number $N$, there is a huge number $K\leq N$ satisfying $P(K)$.
\end{description}
\end{Theorem}

\subsection{Concrete Sequences}
Let $A$ be a class. A function $a:\nat_{acc}\rightarrow A$
is called a \textit{concrete
  sequence on} $A$, which is often written as $a=(\seqinf{a})$.
A function $[1..n]\rightarrow A$ is called a \textit{huge} sequence
if $n$ is huge. The Axiom \ref{axiom:extension} implies the 
following extension property of concrete sequences.

\index{concrete sequence@concrete sequence}
\begin{Theorem}
\label{Theorem:extension-of-sequence}
 A concrete sequence in a class can be extended to a huge sequence in it.
More precisely, if $a=(\seqinf{a})$ is a concrete sequence in a class $A$, 
then there is a huge $N$ and a function $f:[1..N]\rightarrow A$ 
satisfying $f(i)=a_i$ if  $i$ is \accessible{}.
\end{Theorem}
\begin{eproof}
By Axiom~\ref{axiom:extension}, there is a subset $b$ such that 
$\nat_{acc}\sqsubset b\subset [1..N]$ and a function
$g:b\rightarrow A$ which restricts to $a$ on $\nat_{acc}$. 
The condition on $n$ that $[1..n]\subset b$ is 
obviously objective and definite. Moreover it is satisfied by all
$x\in \nat_{acc}$, whence by Theorem \ref{overspillaxiom},
there is a huge $N\gg 1$ such that $[1..N]\subset b$, whence $f=g|[1..N]$ is an extention
with the desired properties.
\end{eproof}

Although extensions of a concrete sequence to huge sequences are not unique,
the ``germ'' of the extensions is unique in the following sense.
\begin{Proposition}
If $N_i$ \itwo are huge numbers and maps
$$f_i:[1..N_i]\rightarrow A, \quad i=1,2$$
satisfy $f_i(k)=a_k$ for accessible $k$ \itwo.  Then there is a huge
$K\leq \min\seti{N_1,N_2}$ such that $f_1(j)=f_2(j)$ for $j\leq K$.
\end{Proposition}
\begin{eproof}
Let $N=\min\seti{N_1,N_2}$. Since the condition on natural numbers $n$ that
\begin{equation}
\label{eq:3-827}
n\leq N \mbox{ and } f_1(n)=f_2(n)
\end{equation}
hold for every \accessible{} $n$, there is a huge $K$ such that \eqref{eq:3-827}
holds for any $n\leq K$. From  \eqref{eq:3-827} for $n=K$, we have
$K\leq N$.
\end{eproof}

We say a class $A$ a \textit{quasi-set} if every concrete sequence in 
$A$ can be extended to a huge sequence in $A$.
\index{quasi-set}

We use often the following lemma.
\begin{prop}
\label{prop:minimum-huge-number}
Suppose, for $f:\nat_{acc}\rightarrow \nat$ is a function such that
$f(i)$ is huge for all $i$. Then there is a huge number $N$ satisfying 
$N\leq f(i)$ for all \accessible{}{} $i$.
\end{prop}
\begin{eproof}
  By Theorem \ref{Theorem:extension-of-sequence}, the concrete
sequence $f$ can be extended to a huge sequence $\tilde{f}:[1..M]\rightarrow \nat$ 
for some huge $M$. Define a mapping $g:[1..M]\rightarrow \nat$
by 
   $$g(i):=\min_{1\leq j\leq i}f(j).$$
Then, for accessible $i$, the following holds  
\begin{equation}
\label{eq:1-902-ff}
i<g(i),
\end{equation}
\begin{equation}
\label{eq:2-902-aa}
g(i)\leq g(j) \mbox{ for all }j\leq i.
\end{equation}
Since these conditions are definite, there is a huge number $K$ for which
the conditions \eqref{eq:1-902-ff} and ~\eqref{eq:2-902-aa} hold for $i=K$.
Take then $ N=g(K)$.
\end{eproof}

If $P$ is definite, Theorem \ref{overspillaxiom} implies the following.
\begin{Proposition}
  Let $P$ be a definite objective weakly inductive condition. 
Then there is a huge $M$ such that $P(i)$ holds for all $i\leq M$.
\end{Proposition}
\begin{eproof}
By the weak induction axiom, $P(n)$ is true for all \accessible{}{} $n$, 
hence by the overspill principle of Theorem \ref{overspillaxiom},
there is a huge $M$ such that $P(n)$ is true for $n\leq M$.
\end{eproof}

\newpage
\section{Continuum}
\def\runningtitle{Continuum}
\label{sec:basic-continuum}
\subsection{Sorites Relations}
A subclass $R$ of the product class $X\times X$ of a class $X$ 
is called a binary relation on $X$. 

A basic example is the binary relation $x\approx y$ of indistinguishablity
between rationals $x,y$. 
This relation is not objective but is definite 
since it can be expressed by bounded quantifier
$$ \forall n\leq N. \mbox{if $n$ is accessible then $|x-y|<\frac1n$},$$
using any huge $N$.

Usual notions of symmetry, anti-symmetry,
reflexivity, transitivity have meaning for $R$, 
with the proviso that the validity of unbounded $\forall$-statements is
understood proof-theoretically as in \S\ref{sss:subclasses}. 

A sequence $(\seq{N}{x})$ is called an \textit{$R$-chain} 
\index{R-chain@$R$-chain}
if $x_iRx_{i+1}$ for $i\in [1..N-1]$. 
\begin{Proposition}
If a transitive binary relation $R$ is objective, the following condition holds
\begin{equation}
\label{eq:2012-02-14}
\mbox{If  $(\seq{n}{x})$ is an $R$-chain on $X$, then $x_1Rx_n$}
\end{equation}
for all $n\in \nat$.
\end{Proposition}
\begin{eproof}
Suppose $(\seq{n}{x})$ is an $R$-chain.
Since $R$ is objective, the class $A:=\setii{i}{x_1Rx_i}\subset[1..n]$ 
is a set. Suppose its greatest element $m$ is less than $n$. Then
$x_1Rx_{m}$ and $x_{m}Rx_{m+1}$ but not $x_1{R}x_{m+1}$, which contradicts to the transitivity of $R$. Hence $m=n$ and we have $x_1Rx_n$. 
\end{eproof}

Note that if $R$ is not objective, the subclass $A$ in the proof may
be proper and have no greatest element and the above arguments fail.
Hence (\ref{eq:2012-02-14}) might not hold for huge $n$ although it
holds for \accessible{} $n$ by the weak induction axiom.

The relation $R$ is called \textit{strictly transitive} if
(\ref{eq:2012-02-14}) holds for every $n$ and an equivalence relation
is called \textit{strict} if it is strictly transitive.  An
equivalence relation which is not strict is called a \textit{sorites
  relation}. An $R$-chain $(\seq{N}{x})$ without the validity of $x_1{R}x_{N}$ is
called a \textit{sorites sequence}.

For example, the equivalence relation $\approx$ on 
the set $[1..N]$ defined by $i\approx j$ if and only if the
rational number $\frac{i-j}N$ is infinitesimal
is a sorites relation 
since $(1,2,\cdots,N)$ is a sorites sequence.
In fact $i\approx i+1$ for $i\in [1..N-1]$ but $1\not\approx N$.

\index{sorites sequence} 
\index{sorites relation} 
\index{transitive, weakly} 
\index{transitive, strongly} 
\index{equivalence relation, weak} 
\index{equivalence relation, strong}

\subsection{Continuum}
Continuum is a pair $C=(\support{C},\approx_C)$ of a class $\support{C}$ 
and an equivalence
relation $\approx$ on $\support{C}$ which might not be strict. 
The class $\support{C}$ is called the \textit{support}
\index{continuum, support}
 of the continuum $C$.

Elements of $\support{C}$ are called the \textit{positions}
\index{continuum, position}
 of $C$. The relation $\approx_C$ is called the
\textit{indistinguishability relation}
\index{continuum, indistinguishability relation}
\index{$p\approx q$ on continuum}
 of $C$. We say two positions $p,q$ are 
indistinguishable if $p\approx_C q$. 

For a position $a\in C$, \textit{the point determined by} $a$, 
is defined as the subclass
$$
[a]:=\setii{x\in \support{C}}{x\approx a},
$$
which is proper in most cases. 
A point of $C$ is the point determined by some position of $C$.
The notation $a\in C$ stand for the phrase that $a$ is a point of $C$.
\index{continuum, point}
For $p\in C$, a position $x\in |C|$ such that $x\in p$
is said to \textit{represent the point} $p$ and the point $p$ is
represented by the position $x$.

\begin{Remark}
The collection of points do not form a class in most cases, 
since the condition of equality of proper classes is not definite.
This is reasonable since if it formed a class, then we would have 
a paradox ``$[x_1]=[x_2]=\cdots=[x_N]$ but $[x_1]\neq [x_N]$'' 
if $(\seq{N}{x})$ is a sorites sequence.

This conforms to the view that the ``points'' of a continuum have
inevitable indefiniteness which however is not perceived by any observation 
however accurate it may be.
\end{Remark}

A continuum $C$ is called a \textit{mesh continuum} if $\support{C}$ is a semiset
and a \textit{rigid mesh continuum} if $\support{C}$ is a set.
\index{continuum, mesh}
\index{continuum, rigid mesh}

\subsubsection{Examples}
\label{subsec:examples-of-continuum}
\paragraph{Linear continuum}
A basic example is the continuum $(\rational,\approx)$, called the linear continuum 
denoted by $\real$.
\index{continuum, linear}

\paragraph{Metric continuum}
A \textit{metric class} $(X,d)$ is a class with a function $d:X\times
X\rightarrow \rational$ satisfying the usual property of distance
function. The relation $x\approx_dy$ defined by $d(x,y)\approx 0$
is an equivalence relation which might not be strict.
We call $(X,\approx_d)$ the \textit{metric continuum} defined by $(X,d)$.
\index{continuum, metric}
If $X$ is a semiset, $(X,d)$ is called a metric semispace and
$(X,\approx_d)$ is a mesh continuum.
If $X$ is a set, $(X,d)$ is called a metric space and the continuum 
$(X,\approx_d)$ is a rigid mesh continuum. 

\paragraph{Subcontinuum}

If $C$ is a continuum class, a subclass $Y\subset \support{C}$ defines
a continuum $(Y,\approx|Y)$ called the \textit{subcontinuum} of $C$ with 
support $Y$. 
\index{continuum@subcontinuum}
\index{metric space}
\index{metric class}
\index{metric semispace}

\paragraph{Interval continuum}
Let $a\in \rational$. Then the definite conditions $|x|<\infty$,
$a\prec x<\infty$, $a\leq x<\infty$, $-\infty<x\prec a$ and $-\infty<x\leq a$ 
define respectively the subclasses
$(-\infty,\infty)_{\rational}$, $(a,\infty)_{\rational}$,
$[a,\infty)_{\rational}$, $(-\infty,a)_{\rational}$ and
$(-\infty,a]_{\rational}$. 

Let $a,b\in \rational$ with $a\prec b$. Then the definite condition
$a\prec x \prec b$ defines the subclass $(a,b)_{\rational} \subset \rational$. 
Similarly the conditions $a\leq  x \prec b$, $a\prec x\leq b$ and 
$a\leq x\leq b$ define
respectively the objective subclasses $(a,b)_{\rational}, [a,b)_{\rational}$,$(a,b]_{\rational}$ and $[a,b]_{\rational}$. 
Note that only $[a,b]_{\rational}$ is objective subclass.

The nine strings $(-\infty,\infty)$, $(a,\infty)$,$[a,\infty)$, $(-\infty,a)$,$(-\infty,a]$,$(a,b)$,$[a,b)$, $(a,b]$,$[a,b]$ are called \textit{interval symbols}. 
\index{interval symbols}
The interval symbols without $\infty$ is called \textit{finite interval symbols}.
\index{interval symbols, finite}
If $I$ is an interval symbol, then the subclass $I_{\rational}\subset\rational$ 
defines a subcontinuum denoted by $I$. For example $[0,1]$ denotes the subcontinuum
$([0,1]_{\rational},\approx)$. 

Let $r$ be a nonzero rational number. We write by $r\integer$ the class of rationals which can be written as $nr$ with $n\in \integer$. 
For an interval symbol $I$, we write $I_{r}:=I_{\rational}\bigcap r\integer$.
\begin{Lemma}
If $I$ is an interval symbol then $I_{r}$ is a semiset. If $I=[a,b]$, then
$I_r$ is a set.
\end{Lemma}
\begin{eproof}
Let $I$ be an interval symbol. Let $M\gg 1$. 
Then $I_r\subset [-M,M]_r$. Take $K$ such that $rK>M$, then
$$I_r\subset [-M,M]_r\subset\setii{xr}{x\in [-K..K]}$$
whence $I_r$ is a semiset.

Suppose $I=[a,b]$. Since $a\leq xr\leq b$ means $\frac{a}{r}\leq x\leq \frac{b}{r}$,
we have $I_r=\setii{xr}{x\in [a'+1,b']}$ where
$a'$ and $b'$ are the integer parts of $\frac{a}{r}$ and $\frac{b}{r}$ respectively.
Hence $I_r$ is a set.
\end{eproof}

For each rational $r\neq 0$ and interval symbol $I$, we obtain a
subcontinuum with support $I_r$ denoted also by $I_r$, which are mesh
continuum by the above lemma, among which $[a,b]_r$ is rigid.

\paragraph{Euclidean continuum}
\index{continuum, Euclidean}
Since $\rational$ is $\sigma$-finite, every number $n$  defines
the product $\rational^n$ with the metric function $d_{\infty}$ 
defined by 
$$  d_{\infty}(x,y):=\max_{1\leq i\leq n}|x_i-y_i|, $$
where $x_i$ is the $i$-th coordinate of $x\in \rational^n$. 
The continuum $(\rational^n,\approx_{d_{\infty}})$ is called the \textit{$n$-dimensional Euclidean continuum}.
\index{continuum, $n$-dimensional Euclidean}

The subcontinuum defined by the subclass $[0,1]_{\rational}^n\subset \rational^n$
is called \textit{the unit $n$-hypercube}.
\index{unit $n$-hypercube}

\paragraph{Product continuum}
\index{continuum}
Let $C_i$ ($i\in [1..N]$) be rigid mesh continua. 
Then the product set $\prod_{i\in [1..N]}\support{C_i}$
has an equivalence relation $x\approx y$ 
defined by $x_i\approx_{C_i} y_i$ for all $i$.
We call the rigid mesh continuum $(\prod_{i\in [1..N]}\support{C_i},\approx)$ 
the \textit{product continuum} of the family $\setii{C_i}{i\in [1..N]}$
and denote it by $\prod_{i\in [1..N]}C_i$
\index{continuum, product}.

Let $C_i$ ($i\in [1..k]$) be mesh continua with $k$ a very small number
so that we can write $1,2,\cdots,k$ without the ellipsis.
Then we have the product semiset $\prod_{i\in [1..k]}\support{C_i}$ and an
equivalence relation $\approx$ defined by
$$x\approx y \mbox{ if and only if } x_i\approx y_i \mbox{ for all $i$}.$$

\paragraph{Graphs as rigid mesh continua}

Let $G=(V,E)$ be a connected symmetric graph with $V,E\subset X\times X$ 
being sets. 
Let $w:E\rightarrow \rational$ be a symmetric positive valued
function.  Define the length of a path $\gamma=(\seq{n}{x})$ by
$\ell_w(\gamma)=\sum_{1\leq i< n}w((x_i,x_{i+1}))$.  Let $d_{w}(p,q)$
be the minimum of the length of paths connecting $p$ and $q$. Then
$d_{G,w}$ is a rational valued metric function on $V$.  The rigid mesh
continuum defined by $(V,d_w)$ is called the \textit{rigid mesh
  continuum generated by the graph $G$ with the edge length function
  $w$.}  \index{continuum generated by a graph}

\begin{Example}
Fix a huge $\Omega\gg1$ and put $V_{\Omega}:=\coprod_{i\leq \Omega}\seti{0,1}^i$, the set of
finite $01$-words of length less than or equal to $\Omega$. Define
$$E=\setii{\pair{w,wi},\pair{wi,w}}{w\in V_{\Omega-1},i=0,1}.$$
Then the graph $G=(V,E)$ is the binary tree of depth $\Omega$. 
If we give uniform infinitesimal length $\frac1\Omega$ to edges, we 
obtain a continuum of hyperbolic type. If the length of the edge $\pair{w,wi}$  
is given $2^{-|w|}$ where $|w|$ denotes the length of the word $w$, then
the continuum induced from $(G,w)$ is the Cantor spaces. In \S\ref{subsec:binary-tree},
we study the topological properties of these continua. 
\end{Example}
\newpage

\subsection{Morphism}
\label{sec:morphisms-continua}
Let $C_i$ \itwo be continua. 
A function  $f:|C_1|\rightarrow |C_2|$
is called \textit{continuous}
\index{continuous function}
if $x\approx y$ implies $f(x)\approx f(y)$. 

We call two continuous functions $f,g:|C_1|\rightarrow |C_2|$ 
\textit{indistinguishable} 
and write $f\approx g$ if we can prove $f(x)\approx g(x)$ for all $x\in |C_1|$.

The collection of continuous function indistinguishable from $f$  is not
generally a class. 
So we formally introduce
a symbol $[f]$ and use it as if it were a class as follows. 
\begin{enumerate}
\item A \textit{morphism} from $C_1$ to $C_2$ is a symbol of the form $[f]$ 
for some continuous function from $\support{C_{1}}$ to $\support{C_2}$. 
\item If $C_i$ \itwo are continua, the notation $F:C_1\rightarrow C_2$ means that 
$F$ is a morphism from $C_1$ to $C_2$.
\item If $F$ is a morphism, then the expression $g\in F$ means $g\approx f$ if $F=[f]$.
If $g\in F$, we say that the morphism $F$ is \textit{represented by} $g$ and
\textit{$g$ represents $F$}.
\item If $F$ and $G$ are morphisms represented respectively by $f$ and $g$, 
then the expression $F=G$ means $f\approx g$. 
\item A condition on continuous functions is called a \textit{condition on
morphisms} if it is $\approx$-invariant.
\end{enumerate}
\index{continuum morphism}

The \textit{identity morphism} $id_C$ is represented by the identity map $id_{\support{C}}$.

The \textit{value of a morphism} $F$ \textit{at a point} $p$ of 
$C_1$ is defined to be the point  $[f(t)]$ for $f\in F$ and $t\in p$.
This does not depend on the choice of representations.

If $F_1:C_1\rightarrow F_2$ and $F_2:C_2\rightarrow C_3$, then the composition
$F_2\circ F_1:C_1\rightarrow C_3$ is the morphism, 
represented by $f_2\circ f_1$,
where $f_i\in [F_i]$ \itwo, which does not depend on the choices of $f_i$ \itwo.
We may write the definition symbolically by
$$(F_2\circ F_1)(p):=F_2(F_1(p)) \mbox{ for } p\in C_1,$$
with proviso that the precise meaning is understood as above since
a morphism cannot be defined as a correspondence which maps
points to arbitrary points.

When $C_i$ \itwo are rigid mesh continua, 
we will construct in \S~\ref{subsec:continuum-functions} a continuum $C(C_1,C_2)$
whose points are precisely morphisms from $C_1$  to $C_2$.

\begin{Example}[Morphism defined by $\frac1x$]
  Taking the inverses of nonzero rationals define a function
  $f:(0,1]_{\rational}\rightarrow [1,\infty)_{\rational}$. 
  In fact, if $x\in (0,1]_{\rational}$, namely,  $0\prec x\leq 1$,
 then $ 1\leq \frac1x<\infty$ and $\frac1x\in [1,\infty)_{\rational}$.
  Moreover if $r,s\in (0,1]_{\rational}$ and $r\approx s$, 
  then since $r,s\not\approx0$, we have
  $\frac1{r}\approx \frac1{s}$ as will be shown in \S\ref{sec:algebr-oper-real}. 
  Thus  $[f]$ is a morphism from $(0,1]$ to $[1,\infty)$.
\end{Example}

\begin{Example}[Morphism to product continuum mesh ]
Let $C$ and $C_i$ ($i\in [1..,n]$) be continua.
An $n$-tupple $(\seq{n}{F})$ of morphisms 
$F_i:C\rightarrow C_i$
defines a morphism $F:C\rightarrow \prod_{1\leq i\leq n}C_i$ 
which is represented by $f$ which assigns $a\in C$ to 
$f(a):=(f_{i}(a))_{i\in [1..n]}$, where $f_i\in F_i$ ($i\in [1..n]$).

Conversely a morphism $F:C\rightarrow  \prod_{i=1}^nC_i$ defines
morphisms $F_i=\pi_i\circ F:C\rightarrow C_i$ ($i\in [1..n]$),
where $\pi_i:\prod_{i=1}^nC_i\rightarrow C_i$ is the projection morphism
represented by $\pi_i:(\seq{n}{x})\mapsto x_i$.
The morphism $F_i$ is called the \textit{$i$-th component} of $F$.
\end{Example}

\subsection{Equivalence}
\label{sec:equivalence-continua}
A morphism $F:C_1\rightarrow C_2$ is called \textit{injective} and
\textit{surjective} if it is represented by a continuous map 
$f:\support{C_1}\rightarrow \support{C_2}$
satisfying respectively
\begin{equation}
\label{eq:903-1}
\mbox{$f(x)\approx f(y)$ implies $x \approx y$ for all $x,y\in C_1$,}
\end{equation}
and
\begin{equation}
\label{eq:903-2}
\mbox{for every $x_2\in C_2$, there is an $x_1\in C_1$ with $f(x_1)\approx x_2$.}
\end{equation}
Note that unless $C_i$ \itwo are rigid mesh continua, these conditions are not definite.
Note also that the above conditions are independent of the choice of $f\in F$.

A morphism $F:C_1\rightarrow C_2$ is an \textit{equivalence} 
\index{continuum, equivalent}
if there is a morphism  $G:C_2\rightarrow C_1$ satisfying 
$$G\circ F=id_{C_1} \quad \mbox{ and } F\circ G=id_{{C_2}}.$$
The morphism $G$ is uniquely determined by $F$ and is called the \textit{inverse} of $F$ and is denoted by $F^{-1}$.

If and equivalence $F$ and its inverse $F^{-1}$ are represented respectively by 
$f:\support{C_1}\rightarrow \support{C_2}$
and 
$g:\support{C_2}\rightarrow \support{C_1}$,
then 
$$g\circ f\approx id_{\support{C_1}} \quad \mbox{ and } f\circ g\approx id_{\support{C_2}}.$$
Such $g$ is uniquely determined by $f$ up to indistinguishability 
and is called an \textit{almost inverse} of $f$.
\index{almost inverse}

\begin{Proposition} Suppose $C_i$ \itwo are rigid mesh continuum. Then
a morphism $F:C_1\rightarrow C_2$ is an equivalence 
if and only if it is injective and surjective.
\end{Proposition}
\begin{eproof}
Suppose $F:C_1\rightarrow C_2$ is an equivalence and let $f\in F$
and $g\in F^{-1}$.
Then $f(x)\approx f(y)$ implies 
$$x\approx g(f(x))\approx g(f(y))\approx y.$$
Furthermore, for every $x_2\in C_2$, we have $f(x_1)\approx x_2$ if we put $x_1=g(x_2)$. 

Conversely suppose $f\in F$ satisfyies the conditions (\ref{eq:903-1}) 
and (\ref{eq:903-2}). 
For any $x_2\in \support{C_2}$, we can choose by (\ref{eq:903-2}) an $x_1\in \support{C_1}$ 
satisfying $f(x_1)\approx x_2$. Define $g(x_2):=x_1$. Then the map $g:\support{C_2}\rightarrow \support{C_1}$ 
is continuous by the condition (\ref{eq:903-1}). 

By definition $f(g(x_2))\approx x_2$ hold for $x_2\in \support{C_2}$.  
For $x_1\in \support{C_1}$, we have $f(g(f(x_1)))\approx f(x_1)$ by definition
of $g$, whence by the condition (\ref{eq:903-1}) we obtain $g(f(x_1))\approx x_1$.
Hence $g$ is an almost inverse to $f$.
\end{eproof}

If there is an equivalence $F:C_1\rightarrow C_2$, we say
the continuum  $C_1$ is \textit{equivalent} to $C_2$ and write
$C_1\simeq C_2$.  Then $\simeq $ satisfies the axiom of 
equivalence relations. 

Let $C$ be a continuum and $C_i\subset C$ \itwo be subcontinua.
An equivalence $\alpha:C_1\simeq C_2$ is called a \textit{quasi-identity}
\index{quasi-identity}
if it satisfies $\alpha(x)=x$ for all $x\in C_1$. 
By definition, a quasi-identity is uniquely determined if it exsits.

\subsubsection{Examples}
\label{example:interval-equivalence}
\begin{Proposition}
\label{prop:r-int-equiv-rat} 
Let $r$ be a nonzero infinitesimal rational number. The inclusion
$$  \imath_r:r\integer \rightarrow \rational$$
represents a quasi-identity, for which the function 
$\kappa_r:\rational \rightarrow r\integer $ 
defined by
$$ \kappa_r(s):=[sr^{-1}]r,$$
where $[x]$ denotes the integer part of $x$,  
gives an almost inverse of $\imath_{r}$.
\end{Proposition}
\begin{eproof}
Obviously $\kappa_r(\imath(nr))=nr$ for $n\in \integer$. 

From $x\leq [x]< x+1$ it follows 
$$s=(sr^{-1})r< [sr^{-1}]r < (sr^{-1}+1)r=s+r$$
whence $[sr^{-1}]r\approx s$ since $r\approx 0$. Hence
$$\imath_r\circ \kappa_r\approx id.$$

Hence $\kappa_r$ is an almost inverse to $\imath_r$.
\end{eproof}

If $I$ is an interval symbol, the inclusion
$$\iota_r|I_{r}:I_r\rightarrow I_{\rational}$$
represents a quasi-identity $(I_r,\approx)\rightarrow I$.
\begin{Corollary}
If $I$ is an interval symbol, then
$$(I_r,\approx) \simeq I$$
for every nonzero infinitesimal rational $r$.
\end{Corollary}
The mesh continuum $(I_r,\approx)$ 
is called a \textit{representation of the continuum} $I$.
\index{continuum, representation}
Obviously representations are unique up to equivalences.

For example, let $r,s$ be nonzero infinitesimal rationals and $I$ an interval symbol.
Then the representations $(I_r,\approx)$ and $(I_s,\approx)$ of $I$ 
are equivalent by the morphism given by 
$$
h(nr)=[nr/s]s,
$$
which is the restriction $h:=g_{s}\circ \iota_r:r\integer \rightarrow s\integer $
on $I_r$.

\begin{Example}
Let $a\prec b$ be finite rationals. Then $[a,b]\approx[0,1]$.
In fact define $f:[a,b]_{\rational}\rightarrow [0,1]_{\rational}$
by $f(x)=\frac{x-a}{b-a}$. Then it has an inverse
$g:[0,1]_{\rational}\rightarrow [a,b]_{\rational}$
defined by $g(x)=(b-a)x + a $.

Since $([a,b]_r,\approx)\simeq [a,b]$ for nonzero infinitesimal $r$,
we have equivalence between rigid mesh continua $f_{rs}:[a,b]_r\approx [0,1]_s$  for each nonzero infinitesimals $r,s$. This equivalence is given by
$$
f_{rs}(nr)=\left[\frac{nr-a}{s(b-a)}\right]s.
$$
\end{Example}

\subsection{Saturation}
Let $C$ be a continuum. A subset $a\subset \support{C}$ 
defines a subclass $\overline{a}\subset \support{C}$ 
\index{saturation}
\index{$\overline{a}$}
of such elements $y$ as satisfying the definite condition 
that there is an $x\in a$ with $y\approx x$. 

A subset $a\subset \support{C}$ is called \textit{dense in $C$} \index{dense} 
if $\overline{a}=\support{C}$. For example $[a,b]_{r}\subset [a,b]_{\rational}$ 
is dense if $r$ is nonzero infinitesimal.  

Obviously the correspondence $a\mapsto \overline{a}$ satisfies the following 
conditions of closure operators
\begin{Proposition}
\begin{enumerate}
\item $a\subset \overline{a}$,
\item $a\subset b$ implies $\overline{a}\subset \overline{b}$,
\item If $\overline{a}$ is a set, then $\overline{\overline{a}}=\overline{a}$.
\end{enumerate}

Moreover the following holds.
\begin{equation}
\label{eq:appearance-closure-union}
\overline{a\bigcup b}=\overline{a}\bigcup \overline{b}
\end{equation}
\end{Proposition}

A subset $a\subset \support{C}$ is called \textit{saturated} if $\overline{a}=a$. 
The saturated subsets are rare. In fact we have
\begin{Proposition}
Suppose  $C$ is connected, namely, every two elements are connected by 
an $\approx$-chain. Then a saturated subset is either $\emptyset$ or $\support{C}$.
\end{Proposition}
\begin{eproof}
Suppose $a=\overline{a}$ and $a\neq\emptyset$ and $a\neq \support{C}$. 
Let $x\in a$ and $y\nin a$. Let $(\seq{N}{x})$ be a $\approx$-chain 
such that $x_1=x$ and $x_N=y$. Let $i\in [1..N]$ be the minimum 
satisfying $x_i\nin a$. Then $x_{i-1}\in a$ and $x_i\nin a$.
The former condition and $x_i\approx x_{i-1}$ implies $x_{i}\in \overline{a}=a$
which contradicts to the latter condition.
\end{eproof}

\newpage
\section{Topology of Continuum}
\label{sec:topology-continua}
\def\runningtitle{Topology of continua}

\subsection{Convergence of Sequences}
Let $C$ be a continuum.
\label{subsec:conv-sequ}
A sequence $a=(\seq{N}{a})$ in $\support{C}$ is called a sequence in $C$.
A sequence $a$ in $C$ \textit{converges}
\index{convergence@congergence, finite sequences}
to $c\in \support{C}$ if there is a
huge $K\leq N$ such that for every huge $I\leq K$, $a_I\approx c$.
The limit $c$ is uniquely defined up to indistinguishability.

Let $F:C_1\rightarrow C_2$ be a morphism and $f\in F$. Then a sequence
$a$ in $C_1$ defines a sequence $f(a)=(f(a_{1}),\cdots,f(a_{N}))$ 
in ${C_2}$ and if $c\in \support{C_1}$ is a limit
of the sequence $a$ then $f(c)$  is a limit of the sequence $f(a)$.

Let $B\subset \support{C}$ be a subset and $c\in C$.  
We call $c$ \textit{an accumulation  point of $B$}
\index{accumulation point of subset} 
if 
there is a huge subset $B'\subset B$ whose elements are indistinguishable from $c$. 
We call $c\in B$ an \textit{isolated point} 
if $c\neq b\in B$ implies $c\not\approx b$. 
Note that there can be an element $c\in B$ which is neither isolated nor 
accumulation point of $B$ since it is possible
that $\overline{\seti{c}}\bigcap \left(B\setminus\seti{c}\right)$ 
contains only one element.

A continuum $C$ is called \textit{perfect} if every $c\in C$ is
an accumulation point of $\support{C}\setminus\seti{c}$. 
\index{perfect}
Basic Euclidean continua are obviously perfect. 

\subsection{Compactness}
A continuum $C$ is called \textit{compact} 
\index{continuum, compact}
if for every huge $N$ there is a dense subset $A\subset \support{C}$ with $\kosuu{A}\leq N$.

\begin{Proposition}
\label{prop:compactness-of-[a,b]}
For finite rationals $a,b$ with $a\prec b$, the continuum $[a,b]$ is compact.
\end{Proposition}
\begin{eproof}
For every huge $N$, the subset $[a,b]_{\frac1{N-1}}\subset [a,b]_{\rational}$
is dense and has $N$ elements.
\end{eproof}

A continuum equivalent to a compact continuum is compact. In fact we have
the following.
\begin{Proposition}
If there is a surjective morphism $F:C_1\rightarrow C_2$ and $C_1$ is compact 
then $C_2$ is also compact.
\end{Proposition}
\begin{eproof}
Let $f\in F$. Let $N$ be a huge number. Take a dense subset 
$A\subset \support{C_1}$ 
with $\kosuu{A}\leq N$. Then $f(A)\subset \support{C_{2}}$ is dense since for every 
$x\in \support{C_2}$,
there is an $y\in \support{C_1}$ with $f(y)\approx x$. Take $a\in A$ such that $a\approx y$.
Then $x\approx f(y)\approx f(a)\in f(A)$. Hence $f(A)$ is dense in $\support{C_2}$.
\end{eproof}

\begin{Theorem}
\label{th:existence-of-accumulation-in-compactspace}
If a continuum $C$ is compact, then every huge subset $B\subset \support{C}$  
has an accumulation point.
\end{Theorem}
\begin{eproof}
  Suppose $C$ is compact and $B\subset \support{C}$ a huge set.  
  Take a huge $N$ with $N^2<\kosuu{B}$. Let $E\subset \support{C}$ be a dense
  subset with $\kosuu{E}\leq N$. Define a map $f:B\rightarrow E$ which
  carries $b\in B$ to an $e\in E$ such that $e\approx b$. 
  If $\kosuu{f^{-1}e}<N$ for all $e\in E$, then 
  $$\kosuu{B}=\sum_{e\in E}\kosuu{f^{-1}e}< N\;\kosuu{E}\leq N^2<\kosuu{B}  $$
  a contradiction. Hence there is an $e\in E$ with $\kosuu{f^{-1}e}\geq \kosuu{E}$. 
  Let $b\in B\bigcap f^{-1}e$. Then $x\in f^{-1}e$ implies
  $x\approx e \approx b$ whence every element of $f^{-1}e$ 
  is indistinguishable from $b$, whence
  $b$ is an accumulation point of $B$.
\end{eproof}

\begin{Corollary}
\label{cor:compactness}
If $C$ is compact, every huge subset of $\support{C}$ has a pair of indistinguishable elements.
\end{Corollary}

A subset  $a\subset \support{C}$ is called
\textit{discrete}
\index{discrete subset}
if $x,y\in a$ and $x\approx y$ implies $x=y$. 
For example $[-M..M]$ is
a discrete subset of $[-M,M]_{\rational}$. 

\begin{Proposition}
\label{prop:noncompactness}
A continuum with a discrete huge subset is not compact.
\end{Proposition}
\begin{eproof} 
Suppose a continuum $C$ is compact with a huge discrete subset $a\subset \support{C}$.
Let $b\subset \support{C}$ be a dense subset with $\kosuu{b}<\kosuu{a}$.
Define a map $f:a\rightarrow b$ by assigning $x\in a$ to an $y\in b$ such 
that $x\approx y$. Since $\kosuu{b}<\kosuu{a}$, there must be $x,y\in a$ 
with $f(x)=f(y)$ but $x\neq y$. Then
$x\approx f(x)=f(y)\approx y$ contradicting to the discreteness of $a$.
\end{eproof}

Let $C$ be a continuum. 
An objective subclass $R\subset \support{C}\times \support{C}$ is called an \textit{objective
  discrimination of $C$}
\index{objective discrimination@objective discrimination}
if $R(x,y)$ implies $x\not\approx y$.  A subclass 
$A\subset X$ is called $R$-discrete if $x\neq y$ implies $R(x,y)$
for $x,y\in A$. For
example, in $\rational$, if $k$ is \accessible{} then the relation
$R_k:=\setii{(x,y)}{|x-y|>\frac1k}$ is an objective discrimination 
and the subclass $\integer\subset \rational$ is $R_{2}$-discrete.

\begin{Proposition}
  A continuum $C$ is not compact if it has an objective
  discrimination $R$ and a subset $a\subset \support{C}$ 
  such that for each \accessible{} number $k$, 
  there is an $R$-discrete subset of $a$ of size greater than $k$.
\end{Proposition}
\begin{eproof}
  Let $R$ be an objective discrimination of $C$ and suppose that for each \accessible{} $k$ the condition
\begin{equation}
\label{eq:20110223-1}
\mbox{there is an $R$-discrete subset of $a$ with at least $k$ elements}
\end{equation}
is satisfied. Since the condition~(\ref{eq:20110223-1}) is objective
and definite and satisfied by all \accessible{} $k$, it is satisfied also
by a huge $k$. Hence by Proposition~\ref{prop:noncompactness}, $C$ is not
compact.
\end{eproof}

\begin{Example}
\begin{enumerate}
\item If $m-n>0$ is huge, the continuum $[n,m]$ is not compact since it has the huge 
discrete subset $[n..m]$.

\item The continuum $(-\infty,\infty)$ is not  compact. 
In fact for any concrete $k$, 
the subset $[1..k]\subset (-\infty,\infty)_{\rational}$ has $k$ elements and $R$-discrete
for the objective discrimination  
  $R=\setii{(x,y)}{|x-y|>\frac12}$.
\end{enumerate}
\end{Example}

\subsection{Connectedness}
Let $C$ be a continuum.  We say $x\in C$ is \textit{connected to}
$y\in C$ and write $x\smile y$ if there is an $\approx$-chain
connecting $x$ and $y$. If $x\smile y$ for every $x,y\in C$, the
continuum is called \textit{connected}. Note that this condition is not definite
generally.
\index{continuum, connected} 

Obviously, if $f:C_1\rightarrow C_2$ is a surjective morphism 
and $C_1$ is connected then $C_2$ is also connected.
In particular connectedness is equivalence invariant.  

Suppose $C$ is a rigid mesh continuum. 
Then the binary relation $x\smile y$ is definite and 
for each $x\in \support{C}$, we have the equivalence class
$$[x]:=\setii{y}{x\smile y}$$
which is a semiset called the \textit{connected component} containing $x$. 
\index{connected component}
Each $x\in \support{C}$ belongs to the connected component $[x]$.
A rigid mesh continuum $C$ is called \textit{totally disconnected}
if $[x]=\overline{\seti{x}}$ for all $x\in \support{C}$.
\index{totally disconnected}

\begin{Remark}
\begin{enumerate}
\item  
One may think that the terminology ``arcwise connected'' conforms
with usual mathematics. However, the popular example of connected space
which is not arcwise connected in the usual mathematics turns out to be connected in our
sense. A continuum corresponding to it is the subcontinuum $H$ of $[0,1]^{2}$ 
defined by 
$$
H=(0,1]_{\rational}\times \seti{0} 
\bigcup 
\seti{0}\times (0,1]_{\rational}
\bigcup 
\left(\bigcup \setii{\frac1i}{i\leq N}\right)\times [0,1]_{\rational}
$$
where $N$ is a huge number. Then $(0,y)\approx (\frac1N,y)$ and
hence $(0,y)$ can be connected by a sorites sequence to any other point.

\item
Note that the usual definition of connectedness asserts that there is
a nontrivial disjoint decomposition $X=A_1\bigcup A_2$ with $A_i$ \itwo being
open and closed, which cannot be used since every rigid mesh continuum $C$ is 
totally disconnected with respect to the ``$S$-topology''.
\end{enumerate}
\end{Remark}

\subsection{Topology of Metric Continuum}
\label{subsec:topol-metr-spac}
The distance function gives refined statements on the topology of contina.

A \textit{metric continuum} 
is a triple $(X,d,\approx) $ where $(X,d)$ is a metric class 
and $(X,\approx)$ is a continuum defined by $x\approx y$ 
if and only if $d(x,y)\approx 0$.  
\index{continuum, metric}
A metric continuum $(X,d,\approx)$ is called \textit{metric mesh continuum}
if $X$ is a semiset and \textit{rigid metric mesh continuum} if
$X$ is a set.

\index{visible ball}

\subsubsection{Completeness}

Let $(X,d,\approx)$ be a metric continuum. 
A concrete sequence $a=(\seqinf{a})$ in $X$ \textit{converges to
  $c\in X$} if for each \accessible{} number $k$ there is an \accessible{}  number
$\ell$ such that for every \accessible{}  $i\geq \ell$ 
$$ d(a_i,c)<\frac1k.  $$

A concrete sequence $a=(\seqinf{a})$ is a \textit{Cauchy
  sequence} if for each \accessible{}  number $k$ there is an \accessible{} 
number $\ell$ such that for every \accessible{}  $i,j\geq \ell$ 
$$  d(a_i,a_j)<\frac1k.$$

\begin{Proposition}
  \label{prop:convergence}
  Let $a=(\seqinf{a})$ be a concrete sequence on a rigid metric
  mesh continuum $(X,d,\approx)$. 
  Let $\tilde{a}=(\seq{N}{a})$ be an extension of it to a huge sequence.  
  Then
\begin{enumerate}
  \item For $c\in X$, the \accessible{}  sequence $a$ converges to $c$ 
  if and only if the extended $\tilde{a}$ converges to $c$.
  \item The sequence $a$ is a Cauchy sequence if and only if the extended $\tilde{a}$ is convergent.
  \end{enumerate}
  Hence every \accessible{}  sequence is convergent if and only if it is a Cauchy sequence.
\end{Proposition}
\begin{eproof}
  Let $a$ be a concrete sequence with a huge extension $\tilde{a}=(\seq{N}{a})$. 
  Suppose $\tilde{a}$ converges to $c$. There is a huge $M$ such that 
  $d(a_L,c)\approx0$ for every huge $L\leq M$. 
  Let $k$ be an arbitrary \accessible{}  number. 
  Since the objective condition
  \begin{equation}
    \label{eq:822-1}
    d(a_i,c)<\frac1k
  \end{equation}
  is satisfied by every huge number $i\leq M$, 
  there is an \accessible{}  number $\ell$
  such that every \accessible{}  $i>\ell$ satisfies \eqref{eq:822-1}
  by Theorem \ref{th:overspill-equivalence}. Hence the concrete
  sequence $a$ converges to $c$.

Conversely suppose that the concrete sequence $a$ converges to $c$. 
Let $k$ be an \accessible{}  number. 
There is an \accessible{} $\ell_k$ such that \eqref{eq:822-1} holds 
for every \accessible{} $i\geq \ell_k$, whence there is a huge $M_k$ such that 
\eqref{eq:822-1} holds for $i\in [\ell_k..M_k]$. 
By Proposition~\ref{prop:minimum-huge-number},
There is a huge $M$ satisfying $M\leq M_k$ 
for every \accessible{}  $k$.
If $I\leq M$ is huge, 
then \eqref{eq:822-1} with $i=I$ holds for each \accessible{} $k$ 
since $I\leq M_k$.
Hence $a_I\approx c$, which means that $\tilde{a}$  converges to $c$.

Suppose now that $\tilde{a}$ converges to $c$. Then $a$ converges 
and hence it is a Cauchy sequence by the usual arguments.

Conversely  suppose  that  $a$  is  a  Cauchy  sequence.  Let  $k$  be 
an \accessible{} number. Then there is an \accessible{}  number $\ell_k$ such that
\begin{equation}
  \label{eq:822-3}
  d(a_i,a_j)<\frac1k
\end{equation}
holds for every \accessible{}  $i,j\geq \ell_k$. Hence for every \accessible{}  $p$,
(\ref{eq:822-3}) holds for all $i,j\in [\ell_k..p]$, 
whence there is a huge
$M_k$ such that \eqref{eq:822-3} holds for every $i,j\in [\ell_k..M_k]$.
Let $M$ be a huge number satisfying $M\leq M_k$ for every \accessible{}  $k$.
Then for every huge $I\leq M$, \eqref{eq:822-3} for $i=I, j=M$ 
holds for every \accessible{}  $k$ since $I,M\in [\ell_k..M_k]$, 
whence $a_I\approx a_M$. Hence $\tilde{a}$ converges to $a_M$.
\end{eproof}

A metric class is called \textit{complete} if every concrete Cauchy sequence converges.

By Proposition~\ref{prop:convergence} we have the following.
\begin{Theorem}
\label{thm:completeness-1}
A metric space $(X,d)$ is complete.
\end{Theorem}

Moreover we have the following.
\begin{Theorem}
\label{thm:completeness-2}
Suppose $(X,d,\approx)$ is a metric class and $A\sqsubset X$ is a quasi-subset,
\index{quasi-subset}
namely, every concrete sequence of $A$ can be extended 
to a huge sequence in $A$. 
Then the metric mesh continuum $(A,d,\approx~)$ is complete.
In particular, for a subset $b\subset X$, the metric class $(\overline{b},d)$ 
is complete.
\end{Theorem}
\begin{eproof}
Suppose $A\subset X$ is a quasi-subset. Let $\seqinf{a}$ be a concrete Cauchy 
sequence of $A$. Extend it to a huge sequence $\tilde{a}$ in $A$. 
By Proposition~\ref{prop:convergence}, $\tilde{a}$ converges and 
hence $a$ converges.

Let $b$ be a subset of $X$. It suffices to show that $A:=\overline{b}$ 
is a quasi-subset.  Let $a$ be a concrete sequence in $A$.
Then there is a concrete sequence $x$ in $b$ with $a_i\approx x_i$
for all \accessible{} $i$. Extend $a$ to huge sequence $\tilde{a}=(\seq{N_1}{a})$ 
in $X$ and $x$ to  huge sequence $\tilde{x}=(\seq{N_2}{x})$. Since
for \accessible{}  $i$, we have $d(a_i,x_i)<\frac1i$, there is a huge 
$M$ with $d(a_I,x_I)<\frac1I$ for all $I\leq M$. Hence 
$a_I\approx x_I$ for huge $I\leq M$
and $a_I\in \overline{b}=A$. Namely $(\seq{M}{a})$ 
is an extention of $a$ in $A$.
\end{eproof}

\subsubsection{Compactness}
Let $(X,d,\approx)$ be a metric continuum. 
For $x\in X$ and a positive rational number  $r$,  
define the \textit{$r$-ball with center $x$} by
$$
B_r(x):=\setii{y\in X}{d(x,y)\leq r},
$$
which is a set if $X$ is a set. 
If $r\succ 0$, then $B_r(x)$ is called \textit{a visible ball}. 

We say the continuum is \textit{precompact} if 
for each \accessible{}  $k$, there is an \accessible{} 
number of points $\seti{\seq{\ell}{x}}$ such that
$$X=\bigcup_{1\leq i\leq \ell}B_{\frac1k}(x_i).$$
\index{precompact}

\begin{Proposition}
An rigid mesh metric continuum is compact if and only if it is precompact.
\end{Proposition}
\begin{eproof}
Let $(X,d,\approx)$ be a rigid mesh metric continuum.

Assume $(X,d)$ is precompact. 
Let $K$ be a huge number and $I$ the set of numbers $n$ 
such that there is a subset $Y\subset X$ satisfying 
$\kosuu{Y}\leq K$ and $X=\bigcup_{y\in Y}B_{1/n}(y)$. 
By assumption $I$ contains every \accessible{} number 
and hence a huge number $N$.
Then there is a subset with $\kosuu{Y}\leq K$ such that $X=\bigcup_{y\in Y}B_{1/N}(y)$.
Hence $\overline{Y}=X$. This means that $(X,\approx)$ is compact.

Conversely suppose $(X,\approx)$ is compact. 
Let $n$ be an \accessible{} number. Let $I$ be the set of numbers
$K$ such that there is a subset $Y\subset X$ 
satisfying $\kosuu{Y}\leq K$ and $\bigcup_{y\in Y}B_{1/n}(y)=X$. 
If $K$ is huge then there is a dense subset $Y\subset X$ with
$\kosuu{Y}\leq K$, whence $K\in I$. Thus $I$ contains all huge 
numbers and hence an \accessible{} number $k$.  
Namely $X$ is covered by an \accessible{}  number of balls of radius $\frac1n$.
Hence $(X,d)$ is precompact.
\end{eproof}

\begin{Corollary}
\label{cor:subcontinuum-compactness}
A rigid mesh subcontinuum of a compact rigid metric continuum is compact.
\end{Corollary}
\begin{eproof}
Suppose $C$ is a compact rigid metric continuum and $X\subset \support{C}$.
Since $C$ is precompact, the metric continuum $(X,d,\approx)$ with the 
restricted distance function is precompact and hence is compact.
\end{eproof}

An $x\in X$ is an \textit{accumulation point} 
of a huge sequence $a=(\seq{N}{a})$ in $X$ 
\index{accumulation point of sequences}
if either there is a huge number of $j$ satisfying $a_j=x$
or the \textit{support} $\setii{a_i}{i\in [1..N]}$ is huge and has $x$ 
as its accumulation point.

An element $x\in X$ is an \textit{accumulation point} of 
a concrete sequence $(\seqinf{a})$ 
if for each \accessible{} number $k$ 
there is an \accessible{}  $i\geq k$ with $d(x,a_i)<\frac1k$.

\begin{Proposition}
\label{prop:characterization-of-accumlation-point}
A concrete sequence on 
a rigid mesh metric continuum $(X,d,\approx)$ 
has an accumulation point
if and only if every huge extension of it has an accumulation point.
\end{Proposition}
\begin{eproof}
Let $a=(\seqinf{a})$ be a concrete sequence in $X$. 
If $a$ is contained in a set with \accessible{} number of points, 
then the assertion is obvious.
Otherwise every extended sequence has the huge support.

Suppose $x$ is an accumulation point of the sequence $a$
and $\tilde{a}=(\seq{N}{a})$ is a huge sequence extending it
such that $\kosuu{\setii{i}{a_i=x}}$ is accessible. 
We show that $x$ is an accumulation point of $A:=\setii{a_i}{i\in [1..N]}$. 
For each \accessible{} $k$, the number of elements $B_{\frac1k}(x)\bigcap A$
is huge hence there is a huge $M$ such that
\begin{equation}
\label{eq:20120323}
\kosuu{B_{\frac1k}(x)\bigcap A}\geq M
\end{equation}
for all \accessible{} $k$ by Proposition~\ref{prop:minimum-huge-number}.
Hence there is a huge $K$ such that \eqref{eq:20120323} holds for
$k=K$, which means $x$ is an accumulation point of $A$.

Suppose every huge extension of $a$ has an accumulation point and
suppose $a$ has no accumulation point. Then for every $x\in X$ there is
an \accessible{}  $k_x$ such that for \accessible{}  $i\geq k_x$
$$d(x,a_i)\geq \frac1{k_x}.$$

Put $k=\max_{x\in X}k_x$. 
Then for all $x\in X$ and for \accessible{} $i\geq k$,
$$
d(x,a_i)\geq \frac1k.
$$
In particular, for every \accessible{}  $i,j\geq k$ we have
\begin{equation}
\label{eq:2012-2-13-1}
d(a_i,a_j)\geq \frac1k.
\end{equation}
If $(\seq{N}{a})$ is a huge extension of $a$, 
then there is a huge $M\leq N$ such that (\ref{eq:2012-2-13-1}) holds 
for every  $i,j\in [k..M]$. 
This means the extended sequence $(\seq{M}{a})$
has no accumulation point, a contradiction.
\end{eproof}

Hence by virtue of Theorem \ref{th:existence-of-accumulation-in-compactspace} we
have
\begin{Corollary}
\label{corollary:accumulation_point_in_compact_space}
A concrete sequence in a compact rigid mesh metric continuum 
has an accumulation point.
\end{Corollary}

\newpage
\section{Continua of Binary Words}
\def\runningtitle{Continua of huge binary words}
\label{sec:binary-words}
The concept of continuum makes it possible to construct
continuum directly from syntactic objects. As an illustration
we examine topological properties of four metric continua 
of huge binary words.

\newcommand{\hugeword}[1]{\seti{0,1}^{#1}}
\newcommand{\hugewordii}[1]{\seti{0,1}^{\leq #1}}
Denote by 
$\hugewordii{N}$
the set of words on $0,1$ of length less than or equal to $N$ and
$\hugeword{N}$
the subset consisting of words of length $N$.

The following interpretations with appropriate distance functions 
give four rigid mesh continua with
topolocical properties different from one another.
\begin{enumerate}
\item $\hugewordii{N}$ is the vertex set of binary trees of depth 
$N$ and $\hugeword{N}$ is the set of its leaves. 
\item $\hugeword{N}$ is the vertex set of the $N$-dimensional hypercube.
\item $\hugeword{N}$ is the set of the characteristic functions of subsets in $[0..1]_{\frac1N}$.
\end{enumerate}

\subsection{Binary Trees}
\label{subsec:binary-tree}
\index{binary tree@binary tree}
Consider the symmetric graph $BTree_{N}$ with the vertext set $\hugewordii{N}$ and
the edges are $\seti{w,w0},\seti{w,w1}$ for $w\in \hugewordii{N-1}$. 
Let $d_0(p,q)$ 
be the path distance, namely the length of the shortest path joining $p$ and 
$q$  where every edges are given unit length.
\begin{Lemma}
\label{lemma:path-distance-binary-words}
$$d_0(x,y)=|x|+|y|-2m(x,y),$$
where 
$$m(x,y):=\min\setii{i}{x_i\neq y_i}.$$
\end{Lemma}
\begin{eproof}
Denote by $x_i$ the $i$-th character of the word $x$.
Put $x=ux',y=uy'$ with $x'_1\neq y'_1$ if both $x'$ and $y'$ are not empty word.
Note that this decomposition is unique. 
Then the shortest path joining $x$ and $y$ is composed the
path of length $|x'|$ from $x$ to $u$ and the path of length $|y'|$
joining $u$ to $y$. Hence, noting $|u|=m(x,y)$ we have 
$$ d_{0}(x,y)=|x'|+|y'|=|x|-|u|+|y|-|u|=|x|+|y|-2m(x,y).$$
\end{eproof}

\begin{Lemma}
\label{lemma:ultra-triangular-inequality}
$$m(x,z) \geq \min\seti{m(x,y),m(y,z)}. $$
\end{Lemma}
\begin{eproof}
Suppose $m(x,y)=m(y,z)$. Then 
$$x=ux',y=uy',z=uz'$$
with $x'_1\neq y_{1}'$ if both $x'$ and $y'$ are nonempty words
and $y_1'\neq z_1'$ if both $y'$ and $z'$ are not empty words.
Hence $m(x,z)\geq |u|=m(x,y)=m(y,z)$.

Suppose $m(x,y)\neq m(y,z)$. We may assume $m(x,y)<m(y,z)$. 
Then we can write
$$x=ux',y=uvy',z=uvz'$$
with $x_{1}'\neq v_1$ and $y'_1\neq z'_1$ if $y'$ and $z'$ are nonempty.
Hence $m(x,z)=|u|=m(x,y)\leq \min\seti{m(x,y),m(y,z)}$.
\end{eproof}

From the function $m$ we obtain various ultrametrics.
\begin{Lemma}
If $f$ is a positive descreasing function on positive rationals, 
then the function $d_f$ defined by $d_f(x,x)=0$ and 
$d_f(x,y):=f(m(x,y))$ if $x\neq y$ is an ultrametric.
\end{Lemma}
\begin{eproof}
Obviously $d_f$ is symmetric and reflexive and if $x\neq y$
then $m(x,y)>0$, whence $d_f(x,y)>0$.
The ultrametric triangle relation
$$
   d_{f}(x,z)\leq \max\seti{d_{f}(x,y),d_f(y,z)}.
$$
follows from  Lemma\ref{lemma:ultra-triangular-inequality}.
\end{eproof}

\begin{Lemma}
\label{lemma:ultrametric-totally-disconnected}
Suppose $C$ is a metric continuum with 
an ultrametric $d$. Then there are no sorites sequences.
If $C$ is a rigid mesh continuum, then it is totally disconnected.
\end{Lemma}
\begin{eproof}
Let $(\seq{N}{x})$ be an $\approx$-chain. 
Put $\varepsilon=\max\setii{d(x_i,x_{i+1})}{i\in [1..N-1]}\approx0$. 
Suppose $d(x_1,x_N)> \varepsilon$. Then
$$k:=\min\setii{i}{d(x_1,d_i)> \varepsilon}\leq N.$$ 
Since $d(x_{k-1},x_k)\leq \varepsilon$, we have
$$\varepsilon< d(x_1,x_k)\leq \max\seti{d(x_1,x_{k-1}),d(x_{k-1},x_k)}\leq \varepsilon$$
a contradition. Hence we have $d(x_1,x_N)\leq \varepsilon$ and 
hence $x_1\approx x_{N}$.
\end{eproof}

\subsubsection{Hyperbolic Space}
\index{binary tree@binary tree, hyperbolic distance}
Define on $\hugewordii{\Omega}$
$$d_{hyp}(x,y)=\frac{|x|+|y|-2m(x,y)}{2\Omega}.$$ 
Restrited on $\hugeword{\Omega}$ we have
$$d_{hyp}(x,y)=1-\Frac{m(x,y)}\Omega,$$
which is an ultrametric.

\begin{figure}[htbp]
\begin{center}
\includegraphics[width=0.4\textwidth]{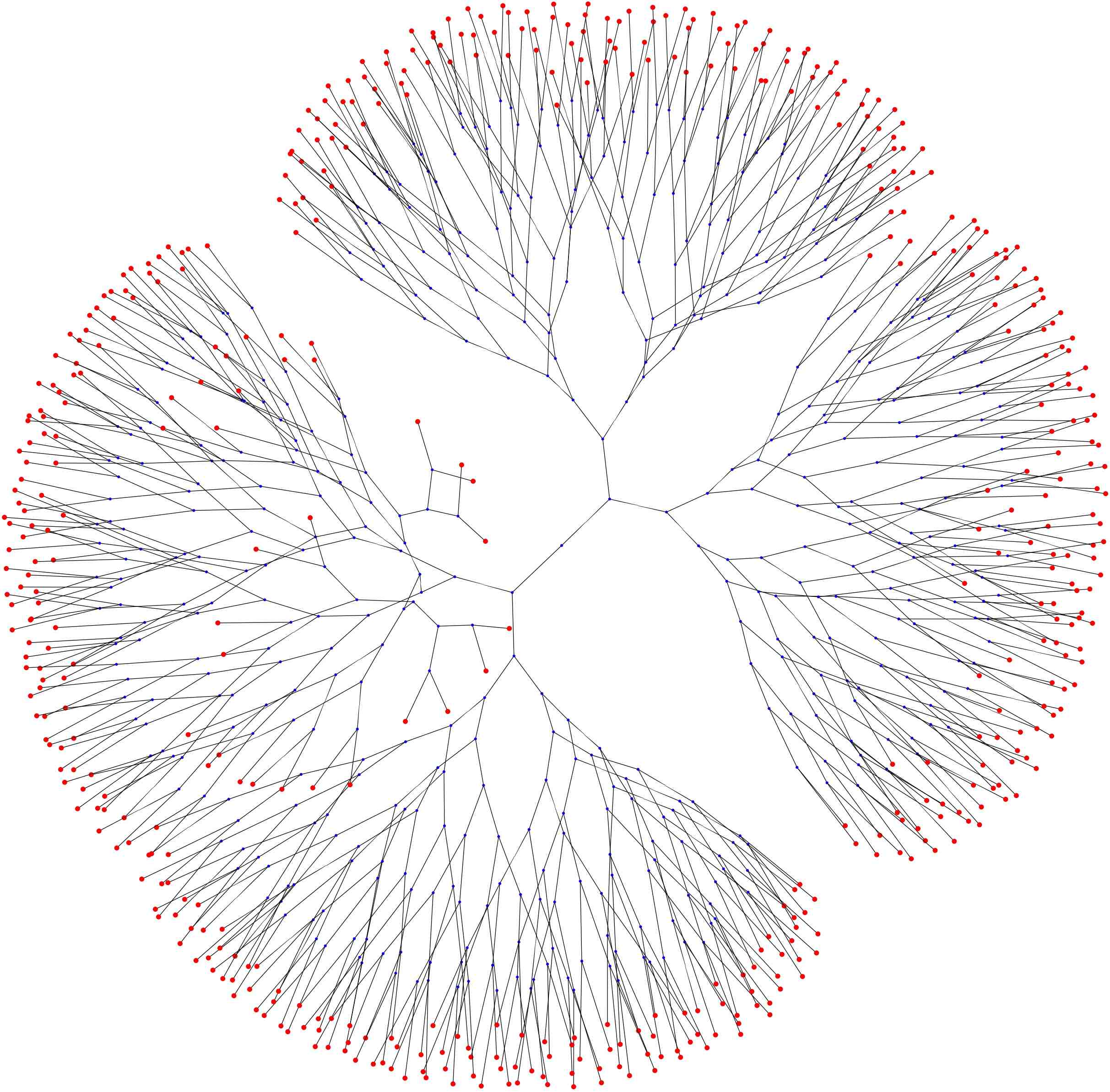}
\includegraphics[width=0.4\textwidth]{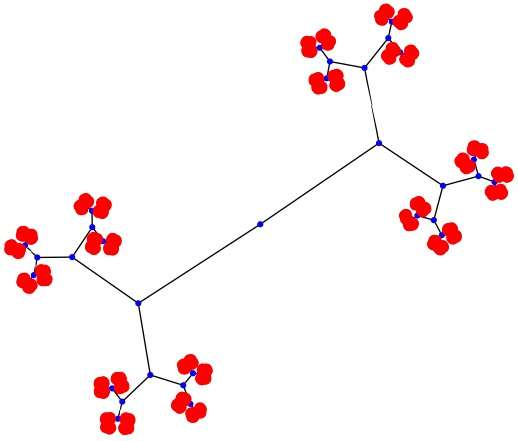}
\caption{Binary tree with $256$ leaves marked by red dots.
In the left graph edges are given uniform length 
whereas in the right the edges connecting the $k$-th level vertices
to its children is given the length $2^{-k}$
The subspace of red dots of the left graph ``approximates'' the hyperbolic space and
that on the right the Cantor space.}
\label{fig:hyperbolic}
\end{center}
\end{figure}

The metric continuum mesh $Hyp_{\Omega}:=(\hugeword{\Omega},d_{hyp},\approx)$ is 
called \textit{the hyperbolic continuum of binary words of length $\Omega$}. 
See the left graph of Fig.~\ref{fig:hyperbolic}.

A metric continuum is \textit{locally compact} if for
every virtual point $x$ there is a rational $r\succ 0$ such that $B_r(x)$ is
compact.
\index{continuum, locally compact metric}
\begin{Proposition}
The hyperbolic continuum mesh $Hyp_{\Omega}$
is perfect but is neither connected nor locally compact.
\end{Proposition}
\begin{eproof}
Since for every $w\in \seti{0,1}^{\Omega}$, 
the ball $B_r(w)$ with $0\neq r\approx0$ is a huge set, 
the continuum $Hyp_{\Omega}$ is perfect.

Since $d_{hyp}$ is ultrametric, the continuum $Hyp_{\Omega}$ is totally
disconnected by Lemma\ref{lemma:ultrametric-totally-disconnected}.

To show $Hyp_{\Omega}$ is not locally compact, 
let $w$ be any word of length $\Omega$.
Let $r\succ 0$ be a rational. 
Let $K$ be the integer part of $r\Omega$ so that $\frac{K}{\Omega}<r$ and
$\frac{K}{\Omega}\approx r$. Decompose as $w=w_{1}w_{2}$ with $|w_2|=K$. 
Let $L$ be the integer part of $\frac{K}2+1$ so that $L\geq \frac{K}2$.
Then 
$$B_r(w)\supset B_{\frac{K}{\Omega}}(w)=\setii{w_1u}{|u|=K}\supset P$$
where 
$$P:=\setii{w_1u_11^{L}}{|u_1|=K-L}.$$
The set $P$ is huge with $2^{K-L}$ elements and has no accumulation points.
In fact if $x_i=w_1u_i1^{L}\in P$ \itwo, then  
$$m(x_1,x_2)\leq |w_1u_1|=\Omega-K+(K-L)=\Omega-L,$$
whence 
$$
d_{hyp}(x_1,x_2)\geq \frac{L}{\Omega}\geq \frac{K}{2\Omega}\approx\frac{r}2.
$$
Hence $B_r(w)$ is not compact for every $r\succ 0$.
\end{eproof}

\subsubsection{Cantor Space}
\index{binary tree@binary tree, Cantor continuum}
\index{Cantor continuum@Cantor continuum}
Let $d_C$ be the distance function on the graph $BTree$ when 
the edge of level $n$ is given the length $2^{-n}$.

\begin{Lemma}
The distance function $d_C$ is given by
$$d_C(x,y):=2^{-m(x,y)+1}-2^{-|x|}-2^{-|y|}$$
for $x,y\in \hugewordii{\Omega}$ and  
if $x,y\in \hugeword{\Omega}$, 
$$d_C(x,y)=2(2^{-m(x,y)}-2^{-\Omega}).$$
Hence $d_C$ on $\hugeword{\Omega}$ is an ultrametric.
\end{Lemma}
\begin{eproof}
The path which connects $x$ to the empty word is
$$\sum_{i=1}^{|x|}2^{-i}=1-2^{-|x|}.$$
The length of the path connecting $x=ux'$ and $y=uy'$ with $|u|=m(x,y)$
is the sum of the length of the paths from $x$ to $u$ and from $u$ to $y$,
whence
\begin{eqnarray*}
d_C(x,y)&=&d_C(x,\lambda)+d_C(y,\lambda)-2d_C(u,\lambda)  \\ 
  &=&
(1-2^{-|x|})+(1-2^{-|y|})-2(1-2^{-m(x,y)})=2^{-m(x,y)+1}-2^{-|x|}-2^{-|y|}.
\end{eqnarray*}
Hence if $|x|=|y|=\Omega$, then
$$
d_C(x,y)=2(2^{-m(x,y)}-2^{-\Omega}).
$$
Hence by Lemma~\ref{lemma:ultra-triangular-inequality}, $d_C$ 
is an ultrametric.
\end{eproof}

The rigid mesh metric continuum $(\seti{0,1}^{\Omega},d_C,\approx)$ is called the 
\textit{Cantor space}.
See the right graph of Fig.~\ref{fig:hyperbolic}.

\begin{Proposition}
\label{prop:Cantor-space}
Let $\Omega$ be a huge number. 
\begin{enumerate}
\item
The rigid mesh continuum 
$(\hugewordii{\Omega},d_{C},\approx)$
is compact but is not perfect nor connected. 
\item The Cantor spacde $(\hugeword{\Omega},d_{C},\approx)$
is compact and perfect but is not connected.
\end{enumerate}
\end{Proposition}
\begin{eproof}
Let $N$ be a huge number 
and take a huge $K< \Omega $ with $2^{K+1}<N$.  
Denote by $Y$ the set of all the words of length $\leq K$. Then
$$\kosuu{Y}=2^{K+1}< N.$$

Every word of length $\leq K$ is in $Y$. Words $w$ 
of length in $[K+1..\Omega]$ is decomposed as
$w=uv$ with $|u|=K$ and   
$$ d(w,u)=2^{-|u|+1}-2^{-|w|}-2^{-|u|}\leq 2^{-|u|}=2^{-K}\approx 0,$$
whence  $w\approx u\in Y$, namely $w\in \overline{Y}$.
As a result $\overline{Y}=\seti{0,1}^{\leq\Omega}$.  
Thus we can make the sizes of dense subsets as small as possible within huge numbers.
Hence $(\hugewordii{\Omega},d_{C},\approx)$ is compact.
It is not perfect nor connected 
since the words of \accessible{}  length are isolated points.

By Corollary~\ref{cor:subcontinuum-compactness}, 
the Cantor space is compact
since it is a rigid mesh subcontinuum of the compact
rigid mesh continuum $(\hugewordii{\Omega},d_{C},\approx)$.

The Cantor space is not connected because $d_{C}$ is an ultrametric.
However it is perfect. In fact, let $x\in \seti{0,1}^{\Omega}$. 
Put $x=x_1x_2$ with $|x_1|=\Omega//2$. Then for any $w$ with $|w|=|x_2|$,
$d(x,x_1w)=2(2^{-|x_1|}-2^{\Omega})\approx 2^{\Omega}\approx0$. 
Hence $x$ is an accumulation.
\end{eproof}

\subsection{Power Set}
\index{binary tree power set@binary tree, power set continuum}

Let $C$ be a rigid mesh continuum.
We can regard $\seti{0,1}^{\support{C}}$ as the power set of 
$\support{C}$ identifying $\chi$ with the subset 
$
\setii{x}{\chi(x)=1}\subset \support{C}.
$
Let $\pow^{+}(|C|)$ be the set of nonempty subset of $|C|$.

Define for nonempty subsets $A,B\subset \support{C}$,
$$d_p(A,B):=\max\seti{\max_{a\in A}d(a,B),\max_{b\in B}d(A,b)},$$
where $d(a,B):=\min\setii{|a-b|}{b\in B}$. Obviously $d_p$ is a metric function.

\begin{Lemma}
$d_p(A,B)\approx 0$ if and only if $\overline{A}=\overline{B}$.
\end{Lemma}
\begin{eproof}
Suppose $d_p(A,B)\approx 0$. Then $d(a,B)\approx 0$ for every
$a\in A$, which means $a\approx b$ for some $b\in B$. Hence $A\subset \overline{B}$.
Similarly $B\subset \overline{A}$, whence $\overline{A}=\overline{B}$.

Suppose $\overline{A}=\overline{B}$. Then for every $a\in A$, there is a $b\in B$
with $d(a,b)\approx 0$. Hence $d(a,B)\approx0$ for every $a\in A$. Hence
$\max_{a\in A}d(a,B)\approx 0$. Similarly $\max_{b\in B}d(A,b)\approx 0$. Hence
$d_p(A,B)\approx 0$.
\end{eproof}

The rigid mesh metric continuum $\pow^{+}(C)
=(\seti{0,1}^{\support{C}}\setminus\seti{\bfzero},d_p,\approx)$ 
is called the \textit{power continuum of $C$}. 
\index{continuum, power}
\begin{Proposition}
\label{prop:power-continuum}
The power continuum $\pow^{+}(C)$ of a rigid mesh continuum $C$ 
is connected if $C$ is connected, perfect if $C$ is perfect and 
compact if $C$ is compact.
\end{Proposition}
\begin{eproof}
Suppose $C$ is connected. Let $A\subset \support{C}$.  
Put $A=\seti{\seq{K}{a}}$. Fix $c\in |C|$ and for each $i\in [1..K]$,
Let $a_i=a_{i0},a_{i1},\cdots,a_{iL_i}=c$ be a sorites sequence connecting
$a_i$ to $c$. Let $L=\max\setii{L_i}{i\in [1..K]}$ and put
$  a_{ip}=c $ for $p>L_i$. Define for $i\in [1..K]$,
$$ A_{\ell}=\setii{a_{ij}}{i\in [1..K],j\in [0..\ell]}. $$
Then $A=A_{0},A_1,\cdots,A_L$ is a sorites sequence.  Define now
$$
B_{p}=\setii{a_{ij}}{i\in [1..K],j\in [p..L]}.
$$
Then 
$B_0=A_L,B_1,\cdots,B_{L}=\seti{c}$ is a sorites sequence. Hence every
subset is connected by a sorites sequence to $\seti{c}$ whence the 
power continuum $\pow^{+}(C)$ is connected.

Suppose now $C$ is perfect and $A\subset \support{C}$ be nonempty and 
$a\in A$. Since $C$ is perfect, there is a huge subset $B$ whose elements
are indistinguishable from $a$. For each $b\in B$, we have $A\approx A_b$
where 
$$A_b:=A\triangle \seti{b}, $$
$\triangle$ denoting the symmetric difference.
Hence there are huge number of subsets indistinguishable
from $A$. Hence the power continuum $\pow^{+}(C)$ is perfect.

Suppose now $C$ is compact.
Let $N$ be a huge number.
Let $M$ be the number satisfying
$$ 2^{M}\leq N< 2^{M+1}.$$
Then $M$ is huge.
Since $X$ is compact, there is a subset $A\subset X$ with $\kosuu{A}\leq M$
and $\overline{A}=X$.
Let $B\subset X$. Then for each $b \in B$, there is an $a_{b}\in A$
with $a_{b}\approx b$. Define
$$\widetilde{B}:=\setii{a_{b}}{b \in B}\subset A,$$
then $B\approx \widetilde{B}$. Hence
$\overline{\pow(A)}=\pow(X)$, with $\kosuu{\pow(A)}=2^{\kosuu{A}}\leq 2^{M}\leq
N$. Since $N$ is an arbitrarily huge number, the power continuum is
compact. 
\end{eproof}

\subsection{Hypercube}
\label{subsec:hypercube}
\index{binary tree hamming@binary tree, hamming distance}
\index{hypercube@hypercube}
As in the previous subsection, we consider
  $\hugeword{\Omega}$ as the power set of
  $X=[0,1]_{\frac1\Omega}\setminus\seti{0}$.  
Let $d_{h}$ be the distance function
  of the graph $Hyper_{\Omega}$ whose nodes are subsets of $X$ and the edges are
  $\seti{A,A\bigcup\seti{b}}$ ($b\nin A$) with length $\frac1\Omega$. 
Obviously we have 
$$  d_{h}(A,B)=\frac{\kosuu{A\triangle B}}{\Omega}.$$
The rigid mesh metric continuum $(\seti{0,1}^{\Omega},d_{h},\approx)$
is called \textit{the hypercube continuum of size $\Omega$}.
\index{continuum, hyperube} 
We remark that this continuum seems
essentially the same as a metric space constructed
in \cite{MR2435141} from finite hypercubes by limiting process.

\begin{figure}[htbp]
\includegraphics[width=0.3\textwidth]{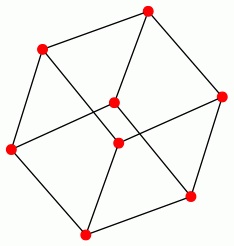}
\includegraphics[width=0.3\textwidth]{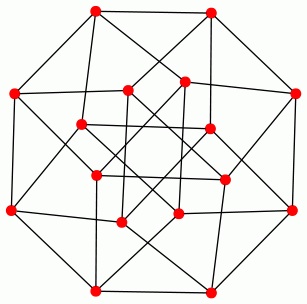}
\includegraphics[width=0.3\textwidth]{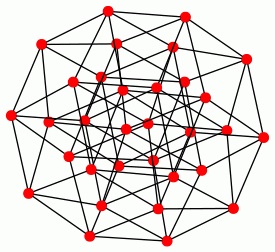}
\\
\includegraphics[width=0.5\textwidth]{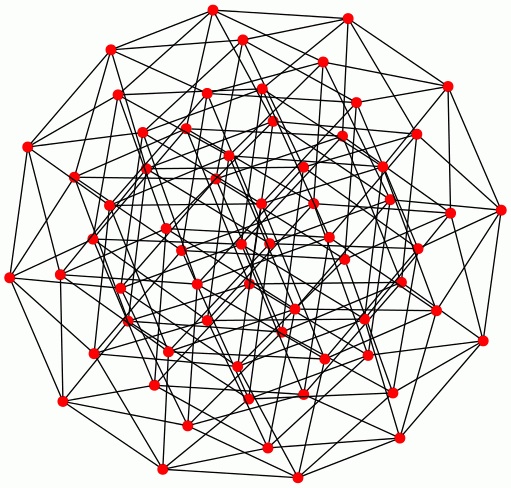}
\includegraphics[width=0.5\textwidth]{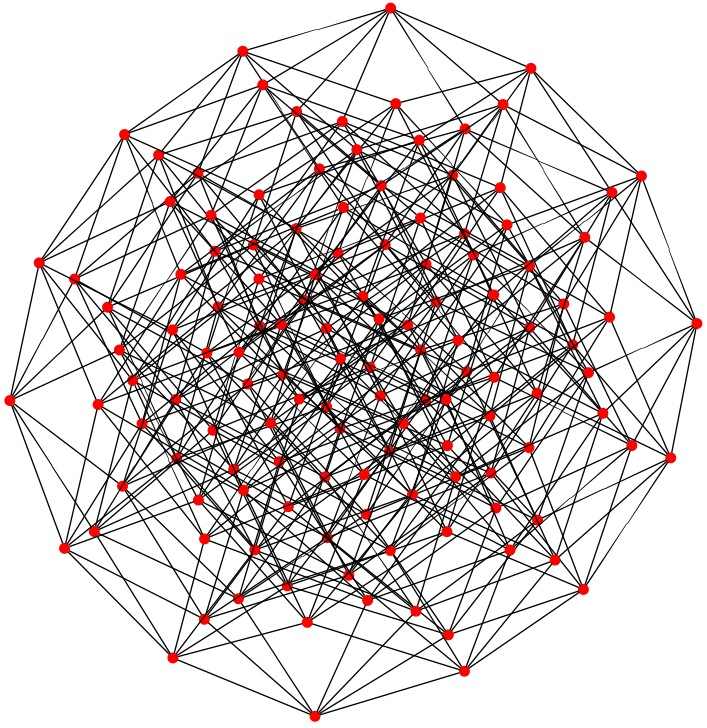}
\caption{Hypercubes on $3,4,5,6,7$ nodes}
\label{fig:hypercube}       
\end{figure}

\begin{Proposition}
The hypercube continuum is connected and perfect but is not locally compact.
\end{Proposition}
\begin{eproof}
By similar arguments as in the proof of Proposition~\ref{prop:power-continuum},
the hypercube continuum is connected and perfect.

However the hypercube continuum is not compact. In fact we show that
there are huge number of words with mutual distance greater than
$\frac12$. Choose a huge $M$ with $2^{M+1}<\Omega\leq 2^{M+2}$.  For
$i\leq M$, Let $A_i$ be the set of integers less than $2^{M+1}$ whose
binary expansion have $1$ on the $i$-th position. Then $\kosuu{A_i}=2^{M-1}$
and $\kosuu{A_i\triangle A_j}=2^{M-2}$ for $i\neq j$. Hence, since
$2^{M+2}\geq\Omega$,
$$
d_h(A_i,A_j)=\Frac{2^{M-1}}{\Omega}\geq\frac18.
$$

Similar arguments show that the hypercube is not locally compact.
\end{eproof}

If we give uniform probability density on $X$, then $A\approx B$ if
and only if $m(A\Delta B)\approx 0$. Hence the hypercube continuum of
binary words is a special case of the continuum of the powerset of
a probability space with this distance function. See \S\ref{sec:measure}.

\begin{Remark}
The metrics $d_p$ and $d_{h}$ on the power set $\seti{0,1}^{\Omega}$
are not comparable. For example, for
  $$A=\setii{\Frac{i}\Omega}{2i\leq \Omega},\quad B=A\bigcup\seti{1},$$
$A\approx_{d_h}B$ but $d_p(A,B)=\frac12$. On the other hand, for
  $$C=\setii{\frac{2i}{\Omega}}{2i \leq \Omega},
\quad D=\setii{\frac{2i+1}{\Omega}}{2i+1 \leq \Omega},$$
$C\approx_{d_p}D$ but $d_h(C,D)=1$ since $C\Delta D=X$.
\end{Remark}

\newpage
\section{Continuum of Morphisms}
\label{sec:continuous-maps}
\def\runningtitle{Continuum of morphisms} 
One might think that this alternative mathematics cannot treat function spaces, 
for which ``infinite sets'' are indispensable. 
However every compact continua are 
represented by rigid mesh continua whose supports are sets 
and any class of morphisms between continua is 
represented by maps between sets, which forms a set.

In this section, we show how to formulate the continuum of 
morphisms between two continua and show an Ascoli-Arzela type
theorem as an illustration showing the usability of our framework 
for ``usual mathematics'' involving infinite sets.

\subsection{Continuum of Functions}
\label{subsec:continuum-functions}
Let $C_i$ \itwo be continua. 
We call two functions
$f_i:\support{C_1}\rightarrow \support{C_2}$ \itwo 
are \textit{indistinguishable}
\index{functions@indistinguishable functions} 
and write $f\approx g$
if $f(x)\approx g(x)$ for all $x\in C_1$.  

If $C_1$ is a mesh continuum and $|C_2|$ is set-like, then we have a continuum 
$$Fun(C_1,C_2)=(Fun(\support{C_1},\support{C_2}),\approx),$$ 
called the \textit{continuum of functions} by Proposition~\ref{prop:continuum-functions}. 

\index{continuum of functions}
\index{$Fun(C_1,C_2)$}

Let $C_i$ \itwo be rigid mesh continua. 
Then the continuum
$
Fun(C_1,C_2)
$ 
is a rigid mesh continuum with $\kosuu{\support{C_2}}^{\kosuu{\support{C_1}}}$ 
virtual points.

The subclass of the continuous functions 
in $Fun(\support{C_1},\support{C_2})$ 
forms a mesh continuum, 
called the \textit{continuum of morphisms from $C_1$ to $C_2$ } 
and written $C(C_1,C_2)$, which is a subcontinuum of $Fun(C_1,C_2)$.
\index{continuum of morphisms}
Note that a point of $C(C_1,C_2)$ is the collection of continuous functions
indistinguishable from a fixed continuous function and hence is exactly 
a morphism from $C_1$ to $C_2$ introduced in \S~\ref{sec:morphisms-continua}.

Even if  $C_i\approx C_i'$ \itwo, the continuum $Fun(C_1,C_2)$
is not necessarily equivalent to $Fun(C_1',C_2')$ but 
the continua of morphisms are equivalent as is seen as follows.

\begin{Proposition} 
Suppose $C_i$ \itwo are mesh continua and
$C'_i$ \itwo are continua. 
Let $g_i:\support{C_i}\rightarrow \support{C_i'}$ \itwo be 
representations of equivalences with almost inverses 
$g_i^{-1}:\support{C_i'}\rightarrow \support{C_i}$ \itwo. 
Define 
$$\alpha:Fun(C_1,C_2)\rightarrow Fun(C_1',C_2')$$
by 
$\alpha(f)=g_2\circ f\circ g_1^{-1}$, and
$$\beta:Fun(C_1',C_2')\rightarrow Fun(C_1,C_2)$$
by 
$$\beta(f')=g_2^{-1}\circ f'\circ g_1.$$
\begin{center}
$$\xymatrix{ C_1 \ar[r]^{f} &
      C_2\ar[d]^{g_2}_{\simeq} \\
      C_1' \ar[u]^{g_1^{-1}}_{\simeq}\ar[r]_{\alpha(f)} & C_2'
   } 
\qquad
    \xymatrix{ C_1 \ar[r]^{\beta{(f')}}\ar[d]_{g_1}^{\simeq} &
      C_2 \\
      C_1' \ar[r]_{f'} & C_2'\ar[u]_{g_2^{-1}}^{\simeq}
   } $$
\end{center}

Then 
$$\alpha:\support{C(C_1,C_2)}\rightarrow \support{C(C_1',C_2')}$$
represents an equivalence with an almost inverse $\beta$.

\end{Proposition}
\begin{eproof}
To show $\alpha$ is continuous, suppose $f\approx f'$.
Then $f(g_1^{-1}(x))\approx f'(g_1^{-1}(x))$ whence
$$\alpha(f)(x)=g_2(f(g_1(x)))\approx g_2(f'(g_1(x)))\approx \alpha(f')(x).$$
Similarly $\beta$ is continuous.

To show that $\beta$  is an almost inverse of $\alpha$
we need the continuity of $f$ and $f'$.
Since $f$ is continuous, so is $\alpha(f)$. Hence 
$$\beta(\alpha(f))\approx f$$
since
$$\beta(\alpha(f))(x)=g_2^{-1}g_2((f(g_1(g_1^{-1}(x)))))
\approx g_2^{-1}(g_2(f(x)))\approx  f(x).$$
Similarly $\alpha(\beta(f'))\approx f'$ if  $f'$ is continuous.
\end{eproof}

The following is the well-known lemma which plays important roles everywhere.
\index{Robinson lemma@Robinson's lemma}
\begin{Lemma}[Robinson]
\label{prop:robinson}
Let $(X,d,\approx)$ be a metric continuum.
Let $(\seq{M}{a})$,$(\seq{M}{b})$ be huge sequences in $X$.
If $a_i\approx b_i$ for \accessible{} $i$, then it holds for 
$i\in [1..K]$ for some huge $K\leq M$.
\end{Lemma}
\begin{eproof}
  Since the objective condition $d(a_i,b_i)<\frac1i$ holds for all
  \accessible{}  $i$, it holds for $i\leq K$ for some huge $K\leq M$. 
For huge $I\leq K$, $d(a_I,b_I)<\frac1I$ implies $a_I\approx b_I$.
\end{eproof}

The continuity can be rephrased by an $\varepsilon-\delta$-like condition.
\begin{Proposition}
\label{prop:rephrase-of-continuity}
  If $C_i=(\support{C_i},d_i,\approx_i)$ \itwo are rigid mesh metric continua,
  then a map $f:\support{C_{1}}\rightarrow \support{C_2}$ is continuous 
  if and only if for every \accessible{} 
  $k$, there is an \accessible{} $\ell$ such that for all $x,y\in \support{C_{1}}$,
  $d(x,y)<\frac1\ell$ implies $d(f(x),f(y))<\frac1k$.
\end{Proposition}
\begin{eproof}
  Suppose $f$ is continuous. Let $k$ be an \accessible{}  number. For every huge number
$N$, the condition $d(x,y)<\frac1N$  implies $x\approx y$ , $f(x)\approx f(y)$ and 
hence $d(f(x),f(y))<\frac1k$.
Since the condition on $i$ that 
\begin{equation}
\label{eq:2012-2-13-2}
\mbox{$d(x,y)<\frac1i$ implies $d(f(x),f(y))<\frac1k$}
\end{equation}
is objective and satisfied by every huge $i$, 
we have an \accessible{}  $i$ for which (\ref{eq:2012-2-13-2}) holds.

  Conversely suppose that for every \accessible{}  $k$, there is an \accessible{} 
  $\ell_k$ such that (\ref{eq:2012-2-13-2}) holds for $i=\ell_k$. Then,
for every \accessible{}  $k$, $x\approx y$ implies $d(f(x),f(y))<\frac1k$
since $d(x,y)<\frac1{\ell_k}$. Hence $f(x)\approx f(y)$.
\end{eproof}

\subsection{Ascoli-Arzela Theorem}
A map $\kappa:[1..N]\rightarrow [1..M]$ 
with huge $N,M$ is called a \textit{scale of approximation}
\index{scale of approximation@scale of approximation}
if it satisfies the condition that 
$i$ is \accessible{} if and only if $\kappa(i)$ is \accessible{} . 

\begin{Lemma} 
\label{lemma:equicontinuity}
Suppose $(X_i,d_i,\approx_i)$ \itwo are rigid mesh metric continua,
Then a map
$f:X_1\rightarrow X_2$ is continuous if 
there is a scale of approximation $\kappa:[1..N_2]\rightarrow [1..N_1]$
and a huge number $K\leq N_2$ such that 
\begin{equation}
\label{eq:1-902}
d(x,y)<\frac1{\kappa(i)} \mbox{ implies } d(f(x),f(y))<\frac1i
\end{equation}
holds for every $i\leq K$ and $x,y\in X_{1}$. 
\end{Lemma}
\begin{eproof}
Suppose there is a scale of approximation $\kappa$ such that
(\ref{eq:1-902}) holds for every $i\leq K$ with a huge $K$. Suppose
$x\approx y$. Let $k$ be an \accessible{} number. 

Since $d(x,y) < \frac1{\kappa(k)}$, 
we have $d(f(x),f(y))<\frac1k$. Hence $f(x)\approx f(y)$.
\end{eproof}

Let $C_i=(X_i,d_i,\approx)$ \itwo be rigid mesh metric continua.
For a scale of approximation $\kappa:[1..N_2]\rightarrow[1..N_1]$, 
denote by $C_{\kappa}(C_1,C_2)$ the set of functions 
$f:\support{C_1}\rightarrow \support{C_{2}}$
satisfying 
$$
\mbox{ $d(x,y)<\frac1{\kappa(i)} $ 
implies $d(f(x),f(y))<\frac1i$ for all $i\in [1..N_2]$ and $x,y\in \support{C_1}$}. 
$$
By Lemma~\ref{lemma:equicontinuity},$C_{\kappa}(C_1,C_2)$ is a subclass of $C(C_1,C_2)$
and in fact is a subset since the condition of membership is objective.

A class $F$ of morphisms from $C_1$ to $C_2$ is called
\textit{equicontinuous}
\index{equicontinuous@equicontinuous}
if there is a scale of approximation $\kappa$ such that
$F\subset C_{\kappa}(C_1,C_2)$.

The following shows that a set of morphisms from $C_1$ to $C_2$ is 
necessarily equicontinuous. 
\begin{Proposition} 
Suppose $C_i=(X_i,d_i,\approx_i)$ \itwo are rigid mesh metric continua.
Then for every subset $F\subset |C(C_1,C_2)|$, there is
a scale of approximation such that $F\subset C_{\kappa}(C_1,C_2)$.
\end{Proposition}

\begin{eproof}
 For each number $i$, define
$$
r_{F}(i):=\min\setii{d(x,y)}%
{
d(f(x),f(y))\geq \frac1i \mbox{ for some $f\in F$} 
}.
$$
If $d(x,y)<r_F(i)$ then $d(f(x),f(y))<\frac1i$ for all $f\in F$.
Define 
$$\kappa(i):=\max\seti{\left[\frac1{r_F(i)}\right]+1,i}.$$
Then $\kappa(i)\geq i$ and $\kappa(i)>\frac1{r_F(i)}$, whence $\frac1{\kappa(i)}<r_F(i)$.
Hence 
\begin{equation}
\label{eq:1-902a}
d(x,y)< \frac1{\kappa(i)} \mbox{ implies } d(f(x),f(y))<\frac1i
\end{equation}
holds for all $f\in F$. 
By Proposition~\ref{prop:rephrase-of-continuity}, 
$\kappa(i)$ is \accessible{} if $i$ is \accessible{} and 
$\kappa$ is a scale of approximation.
Hence $F\subset C_{\kappa}(C_1,C_2)$
\end{eproof}

\index{ascoli-arzela@Ascoli Arzela Theorem}
\begin{Theorem}[Ascoli-Arzela]
\label{theorem:ascoli-arzela}
  If continua $C_i$ \itwo are compact rigid mesh metric continua, 
  then every subcontinuum of $C(C_1,C_2)$ 
  with set support is compact. 
  In particular, if $\kappa$ is a scale of approximation, 
  then $C_{\kappa}(C_1,C_2)$ is compact. 
\end{Theorem}

\begin{eproof}
Let $K$ be a huge number. Select a huge $L$ with $L^{L}\leq K$.  
Since $C_i$ \itwo are compact, there are dense subsets $A_i\subset \support{C_i}$ 
with $\kosuu{A_{i}}\leq L$ \itwo. 
Let $F\subset |C(C_1,C_2)|$ be a subset.

For each $f\in F$ and $a\in A_1$, choose an element $b\in A_{2}$ such that
$f(a)\approx b$ and put $\tilde{f}(a):=b$. 
Then $\tilde{f}\in Fun(A_1,A_2)$.  
Define a map $\alpha: F\rightarrow Fun(A_1,A_2)$ by $\alpha(f):=\tilde{f}$. 

Then $\alpha$ is injective, namely,
$\alpha(f)\approx \alpha(g)$ implies $f\approx g$. 
In fact, suppose $\alpha(f)\approx \alpha(g)$ and $x\in \support{C_1}$. 
Choose $y\in A_1$ such that $x\approx y$. Then 
$$f(x)\approx f(y)\approx \alpha(f)(y)\approx  \alpha(g)(y)\approx g(y) \approx g(x).$$
Hence $f\approx g$.

Let $\beta:\alpha(F)\rightarrow |C(C_1,C_2)|$ be a right inverse of
$\alpha$, namely, $\alpha(\beta(g))=g$. Then $Im(\beta)$ is dense in $F$.
In fact, for each $f\in F$, $\alpha(\beta(\alpha(f)))=\alpha(f)$
implies $f\approx \beta(\alpha(f))$. 

Since 
$$\kosuu{Im(\beta)}\leq \kosuu{\alpha(F)} \leq \kosuu{Fun(A_1,A_2)}\leq L^L<K,$$
we have a dense subset of $F$ with the number of elements less than $K$.
Hence $F$ is compact. The latter assertion is the special case of the former
since $C_{\kappa}(C_1,C_2)$ is a set.
\end{eproof}

\begin{Corollary}[Ascoli-Arzela]
An equicontinuous concrete sequence of morphisms between compact rigid mesh metric
continua has an accumulation point.
\end{Corollary}
\begin{eproof}
Suppose rigid mesh metric continua $C_i$ \itwo are compact.
Let $f=(\seqinf{f})$ be an equicontinuous concrete sequence of morphisms
from $C_1$ to $C_2$. Then $f$ is a concrete sequence
in the continuum $C_{\kappa}(C_1,C_2)$ for some scale $\kappa$ of approximation,
which is compact by Theorem~\ref{theorem:ascoli-arzela},
whence has an accumulation point by Corollary~\ref{corollary:accumulation_point_in_compact_space}.
\end{eproof}
\index{ascoli-arzela@Ascoli Arzela Theorem}

\newpage
\section{Real Numbers}
\label{sec:real-numbers}
\def\runningtitle{Real Numbers} 
\subsection{Real Numbers}
A point of the continuum $\real=(\rational,\approx)$ is called 
a \textit{real number}, namely a real number is a class 
$$[r]:=\setii{s\in \rational}{s\approx r}$$
for some rational number $r\in \rational$. 
A real number $a$ is said to be represented by a rational number $r$
if $a=[r]$. 
\index{real number}

By Proposition~\ref{example:interval-equivalence}, we have the following.
\begin{Lemma}
If $\varepsilon>0$ is an infinitesimal, we have
$r\approx [\frac{r}{\varepsilon}]\varepsilon$ for every $r\in \rational$. 
In particular, every real number is represented 
by a rational of the form $m\varepsilon$ with $m\in \integer$. 
\end{Lemma}
A representation of a real number by a rational in $\varepsilon\integer$
is called \textit{$\varepsilon$-separate}.
\index{real number,$\varepsilon$-separate representation}

Since each real number is a proper class, 
the equality of real numbers is not a definite condition
and the collection of real numbers do not form a class.
In other words, the ``a real number'' should not be regarded 
as a definite object. 
They have uneliminable indefiniteness indicated by the sorites paradox
that $x_1=x_2=\cdots=x_N$ but $x_1\neq x_N$ if the equality of real numbers
had definite meaning.

We defined the relations $r\approx s$, $r\prec s$ and $r\preceq s$ 
for rationals $r,s$ in \S\ref{subsub:infinitesimal-rational}, 
which induce relations of real numbers $p=q,p<q$ and $p\leq q$ respectively 
owing to Proposition~\ref{prop:daishou}.
A real number $p$ is called positive and negative
respectively when $p>0$ and $p<0$.

The absolute value $|p|$ of a real number $p=[r]$ is defined by 
$|p|:=[|r|]$. 

The following is obvious but shows that there are no nonzero infinitesimal reals.
\begin{Proposition}
\label{prop:no-infinitesimal-reals}
If a real number $p$ satisfies 
$|p|\leq \frac1k$ for every \accessible{} number $k$,then $p=0$.
\end{Proposition}

We also defined the notion of finiteness of rationals 
in \S\ref{subsub:accessibility}, which induces finiteness of real numbers.
A real number is called \textit{commensurable} 
\index{real number,commensurable}
\index{real number, finite}
if it is represented by an \accessible{} rational number. 
Since for accessible rational numbers $r,s$,
the indistinguishability implies equality, every commensurable real number
is represented by a unique \accessible{} rational number.
We identify each \accessible{} rational number with
the commensurable real number represented by it. 
For example, the \accessible{} rational number $\frac12$ denotes 
also the commensurable real number $[\frac12]$.

In the following, we define following operations and functions of real numbers.
\begin{enumerate}
\item Addition and multiplication of finite real numbers,
\item For \accessible{} number $n$, the $n$-power of finite real numbers,
and $n$-th root of finite non-negative real numbers.
\item Exponentiation of finite real numbers and logarithm of positive real numbers,
\item Power of finite positive real numbers to finite real numbers.
\end{enumerate}

\subsection{Arithmetic Operations}
\label{sec:algebr-oper-real}
The arithmetic operations on rational numbers induce those on
real continuum.
\begin{Lemma}
\label{Lemma:arithm-oper-real}
Suppose $r,s,r_i,s_i\in \rational$ satisfy $r\approx s$ and $r_i\approx s_i$ \itwo. 
\begin{enumerate}
\def\labelenumi{(\theenumi)}
\item $r_1+r_2\approx s_1+ s_2$.

\item If $r_i$ \itwo are finite then $s_i$ \itwo are finite and
$$ r_1r_2\approx s_1s_2.$$
\item If $r$ is not infinitesimal, then $s$ is neither infinitesimal and
both $\frac1r$ and $\frac1s$ are finite and satisfies
$$ \frac1{r}\approx \frac1{s}. $$
\item\label{item:3} 
If $r$ is finite and $n$ is an \accessible{} number then
$$  r^n\approx s^n.$$
\end{enumerate}
\end{Lemma}
\begin{eproof}
Suppose $r_i$ \itwo are finite. Since a rational indistinguishable from
a finite rational is finite, $s_i$ \itwo are finite. 
From
$$
|r_1r_2-s_1s_2|\leq |r_1||r_2-s_2|+|s_2||r_1-r_2|\leq k(|r_2-s_2|+|r_1-r_2|)
$$
where $k=\max\seti{|r_1|,|s_2|}<\infty$ it follows $r_1r_2\approx s_1s_2$. 

Let $inv:\rational\setminus\seti{0}\rightarrow \rational$ be the function
$inv(r)=\frac1r$.
If $r\not\approx 0$, then $s\not\approx 0$ hence there is an \accessible{}
 $k$ such that $|r|,|s|>\frac1k$. Hence 
$$
\left|\frac1{r}-\frac1{s}\right|=\frac{|r-s|}{|rs|} < k|r-s|\approx 0.
$$
On the other hand the condition $|r|,|s|<\infty$ implies that
$|\frac1r|,|\frac1s|\not\approx 0$.
Hence $inv$ represents a morphpism from $(-\infty,0)_{\rational}\bigcup(0,\infty)_{\rational}$ to itselft.

Put $K=n(\max\seti{|r|,|s|})^n$, then $K$ is finite and 
$$  |r^n-s^n|\leq K|r-s|,$$
hence $r^n\approx s^n$.
\end{eproof}

Hence we have the following morphisms
\begin{Theorem}
\label{th:arithm-oper-real}
\begin{enumerate}
\def\labelenumi{(\theenumi)}
\item The addition defines a morphism 
$$+:\real^2 \rightarrow \real.$$
\item The multiplication defines a morphism 
$$\times:(-\infty,\infty)^2\rightarrow (-\infty,\infty).$$
\item The inverse defines an equivalence 
$$
(-)^{-1}:(-\infty,0)\bigcup (0,\infty) \rightarrow (-\infty,0)\bigcup (0,\infty).
$$
\item The powers $r\mapsto r^{n}$ defines 
morphisms 
$$pow_n:\virtualline\rightarrow \virtualline$$ 
for each \accessible{} $n$.
\end{enumerate}
\end{Theorem}

We express this symbolically by the following point wise ``definition'' 
on real numbers. 
\begin{Definition}
Let $p=[r]$ and $q=[s]$.
\begin{enumerate}
\item $p+q:=[r+s]$,
\item $pq:=[rs]$, when $p,q$ are finite, 
\item $\frac1p=[\frac1r]$, when $p\neq 0$, 
\item $\frac{p}{q}:=p\frac1q$, when $q\neq 0$,
\item $p^n:=[r^n]$, when $p$ is finite and $n$ is \accessible.
\end{enumerate}
\end{Definition}

It should be noted that since the real numbers are vague
objects without definite identity, precise meaning of this definition
is given by the above Theorem~\ref{th:arithm-oper-real}.

The usual axiom of field is satisfied by these operations. 
Let $0$ and $1$ denotes the commensurable real numbers represented
by the rational $0$ and $1$ respectively. Define $-[r]:=[-r]$.
\begin{Proposition}
\label{prop:field-axiom-reals}
Let $p,q,r$ be finite real numbers. Then
\begin{enumerate}
\item The addition and multiplication are associative and commutative.
\item $0+p=p$,$1\times p=p$,
\item $p+(-p)=0$,
\item if $p\neq 0$ then $p\times\frac1p=1$,
\item $p\times(q+r)=p\times q+p\times r$.
\end{enumerate}
\end{Proposition}

This implies for example the following.
\begin{Lemma}
\label{lemma:indistinguishbility-by-quotient}
If $a,b$ are finite rationals and $a$ is not infinitesimal,
then $a\approx b$ if and only if $\frac{a}{b}\approx 1$.
\end{Lemma}
\begin{eproof}
Put $p=[a]$ and $q=[b]$. Then the statement means
$p=q$ if and only if $\frac{p}{q}=1$, which follows from
Proposition~\ref{prop:field-axiom-reals}.
\end{eproof}

Let $\varepsilon>0$ be an infinitesimal. Then the above operations can
be realized by those on the $\varepsilon$-separate representations.
These representations of operations of reals are considered
fundamental in the computational treatments of real numbers. See
\cite{Reveilles-Richard1996,chollet2009insight} for example.

\begin{Proposition}
\begin{enumerate}
\def\labelenumi{(\theenumi)}
\item The addition defines a morphism 
$$+:(\varepsilon\integer,\approx)^2 \rightarrow (\varepsilon\integer,\approx).$$
\item The multiplication is represented by the morphism 
$$mult_\varepsilon:((-\infty,\infty)_\varepsilon,\approx)^2\rightarrow ((-\infty,\infty)_\varepsilon,\approx)$$
defined by 
$$ mult_\varepsilon(n\varepsilon,m\varepsilon)
:=[mn\varepsilon]\varepsilon\quad m,n \in \integer$$
\item The inverse is represented by the morphism
$$
inv_\varepsilon:((-\infty,0)_\varepsilon\bigcup (0,\infty)_\varepsilon,\approx) 
\rightarrow ((-\infty,0)_\varepsilon\bigcup (0,\infty)_\varepsilon,\approx)
$$
defined by 
$$
inv_\varepsilon(n\varepsilon)=\left[\frac1{n\varepsilon^2}\right]\varepsilon
$$
\item If $n$ is an \accessible{} number, the the 
power morphism $r\mapsto r^{n}$ is represented by the morphism
$$ pow_{n,\varepsilon}: 
(\virtualline_\varepsilon,\approx)\rightarrow (\virtualline_\varepsilon,\approx)$$
defined by 
$$pow_{n,\varepsilon}(k\varepsilon):=[k^n\varepsilon^{n-1}]\varepsilon.$$
\end{enumerate}

\end{Proposition}

As for the root operation, we can define it only through representations.
\begin{Lemma}
\label{lemma:k-th-root}
Let $k$ be an \accessible{}  number and $\varepsilon$  a positive infinitesimal.
Then for each finite positive rational number $x$, 
$$ root_{k,\varepsilon}(x)^{k}\approx x,  $$
where 
\begin{equation}
\label{eq:20120225-1}
root_{k,\varepsilon}(x):=\varepsilon\max\setii{m\in\nat}{(m\varepsilon)^{k}\leq x}.
\end{equation}
Moreover $root_{k,\varepsilon}(x)=0$ if and only if $x=0$.
\end{Lemma}

\begin{eproof}
Put $m:=\max\setii{m\in\nat}{(m\varepsilon)^{k}\leq x}$. Then
$$(m\varepsilon)^k\leq x< ((m+1)\varepsilon)^k.$$
Since $m\varepsilon\approx (m+1)\varepsilon$, 
we have $(m\varepsilon)^k\approx ((m+1)\varepsilon)^k$ 
by the accessibility of $k$ 
and hence $x\approx (m\varepsilon)^k$. 
\end{eproof}

\begin{Lemma}
\label{lemma:continuity-of-root}
If finite nonnegative rationals $u,v$ satisfy $u^k\approx v^k$ for an \accessible{} $k$,
then $u\approx v$. In particular, the function $root_{k,\varepsilon}$ is continuous for
infinitesimal $\varepsilon>0$.
\end{Lemma}
\begin{eproof}
Suppose $u\not\approx v$ but $u^k\approx v^k$.
We may assume $u\prec v$. 
Then $u+\frac1k<v$ for some \accessible{} $k$.
We may assume $0\prec u$.
Then $\frac{nu^{k-1}}k\succ 0$, whence
$$u^n \prec u^n+ \frac{nu^{k-1}}{k}\leq (u+\frac1k)^n <v^n$$
which contradicts $u^n\approx x\approx y\approx v^n$.
Hence $u\approx v$.

Suppose rational numbers $a,b$ satisfy $a\approx b$.
Then 
$$
(root_{k,\varepsilon}(a))^k\approx a\approx b\approx (root_{k,\varepsilon}(b))^k
$$
and the above conclusion implies 
$
root_{k,\varepsilon}(a)\approx root_{k,\varepsilon}(b).
$
\end{eproof}

Hence we have proved
\begin{Theorem}
If $k$ is \accessible{} 
then the power operator 
$$pow_{k}:[0,\infty)\rightarrow[0,\infty)$$ 
represents an equivalence. For each nonzero infinitesimal $\varepsilon$, the function
$$x\mapsto root_{k,\varepsilon}([x/\varepsilon]\varepsilon)$$ 
is an almost inverse of the power operator $pow_k$. 
\end{Theorem}

From this we define the $k$-th root $p^{\frac1k}$
of a finite nonnegative real number $p=[r]$ by 
$$p^{\frac1k}:=[root_{k,\varepsilon}(r)],$$
where $\varepsilon>0$ is an infinitesimal. 
The above theorem shows that this does not depend on the choice of
$r$ and $\varepsilon$ and
the following holds:
$$(p^{\frac1k})^k=p,
(p^k)^{\frac1k}=p.
$$

If $s=\frac{\ell}{k}$ is \accessible{}, 
namely, $\ell,k$ are \accessible{}  numbers, 
we define for a finite positive real number $p$
$$p^{\frac{\ell}{k}}:=(p^{\ell})^{\frac1k}.$$
\begin{Remark}
The exponentiation $x^{y}$ for general finite $x,y$ will 
be defined as $\exp(y \log x)$ after defining the exponentiation
$\exp(x)$ and the logarithm function $\log$ as the inverse of $\exp$.
\end{Remark}

\subsection{Sequence}
\label{subsec:sequ-real-numb}
In \S\ref{subsec:topol-metr-spac} we defined convergence of concrete
sequences and Cauchy sequences on on metric spaces. 

Two concrete sequence of rational numbers $(\seqinf{a})$ and
$(\seqinf{b})$ are indistinguishable if $a_i\approx b_i$ for
all $i$. A \textit{concrete sequence of real numbers} is 
the collection of the concrete sequence of rational numbers 
indistinguishable with one such $(\seqinf{a})$. This is not
a class but we use the symbol $[a]=([a_1],[a_2],\cdots)$
to denote this collection.

If $p=(\seqinf{p})$ is a concrete sequence of real numbers,
then a concrete sequence of rational numbers $a=(\seqinf{a})$
is said to represent $p$ if $a_i\in p_i$ for all $i$.
Note that we cannot form a representation by arbitrarily
choosing elements of each $p_i$.

A concrete sequence of real numbers $(\seqinf{p})$
converges to a real number $q$ if 
for each \accessible{}  number $k$ 
there is an \accessible{} number $\ell$ 
such that for every \accessible{} $i\geq \ell$ we have
$$ |p_i-q|<\frac1k. $$
By Proposition \ref{prop:no-infinitesimal-reals}, such $q$ is
uniquely determined and is called the limit of the sequence $p$ 
and is denoted by $\lim_{i\rightarrow \infty}p_i$.

We say that a concrete sequence of real numbers $p=(\seqinf{p})$ 
is a \textit{Cauchy sequence} 
if for every \accessible{}  $k$, there is an \accessible{}  $\ell$ 
such that for every \accessible{}  $i,j\geq \ell$ 
we have
$$  |p_i-p_j|<\frac1k.$$
This means that $p$ is represented by a Cauchy concrete sequence of
rational numbers. 

By Proposition~\ref{prop:convergence}, a concrete sequence of rational
numbers $a$ converges if and only if it is a Cauchy sequence
whence we have the following ``completeness'' of the metric continuum $\real$.
\begin{Theorem}
\label{th:completeness-of-real-numb}
Concrete Cauchy sequences of real numbers converge.
\end{Theorem}

We have also
\begin{Theorem}
\label{th:convergence-of-increasing-sequence}
An increasing concrete sequence of real numbers bounded from above converges.
\end{Theorem}
\begin{eproof}
Let $p=(\seqinf{p})$ be a concrete sequence of real numbers 
such that $p_i\leq p_{i+1}$ for all $i$ and, 
for some some \accessible{} number $k$,  $p_i\leq k$ for all $i$.

Let $\seqinf{a}$ be a concrete sequence of rational numbers representing $p$.
Then $a_i\preceq a_{i+1}$ and $a_i\preceq k$ for all $i$.

Since the continuum $[a_1,k]$ is compact by
Proposition\ref{prop:compactness-of-[a,b]}, 
the sequence has an accumulation point $c$. 
Hence for every \accessible{}  $\ell$, 
the numbers $i\geq \ell$ satisfying
\begin{equation}
\label{eq:9-831}
|a_i-c|<\frac1\ell.
\end{equation}
is not finite. Let $i_0$ be one such number. 
If there is an \accessible{}  $j>i_0$ with $c+\frac1\ell\leq a_j$,
then $j<m$ implies 
$$c+\frac1{2\ell}\prec c+\frac1\ell\leq a_j\preceq a_m$$
and hence the number of $i$ satisfying ~\eqref{eq:9-831}
with $\ell$ replaced by $2\ell$
is less than or equal to $j$, a contradiction. 
Hence $i_0<i$ implies ~\eqref{eq:9-831}, 
which means that $(\seqinf{a})$ converges to $c$.
\end{eproof}

\subsection{Series}
The addition of rationals can be extended to a function

$$ \rational^{N}\ni (\seq{N}{a})\mapsto \sum_{i\in [1..N]}a_i\in \rational.$$

However this does not define a morphism $\real^N\rightarrow \real$ since
generally $\sum_{i=1}^{N}a_i\not\approx \sum_{i=1}^{N}b_i$ even if
 $a_i\approx b_i$ ($i\in [1..N]$).

We call that the sum \textit{$\sum_{i=1}^{N}a_i$ converges} 
\index{convergence series@convergence of series}
if the following holds.
\begin{equation}
\label{eq:20120226-3}
\mbox{$\sum_{i=I}^{N}a_i\approx 0$ for every huge $I\leq N$}.
\end{equation}

\index{convergence series@convergence of series, absolute}
Similarly the \textit{huge sum $\sum_{i=1}^{N}a_i$ converges absolutely} 
if $\sum_{i=I}^{N}|a_i|\approx 0$ for every huge $I\leq N$.

Note that if we define 
$$  S_{k}:=\sum_{i=1}^{k}a_i, \quad k\in [1..N],$$
the condition \eqref{eq:20120226-3} is equivalent to the convergence of 
the sequence $(S_1,\cdots,S_{N})$.

The following can be easily proved.
\begin{Lemma}
\label{lemma:convergense-series-real-numbers}
  If $a_i\approx b_i$ for $i\in [1..K]$ and $\sum_{i\in[1..K]}a_i$
  converges then for some huge $L\leq K$, $\sum_{i\in [1..L]}b_i$
  converges and their limits coincide up to indistinguishability.
\end{Lemma}
\begin{eproof}
Put $\varepsilon=\max\setii{|a_i-b_i|}{i\in [1..K]}\approx 0$.
Since $n\varepsilon\approx0$ for every \accessible{} $n$, we 
can choose a huge $L$ such that $L\varepsilon \approx 0$. 
Then for huge $I\leq L$,
\begin{eqnarray*}
\sum_{i\in [I..L]}|b_i| & \leq & 
\sum_{i\in [I..L]}|b_i-a_i| + \sum_{i\in [I..L]}|a_i|
\\
&\leq& \varepsilon L+\sum_{i\in [I..L]}|a_i|\approx0.
\end{eqnarray*}
hence the sum $\sum_{1\leq i\leq L}b_i$ converges.
Moreover
$$
\left|\sum_{i\in [1..L]}a_i-\sum_{i\in [1..L]}b_i\right|
\leq \sum_{i\in [1..L]}|a_i-b_i|\leq \varepsilon L\approx 0.
$$
\end{eproof}

A point of $\real^N$ is called a sequence of real numbers
and is denoted by $p=(\seq{N}{p})$. It is represented
by a sequence or rational number $a=(\seq{N}{a})\in \rational^{N}$.
We say that $p$ converges if its representation converges and
define the sum $\sum_{i}p_i=[\sum_ia_i]$.
By Lemma~\ref{lemma:convergense-series-real-numbers}, the 
condition of convergence and the value of the sum are 
independent of the choice of representations.

\newpage

\section{Real Functions on Continua} 
\label{sec:real-valued-functions}
\def\runningtitle{real functions on continua}
\subsection{Real Functions}
\label{sec:real-functions}

Let $C$ be a continuum. A morphisms from $C$ to $\real$ is called a real function
on $C$. Recall it is a formal symbol $[f]$ where $f$ is a rational valued continuous
function on $\support{C}$. See \S~\ref{sec:morphisms-continua}.

The \textit{value of a real function} $F$ \textit{at a point} $p$ of 
$C$ is defined to be the real number $[f(t)]$ for $f\in F$ and $t\in p$.
This does not depend on the choice of representations.

We saw in \S\ref{subsec:continuum-functions} that
if $C$ is a mesh continuum, 
the continuous rational valued  functions form a subcontinuum
$$C(C,\real)\subset Fun(\support{C},\rational)$$
and the indistinguishablity condition is definite. 
Hence in this case, the symbol $[f]$ can be interpreted by 
the class $\setii{g\in C(C,\real)}{g\approx f}$ and the 
above definition of the symbol of real function 
conforms to this interpretation.

Suppose $\alpha:C_1\rightarrow C_2$ is a morphism 
between continua $C_i$ \itwo. 

If $F$ is a real function on $C_2$ represented by $f$,
then the real function $[f\circ\alpha]$ 
does not depend on $f$
since $f\approx g$ implies $f\circ\alpha\approx g\circ\alpha$.
The real function $[f\circ\alpha]$ on $C_1$ 
is called the \textit{pull back of $F$ by $\alpha$}
and denoted by $F\circ \alpha$.
Note that if $\alpha\approx\alpha'$, then
$F\circ\alpha=F\circ\alpha'$ since $f\circ\alpha\approx f\circ\alpha'$.

\textit{A representation of a real function} $F$ on $C$ is defined to be 
a pair $(f,\alpha)$, where $\alpha:C\rightarrow C_1$ is an
equivalence of continua and $f$ is a rational valued continuous function 
on $C_1$ such that $f\circ\alpha\in F$.
\index{real function, representation}
Obviously we have the following.
\begin{Proposition}
\label{lem:existence-of-representation}
Suppose $\alpha:C\rightarrow C_1$ is an equivalence of continua.
The assignment $F\mapsto F\circ\alpha$ 
defines a one-to-one correspondence from the collection of real functions on $C_1$ 
onto those on $C$. In particular, 
every real function is represented as $(f,\alpha)$ for some 
rational valued continuous function $f$ on $C_1$. 

\end{Proposition}

Let $D$ be a subcontinuum of the linear continuum $\real$. 
A continuous rational valued  function $f$ on $C$ is called \textit{$D$-valued} 
if  $f(x)\in D$ for all $x\in C$. 
\begin{Proposition}
If $C$ is a mesh continuum, the condition of being $D$-valued is definite.
In particular, the $D$-valued continuous rational valued  functions on $C$
forms a subcontinuum $C(C,D)\subset C(C,\real)$.
\end{Proposition}
\begin{eproof} Let $f$ be a real function on $C$. 
Let $\tilde{f}:b\rightarrow \rational$ be an extension of it. 
Then $f$ is $D$-valued if and only if it satisfies the bounded condition
$$\exists b'\subset b\; \forall x\in b' \;
[ \mbox{ $x\in \support{C}$ implies $\tilde{f}(x)\in D$} \;]. $$
\end{eproof}
Note that $D$-valuedness is not objective condition in general but 
if $D$ is an objective subclass and $C$ is rigid then it is objective.

Let $I$ be an interval symbol 
defined in \S~\ref{subsec:examples-of-continuum}
and $C$ is a mesh continuum 
then $I$-valued continuous rational valued  function on $C$ defines a subcontinuum 
denoted by $C(C,I)\subset C(C,\real)$. 
A continuous rational valued  function $f$ is called \textit{finite} 
if $f$ is $\virtualline$-valued.
\index{real function@finite real function}
\index{real function@$D$-valued real function}

Let $C$ be a continuum. A real function $F$ on $C$ is called 
\emph{$D$-valued} if it is represented by a $D$-valued continuous rational valued  function.
In particular $F$ is called finite if it is represented by a finite continuous
rational valued  function.

Let $C$ be a mesh continuum.
Let $D_{i}$ \itwo be subcontinua of $\real$ and 
$\beta:D_1\rightarrow D_2$ be a quasi-identity 
in the sense explained in \S~\ref{sec:equivalence-continua}. 
If $f$ is a $D_1$-valued continuous rational valued function, 
then $\beta\circ f$ is $D_2$-valued and
if $f_1\approx f_2$, then $\beta\circ f_1\approx \beta\circ f_2$
whence $\beta$ induces a morphism
$$\beta_{*}:C(C,D_1)\rightarrow C(C,D_2)$$
which is an equivalence 
since $\gamma_{*}$ is an almost inverse whenever
$\gamma$ is an almost inverse of $\beta$.

Note that a real function $F$ is $D_{1}$-valued if and only if
$D_2$-valued, since if $f\in F$ is $D_1$-valued then $f\approx
\beta\circ f$ is $D_{2}$-valued and hence $F$ is also represented by
$D_2$-valued function.

For example if $\varepsilon>0$ is infinitesimal, 
the inclusion function
$$\imath_\varepsilon:\varepsilon\integer\rightarrow \rational$$
defines an equivalence
$$
\imath_{\varepsilon*}:C(C,\varepsilon\integer)\rightarrow C(C,\real)
$$
with the almost inverse given by $\kappa_{\varepsilon*}$
where $\kappa_\varepsilon:\rational \rightarrow \varepsilon\integer$ 
is the almost inverse defined in Proposition~\ref{prop:r-int-equiv-rat}.

Similarly, for every interval symbol $I$ and an infinitesimal $\varepsilon>0$, 
we have an equivalence 
$$
\imath_{\varepsilon*}:C(C,I_\varepsilon)\rightarrow C(C,I).
$$

Note that even if $C$ is not a mesh continuum, every real function $F=[f]$ 
on $C$ is represented by an $\varepsilon\integer$-valued continuous function
such as $\kappa_{\varepsilon}\circ f$.
Similarly every $I$-valued real function $F=[f]$ on $C$ 
is represented by an $I_\varepsilon$-valued function 
$\kappa_\varepsilon\circ f$.

\begin{Proposition}
Suppose a continuum $C$ has a dense subcontinuum $M$ with an almost
inverse $\kappa$ for the inclusion morphism $\imath:M\rightarrow
C$. Then the assignments $\imath^{*}:F\mapsto F\circ\imath$ 
and $\kappa^{*}:G\mapsto G\circ\kappa$ are inverse to one another and 
defines a one-to-one correspondence between the real functions on $C$
and those on $M$.
\end{Proposition}
\begin{eproof}
Since $\kappa\circ\imath=id_M$ and $\imath\circ\kappa\approx id_{C}$,
$F\circ\imath\circ\kappa=F$ and $G\circ\kappa\circ\imath=G$.
\end{eproof}

Thus if $C$ and $D\subset \real$  are continua and
there are dense mesh subcontinua 
$$\imath:C_0\subset C \mbox{ and } \jmath:D_0\subset D$$
whose inclusions morphisms have almost inverses 
$$\lambda:C\rightarrow C_0 \mbox{ and }\kappa:D\rightarrow D_0,$$
then the correspondence $F\leftrightarrow \kappa\circ F\circ\imath$
defines one-to-one correspondence between $D$-valued real functions on $C$
and $D_0$-valued real functions on $C_0$.
Hence although there is no such continuum as ``$C(C,D)$'', we can
treat $D$-valued real functions on $C$ via mesh continua such as $C(C_0,D_0)$.

\paragraph{Composition}
Suppose $C$ is a mesh continuum and $D\subset \real$. 
Let $F$ be a $D$-valued real function on $C$ and
$G$ be a real function on $D$. Then a real function $G\circ F$
is defined by
$$ G\circ F:=[g\circ f]$$
with $f\in F$ and $g\in G$. This is well-defined since 
$f\approx f'$ and $g\approx g'$ implies
$ g\circ f\approx g'\circ f' $.

However this ``point wise definition'' cannot be given precise meaning
as a morphism 
$$C(C,D)\times C(D,\real)\rightarrow C(C,\real) $$
since there are no such continuum as ``$C(D,\real)$''.
However if $D'\subset D$ is a dense mesh subcontinuum with
an almost inverse $\kappa:D\rightarrow D'$. Then 
the real functions on $D$ corresponds to those on $D'$ 
in bijective way and we can take $C(D',\real)$ as one 
realization of the phantom ``$C(D,\real)$''.

Then the composition $(F,G)\mapsto G\circ F$ is realized
by the morphism 
$$ \gamma: C(C,D)\times C(D',\real)\rightarrow C(C,\real) $$
defined by $(f,g)\mapsto g\circ\kappa\circ f$.

This does not depend on the choice of $D'$ in the sense that
the following diagram commutes up to indistinguishability, 
whenever $\kappa_i:D\rightarrow D'_i$ \itwo are almost inverse
of the inclusions and $\beta$ is the restriction of $\kappa_2$ on $D'_2$.
\begin{center}
  $\xymatrix{ 
   C(C,D)\times C(D_2',\real)\ar[rd]_{1\times\kappa_1^{*}}\ar[rrd]^{\gamma_1}\ar[dd]_{1\times \beta^{*}}  && \\
   & C(C,D)\times "C(D,\real)"\ar[r]    & C(C,\real)   \\
   C(C,D)\times C(D_1',\real)\ar[ru]^{1\times\kappa_2^{*}}\ar[rru]_{\gamma_2}  &&   
}$
\end{center}

\subsection{Examples}

\subsubsection{Polynomial Functions}

Let $f(x)$ be a polynomial
$$  f(x)=\sum_{i=0}^{n}a_ix^{i}$$
with $n$ \accessible{}  and $a_{i}\in \virtualline_{\rational}$. Then
the function $r\mapsto f(r)$ is continuous on 
$\virtualline_{\rational}$ and defines a
real function $\lambda x.P(x)$ on $\virtualline$, 
called the polynomial functions defined by $f$. 
It is denoted by 
$$  F(x):=\sum_{i=0}^{n}p_ix^{i},$$
where $p_i=[a_i]$ and is called the real polynomial of degree $n$ if $p_n\neq 0$. 
For a finite real number $t$, its value is 
$$ F(t):=\sum_{i=0}^{n}p_it^{i}.$$

\subsubsection{Exponential}
\index{exponential morphism}
For huge $T$ and rational $r$, define a rational number by
  $$exp(r,T):=\sum_{i=0}^{T}\Frac{r^{i}}{i!}.$$
\begin{Proposition}
\label{prop:exp-independence-of-T}
If $r\in \virtualline_{\rational}$, the sum $exp(r,T)$ 
converges. In particular, if $T,S$ are huge, then 
$$exp(r,T)\approx exp(r,S).$$
\end{Proposition}
\def\finite{ finite }
\def\huge{ huge }

The proposition follows directly from the following
Lemmas~\ref{lem:finiteness},~\ref{lem:estimate-of-tail}.

\begin{Lemma}
\label{lem:finiteness}
If $r>0$ is a \finite positive rational, and $T$ a huge 
  number, then $\sum_{i=0}^T \Frac{r^i}{i!}$ is \finite.
\end{Lemma}
\begin{eproof}
  Take an \accessible{}  number $k$ satisfying $2r<k$. The sum
  $\sum_{i=0}^{k-1}\Frac{r^{i}}{i!}$ is \finite{} being 
  the sum of an \accessible{}  number of bounded rational. 
  Hence it suffices to show that
  $$\sum_{i=k}^{T}\Frac{r^{i}}{i!}$$
  is bounded.
  \begin{eqnarray*}
    \sum_{i=k}^T \Frac{r^i}{i!} 
    &=& \Frac{r^k}{k!}\left(1+\sum_{i=1}^{T-k} \Frac{r^i}{(k+i)(k+i-1)\cdots(k+1)}\right) \\
    &<& \Frac{r^k}{k!}\left(1+\sum_{i=1}^{T-k} \left(\Frac{r}{k}\right)^i\right) \\
    &<& \Frac{r^k}{k!}\Frac{1}{1-\frac{r}k}< 2\Frac{k^k}{2^kk!} \\
  \end{eqnarray*}
\end{eproof}

\begin{Lemma}
\label{lem:estimate-of-tail}
  If $r>0$ is a \finite{}  rational, and $T,N$ are \huge numbers, then
  $$\sum_{T\leq i\leq T+N}\Frac{r^i}{i!}\approx 0.$$
\end{Lemma}

\begin{eproof}
  \begin{eqnarray*}
    \sum_{i=T}^{T+N}\Frac{r^i}{i!}
    &=&
    \Frac{r^{T}}{T!}\left(\;
      1+\sum_{i=1}^N \Frac{r^i}{(T+1)(T+2)\cdots (T+i)}
      \;\right)
    \\
    &\leq&
    \Frac{r^{T}}{T!}\sum_{i=0}^N \Frac{r^i}{i!}
  \end{eqnarray*} 

  Take an \accessible{}  $k$ with $r<k$. Since $\frac{r}{i}<1$  for $i>k$, we have
  $$ \Frac{r^T}{T!}= \Frac{r^k}{k!}\Frac{r^{T-k}}{(k+1)(k+2)\cdots T}
  =\Frac{r^k}{k!}\frac{r}{k+1}\; \frac{r}{k+2}\; \frac{r}{k+3}\; \cdots
  \;\frac{r}{T} \leq \Frac{r^k}{k!}\Frac{r}{T} \approx 0. $$
  By Lemma \ref{lem:finiteness}, $\sum_{i=0}^N \Frac{r^i}{i!}$ is \finite,
  whence
  $$\Frac{r^{T}}{T!}\sum_{i=0}^N \Frac{r^i}{i!}\approx0.$$
\end{eproof}

The following is a more primitive expression for the exponential
function.  This gives an example of huge number of product of
rationals indistinguishable from $1$ gives a number $\succ 1$.
\begin{Proposition} If $r$ is \finite and $T$ is huge, then
\label{prop:product-formula-exp}
$$
  exp(r,T)\approx \left(\;
1+\Frac{r}{T} \;\right)^{T}.
$$
\end{Proposition}
\begin{eproof}
  \begin{eqnarray*}
    \left(\;
      1+\Frac{r}{T}
      \;\right)^{T}
    &=& 
    \sum_{i=0}^T \binom{T}{i}\Frac{r^i}{T^i}
    \\  &=& 
    \sum_{i=0}^T a_i\Frac{r^i}{i!} 
  \end{eqnarray*}
where
$$a_i:=\left(1-\Frac1T\right)\left(1-\Frac2T\right)
\cdots\left(1-\Frac{i-1}{T}\right).$$
If $i$ is \accessible{}  $a_i\approx 1$ hence 
\begin{equation}
\label{eq:111013}
\sum_{i=0}^n \Frac{a_ix^i}{i!}\approx \sum_{i=1}^n\Frac{x^i}{i!}
\end{equation}
holds for \accessible{}  $n$. Hence by Robinson's lemma,  there is a \huge $N$ such that
(\ref{eq:111013}) holds for $n\leq N$. If $T\leq N$ we have nothing more to show.
Suppose $T> N$. 
$$
\left(\;1+\Frac{r}{T}\;\right)^{T}\approx \sum_{i=0}^N
\Frac{r^i}{i!}+\sum_{i=N+1}^T \Frac{a_ir^i}{i!}.
$$
Since $a_i<1$, the second term is infinitesimal by Lemma~\ref{lem:estimate-of-tail}.
\end{eproof}

\begin{Proposition}
\label{prop:exp-continuity} 
  If $p$ is a finite real number and $r\in p$ and $T\gg1$. Then
the real number $[exp(r,T)]$ does not depend on the representation $r$ and
$T$.
\end{Proposition}
\begin{eproof}
  If $r,r'\in p$ then  
  \begin{equation}
    \label{eq:3-831}
    \sum_{i=0}^n \Frac{r^i}{i!}  \approx \sum_{i=0}^n \Frac{r'^i}{i!},
  \end{equation}
for \accessible{} $n$. 
  Hence by the Robinson's lemma \ref{prop:robinson}, there is a huge $N$
  such that \eqref{eq:3-831} holds for $n\leq N$.
Hence, by Proposition~\ref{prop:exp-independence-of-T}, 
$$
exp(r,T)\approx exp(r,N)\approx exp(r',N)\approx exp(r',T').
$$
\end{eproof}

Hence the function 
$$exp_T:\virtualline_{\rational}\rightarrow (0,\infty)_{\rational}$$
defined by $exp_T(r):=exp(r,T)$
is continuous by Proposition~\ref{prop:exp-continuity} and 
represents a real function, called \textit{exponential function}:
$$
exp:\virtualline\rightarrow (0,\infty),
$$
which is independent of $T$.
Its value at a real number $p=[r]$ is given by 
$$exp(p):=[exp(r,T)]$$
by any $T\gg 1$.

\begin{Proposition}
  \label{prop:additivity} If $x_i$ \itwo are finite real numbers, then
  $$ exp(x_1+x_2)=exp(x_1)exp(x_2). $$
In other words, if $r_i$ \itwo are \finite rationals 
then $\exp(r_1+r_2,T)\simeq \exp(r_1,T)\exp(r_2,T)$.
\end{Proposition}
\begin{eproof}
  If $N$ is \huge, 
  \begin{eqnarray*}\exp(r_1,N)\exp(r_2,N)&=&
    \left(\;\sum_{i=0}^{N}\Frac{r_1^{i}}{i!}\;\right)
    \left(\;\sum_{j=0}^{N}\Frac{r_2^{j}}{j!}\;\right)
    =
    \sum_{k=0}^{2N}\sum_{i+j=k,i\leq N,j\leq N}\Frac{r_1^{i}r_{2}^{j}}{i!j!}
    \\
    &=&
    \sum_{k=0}^{N}\Frac{(r_1+r_{2})^{k}}{k!} + U =\exp(r_1+r_2,N)+U,
  \end{eqnarray*}
  where
  $$U:=
  \sum_{k=N+1}^{2N}\sum_{i+j=k,i\leq N,j\leq
    N}\Frac{r_1^{i}r_{2}^{j}}{i!j!}.$$
  By Lemma~\ref{lem:estimate-of-tail},
  \begin{eqnarray}
    \mylabeleq{eq:2-20120607}
    |U| & \leq &
    \sum_{k=N+1}^{2N}\sum_{i+j=k,i\leq N,j\leq
      N}\Frac{|r_1|^{i}|r_{2}|^{j}}{i!j!}
    \\
    &\leq &
    \sum_{k=N+1}^{2N}\sum_{i+j=k}\Frac{|r_1|^{i}|r_{2}|^{j}}{i!j!}
    =\sum_{k=N+1}^{2N}\Frac{(|r_1|+|r_2|)^{k}}{k!}\approx 0
  \end{eqnarray}
\end{eproof}

\begin{Corollary} If $x$ is a finite real number then
$exp(x)exp(-x)=1$. 
In other words, for a \finite{} rational $r$ and huge $T$,
  $$\exp(r,T)\exp(-r,T)\approx 1.$$
\end{Corollary}

\begin{Lemma}
\label{lemma:exp-estimate-below}
 If $x$ is a nonzero finite real number, then
$$  exp(x) > 1+x. $$
In other words, if $r\in (-\infty,\infty)_{\rational}$ satisfies
$r\not\approx 0$ then 
\begin{equation}
\label{eq:20120226-10}
\exp(r,T)\succ 1+r.
\end{equation}
  \mylabel{lem:exp-estimate}
\end{Lemma}
\begin{eproof}

We may assume $r\succeq -1$ since otherwise the right hand is non positive.

Suppose $0\prec r$. Then 
  $$\exp(r,T)= 1+r+\sum_{i=2}^T\Frac{r^i}{i!} \geq 1+r+\Frac{r^2}2 \succ 1+r,$$

Suppose now $-1\prec r\prec 0$. Then $-r\succ 0$ and
  $$\exp(-r,T)= 1+(-r)+\sum_{i=2}^T\Frac{(-r)^i}{i!} 
\leq -\frac{r^2}{2}+ \sum_{i=0}^{T}(-r)^{i} 
\prec \frac{1-(-r)^{T+1}}{1+r}\leq \frac1{1+r}.
$$
Hence $exp(-r,T)\prec \frac1{1+r}$ namely
$exp(r,T)\succ 1+r$ and 
\eqref{eq:20120226-10} holds also for $-1\prec r\prec 0$.
\end{eproof}

\begin{Corollary}
\label{cor:exp(-inf)}
If $T,S\gg1$ then $exp(-T,S)\approx 0$ and $exp(T,S)\gg 1$.
\end{Corollary}
\begin{eproof}
By Lemma~\ref{lemma:exp-estimate-below}, for \accessible{} $n$, 
$$ exp(n,S)> 1+n,$$
whence for some huge $T$, 
$$ exp(T,S)> S+T\gg1.$$
If $T_1>T$, then
$$ exp(T_1,S)> exp(T,S)\gg1.$$

By Proposition~\ref{prop:product-formula-exp} and
Lemma~\ref{lemma:exp-estimate-below}, for \accessible{} $n$, 
$$exp(-n,S)\approx \frac{1}{exp(n,S)}
\prec \frac{1}{1+n},$$
whence
$$exp(-n,S)<\frac1{1+n}$$
for all \accessible{} $n$, whence there is a huge $M$ satisfying such that
for all $T\leq M$
$$exp(-T,S)<\frac1{1+T}\approx 0,$$
whence $exp(-T,S)\approx 0$. If $T_1>T$, then
$$
exp(-T_1,S)\approx \frac{1}{exp(T_1,S)}
\leq
\frac{1}{exp(T,S)}\approx exp(-T,S)\approx 0.
$$
\end{eproof}

\begin{Proposition}
\label{prop:exponential-preserving-prec}
The exponential is injective and order preserving. Namely, 
for finite real numbers $p,q$, $exp(p)=exp(q)$ if and only if $p=q$ and
$exp(p)<exp(q)$ if and only if  $p<q$.
In other words, for \finite rationals $x,y$ and huge $N$, 
the condition $\exp(x,N)\approx  \exp(y,N)$ implies $x\approx y$, 
and the condition $\exp(x,N)\prec \exp(y,N)$ implies $ x\prec y$.
\end{Proposition}
\begin{eproof}
  If $p> 0$, then by Lemma~\ref{lem:exp-estimate},
  $\exp(p)> 1+p> 1$.  Hence if $p<q$, then
  $\exp(q-p)>1 $, which implies $\exp(p)<\exp(q)$
when multiplied by $exp(p)$.
  The other assertions follow from observing that 
  the mutually disjoint and exhausting conditions
  $$exp(p)<exp(q),\quad \exp(p)=exp(q),\quad \exp(p)> \exp(q)$$
  hold according respectively to 
the mutually disjoint and exhausting conditions $p<q$,$p=q$,$p>q$.

\end{eproof}

\subsubsection{Logarithm}
\index{logarithm@logarithm}\
\begin{Lemma}
\label{lemma:logarithm}
  \mylabel{prop:logdef} Let $r\in (0,\infty)_{\rational}$ and $T\gg1$. 
  Define
$$
\log(r,T):=\frac1T\max\setii{k\in [-T^2..T^2]}{\exp(\frac{k}{T},T)\leq r}.
$$
Then for $s\in \virtualline_{\rational},r,r'\in (0,\infty)_{\rational}$
\begin{enumerate}
\def\labelenumi{(\theenumi)}
\item\label{item:20120226-a21}  
$\exp(\log(r,T),T)\approx r$,
\item\label{item:20120226-a22}  
$\log(\exp(s,T),T)\approx s$,
  \item\label{item:20120226-a23} 
$\log(r,T)\approx s$ if and only if $r\approx \exp(s,T)$,
  \item\label{item:20120226-a23-2} 
$\log(rr',T)\approx \log(r,T)+\log(r',T)$,
  \item\label{item:20120226-a24} 
 $\log(r,T)$ is monotone increasing,
\item\label{item:20120226-a25} 
$\log(r,T)\succ 0$, $\log(r,T){\approx}0$, and $\log(r,T)\prec 0$ according 
respectively to     $r\succ 1$, $r\approx 1$, and $r\prec 1$.
\item\label{item:20120226-a26} 
If $r\approx r'$, then $\log(r,T)\approx \log(r',T)$.
\end{enumerate}
Hence the function $r\mapsto\log(r,T)$ defines a morphism
$$\log_T:(0,\infty)\rightarrow \virtualline, $$
which is an almost inverse of $\exp_T$. 
If $T'$ is another huge number then the morphisms $\log_{T}$ and $\log_{T'}$ are indistinguishable.

\end{Lemma}

\begin{eproof}
First, we show that $\log(r,T)$ is well-defined. By
Corollary~\ref{cor:exp(-inf)}, $exp(-T,T)\approx 0$ whence
the set $\setii{t\in [-T..T]_{\frac1T}}{\exp(t,T)\leq r}$
is not empty and its maximum $t_0$ is defined, which is finite. 
In fact if $r\geq 1$ then $t_0<r-1$
and $t_0\geq 0$ since $exp(0,T)=1$ 
Suppose $r<1$. Obviously $t_0<0$. 
If $t_1:=\frac2{r}-1$ then
$$
exp(-t_1,T)
\approx
\frac1{exp(t_1,T)}\leq \frac1{1+t_1}=\frac{r}2\prec r,
$$
whence $ t_0>-t_1>-\infty$.

Then 
$$\exp(\frac{k}{T},T)\leq r< \exp(\frac{k+1}{T},T)),$$
Hence $\log(r,T)=\frac{k}T$ satisfies
$$\exp(\log(r,T),T)\approx r.$$

Let $s\in \virtualline_{\rational}$. Then 
$\exp(\frac{k}{T},T)\leq \exp(s,T) $
is equivalent to $\frac{k}{T}\leq s$, hence 
$$\log(\exp(s,T),T)= \frac{[sT]}{T} \approx s.$$

The third assertions follow from the first and Proposition \ref{prop:exponential-preserving-prec}.

The fourth assertion can be verified as follows. Let $r,r'\in (0,\infty)_{\rational}$.
Then
\begin{eqnarray*}
\exp(\log(rr',T),T)&\approx& 
rr'=\exp(\log(r,T),T)\exp(\log(r',T),T)\\
&\approx&
\exp(\log(r,T)+\log(r',T),T),
\end{eqnarray*}
whence the assertion follows by Proposition~\ref{prop:exponential-preserving-prec}. 
\end{eproof}

The real function $\log:=[\log_T]$ on the continuum $(0,\infty)$ 
is called the \textit{natural logarithm function}. Lemma~\ref{lemma:logarithm}
can be rephrased as follows.
\begin{Proposition} Let $x,x'>0$ and $y$ be finite real numbers. Then
\label{prop:logarithm}
\begin{enumerate}
\def\labelenumi{(\theenumi)}
\item\label{item:20120226-b21}  
$\exp(\log(x))=x$,
\item\label{item:20120226-b22}  
$\log(\exp(y))=y$,
  \item\label{item:20120226-b23} 
$\log(x)=y$ if and only if $x=\exp(y)$,
  \item\label{item:20120226-b23-2} 
$\log(xx')=\log(x)+\log(x')$,
  \item\label{item:20120226-b24} 
 $\log(x,T)$ is monotone increasing,
\item\label{item:20120226-b25} 
$\log(x)>0$, $\log(x)=0$, and $\log(x)<0$ according 
respectively to $x>1$, $x=1$, and $x<1$.
\end{enumerate}
\end{Proposition}

For finite real numbers  $x,y$ with $x>0$, we define 
$$
x^y:=\exp(y\log(x)).
$$
It is easily verified the following.
\begin{Proposition}
  \begin{description}
  \item[(1)] If $y$ is commensurable, then $x^{y}$ 
coincides with the $x^{y}$ defined in \S \ref{sec:algebr-oper-real}.
  \item[(2)] $x^{y+z}=x^{y}x^{z}$,
  \item[(3)] $x^{yz}=(x^{y})^{z}$.
  \end{description}
\end{Proposition}

\subsection{Mean Value Theorem}

\begin{Theorem}[Mean value theorem]

\label{th:mean-value-theorem}
Let $p,q$ be real numbers with $p<q$ represented respectively by $a,b\in \rational$.
If a real function $F$ on $[a,b]$ satisfies
$F(p)<t<F(q)$ for a real number $t$, then there is a real number $s$ satisfying
$F(s)=r$ and $p<s<q$.

In other words, let $a,b\in \rational$ with $a\prec b$.
If a continuous rational valued function $f$ on $[a,b]$ satisfies  
$f(a)\prec c\prec f(b)$ for a rational number $c$, 
then there is an $r\in [a,b]_{\rational}$  satisfying $f(r)\approx c$.
\end{Theorem}
\begin{eproof}
Let $\varepsilon>0$ be an infinitesimal and 
put 
$m_1:=[a/\varepsilon]+1$ and $m_2:=[b/\varepsilon]$.
Then $m_1\varepsilon,m_2\varepsilon\in[a,b]_{\varepsilon}$ and 
$$f(m_1\varepsilon)\approx f(a) \prec c\prec f(b)\approx f(m_2\varepsilon).$$
Hence we can define
$$
x:=\varepsilon \min\setii{i\in [m_1..m_2]}{f(i\varepsilon)\geq c},
$$
which satisfies $f(x-\varepsilon)<c\leq f(x)$.
Since $f(x-\varepsilon)\approx f(x)$,
we have $f(x)\approx c$.

\end{eproof}

Note that in the theorem the real numbers $p$ or $q$ may be infinite.
\subsection{Maximum Principle}
\begin{Theorem}[Maximum principle] 
Let $F$ be a real function on a continuum $C$ with a dense subset $A$.
Then there are points $p,q$ of $C$ such that
$$   F(p)\leq F(x)\leq F(q)$$
for all points $x$ of $C$.
In other words, if $f\in F$ is a continuous rational valued function on $C$,
then there are positions $x_m,x_M\in A$ such that 
\begin{equation}
\label{eq:906-10}
f(x_m) \preceq f(x)\preceq f(x_M)
\end{equation}
for all $x\in X$.
\end{Theorem}
\begin{eproof}
Let $f(x_m)=\min_{x\in A}f(x)$ and $f(x_M)=\max_{x\in A}f(x)$. 
Then for all $x\in A$
\begin{equation}
\label{eq:20120227-1}
f(x_m) \preceq f(x) \preceq  f(x_M).
\end{equation}
Let $y\in C$. Then there is a $z\in A$ with $y\approx z$. 

Since \eqref{eq:20120227-1} is valid for $x=z$ and $f(y)\approx f(z)$,
it is valid also for $x=y$ by Proposition \ref{prop:daishou}.
\end{eproof}

\begin{Corollary}
Every real function on a compact continuum has maximum and minimum values.
\end{Corollary}

\subsection{Behavior of Real Funcrtions on a Point}
\label{sec:landau-symbol}
We introduce a method of describing the infinitesimal behavior of functions
on a point which will be used extensively in the treatment of calculus.

In the following, $C=(\support{C},d,\approx_d)$ denotes a metric continuum and 
$f(x),g(x)$ are rational valued functions on $\support{C}$.

\begin{Definition}
\label{def:Landau}
We write for $a\in \support{C}$ and \accessible{}  $n$
\begin{equation} 
\label{eq:904-20}
f(x)\equiv_n 0 \mbox{ if } x\approx a
\end{equation}
if 
\begin{equation}
\label{eq:911-01}
\Frac{|f(x)|}{d(x,a)^n}\approx0 
\end{equation}
whenever $x\neq a $ and $x\approx a$
\end{Definition}

This relation is not indistinguishability-invariant, namely, $f\approx 0$ does not
necessarily imply (\ref{eq:904-20}). For example, let $\varepsilon>0$ be
an infinitesimal and define a rational valued function $f$ on $[-1,1]_{\rational}$ by

\begin{equation}
  \label{eq:20120302}
f(x):= \left\{
  \begin{array}{cc}
  |x|^n  & \mbox{for $|x|\leq \varepsilon$} \\
  \varepsilon^n  & \mbox{for $|x|>\varepsilon$ } 
  \end{array}
\right.
\end{equation}
Then $f\approx 0$ but
$
\frac{f(x)}{|x|^n}=1
$
for $|x|\leq \varepsilon$ and hence it is not the case that $f(x)\equiv_{n}0$ if $x\approx0$.

However the
following weaker condition is indistinguishability-invariant.
\begin{Definition}
\label{def:Landau-approx}
We write, for $a\in \support{C}$ and \accessible{}  $n$
\begin{equation}
\label{eq:904-21}
f(x)\approx_n 0 \mbox{ if } x\approx a
\end{equation}
if there is an infinitesimal $\varepsilon>0$ such that (\ref{eq:911-01}) holds
whenever $d(x,a)\geq \varepsilon$ and $x\approx a$
\end{Definition}
The following two propositions give basic properties of this condition.
\begin{Proposition}
\label{prop:basic-property-landau-symbol}
If $f,g$ are rational valued functions on $\support{C}$ such that $f\approx g$ 
and $a,b\in \support{C}$ are indistinguishable, then
$$ \mbox{$f(x)\approx_n 0 $ if $x\approx a$}
\Longleftrightarrow 
\mbox{
$g(x)\approx_n 0 $ if $x\approx b$
}. $$
\end{Proposition}
\begin{eproof} Let $n$ be an \accessible{} number and $a\in C$. 

First we show that if $f\approx 0$, then $f(x)\approx_n0$ if $x\approx a$.
Put 
$$\eta:=\max\setii{|f(x)|}{x\in \support{C}}\approx0.$$ 
Let $\varepsilon$ be an infinitesimal such that $\eta<\varepsilon^{n+1}$.
For example take a huge $N$ and the minimal number $k$ such that 
$\eta<\left(\frac{k}{N}\right)^{n+1}$ and put $\varepsilon=\frac{k}{N}$.
Then $\eta<\varepsilon^{n+1}$. 
Moreover $\left(\frac{k-1}{N}\right)^{n+1}\leq \eta$ 
implies $\varepsilon^{n+1}\approx \eta\approx 0$ whence
$\varepsilon\approx 0$ by Lemma \ref{lemma:continuity-of-root}.

If $ \varepsilon<d(x,a)$, then
$$
\left|
\Frac{f(x)}{d(x,a)^{n}}\right|<\Frac{\eta}{\varepsilon^n}<\varepsilon \approx 0,
$$
whence $f(x)\approx_n0$ if $x\approx a$.

Suppose now $a\approx b$ and $f(x)\approx_n0 $ if $x\approx a$.
Then there is an infinitesimal $\varepsilon>0$ such that
\eqref{eq:911-01} holds when $x\approx a$ and $d(x,a)>\varepsilon$.
Put $\varepsilon_1:=\max\seti{\varepsilon+d(a,b),3d(a,b)}\approx 0$ and
suppose $0\approx d(x,b)>\varepsilon_1$. Then
$$d(x,a)>d(x,b)-d(b,a)\geq \varepsilon, $$
whence \eqref{eq:911-01} holds.  On the other hand
$$
d(x,a)\geq d(x,b)-d(a,b)\geq 3d(a,b)-d(a,b)=2d(a,b),
$$
whence
$$
d(x,b)\geq d(x,a)-d(a,b)\geq \frac12d(x,a).
$$
This implies
$$
\Frac{|f(x)|}{d(x,b)^n}\leq
2^n\Frac{|f(x)|}{d(x,a)^n}\approx 0
$$
whence $f(x)\approx_n0$ if $x\approx b$.
\end{eproof}

\begin{Proposition}
\label{prop:equivalence-higher-infinitesimal}
Suppose $C_i$ \itwo are metric continua and
$\alpha:\support{C_1}\rightarrow \support{C_2}$ a function satisfying, 
for $x,y\in \support{C_1}$, 
$$ K_1d_1(x,y)<d_2(\alpha(x),\alpha(y))< K_2d_1(x,y) \mbox{ if $d_1(x,y)>\eta$}$$
with finite rational numbers $K_1,K_2$ and an infinitesimal $\eta>0$.
Suppose $f$ is a continuous rational valued function on $\support{C_2}$ and
$a\in \support{C_1}$. Then
\begin{equation}
\label{eq:20120307-1}
f(y)\approx_n0 \mbox{ if } y\approx \alpha(a) \mbox{ on $C_2$,}
\end{equation}
implies
\begin{equation}
\label{eq:20120307-2}
f(\alpha(x))\approx_n 0 \mbox{ if } x\approx a \mbox{ on $C_1$.}
\end{equation}
\end{Proposition}
\begin{eproof}
Suppose \eqref{eq:20120307-1} holds. Then there is an infinitesimal $\varepsilon>0$
such that if $\varepsilon< d_2(y,\alpha(a))\approx 0$, then
$$
\frac{|f(y)|}{d_2(y,\alpha(a))^n}\approx 0.
$$
If $d_1(x,a)>\max\seti{\eta,\frac{\varepsilon}{K_{1}}}$, then
$d_2(\alpha(x),\alpha(a))>K_1d_1(x,a)>\varepsilon$, whence
$$\frac{|f(\alpha(x))|}{d_1(x,a)^n}
\leq 
K_2^{-n}\frac{|f(\alpha(x))|}{d_2(\alpha(x),\alpha(a))^n}\approx 0.
$$
Hence \eqref{eq:20120307-2} holds.
\end{eproof}

The following shows that the relation $\approx_n$ is invariant under 
quasi-identities.
\begin{Proposition}
\label{prop:invariance-of-landau-by-quasi-identity}
\label{prop:invariance-of-landau-by-normal-equivalence}
Let $\varepsilon,\eta>0$ be infinitesimal and  
$\alpha: X:=[0,1]_{\varepsilon}^{n}\rightarrow Y:=[0,1]_{\eta}^n$ represents
a quasi-identity and $f$ a continuous rational valued function on 
$[0,1]_{\eta}^n$. 
Then, for \accessible{} $n$ and $a\in X$,
\begin{equation}
\label{eq:910-10}
 f(y)\approx_n 0 \mbox{ if $y\approx \alpha(a)$ } 
\end{equation}
if and only if
\begin{equation}
\label{eq:910-11}
f(\alpha(x))\approx_n 0\mbox{ if $x\approx a$.} 
\end{equation}
\end{Proposition}
\begin{eproof}
Let $\beta$ be an almost inverse of $\alpha$. Then there is 
an infinitesimal $\delta>0$ such that 
$$
d(x,\alpha(x)),d(y,\beta(y))<\delta
$$
holds for all $x\in [0,1]_{\varepsilon}^{n}$ and
$y\in [0,1]_{\eta}^{n}$.
Then 
$$  d(x,y)-2\delta < d(\alpha(x),\alpha(y))< d(x,y)+2\delta.  $$
Hence if $d(x,y)>4\delta$, then
$$   \frac12d(x,y)<d(\alpha(x),\alpha(y))<\frac{3}{2}d(x,y).$$
Hence by Proposition~\ref{prop:equivalence-higher-infinitesimal},
the condition  \eqref{eq:910-10} implies \eqref{eq:910-11}.

Assume now \eqref{eq:910-11}. Put $b=\alpha(a)$. Since $a\approx \beta(b)$,
we have
\begin{equation}
\label{eq:910-11-a}
f(\alpha(x))\approx_n 0\mbox{ if $x\approx \beta(b)$.} 
\end{equation}
Since 
$$   \frac12d(x,y)<d(\beta(x),\beta(y))<\frac{3}{2}d(x,y)$$
whenever $d(x,y)>4\delta$, we infer from
Proposition~\ref{prop:equivalence-higher-infinitesimal} 
that \eqref{eq:910-11-a} implies
$$
f(\alpha(\beta(y)))\approx_n 0\mbox{ if $y \approx b$.} 
$$
Since $f(\alpha(\beta(y)))\approx f(y)$ and 
$b\approx \alpha(a)$, we obtain \eqref{eq:910-10} by Proposition~\ref{prop:basic-property-landau-symbol}.
\end{eproof}

Although the relation $\approx_{n}$ is weaker than $\equiv_{n}$, 
the following proposition shows that the former implies the latter 
if the continuum is replaced with a suitable equivalent continuum.
\begin{Proposition}
\label{prop:approx-imply-equiv}
Let $\varepsilon>0$ be an infinitesimal. 
Suppose $C=[r,s]_{\varepsilon}^{n}$ is given the
metric $d_{\infty}$ and a continuous rational valued function $f$ on 
$\support{C}\times \support{C}$ satisfies
$$ f(x,a)\approx_n 0 \mbox{ if }  x\approx a$$
for all $a\in C$. 
Then there is an integer $L>0$ such that $L\varepsilon\approx0$ and
$$g(x,b) \equiv_n0 \mbox{ if } x\approx b,$$
for all $b\in [r,s]_{L\varepsilon}$
where $g=f|X\times X$.
\end{Proposition}
\begin{eproof}
For each $a\in [r,s]_{\varepsilon}^n$ let $\delta_a>0$ be an infinitesimal such that
\begin{equation}
\label{eq:20120302-2}
\Frac{|f(x,a)|}{d(x,a)^n}\approx0 
\end{equation}
holds if $x$ satisfies $0\approx d(x,a)>\delta_a$.  
Put $\delta=\max\setii{\delta_a}{a\in [r,s]_{\varepsilon}^n}\approx 0$.
Then \eqref{eq:20120302-2} holds for all $x,a\in [r,s]_{\varepsilon}^n$
with $d(x,a)>\delta$. 
Let $L$ be the least integer greater than $\frac{\delta}{\varepsilon}$.  
Then $\delta<L\varepsilon\approx 0$ and different $y,z\in [0,1]_{L\varepsilon}^n$
satisfies $d(y,z)>L\varepsilon>\delta$. 
Hence \eqref{eq:20120302-2} holds for different $x,a\in [0,1]_{L\varepsilon}^n$,
whence 
$g(y,b)\approx_n0$ if $y\approx b$ for each $b\in [r,s]_{L\varepsilon}^{n}$.
\end{eproof}

The following proposition shows that the condition \eqref{eq:904-21}
is characterized by the behavior of $f$ for $x\not\approx a$.
\begin{Proposition}
\label{prop:landau-symbol-in-other-words}
Let $C$ be a metric continuum and $f$ is a rational valued function on $\support{C}$.
Then for \accessible{} $n$ and $a \in\support{C}$, 
the following conditions are equivalent. 
\begin{description}
\item[(A)] $f(x)\approx_n 0$ if $x\approx a$,
\item[(B)] for each \accessible{}  number $k$, there is an \accessible{}  number $\ell$, such that
if $0\not\approx d(x,a)<\frac1\ell$ then $|f(x)|\leq \frac1{k} d(x,a)^{n}$.
\end{description}
\end{Proposition}
\begin{eproof}
Suppose the condition (A) is satisfied and there is a positive
$\varepsilon\approx0$ such that (\ref{eq:911-01}) if $0\approx d(x,a)>\varepsilon$.
Let $k$ be an \accessible{}  number. Then for all huge $i$, the objective condition 
\begin{equation}
\label{eq:905-2}
\varepsilon<d(x,a)<\Frac1i \mbox{ implies } \left|\Frac{f(x)}{d(x,a)^n}\right|<\frac1k
\end{equation}
holds, whence it holds also for an \accessible{}  $i=\ell$, whence (B).

Conversely suppose the condition (B) holds. Let $k$ be an \accessible{}  number. Then
there is an \accessible{}  $\ell$ such that for all \accessible{}  $p$, the objective condition 
\begin{equation}
\label{eq:905-3}
 \frac1p < d(x,a)<\frac1\ell \mbox{ implies } \left|\Frac{f(x)}{d(x,a)^n}\right|<\frac1k.
\end{equation}
holds. Hence for some huge $M_k$, the condition~\eqref{eq:905-3} holds for $p=M_k$.
By Proposition~\ref{prop:minimum-huge-number}, there is a huge $M$ such that
$M\leq M_k$ for all $k$. Then if $\frac1M< d(x,a)\approx 0$ then
$\left|\frac{f(x)}{d(x,a)^n}\right|<\frac1k$ for all \accessible{}  $k$, whence
(\ref{eq:911-01}).
\end{eproof}

\begin{Proposition}
\label{prop:landau-integral}
Let $n$ be an \accessible{} number and $a\in [0,1]_{\varepsilon}$. 
If a rational valued continuous function $f$ on $[0,1]_{\varepsilon}$ 
satisfies 
$$  f(x)\equiv_n 0 \mbox{ if }x\approx a$$
then
$$
g(x)\equiv_{n+1}0  \mbox{ if }x\approx a$$
where $g$ is defined by 
$$
g(x) = \left\{
  \begin{array}{cc}
\sum_{a\leq u< x}f(u)\Delta x & \mbox{for $x\geq a$} \\
\sum_{x\leq u < a}f(u)\Delta x    & \mbox{for $x< a$} 
  \end{array}
\right.
$$
with $\Delta x:=\varepsilon$.
\end{Proposition}
\begin{eproof}
By hypothesis if $0\approx |x-a|\neq 0$ then
$$
c:=\max_{x\approx a, x\neq a}\frac{|f(x)|}{|x-a|^n}\approx 0.
$$
If $x>a$
$$
\frac{|\sum_{a\leq u< x}f(u)\Delta x|}{|x-a|^{n+1}}
=
\frac1{|x-a|}\sum_{a\leq u< x}\frac{|f(u)|}{|x-a|^{n}}\Delta x
\leq \frac{c}{|x-a|}\sum_{a\leq u<x}\Delta x=c\approx 0
$$
Together with similar arguments for $x < a$ we conclude
$g(x)\equiv_{n+1}0$ if $x\approx a$.
\end{eproof}

The following lemmas show that the monomials are linearly independent
even within a point.
\begin{Lemma}
\label{lem:係数比較}
Let $k$ be an \accessible{}  number and $a_0,\cdots,a_k$ finite rational numbers.  
Suppose $\sum_{i=0}^ka_ix^i\approx_k0$ if $x\approx0$. Then
$a_i\approx 0$ for $i\in [0..k]$.
\end{Lemma}
\begin{eproof}
By hypothesis, there is an infinitesimal $\varepsilon>0$ such that 
if $0\approx |x|\geq \varepsilon$, then
$$\Frac{\sum_{i=0}^ka_ix^i}{|x|^k}\approx0.$$
Let $N$ be a huge number greater than $\frac1\varepsilon$, then  
$\varepsilon<\frac{j}{N}\approx0$ for $j\in [1..k+1]$. Hence
$$\sum_{i=0}^ka_iN^{k-i}j^{i-k}\approx 0, \mbox{ for } j=1,\cdots,k+1.$$
Define $(k+1,k+1)$ matrix $B$ and $k+1$ vector $\ba$ by 
$$B_{j,i}=((j+1)^{i-k}), \ba_i:=(a_iN^{k-i}),$$
($0\leq j,i\leq k$) then $B\ba\approx0$. Since the Van der Monde matrix $B$ has an inverse whose
components are bounded, we have
$$\ba\approx \bfzero.$$
Hence for each $i$, $ a_iN^{k-i}\approx0$, whence
$a_{i}\approx0$.
\end{eproof}

For multi-index 
$J=(j_1,\cdots,j_n)$, 
put $|J|:=j_1+\cdots+j_n$ and 
$x^{J}:=x_1^{j_1}x_2^{j_2}\cdots x_n^{j_n}$.
\begin{Lemma}
\label{lem:compare-coeff-multivar}
Let $k$ and $n$ be \accessible{} numbers. 
Suppose for each multi-index $J$ with $|J|\leq n$ a
finite rational number $a_J$ is given and satisfies
$$\sum_{|J|\leq n}a_Jx^J\approx_k0 \myif x\approx 0$$
on $\kukan{M}^n$. Then $a_J\approx 0$ for all $J$.
\end{Lemma}

\begin{eproof}
By assumption, there is an infinitesimal $\varepsilon>0$ such that
$0\approx d(0,x)\geq \varepsilon$ implies 
\begin{equation}
\label{eq:911-2012}
\Frac{\sum_{|J|\leq k}a_Jx^J}{d(0,x)^{k}}\approx0.
\end{equation}
Let $N=[\frac1\varepsilon]$. 
If components of $ \alpha\in \integer^n $ are \accessible{} ,
then
$$d(0,\frac{\alpha}N)=\Frac{d(0,\alpha)}{N}>\varepsilon,$$
whence by substituting $x=\frac{\alpha}N$ with $\alpha\in [0..k]^{n}$
in (\ref{eq:911-2012}), 
we obtain
\begin{equation}
\mylabeleq{eq:9-718}
\sum_{|J|\leq k}\alpha^J(a_JN^{k-|J|})\approx0,\quad \alpha\in[0..k]^{n}.
\end{equation}
Since the vectors 
$$\setii{(\alpha^J)_{|J|\leq k}}{\alpha\in [0..k]^{n}}
$$ are linearly independent, there is a subset
$T\subset[0..k]^{n}$ such that 
$A=(\alpha^J)_{|J|\leq k,\alpha\in T}$ is regular and the inverse matrix 
has bounded components. Hence ~(\ref{eq:9-718}) implies
$a_{J}N^{k-|J|}\approx 0$
for all $|J|\leq k$, whence $a_J\approx 0$.
\end{eproof}

\newpage
\section{Differentiation: Single variable }
\label{sec:calculus-single}
\renewcommand{\runningtitle}{Differentiation: Single variable} 
It turns out that
the differentiability of real functions can be defined as
the possibility of choosing good representations with
continuous difference quotients.
One might be puzzled that every real function seems to become
differentiable according to this definition since its representations
can be arbitralily specified within a point.
However the fringe of a point is not seperated from that of other 
``neighboring points'' and the nearer we move to the ``boundary'' of a point, the less
freedom we have for the specification of the behaviour of
representations of a real function.  This feature might be understood
by the fact that for every infinitesimal $\varepsilon$, there is an
huge number $L$ such that $L\varepsilon$ is still infinitesimal, which means
that given any two places within a point, we can magnify 
the point to have visible extent without breaking their indistinguishability.

\subsection{Difference Quotient}
In this section $\varepsilon>0$ is a fixed infinitesimal 
and $C=([0,1]_{\varepsilon},\approx)$ the rigid mesh continuum representing $[0,1]$.\footnote{We consider only the unit interval continuum
$[0,1]$ but everything can be directly generalized to general interval $[a,b]$.
}

Denote by $[0,1]_{\varepsilon}^{-}$ the subset of $[0,1]_{\varepsilon}$ obtained
by removing the greatest element, namely, $[\frac1\varepsilon]\varepsilon$.
and denote the next larger element of $x$ of $[0,1]_\varepsilon^{-}$ by $x^{+}$, 
namely, $x^{+}=x+\Delta x$, where $\Delta x=\varepsilon$.  Note that
$([0,1]_\varepsilon^{-},\approx)$ is also a rigid mesh continuum representing $[0,1]$.

\begin{Definition}[Difference operator]
  For a rational valued function $f$ on $[0,1]_\varepsilon$, define its
  \textit{difference} \index{difference operator@difference operator}
  $\Delta f$ defined for $x\in [0,1]_\varepsilon^{-}$ by
    $$\Delta f(x):=f(x^{+})-f(x).$$
The quotient 
  $$\sabun{f}(x):=\frac{f(x^{+})-f(x)}{\Delta x} \in \rational $$ 
   is called the
    \textit{difference quotient} \index{difference quotient@difference
      quotient} of $f$ at $x$ and the rational valued function $\sabun{f}$
   is called the \textit{difference quotient function} of $f$.
\index{difference quotient}

\end{Definition}

For a rational valued function $f$ on $[0,1]_{\varepsilon}$ 
and $a,b\in [0,1]_{\varepsilon}$ with $a<b$, we define
$$
\sum_a^{b} f:=\sum_{u\in [a,b]_{\varepsilon}}f(u)=\sum_{i=0}^{\frac{b-a}{\varepsilon}}f(a+i\varepsilon)
$$

\begin{Proposition}
If $f$ is a rational valued function on $[0,1]_\varepsilon$ and $a,b\in [0,1]_\varepsilon^{-}$ with $a<b$, then 
\begin{equation}
\label{eq:sum-of-difference}
\sum_a^{b}\sabun{f}\Delta x = f(b^{+})-f(a).
\end{equation}
\end{Proposition}

\begin{eproof}
Since 
$$\sabun{f}(x)\Delta x = \Delta f(x),$$
we have 
$$
\sum_{a}^{b}\sabun{f}\Delta x 
=
\sum_{x\in [a,b]_{\varepsilon}}(f(x^{+})-f(x))
=f(b^{+})-f(a).
$$
\end{eproof}

\begin{Proposition}
\label{prop:Sigmaf}
If rational valued functions $f,g$ on $[0,1]_{\varepsilon}$ 
satisfy $f\approx g$, then 
$$ \sum_a^b f\Delta x\approx \sum_a^b g\Delta x,$$
for $a,b\in [0,1]_{\varepsilon}^{-}$ with $a<b$. 
\end{Proposition}
\begin{eproof}
Put $c:=\max\setii{|f(x)-g(x)|}{x\in [0,1]_\varepsilon}$. Then $c\approx 0$ and 
$$
  \left|\sum_a^b f-\sum_a^b g\right|
\leq
  \sum_a^b|f-g|
\leq
\sum_a^bc\Delta x=c(b^{+}-a)\leq c \approx0.
$$
\end{eproof}

Basic relation between $f$ and its difference quotient $\sabun{f}$ is as follows:
\begin{Proposition}
\label{prop:bibun-hyouka}
Suppose a rational valued function $f$ on $[0,1]_\varepsilon$ is continuous and 
the difference quotient function $\sabun{f}$ 
is finite on $[0,1]_\varepsilon^{-}$. Then, for $x,y\in [0,1]_\varepsilon$ with $x<y$,
\begin{enumerate}
\def\labelenumi{(\theenumi)}
\item $ f(y) \preceq f(x)+M(y-x) $, where 
$$ M=\max\setii{\sabun{f}(u)}{u\in [0,1]_\varepsilon,x\leq u<y}, $$
\item $ f(y) \succeq f(x)+m (y-x) $, where
$$ m=\min\setii{\sabun{f}(u)}{u\in [0,1]_\varepsilon,x\leq u<y}, $$
\item $ |f(y)-f(x)|\preceq M_1 |y-x|$, where
$$ M_1=\max\setii{\left|\sabun{f}(u)\right|}{u\in [0,1]_\varepsilon,x\leq u< y}.$$
\end{enumerate}
\end{Proposition}
\begin{eproof}
$$
  f(y)-f(x) 
=
\sum_x^{y-\varepsilon}\sabun{f}\Delta x
\leq
\sum_x^{y-\varepsilon} M\Delta x=M(y-x).
$$

The other assertions can be shown similarly.
\end{eproof}

\begin{Corollary}
\label{cor:inverse-continuity-of-function-with-positive-sabun}

If $\sabun{f}\succ 0$ on $[0,1]_\varepsilon^{-}$, then $x\prec y$ implies $f(x)\prec f(y)$.
In particular $f(x)\approx f(y)$ implies $x\approx y$.
\end{Corollary}
\begin{eproof}
Suppose $\sabun{f}\succ 0$ on $[0,1]_\varepsilon$. Put 
$c=\min\setii{\sabun{f}(x)}{x\in [0,1]_\varepsilon^{-}}$.
If $y\succ x$, then 
$$f(y)-f(x)\succeq c(y-x)\succ 0.$$
If $f(x)\approx f(y)$, then neither $x\prec y $ nor $y\prec x$ is possible whence
$x\approx y$.
\end{eproof}

\begin{Remark}
  Even if $f$ is continuous, the difference quotient $\sabun{f}$ 
  may be neither continuous nor finite. 
  Moreover indistinguishable functions may have distinguishable difference quotients.
\begin{enumerate}
\item Define rational valued function $f,g$ on $[0,1]_{\varepsilon}$ by 
$$f(i\varepsilon):=i\varepsilon,\quad
g(i\varepsilon):={2[i/2]}\varepsilon.$$
Then  $f\approx g$ but 
$$
\sabun{f}(i\varepsilon)=1, 
\mbox{ but}\quad
\sabun{g}(i\varepsilon)=\left\{
  \begin{array}{cl}
  0 & \mbox{if $i$ is even}\\
 2 & \mbox{otherwise}
\end{array}
\right.
$$
Hence $\sabun{f}\not\approx\sabun{g}$. 
Although both $\sabun{f}$ and $\sabun{g}$ are finite,
$\sabun{g}$ is not continuous.

\item\mylabelitem{item:6} Let $N\gg1$ and $\varepsilon=\frac1{N^{2}}$.
Define a rational valued function $h$ on $[0,1]_{\varepsilon}$
by 
$$
h(i\varepsilon) = \left\{
  \begin{array}{cl}
  \frac{1}{N} & \mbox{ if $i$ is even}\\
 0 & \mbox{ otherwise},\\
\end{array}
\right.
$$
Since $h(x)\leq \frac1N$, $h\approx 0$ and $h$ is continuous 
but the difference quotient is not finite since 
$$ \sabun{h}(2i\varepsilon)=\frac{\frac1N}{\frac1{N^2}}=N.$$
\end{enumerate}
\end{Remark}

\subsection{Differentiability}
Let $\varepsilon$ be a positive infinitesimal. 
We say a real function $F$ on $[0,1]$ is \textit{represented by $(f,[0,1]_{\varepsilon})$},
sometimes simply \textit{represneted by $f$}, 
if $f$ is a continuous rational-valued function $f$ on $[0,1]_{\varepsilon}$
and $f\circ \kappa_{\varepsilon}$ represents $F$ 
where $\kappa_{\varepsilon}$ is the map defined in Proposition~\ref{prop:r-int-equiv-rat}.
\footnote{
This is equivalent to the condition that
$F$ is represented by $(f,[0,1]_{\varepsilon},\kappa_{\varepsilon}|_{[0,1]_{\rational}})$
in the terminology of \S~\ref{sec:real-functions}.
} 

\begin{Definition}
  A real function $F$ on $[0,1]$ is called
  \textit{differentiable} if it is represented by $(f,[0,1]_\varepsilon)$ with
$\varepsilon$ a positive infinitesimal whose difference quotient $\sabun{f}$ is continuous. 
We say $(f,[0,1]_\varepsilon)$ is a \textit{representation of $F$
    with continuous difference quotient}.
\end{Definition}

The real function on $[0,1]$ 
represented by $(\sabun{f},[0,1]^{-}_\varepsilon)$
does not depend on the choice of the representation $(f,[0,1]_\varepsilon,\alpha)$
by Proposition~\ref{prop:independence of difference quotient}. 
It is denoted by $F'$ and is called the \textit{derivative} \index{derivative@derivtive}
\index{$\frac{dF}{dx}$, $F'$} of $F$.
It is also denoted by $\frac{dF(x)}{dx}$.

\begin{Proposition}
\label{prop:independence of difference quotient}
Let $F$ be a real function on $[0,1]$. 
Let $(f_i,[0,1]_{\varepsilon_i})$ \itwo be representations of $F$ 
such that the difference quotients of $f_1,f_2$ are continuous.
Then the real functions on $[0,1]$ represented by the
difference quotients $\left(\sabun{f_i},[0,1]^{-}_{\varepsilon_i}\right)$ \itwo
coincides.
\end{Proposition}
\begin{eproof}
Put $\alpha_i=\kappa_{\varepsilon_i}|_{[0,1]_{\rational}}$ \itwo.
By Theorem~\ref{Th:taylor-formula} of the next section, 
we have for $a\in [0,1]_{\rational}$ 
$$ f_i(y)\approx_1 f_i(\alpha_i(a))+\sabun{f_i}(y)(y-\alpha_i(a)) 
\mbox{ if } y\approx \alpha_i(a)$$
on $[0,1]_{\varepsilon_i}$ for $i=1,2$.

For $i=1,2$, 
by Proposition~\ref{prop:invariance-of-landau-by-quasi-identity}
\begin{equation*}
f_i(\alpha_i(x))\approx_1 f_i(\alpha_i(a))+\sabun{f_i}(\alpha_i(a))(\alpha_i(x)-\alpha_i(a)) \mbox{ if } x\approx a
\end{equation*}
and, since $\alpha_i(x)-\alpha_i(a)\approx x-a$, 
Proposition~\ref{prop:basic-property-landau-symbol} implies 
\begin{equation*}
f_i(\alpha_i(x))\approx_1 f_i(\alpha_i(a))+\sabun{f_i}(\alpha_i(a))(x-a) \mbox{ if } x\approx a.
\end{equation*}

Since $f_1\circ\alpha_1\approx f_2\circ \alpha_2$, 
Proposition~\ref{prop:basic-property-landau-symbol} implies
$$
\sabun{f_{1}}(\alpha_1(a))(x-a)
\approx_1
\sabun{f_2}(\alpha_2(a))(x-a) \mbox{ if } x\approx a.
$$
Hence by Lemma~\ref{lem:係数比較}, we have
$$\sabun{f_1}(\alpha_1(a))\approx \sabun{f_2}(\alpha_2(a)).$$
\end{eproof}

Note that in the special case when $\varepsilon_1=\varepsilon_2$,
the independence can be proved directly as follows.
\begin{Lemma}
\label{lem:表現からの独立性}
Suppose $f,g$ are continuous rational valued functions on $X=[0,1]_\varepsilon$
with continuous difference quotients. If $f\approx g$, then
$$\sabun{f}\approx \sabun{g}.$$
\end{Lemma}
\begin{eproof}
It suffices to show that if $f\approx 0$ and $\sabun{f}$ is continous then
$\sabun{f}\approx 0$

Suppose $\sabun{f}\not\approx 0$. We may suppose that 
the maximum of $\sabun{f}$ is positive finite rational number. Let
$\sabun{f}(a)=r$ be one of the maxima. 
If $x\approx a$ then $\sabun{f}(x)\approx r$, whence $\sabun{f}(x)>\Frac{r}2$.
Therefore, if $k$ is huge
\begin{equation}
\label{eq:904-10}
|x-a|<\Frac1k \mbox{ implies } \sabun{f}(x)>\Frac{r}2.
\end{equation}
Hence there is an \accessible{} $n$ such that \eqref{eq:904-10} holds for $k=n$.
Let $x_1$ and $x_2$ be respectively the minimum and the maximum 

of $[a-\frac1n,a+\frac1n]_{\varepsilon}$. Then
$$x_1\approx a-\Frac1n\quad \mbox{ and}\quad x_2\approx a+\Frac1n.$$
whence by Proposition~\ref{prop:bibun-hyouka}
$$
f(x_2)-f(x_1) \geq \Frac{r}2(x_2-x_1)\approx \Frac{r}2
\Frac{2}{n}=\Frac{r}{n} \succ 0,
$$
which contradicts 
$f\approx0$.
\end{eproof}

\subsection{Infinitesimal Taylor Formula}
We fix a positive infinitesimal $\varepsilon$ in this section.
\begin{Theorem}[First order Infinitesimal Taylor formula]
\label{Th:taylor-formula}
If $f$ is a function on $[0,1]_\varepsilon$ with continuous difference quotients, 
then for $a\in [0,1]_\varepsilon^{-}$,
\begin{equation}
\label{eq:20120331-1}
f(x)\equiv_1 f(a)+\sabun{f}(a)(x-a) \mbox{ if } x\approx a.
\end{equation}
In particular 
$$ f(x)\approx_1 f(a)+\sabun{f}(a)(x-a) \mbox{ if } x\approx a.$$
\end{Theorem}
\begin{eproof}
If $x\in [0,1]_{\varepsilon}$, $x\approx a$ and $a<x$, then
  \begin{eqnarray*}
    f(x)-f(a)
    &=& \sum_{a\leq u< x}\sabun{f}(u)(\mynext{u}-u)
    \\ &=&
    \sum_{a\leq u< x}\sabun{f}(a)(\mynext{u}-u)
    +\sum_{a\leq u< x}\left(\sabun{f}(u)-\sabun{f}(a)\right)(\mynext{u}-u)
    \\ &=&
    \sabun{f}(a)(x-a)
    +\sum_{a\leq u< x}\left(\sabun{f}(u)-\sabun{f}(a)\right)(\mynext{u}-u)
  \end{eqnarray*}
If we put 
$$c:=\max_{a\leq u< x}\left|\sabun{f}(u)-\sabun{f}(a)\right|,$$
then
$$
\left|f(x)-f(a)-\sabun{f}(a)(x-a)\right|
\leq 
\sum_{a\leq u< x}c(\mynext{u}-u)=c|x-a|.
$$
By the continuity of $\sabun{f}$, we have $c\approx0$, which implies
\eqref{eq:20120331-1}.

The proof for the case $x<a$ is similar.
\end{eproof}

\begin{Theorem}
\label{th:bibunkanou-tokuchouzuke}
If a real function $F$ on $[0,1]$ is differentiable then 
for every representation $(f,[0,1]_{\varepsilon})$ 
of $F$ there is a continuous rational valued function $g$ on 
$[0,1]_{\varepsilon}$ satisfying, for each $a\in [0,1]_{\varepsilon}$,
\begin{equation}
\label{eq:907-10}
f(x)\approx_1 f(a)+ g(a)(x-a) \mbox{ if }x\approx a.
\end{equation}

Conversely if a real function $F$ on $[0,1]$ 
has a representation $(f,[0,1]_{\varepsilon})$ 
with a continuous rational valued function $g$ on $[0,1]_{\varepsilon}$
satisfying \eqref{eq:907-10} for each $a\in [0,1]_{\varepsilon}$ then 
$F$ is differentiable.
\end{Theorem}
\begin{eproof}
Suppose $F$ is differentiable and let $(f_1,[0,1]_{\varepsilon_1})$
be a representation of $F$ such that
the difference $\sabun{f_1}$ is continuous. By Theorem~\ref{Th:taylor-formula}, we have
\eqref{eq:907-10} for $a\in [0,1]_{\varepsilon_1}$ with $g=\sabun{f_1}$. 
Let $b\in [0,1]_{\rational}$. Then
\begin{equation}
\label{eq:907-11}
f_1(x)\equiv_1 f_1(\alpha_1(b))+ \sabun{f_1}(\alpha_1(b))(x-\alpha_1(b)) \mbox{ if }x\approx \alpha_1(b),
\end{equation}
whence by Proposition~\ref{prop:basic-property-landau-symbol}
$$
f_1(\alpha_1(x))\approx_1 f_1(\alpha_1(b))
+ \sabun{f_1}(\alpha_1(b))(x-b) \mbox{ if }x\approx b.$$
Since $f_1\circ\alpha_1\in F$, the relation \eqref{eq:907-10} holds
for $f=f_1\circ\alpha_1$ and $g=\sabun{f_1}\circ\alpha_1$. 

Now let $(f_2,[0,1]_{\varepsilon_2})$ be an arbitrary representation of $F$. 
By Proposition~\ref{lem:existence-of-representation}, 
there is a continuous rational valued
function $g_2$ on $[0,1]_{\varepsilon_2}$ such that 
$$g_2\circ\alpha_2\approx \sabun{f_1}\circ\alpha_1,$$
where $\alpha_i=\kappa_{\varepsilon_i}$ \itwo.
Since $f_2\circ \alpha_2\approx f_1\circ\alpha_1$, we have 

$$
f_2(\alpha_2(x))\approx_1 f_2(\alpha_2(b))
+ g_2((\alpha_2(b))(x-b) \mbox{ if }x\approx b.$$
Hence by
Proposition~\ref{prop:invariance-of-landau-by-normal-equivalence}, 
the condition \eqref{eq:907-10} holds for $f=f_2$ and $g=g_2$.

Conversely suppose $(f,[0,1]_{\varepsilon})$ is a representation of $F$ 
and there is a continuous rational valued function $g$ on $[0,1]_{\varepsilon}$ 
satisfying \eqref{eq:907-10} for each $a\in [0,1]_{\varepsilon}$. 
By Proposition~\ref{prop:approx-imply-equiv}, 
there is a subcontinuum $[0,1]_{\varepsilon'}$ of $[0,1]_{\varepsilon}$ 
such that 
\begin{equation}
\label{eq:111015-10}
f_1(y)\equiv_1 f_1(a)+ g_1(a)(y-a) \mbox{ if }y\approx a.
\end{equation}
for each $a\in [0,1]_{\varepsilon'}$, 
where $f_1:=f|[0,1]_{\varepsilon'}$ and $g_1:=g|[0,1]_{\varepsilon'}$. 
Hence 
$$
\sabun{f_1}(a)-g_1(a)=\frac{(f_1(a+\varepsilon')-f_1(a))-g_1(a)\varepsilon'}{\varepsilon'}
\approx 0
$$
for each $a\in [0,1]_{\varepsilon'}$ and $\sabun{f_1}\approx g_{1}$ is continuous. 
Then $(f_1,[0,1]_{\varepsilon'}) $ represents $F$ and have continuous difference quotients.
Hence $F$ is differentiable.
\end{eproof}

The last part of the proof shows the following Corollary which asserts
that an arbitrary representation of a differentiable function has
continuous difference quotients when restricted on a coarser but dense 
rigid mesh subcontinuum. 
\begin{Corollary}
\label{cor:restriction}
Suppose a differentiable function $F$ on $[0,1]$ is represented by
$(f,[0,1]_{\varepsilon})$. 
Then there is 
an infinitesimal $\varepsilon'\in [0,1]_{\varepsilon}$ 
and
a continuous rational valued function $g$ on $[0,1]_{\varepsilon'}$ 
such that for $a\in [0,1]_{\varepsilon'}$
\begin{equation*}
f(x)\equiv_1 f(a)+ g(a)(x-a) \mbox{ if }x\approx a.
\end{equation*}
on $[0,1]_{\varepsilon'}$.  
In particular the difference quotient of $f|[0,1]_{\varepsilon'}$ is continuous.  
\end{Corollary}

This can be rephrased as follows.
\begin{Corollary}
\label{cor:restriction2}
If a differentiable function $F$ on $[0,1]$ is represented by
$(f,[0,1]_{\varepsilon})$, 
then there is an infinitesimal $\varepsilon'\in [0,1]_{\varepsilon}$ 
such that 
the rational valued function
$$g(x):=\Frac{f(x+\varepsilon')-f(x)}{\varepsilon'}$$
on $[0,1]^{-}_{\varepsilon'}$ is continuous and represents $F'$.
\end{Corollary}
\begin{Proposition}
\label{th:char-diff}
Let $\seq{k}{F}$ be an \accessible{} number of differentiable functions on
$[0,1]$ and $\alpha:[0,1]\stackrel{\simeq}{\rightarrow} [0,1]_{\varepsilon}$ 
be a quasi-identity. Then each $F_i$ is represented by a continuous
rational valued function on $[0,1]_\varepsilon$ whose difference quotient is
continuous.
\end{Proposition}
\begin{eproof} First we represent each $F_i$ by a rational valued
continuous function on $[0,1]_{\varepsilon}$.
Using Corollary~\ref{cor:restriction} \accessible{} number of times, we
obtain a dense subset $[0,1]_{\varepsilon'}\subset [0,1]_{\varepsilon}$ 
for which the exact Taylor formula holds for each $f_i$. 
Then we extend them to functions on $[0,1]_{\varepsilon}$ by
linear interpolation. Details are omitted.
\end{eproof}

The next lemma is used in the section of inverse function theorem.
\begin{Lemma}
\label{lemma:lipshitz-condition-for-differentiable-function}
Let $f$ be a rational valued function on $[0,1]_{\varepsilon}$
with continuous difference quotient. If $\sabun{f}(a)\not\approx 0$,
then for some rationals $K_1,K_2,c\succ 0$ 
$$
K_1|x-a|\preceq |f(x)-f(a)|\preceq K_2|x-a|
$$
holds for $|x-a|<c$.
\end{Lemma}
\begin{eproof}
By Theorem~\ref{Th:taylor-formula}, we have
$$
f(x)-f(a)\equiv_1 f'(a)(x-a) \myif x\approx a,
$$
namely if $0<|x-a|\approx 0$ then
$$
\left|\frac{f(x)-f(a)}{x-a}-f'(a)\right|\approx 0.
$$
Hence if $K_1:=\frac{|f'(a)|}{4}\succ 0$, then 
for every $N\gg 1$ the condition $0<|x-a|<\frac1N$ implies 
\begin{equation}
\label{eq:20120307-a}
\left|\frac{f(x)-f(a)}{x-a}\right|>2K_1\succ K_1
\end{equation}
which, by the overspill principle, holds also if $|x-a|<c:=\frac1n$ 
for some \accessible{} $n$. Hence if $|x-a|<c$ then 
$$|f(x)-f(a)|\succ K_1|x-a|.$$
On the other hand Proposition\ref{prop:bibun-hyouka}
implies
$$
|f(x)-f(a)|<K_2|x-a| \mbox{ if } |x-a|<c
$$
if $K_2=\max\setii{\left|\sabun{f}\right|(x)}{|x-a|<\frac1n}$.
\end{eproof}

\subsection{Chain Rule}
\begin{Theorem}
\label{th:chain-rule}
Let $F_i$ \itwo be differentiable real functions on $[0,1]$ such 
that $F_1$ is $[0,1]$-valued. 
Then the composition $F_2\circ F_1$ is differentiable and satisfies
    $$(F_2\circ F_1)'=(F_2'\circ F_1)F_1'.$$
\end{Theorem}

\begin{eproof}
First we represent $F_{2}$ by 
$(f_2,[0,1]_{\varepsilon_2})$ with continuous difference quotient.

Let $(g,[0,1]_{\eta})$ represents $F_1$.
Then $(g_1:=\kappa_{\varepsilon_2}\circ g,[0,1]_{\eta})$ also represents $F$ 
and it is $[0,1]_{\varepsilon_2}$-valued. 
By Corollary~\ref{cor:restriction},
there is a positive  infinitesimal $\varepsilon_1\in [0,1]_{\eta}$ such that
$(f_1:=g_1|[0,1]_{\varepsilon_1},[0,1]_{\varepsilon_1})$
represents $F_1$ and has the continuous difference quotient.

The composition $F_2\circ F_1$ is represented by $f_2\circ f_1$.
Put $\Delta_ix=\varepsilon_i$ \itwo. 
\newcommand{\temprange}{f_1(a)\leq u<f_1(a+\Delta x)}
\newcommand{\temprangei}{x\leq w<\mynext{x}}
\newcommand{\tempfyxi}{f_1(a+\Delta_1 x)-f_1(a)}
\newcommand{\tempfyx}{\tempfyxi}
\newcommand{\sabunfii}{\frac{\Delta f_2}{\Delta_2 x}}
Suppose $f_1(x_1)<f_1(x_1+\Delta_1 x)$. Then $h=f_2\circ f_1$ satisfies
\begin{eqnarray*}
\Delta h(a)&=&h(a+\Delta_1 x)-h(a) = f_2(f_1(a+\Delta_1 x)-f_2(f_1(a))
      \\&= &
      \sum_{\temprange}\sabunfii(u)(\mynext{u}-u)
      \\ &=&
      A+\sabunfii(f_1(a))\sum_{\temprange}(\mynext{u}-u) 
      \\ &=&
      A+\sabunfii(f_1(a))(\Delta f_1(a)),
    \end{eqnarray*}
where
$$
A=\sum_{\temprange}(\sabunfii(u)-\sabunfii(f_1(a)))(\mynext{u}-u).
$$
Since $\sabunfii(u)$ is continuous, 
$\sabunfii(u)-\sabunfii(f_1(a))\approx0$ hence
by Proposition~\ref{prop:landau-integral}, $A\equiv_{1}0 $ if $x\approx a$.
Hence
$$
\frac{\Delta h}{\Delta_1 x}(a)=\frac{A}{\Delta_1 x}+ \sabunfii(f_1(a))\frac{\Delta f_1}{\Delta_1 x}(a)
\approx \sabunfii(f_1(a))\frac{\Delta f_1}{\Delta_1 x}(a),
$$
which implies that $\sabun{h}(a)$ is continuous and the derivative of $F_2\circ F_1$
is represented by $(\sabun{f_2}\circ f_1)\sabun{f_1}$, namely
$(F_2'\circ F_1)F_1'$.
\end{eproof}

\subsection{Inverse Function Theorem} 
\begin{Theorem}[Inverse Function]
\label{thm:inverse-funct-theor}
Let $F$ be a differentiable real function on $[-1,1]$
such that $F(0)=0$ and $F'(0)\neq 0$. 
Then there is a rational $c\succ 0$ and a 
diffentialble real function on $[-c,c]$ with values in $[-1,1]$
such that
\begin{itemize}
\item $F(G(y))=y$ if $y\in [-c,c]$,
\item $G(F(x))= x $ if $F(x)\in [-c,c]$,
\item $G'=\frac1{F'\circ G}$.
\end{itemize}
\end{Theorem}
It suffices to show the following.
\begin{Lemma}
\label{lem:inverse-funct-theor}
Let $f$ be a rational valued continuous function on $[-1,1]_\varepsilon$
such that $f(0)=0$ and have the continuous difference quotient and
$\sabun{f}(0)\not\approx 0$.
Then there is a rational valued function $g$ on $[-c,c]_{\varepsilon}$
for some rational $c\succ 0$ such that 
\begin{itemize}
\item $f(g(y))\approx y$ if $y\in [-c,c]_\varepsilon$,
\item $g(f(x))\approx x$ if $f(x)\in [-c,c]_\varepsilon$.
\end{itemize}
Moreover for $y,b\in [-c,c]_{\varepsilon}$
\begin{equation}
\label{eq:912-200}
g(y)-g(b) \approx_1 \Frac1{\sabun{f}(g(b))}(y-b) \myif y\approx b
\end{equation}
hence the real function $G$ represented by $g$ is differentiable 
and have the derivative $\frac1{F'\circ G}$.
\end{Lemma}

\begin{eproof}
Let $a=\sabun{f}(0)\succ 0$. Since $\sabun{f}$ is continuous, 
there is a constant $a_1\succ0$ such that 
if $|x|<a_1$ then $\sabun{f}(x)>\frac{a}2$. Hence $f$ is strictly increasing
on $[-a_1,a_1]_\varepsilon$ and by Proposition~\ref{prop:bibun-hyouka},
$$f(a_1)>c,\quad f(-a_1)<-c,$$
where $c=\frac{aa_1}2$.

For $y\in [-c,c]_\varepsilon$, define
$$
g(y):= \max\setii{x\in [-1,1]_{\varepsilon}}{ f(x)\leq y  }
$$
Then if $y\in [-c,c]_\varepsilon$
$$ f(g(y))\leq y< f(g(y)+\Delta x)\approx f(g(y)), $$
whence
$$  f(g(y))\approx y.$$
Suppose $y_1,y_2\in [-c,c]_{\varepsilon}$ 
satisfies $y_1\approx y_2$. 
Then
$f(g(y_1))\approx y_1\approx y_2\approx f(g(y_2))$ 
with $g(y_1),g(y_2)\in [-a_1,a_1]$,  
whence by Corollary~\ref{cor:inverse-continuity-of-function-with-positive-sabun}
$g(y_1)\approx g(y_2)$. Thus $g$ is continuous.

On the other hand, since $f$ is strictly increasing on $[-a_1,a_1]$, 
if $x\in [-a_1,a_1]$ satisfies $f(x)\in [-c,c]$ then
$$
f(g(f(x)))\leq f(x)\leq f(g(f(x))+\Delta x)
$$
implies
$$ g(f(x))\leq x\leq g(f(x))+\Delta x,$$
whence
$$g(f(x))\approx x.$$

Let $b_1\in [-c,c]_{\varepsilon}$.
By Lemma~\ref{lemma:lipshitz-condition-for-differentiable-function},
there is a rational $K_1,K_2,c_1\succ 0$, such that
$$
K_1|x-b_1|\preceq |f(x)-f(b_1)|\preceq K_2|x-b_1|
$$
if $|x-b_1|<c_1$. If we put $y:=f(x)$ and $b:=f(b_1)$, 
then $x\approx g(y)$ and $b_1\approx g(b)$, whence
$$
\frac1K_2|y-b|\preceq |g(y)-g(b)|\preceq \frac1{K_1}|y-b|.
$$
Then 
$$
f(x)-f(b_1) \equiv_1 \sabun{f}(b_1)(x-b_1) \myif x\approx b_1
$$
implies by virtue of Proposition~\ref{prop:invariance-of-landau-by-quasi-identity} 
that 
$$
f(g(y))-f(g(b)) \approx_1 \sabun{f}(g(b))(g(y)-g(b)) \myif y\approx b,
$$
whence \eqref{eq:912-200} since $y-b\approx f(g(y))-f(g(b))$.
\end{eproof}

\subsection{Second Order Differentiability}
\begin{Definition}
A real valued function $F$ on $[0,1]$ is called \textit{differentiable up to second order}  
if it is differentiable and its derivative is also differentiable. The derivative of $F'$ is
denoted by $F''$. We often write it also by $F^{(2)}$ or by $\frac{d^2F(x)}{dx^2}$.
\end{Definition}

\begin{Theorem}
\label{th:good-rep-of-2nd-diff-fn}
If a real valued function $F$ on $[0,1]$ is differentiable up to
second order, then for every positive infinitesimal $\varepsilon$,
$F$ has a representation $(f,[0,1]_{\varepsilon})$ with 
continuous difference quotients up to second order, 
namely, not only \sabun{f} but also its difference quotient
$\koukaisabuni{f}{2}:=\sabun{(\sabun{f})}$ is continuous and 
$(\koukaisabuni{f}{2},\alpha)$ represents $F''$.
\end{Theorem}

\begin{eproof}
By Proposition~\ref{th:char-diff}, the differentiable functions
$F$ and $F'$ have representatives $(f,[0,1]_{\varepsilon})$
and $(g,[0,1]_{\varepsilon})$.
Put 
$$\widehat{f}(x)=f(0) + \sum_{0\leq u<x,u\in X}g(u)\Delta x$$
where $\Delta x=\varepsilon$.
Then $\sabun{f}\approx g$ implies 
$$\widehat{f}(x)  \approx  
f(0)+ \sum_{0\leq u<x}\sabun{f}(u)\Delta x=f(x),$$
whence
$\widehat{f}\approx f$. Then $\widehat{f}$ is a representation 
such that both $\sabun{\widehat{f}}=g$ and $\koukaisabuni{\widehat{f}}{2}=\sabun{g}$ 
are continuous.
The last statement is obvious.
\end{eproof}

\begin{Theorem}
\label{thm:second-order-taylor}
Let $F$ be a real valued function on $[0,1]$ differentiable up to second order.
Let $(f,[0,1]_{\varepsilon})$ be a representation of $F$ with continuous 
$\sabun{f}$ and $\koukaisabuni{f}{2}$.  Then for $a\in [0,1]_{\varepsilon}^{--}$
\begin{equation}
\label{eq:908-90}
    f(x)\approx_2
    f(a)
    +\sabun{f}(a)(x-a)
    +\koukaisabuni{f}{2}(a) \frac{(x-a)^2}2 \myif x\approx a.
\end{equation}

\end{Theorem}
\newcommand{\kukanN}{[0,1]_{\varepsilon}}

\begin{eproof}
 By Theorem~\ref{Th:taylor-formula}, we have for $a\in \kukanN$,
\begin{equation}
\label{eq:908-50}
\sabun{f}(u)\equiv_1 \sabun{f}(a)+\koukaisabuni{f}{2}(a)(u-a)
\myif u\approx a 
\end{equation}
Assume $x>a$. By substituting this into 
\begin{equation}
\label{eq:908-30}
 f(x)=f(a)+\sum_{a\leq u<x}\sabun{f}(u)\Delta x,
\end{equation}
with $\Delta x=\varepsilon$, we obtain by Proposition~\ref{prop:landau-integral}
\begin{eqnarray*}
    f(x)&\equiv_2&
    f(a)
    +\sabun{f}(a)(x-a)
    +\koukaisabuni{f}{2}(a) \sum_{a\leq u<x}(u-a)\Delta x 
\myif x\approx a \mbox{ and $x>a$}
\end{eqnarray*}
since $|x-a|\geq |u-a|$ if $a\leq u < x$ 

Put $M:=\Frac{|x-a|}{\Delta x}$. Then $u\in \kukanN$ with $a\leq u<x$ is
written as $u=a+i\Delta x$ with $i=\frac{u-a}{\Delta x}\in [0..M]$. Hence
\begin{eqnarray*}
\sum_{a\leq u<x}(u-a)\Delta x
 &
=
 & 
\sum_{0\leq i\leq M}i(\Delta x)^2=\Frac{M(M+1)}{2}(\Delta x)^2
=\Frac{(x-a)^2}2+|x-a|\frac{\Delta x}2.
\end{eqnarray*}
Suppose $0\approx |x-a|\geq \sqrt{\Delta x}$.
Then 
$$\Frac{\frac{|x-a|\Delta x}2}{|x-a|^2}
\leq \frac12\sqrt{\Delta x}\approx 0.
$$
Hence 
\begin{eqnarray*}
    f(x)&\approx_2&
    f(a)
    +\sabun{f}(a)(x-a)
    +\koukaisabuni{f}{2}(a) \frac{(x-a)^2}2 \myif x\approx a \mbox{ and $x>a$}.
\end{eqnarray*}
The proof in the case when $x<a$ is similar.

\end{eproof}

\begin{Corollary}
\label{cor:second-order-taylor}
Let $F$ be a real valued function on $[0,1]$ differentiable up to second order.
Let $(g,[0,1]_{\varepsilon})$ be its representation. 
Let $g_1$ and $g_2$ be continuous rational valued function on $[0,1]_{\varepsilon}$ 
representing $F'$ and $F''$. Then for $a\in [0,1]_{\varepsilon}$
\begin{equation}
\label{eq:908-100}
 g(x)\approx_2g(a)+g_1(a)(x-a)+g_2(a)\Frac{(x-a)^{2}}{2}\myif x\approx a.
\end{equation}
\end{Corollary}
\begin{eproof}
By Theorem \ref{th:good-rep-of-2nd-diff-fn}, there is a representation
$(f,[0,1]_{\delta})$ of $F$ with continuous $\sabun{f}$ and $\koukaisabuni{f}{2}$. 
By Theorem~\ref{thm:second-order-taylor}, we have \eqref{eq:908-90} on $[0,1]_{\delta}$.
Put $\gamma:=\kappa_{\delta}|[0,1]_{\varepsilon}$ which is a quasi-identity
from $[0,1]_{\varepsilon}$ to $[0,1]_{\delta}$. 
Then by Proposition~\ref{prop:invariance-of-landau-by-normal-equivalence},
we have
\begin{equation}
    f(\gamma(x))\approx_2
    f(\gamma(a))
    +\sabun{f}(\gamma(a))(\gamma(x)-\gamma(a))
    +\koukaisabuni{f}{2}(\gamma(a)) \frac{(\gamma(x)-\gamma(a))^2}2 \myif x\approx a.
\end{equation}
Since $f\circ\gamma \approx g$, $\sabun{f}\circ\gamma\approx g_1$,
and $\koukaisabuni{f}{2}\circ\gamma\approx g_2$, 
we obtain 
we obtain \eqref{eq:908-100} by Proposition~\ref{prop:basic-property-landau-symbol}.
\end{eproof}

\begin{Proposition}
Suppose a real valued function $F$ has representations
$(f_i,[0,1]_{\varepsilon_i})$ \itwo 
with continuous difference quotients up to second order.
Then $\koukaisabuni{f_{1}}{2}$ and $\koukaisabuni{f_{2}}{2}$ 
represent one and the same real function. In particular,
the second derivative $F''$ is represented by 
$\koukaisabuni{f}{2}$ of any representation $(f,[0,1]_{\varepsilon})$ 
of $F$ with continuous difference quotients up to second order.
\end{Proposition}
\begin{eproof}
By Theorem~\ref{thm:second-order-taylor}, 
\begin{equation}
\label{eq:3-715}
f_i(x)-f_i(a)-\sabun{f_i}(a)(x-a)-\koukaibibun{f_i}{2}(a)\Frac{(x-a)^{2}}{2}\approx_20
\myif x\approx a
\end{equation}
on $[0,1]_{\varepsilon_i}$ \itwo. 

Let $\beta:=\kappa_{\varepsilon_2}|[0,1]_{\varepsilon_1}$ 
be the quasi-identity from
$[0,1]_{\varepsilon_1}$ to $[0,1]_{\varepsilon_2}$. 
Then by
Propositions~\ref{prop:basic-property-landau-symbol} and 
Proposition~\ref{prop:invariance-of-landau-by-normal-equivalence},
\begin{equation}
\label{eq:3-715a}
f_2(\beta(x))-f_2(\beta(a))-\sabun{f_2}(\beta(a))(x-a)
-\koukaibibun{f_2}{2}(\beta(a))\Frac{(x-a)^{2}}{2}\approx_20
\myif x\approx a
\end{equation}

Since $f_2\circ \beta\approx f_1$, 
Propositions~\ref{prop:basic-property-landau-symbol}
implies
\begin{equation}
\left(\sabun{f_{2}}(\beta(a))-\sabun{f_1}(a)\right)(x-a)
+
\left(\koukaibibun{f_2}{2}(\beta(a))-\koukaibibun{f_1}{2}(a)\right)
\Frac{(x-a)^{2}}{2}\approx_20
\myif x\approx a.
\end{equation}
Hence by Lemma~\ref{lem:係数比較}, 
$$
\koukaibibun{f_2}{2}(\beta(a))\approx
\koukaibibun{f_1}{2}(a).
$$
Since $\beta\circ\kappa_{\varepsilon_1}\approx \kappa_{\varepsilon_2}$,
$$
\koukaibibun{f_2}{2}\circ\kappa_{\varepsilon_2}
\approx
\koukaibibun{f_2}{2}\circ \beta \circ \kappa_{\varepsilon_1}
\approx
\koukaibibun{f_1}{2}\circ\kappa_{\varepsilon_1}.
$$
\end{eproof}
\begin{Theorem}[Characterizaion of second order differentiability]
\label{th:second-order-diff}
A real function $F$ on $[0,1]$ is differentiable up to second order 
if it has a representation $(f,[0,1]_{\varepsilon})$ 
and there are continuous rational valued functions 
$g_1$ and $g_2$ on $X$ satisfying for each $a\in \kukanN$
\begin{equation}
\mylabeleq{eq:11-909}
 f(x)\approx_2 f(a)+g_1(a)(x-a)+g_2(a)\Frac{(x-a)^{2}}{2} \myif x\approx a.
\end{equation}
Moreover $F'$ and $F''$ are represented respectively by $g_1$ and $g_2$.
\end{Theorem}

\begin{eproof}
By Lemma~\ref{prop:approx-imply-equiv}, 
there is a positive infinitesimal $\delta\in [0,1]_{\varepsilon}$ such that
on $[0,1]_{\delta}$ 
\begin{equation}
\label{eq:12-909-20}
 f(x)\equiv_2 f(a)+g_1(a)(x-a)+g_2(a)\Frac{(x-a)^{2}}{2} \myif x\approx a
\end{equation}
Substituting $x=a+\delta$ 
we have
$$\sabun{f}(a)\approx g_1(a).$$ 
Substituting $x=a+2\delta$ 
and using
$$\koukaisabuni{f}{2}(a)
=\Frac{f(a+2\delta)-2f(a+\delta)+f(a)}{\delta^2},$$
we obtain
$$
\koukaisabuni{f}{2}(a)\approx g_2(a).
$$
\end{eproof}

\subsection{Higher Order Differentiability}
\begin{Definition}
  Let $k\geq 3$ be an \accessible{} number. 
A real function $F$ on $[0,1]$ is
  differentiable up to order $k$ 
  if it is differentiable up to order
  $k-1$ and its $k-1$-th derivative is differentiable. 
  The derivative of its $k-1$-th derivative is called its $k$-th derivative and is denoted by $F^{(k)}$ and $\frac{d^kF(x)}{dx^k}$.
\end{Definition}
\index{differentiable, up to order $k$}
Theorem \ref{th:good-rep-of-2nd-diff-fn} extends to general order.
\begin{Theorem}
\label{th:good-rep-of-kth-diff-fn}
Let $\varepsilon$ be a positive infinitesimal and $k$ an \accessible{} number.
If a real function $F$ on $[0,1]$ is differentiable up to
$k$-th order, then $F$ has a representation $(f,[0,1]_{\varepsilon})$ 
with continuous difference quotients up to $k$-th order, 
namely, the higher order difference quotients defined inductively by
$\koukaisabuni{f}{i}:=\sabun{(\koukaisabuni{f}{i-1})}$ 
is continuous and represents $F^{(i)}$ for $i\leq k$.
\end{Theorem}
\begin{eproof}
By Proposition~\ref{th:char-diff}, the real functions $F^{(i)}$
have representatives $(g_i,[0,1]_{\varepsilon})$ for $i\in[0..k]$.
Define $\widehat{g}_{k-i}$ for $i\in [0..k]$ inductively by $\widehat{g}_{k}=g_k$
and, for $i\geq 1$,
$$
\widehat{g}_{k-i}(x)=g_{k-i}(0)+\sum_{0\leq u< x}\widehat{g}_{k-i+1}(u)\Delta x.
$$
Then $\widehat{g}_i\approx g_i$ for $i\in [0..k]$ and if $i<k$
$$\sabun{\widehat{g}_{i}}=\widehat{g}_{i+1},$$
whence $\koukaisabuni{\widehat{g}_0}{i}=\widehat{g}_i\approx g_i$ is continuous.
\end{eproof}

\begin{Theorem}[Taylor formula]
\label{th:high-order-taylor}
Let $F$ be a real function on $[0,1]$ differentiable up to $k$-th order with
\accessible{} $k$. Let $(f,[0,1]_{\varepsilon})$ 
be a representation with continuous $\koukaisabuni{f}{i}$ for $i\in [1..k]$. 
Then for $a\in X$
\begin{equation}
\label{eq:909-taylor}
    f(x)\approx_k
    f(a)
    +\sum_{i=1}^{k}\koukaisabuni{f}{i}(a)\frac{(x-a)^i}{i!} \myif x\approx a.
\end{equation}
\end{Theorem}
\begin{eproof}
By induction on $\ell\in [1..k]$, we can show
\begin{equation}
\label{eq:909-10}
f(x)-f(a)\approx_k
\sum_{i=1}^{\ell-1}\koukaisabuni{f}{i}(a)\Frac{(x-a)^i}{i!}
+\sum_{a\leq u_\ell\leq \cdots\leq u_2\leq u_1< x}\koukaisabuni{f}{\ell}(u_\ell)(\Delta x)^{\ell}
\myif x\approx a.
\end{equation}
In fact, for $\ell=1$, this is essentially \eqref{eq:sum-of-difference}.
Suppose \eqref{eq:909-10} holds for $\ell\leq t$. 
\begin{eqnarray*}
&&
\sum_{a\leq u_t\leq \cdots\leq u_2\leq u_1< x}\koukaisabuni{f}{t}(u_{t})(\Delta x)^{t}
\\
&=&
\sum_{a\leq u_t\leq \cdots\leq u_2\leq u_1< x}
\left(
\koukaisabuni{f}{t}(a)
+
\sum_{a\leq u_{t+1}\leq u_{t}}\koukaisabuni{f}{t+1}(u_{t+1})\Delta x
\right)
(\Delta x)^{t}
\\
&=&
\koukaisabuni{f}{t}(a)\sum_{a\leq u_t\leq \cdots\leq u_2\leq u_1< x}(\Delta x)^{t}\\
&&+
\sum_{a\leq u_{t+1}\leq u_t\leq \cdots\leq u_2\leq u_1< x}\koukaisabuni{f}{t+1}(u_{t+1})(\Delta x)^{t+1}
\\
&\approx_k& \koukaisabuni{f}{t}(a)\frac{(x-a)^{t}}{t!}
+
\sum_{a\leq u_{t+1}\leq u_t\leq \cdots\leq u_2\leq u_1< x}\koukaisabuni{f}{t+1}(u_{t+1})(\Delta x)^{t+1},
\end{eqnarray*}
by Lemma~\ref{lem:914-f} below. Hence \eqref{eq:909-10} holds for $\ell=t+1$. 

Finally we calculate the last term of \eqref{eq:909-10} with $\ell=k$.  
Since $\koukaisabuni{f}{k}$ is continuous,
\begin{eqnarray*}
\sum_{a\leq u_{k}\leq \cdots\leq u_2\leq u_1< x}\koukaisabuni{f}{k}(u_{k})(\Delta x)^{k}
&\approx&
\sum_{a\leq u_{k}\leq \cdots\leq u_2\leq u_1< x}\koukaisabuni{f}{k}(a)(\Delta x)^{k}\\
&\approx&
\koukaisabuni{f}{k}(a)\Frac{(x-a)^{k}}{k!}
\end{eqnarray*}
by Lemma~\ref{lem:914-f} below.
Hence by Proposition~\ref{prop:basic-property-landau-symbol}
$$
\sum_{a\leq u_{k}\leq \cdots\leq u_2\leq u_1< x}\koukaisabuni{f}{k}(u_{k})(\Delta x)^{k}
\approx_k
\koukaisabuni{f}{k}(a)\Frac{(x-a)^{k}}{k!} \myif x\approx a
$$
\end{eproof}

\begin{Lemma} 
\label{lem:914-f}
Let $\ell$ be an \accessible{} number and $a,x\in \kukanN$ with $a\approx x$. 
If $\frac{|a-x|}{\varepsilon}$ is huge then
\begin{equation}
\label{eq:909-111}
\sum_{a\leq u_\ell\leq \cdots\leq u_2\leq u_1< x}(\Delta x)^{\ell}\approx\Frac{(x-a)^{\ell}}{\ell!}
\end{equation}
where $\Delta x=\varepsilon$.
In particular, by Proposition~\ref{prop:basic-property-landau-symbol}, for any
\accessible{} number $k$
\begin{equation}
\label{eq:909-111-2}
\sum_{a\leq u_\ell\leq \cdots\leq u_2\leq u_1< x}(\Delta x)^{\ell}\approx_k\Frac{(x-a)^{\ell}}{\ell!}
 \myif x\approx a,
\end{equation}
\end{Lemma}
\begin{eproof}
Put $L=\frac{x-a}{\Delta x}$ and assume $L\gg 1$. Then
\begin{eqnarray*}
\sum_{a\leq u_\ell\leq \cdots\leq u_2\leq u_1< x}(\Delta x)^{\ell}
&=&{}^{\#}\setii{(\seq{\ell}{i}}{0\leq i_1\leq \cdots \leq i_{\ell}<L}(\Delta x)^{\ell}
\\
&=& \binom{L}{\ell}(\Delta x)^{\ell}
= \Frac{L^{\ell}}{\ell!}
\left(1-\frac1L\right)
\left(1-\frac2L\right)
\left(1-\frac{\ell-1}L\right)(\Delta x)^{\ell}
\\
&\approx&
\Frac{L^{\ell}}{\ell!}(\Delta x)^{\ell}=\Frac{(x-a)^{\ell}}{\ell!}.
\end{eqnarray*}
\end{eproof}

\begin{Remark}
The infinitesimal Taylor formula gives usual one by virtue
of Proposition~\ref{prop:landau-symbol-in-other-words}.
\end{Remark}

Corollary~\ref{cor:second-order-taylor} extends to general order. The proof 
is similar and omitted.
\begin{Corollary}
\label{cor:higer-order-taylor}
Let $F$ be a real function on $[0,1]$ differentiable up to $k$-th order
with a representation $(g,[0,1]_{\varepsilon}$.
Let $g_i$ ($1\leq i\leq k$) be continuous rational valued functions on $[0,1]_{\varepsilon}$ 
representing $\frac{d^iF}{dx^i}$. Then for $a\in [0,1]_{\varepsilon}$
\begin{equation}
\label{eq:910-100}
 g(x)\approx_k g(a)+\sum_{i=1}^kg_i(a)\Frac{(x-a)^i}{i!}\myif x\approx a.
\end{equation}
\end{Corollary}

Theorem~\ref{th:second-order-diff} also holds for higher order differentiability
and proved similarly.
\begin{Theorem}[Characterizaion of higher order differentiability]
\label{th:highr-order-diff-characterization}
A real function $F$ on $[0,1]$ is differentiable up to $k$-th order if
it has a representation $(f,[0,1]_{\varepsilon})$ 
and there are continuous rational valued functions
$g_i$ ($i\in [1..k]$) on $[0,1]_{\varepsilon}$ 
satisfying \eqref{eq:910-100} for each $a\in [0,1]_{\varepsilon}$.
Then $g_i$ represents $\frac{d^iF}{dx^i}$ for $i\in [1..k]$.
\end{Theorem}

\subsection{Fundamental Theorem of Calculus}
\def\pcont{continuous}
\begin{Definition}
Let $F$ be a real function on $[0,1]$ with
a representation $(f,\kukanN)$.
The rational valued function $\Sigma f\Delta x$ on $\kukanN$ defined by 
$$
(\Sigma f\Delta x)(u):= \sum_{0}^{u}f\Delta x=\sum_{0\leq x< u}f(x)\Delta x,
$$
where $\Delta x=\varepsilon$, is continous and finite. 
The real function represented by $\Sigma f\Delta x$ is called the \textit{indefinite integral of} $F$ and is
written as $\int_0^tF(x)dx$.
\end{Definition}
\index{indefinite integegral}
This definition is legitamate since 
$\Sigma f\Delta x$ is continuous and finite 
by Proposition~\ref{prop:well-definedness-of-indefinite-integral}
and $\int_0^tF(x)dx$ does not depend on the representation of $F$
by Proposition~\ref{prop:indefinite-integra:independence-of-representation}.

\begin{Proposition}
\label{prop:well-definedness-of-indefinite-integral}
If $f$ is a continuous rational valued function on $[0,1]_{\varepsilon}$,
then the rational valued function $\Sigma f\Delta x$ on $[0,1]_{\varepsilon}$
is \pcont{} and finite. 
\end{Proposition}
\begin{eproof}
Let $M$ be the maximum of $|f|$, which is finite by assumption. 
If $a,b\in [0,1]_{\varepsilon}$ and $a<b$,
then
$$
|(\Sigma f\Delta x)(b)-(\Sigma f\Delta x)(a)|
\leq 
\sum_{a\leq  x<  b}|f(x)|\Delta x
\leq M\sum_{a\leq  x< b}\Delta x=M|b-a|.
$$
Hence $\Sigma f\Delta x$ is \pcont{}. It is finite since
$$
|(\Sigma f\Delta x)(u)|\leq M\sum_0^u\Delta x=Mu.
$$

\end{eproof}

\begin{Lemma}
\label{lem:extension} Suppose $\varepsilon_i$ \itwo are positive infinitesimals 
such that $\varepsilon_1 \in [0,1]_{\varepsilon_{2}}$.
For a \pcont{} rational valued function $g$ on $[0,1]_{\varepsilon_1}$, define
$\tilde{g}:=g\circ\kappa_{\varepsilon_1}$ on $[0,1]_{\varepsilon_2}$, 
namely, put
$$\tilde{g}(x):=g([x/\varepsilon_1]\varepsilon_1).$$
Then $\tilde{g}$ is \pcont{} and on $[0,1]_{\varepsilon_2}$ 
\begin{equation}
\label{eq:20120331-10}
\Sigma \tilde{g}\Delta x\approx (\Sigma g\Delta x)\circ\kappa_{\varepsilon_1}.
\end{equation}
In particular $\Sigma g\Delta x$ 
and $\Sigma \tilde{g}\Delta x$ represent 
one and the same real function on $[0,1]$.
\end{Lemma}
\begin{eproof}
The continuity of $\tilde{g}$ is obvious.

For $u\in [0,1]_{\varepsilon_2}$, put $u_{-}:=\kappa_{\varepsilon_1}(u)=[u/\varepsilon_1]\varepsilon_1$.
Then
\begin{eqnarray*}
(\Sigma \tilde{g}\Delta x)(u)
&=&\sum_{0\leq t<u,t\in [0,1]_{\varepsilon_2}}\tilde{g}(t)\varepsilon_{2}\\
&=&\sum_{0\leq t<u_{-},t\in [0,1]_{\varepsilon_2}}\tilde{g}(t)\varepsilon_2+ \sum_{u_{-}\leq t<u,t\in [0,1]_{\varepsilon_2}}\tilde{g}(t)\varepsilon_2\\
&=&\sum_{0\leq t<u_{-},t\in [0,1]_{\varepsilon_1}}g(t)\varepsilon_1+ \sum_{u_{-}\leq t<u}\tilde{g}(t)\varepsilon_2
\approx(\Sigma g\Delta x)(u_{-}),
\end{eqnarray*}
whence \eqref{eq:20120331-10}.

\end{eproof}

\begin{Proposition}
\label{prop:indefinite-integra:independence-of-representation}
If $(f_i,[0,1]_{\varepsilon_i})$ \itwo 
represent a real function $F$ on $[0,1]$,
then $ \Sigma f_1\Delta x$ and $\Sigma f_2\Delta x$ represent one and the same real function on $[0,1]$.
\end{Proposition}
\begin{eproof}
Let $\varepsilon$ be an infinitesimal such that $\varepsilon_i\in [0,1]_{\varepsilon}$
\itwo. For example, if $\varepsilon_i=\frac{p_i}{q_i}$, then
one may take $\varepsilon=\frac1{q_1q_2}$.

Then $F$ is represented by $\tilde{f}_i:= f_i\circ \kappa_{\varepsilon_i}|[0,1]_{\varepsilon}$ \itwo 
and by Proposition~\ref{prop:Sigmaf}
\begin{equation}
\label{eq:20120305}
\Sigma \tilde{f}_1\Delta x\approx \Sigma \tilde{f}_2\Delta x,
\end{equation}
and by Lemma~\ref{lem:extension}, for $i=1,2$, 
\begin{equation}
\label{eq:20120331-20}
\Sigma \tilde{f}_i\Delta x \approx (\Sigma f_i \Delta x)\circ \kappa_{\varepsilon_i}.
\end{equation}
Since $\kappa_{\varepsilon_i}\circ\kappa_{\varepsilon}=\kappa_{\varepsilon_i}$ \itwo, for
$i=1,2$ we have
$$
\Sigma f_i\Delta x\circ\kappa_{\varepsilon_{1}}
=
\Sigma f_i\Delta x\circ\kappa_{\varepsilon_{1}}\circ\kappa_{\epsilon}
\approx
\Sigma \tilde{f}_i\Delta x\circ\kappa_{\varepsilon}.
$$
Hence by \eqref{eq:20120305},
$$\Sigma f_1\Delta x\circ\kappa_{\varepsilon_{1}}
\approx \Sigma f_2\Delta x\circ\kappa_{\varepsilon_{2}},$$
namely
$\Sigma f_i\Delta x$ \itwo represent one and the same real function
on $[0,1]_{\rational}$.

\end{eproof}

Generally, the indefinite integral $\int_{a}^{t}F(x)dx$ is defined 
for a real function $F$ on general interval $[a,b]$.

\begin{Proposition}
Suppose $F$ is a real function on $[0,1]$. Then the real function 
$\int_{0}^tF(x)dx$ on $[0,1]$ is differentiable and its derivative is $F$.
\end{Proposition}
\begin{eproof}
  The indefinite integral $\int_{0}^tF(x)dx$ is represented by 
  $\Sigma f\Delta x$ using a representaiton 
$(f,[0,1]_{\varepsilon})$ of $F$. 
Its difference quotient is 
$$\frac{\Delta(\Sigma f\Delta x)(x)}{\Delta x}
=\Frac{(\Sigma f\Delta x)(x^+)-(\Sigma f\Delta x)(x)}{\Delta x}=f(x),$$
where $\Delta x=\varepsilon$. 
Thus the indefinite integral has the representation $\Sigma f\Delta x$ 
with the continuous difference quotients $f$, 
whence is differentiable and its derivative is $F$.
\end{eproof}

\subsection{Ordinary Differential Equation}

\begin{Theorem}                 
\label{theorem:ODE-solution}
Let $K$ be a finite positive real number and 
$F$ a real function on $A=[0,1]\times [-2K,2K]$ 
satisfying
$$|F(x,y)| < K$$
for all $(x,y)\in A$. 
Further suppose that there is a finite real number $L$ such that
$$
|F(x,y_1)-F(x,y_2)|\leq L|y_1-y_2|
 $$
holds for all $x\in [0,1]$ and $y_i\in [-2K,2K]$  \itwo.
Then there is a unique real function $G(x,a)$ on $[0,1]\times [-K,K]$ satisfying 
$$
G(0,a)=a, \quad \frac{dG(x,a)}{dx}=F(x,G(x,a)). 
$$
\end{Theorem}

This is equivalent to the following statement.
\begin{Proposition}
Let $\varepsilon>0$ be an infinitesimal and $K\prec 0$ a finite rational number.
If $f$ is a rational valued continuous function on 
$X=\kukanN\times \kukaniie{2K}{\varepsilon}$ 
satisfying
$$ |f(x,y)|\prec K,$$
and there is a finite rational number $L$ such that 
$$
|F(x,y_1)-F(x,y_2)|\preceq L|y_1-y_2|
 $$
for all $x\in [0,1]_{\varepsilon}$ and $y_i\in [-2K,2K]_{\varepsilon}$  \itwo.

Then there is a rational valued continuous function $\varphi$ on 
$Y=\kukanN\times \kukaniie{K}{\varepsilon}$ satisfying
$$ \varphi(0,a)\approx a$$
for all $a\in \kukaniie{K}{\varepsilon}$ and 
\begin{equation}
\label{eq:9-721}
\Frac{\varphi(x+\Delta x,a)-\varphi(x,a)}{\Delta x}\approx f(x,\varphi(x,a))  
\end{equation}
for all $(x,a)\in \kukanN\times \kukaniie{K}{\varepsilon}$.

Furthermore if another continuous rational valued function $\psi$ on
$\kukanN\times \kukaniie{K}{\varepsilon}$
satisfies (\ref{eq:9-721}) with $\varphi$ replaced by $\psi$, then
$\varphi\approx \psi$.
\end{Proposition}

\begin{eproof}
For $x\in \kukanN$ and $a\in \kukaniie{K}{\varepsilon}$, 
we can define $\varphi(x,a)$ by ``induction on $x$'' as follows.
$$\varphi(0,a):=a,$$
$$ \varphi(x+\Delta x,a):=\varphi(x,a)+f(x,\varphi(x,a))\Delta x,$$
where $\Delta x=\varepsilon$.

In fact, suppose we have defined $\varphi(x,a)$ for $x\leq b$
for some $b<1$. Then for $0\leq x\leq b$
$$\left|\Frac{\Delta_x \varphi(x,a)}{\Delta x}\right|=|f(x,\varphi(x,a))|\prec K,$$
where $\Delta_x \varphi(x,a)=\varphi(x+\Delta x,a)-\varphi(x,a)$.
Hence by Proposition~\ref{prop:bibun-hyouka}
$$ |\varphi(b,a)-\varphi(0,a)|\prec bK < K. $$
Hence 
$$ |\varphi(b,a)|\prec K+a< 2K$$
and $f(b,\varphi(b,a))$ has value and $\varphi(b+\Delta x,a)$ is defined.

From $\left|\Frac{\Delta_x\varphi}{\Delta x}\right|\prec K$, it follows
$$   |\varphi(x,a)-\varphi(y,a)|\prec K|x-y|.$$
Hence $\varphi$ is \pcont{}  with respect to the first variable and
the partial difference quotient with respect to $x$
$$
\Frac{\Delta_x \varphi(x,a)}{\Delta x}=f(x,\varphi(x,a))
$$
is also continuous with respect to $x$. 

To show the continuity of $\varphi$ with respect to the second variable, 
put $h(x):=\varphi(x,a)-\varphi(x,b)$. 
Then 
$$\left|
\Frac{h(x+\Delta x)-h(x)}{\Delta x}
 \right|
= |f(x,\varphi(x,a))-f(x,\varphi(x,b))|
\preceq L|h(x)|.
$$
Hence by Lemma~\ref{lemma:ode-unique} below,
$$h(x)\preceq h(0)\exp(Lx)=|a-b|\exp(Lx).$$
Hence $a\approx b$ implies $\varphi(x,a)\approx \varphi(x,b)$.
Thus it is verified that $\varphi$ is continuous.

Let $\psi$ be a rational valued function on
$\kukanN\times \kukaniie{K}{\varepsilon}$
satisfying (\ref{eq:9-721}) with $\varphi$ replaced by $\psi$.

Fix $a\in \kukaniie{K}{\varepsilon}$ and 
put $h(x)=\varphi(x,a)-\psi(x,a)$. 
Then we can show similarly
\begin{equation}
\label{eq:1105-1}
\left|\;
\frac{\Delta h(x)}{\Delta x}
\;\right|
\preceq L|h(x)|
\quad
\mbox{for all $x\in \kukanN$},
\end{equation}
whence by the following Lemma~\ref{lemma:ode-unique}, 
\begin{equation}
\label{eq:1105-2}
|h(x)|\preceq |h(0)|\exp(Lx).
\end{equation}
Since $|h(0)|\approx 0$ and $\exp(xL)$ is bounded, we obtain
$|h(x)|\approx 0$ for all $x\in \kukanN$.
\end{eproof}

\begin{Lemma}
\label{lemma:ode-unique}
If a rational valued continuous function $h$ on $[0,1]$ satisfies
\eqref{eq:1105-1}, then \eqref{eq:1105-2} holds.
\end{Lemma}
\begin{eproof}
Put
$$
\delta:=\max\seti{
0,
\max\setii{
\left|\;\frac{\Delta h(x)}{\Delta x}
\;\right|-L|h(x)|}
{x\in \kukanN}
}\approx 0.
$$
Then
$$
\left|\;
\frac{h(x+\Delta x)-h(x)}{\Delta x}
\;\right|
\leq L|h(x)|+\delta.
$$
Hence
$$
|h(x+\Delta x)|\leq (1+L\Delta x)|h(x)|+\delta\Delta x,
$$
which is equivalent to
$$
|h(x+\Delta x)|+\frac\delta{L} \leq 
\left(\;
|h(x)|+\frac\delta{L}\;\right)(1+L\Delta x).
$$
Hence 
\begin{eqnarray*}
|h(x)|+\frac\delta{L}\leq 
\left(\;h(0)+\frac\delta{L}\;\right)
(1+L\Delta x)^{\frac{x}{\Delta x}}
\end{eqnarray*}
If we put $T=\frac{x}{\Delta x}$, then $\frac1T\approx 0$ whenever $x\not\approx 0$ 
and by Proposition\ref{prop:product-formula-exp}, we have
$$
(1+L\Delta x)^{\frac{x}{\Delta x}}=
\left(\;1+\frac{xL}{T}\;\right)^{T}
\approx 
\exp(xL).
$$
Hence 
$$
|h(x)|\approx |h(x)|+\frac\delta{L}
\preceq
\left(\;
|h(0)|+\frac\delta{L}
\;\right)
\exp(xL)\approx |h(0)|\exp(xL)
$$
when $x\succ 0$, which implies (\ref{eq:1105-2}) by Proposition~\ref{prop:daishou}.
\end{eproof}

\newpage

\section{Differentiation: Multiple Variables}
\label{sec:calculus-multiple}
\def\runningtitle{Differentiation: Multiple variables}
Let $n$ be an \accessible{} number. 
We consider real functions only on the continuum $[0,1]^n$ for simplicity 
but nothing changes essentially for general continuum of the form $\prod_i[a_i,b_i]$.

We say a real function $F$ on $[0,1]^n$ is represented by a 
continuous rational valued function $f:[0,1]_{\varepsilon}^n\rightarrow \rational$
if $f\circ\kappa_{\varepsilon}\in F$. We say then that $(f,[0,1]_{\varepsilon})$ is
a representation of $F$.

\subsection{Partial Difference Quotients}
Let $f$ be a rational valued function on $X=\kukanN^{n}$ with \accessible{} $n$.
Define the partial differences of $f$ by 
$$ \Delta_if(x):= f(x+\Delta x\;\be_i)-f(x)$$
for $x_i+\Delta x\leq 1$ and $\Delta_if(x):=\Delta_i f(x-\Delta x \be_i)$ otherwise.
Here $\Delta x:= \varepsilon$ 
and $\be_i$ the is the $n$-vector with the $i$-th component $1$ and other components $0$.
The quotient $\hensabun{f}{i}(a):=\Delta_if(a)/\Delta x$ is called 
\textit{the $i$-th partial difference quotient at} $a$ and 
$\hensabun{f}{i}$ is a rational valued function on $X$

\begin{Lemma}
\label{lemma:differene-quotient}
If $f$ is a rational valued function on $X=\kukanN^{n}$ and $a\in X$, then
$$
f(x)-f(a)=\varepsilon\sum_{i=1}^n\sum_{u\in I_i(a,x)}sgn(x_i-a_i)\hensabun{f}{i}(x[i-1]+u\be_i),
$$
where
$$x[i]=(x_1,\cdots,x_i,a_{i+1},\cdots,a_n)$$
with $x[0]=a$,
$$
I_i(a,x)= \left\{
  \begin{array}{cc}
\setii{u}{a_i\leq u< x_i}    & \mbox{if $a_i\leq x_i$ } \\
\setii{u}{x_i\leq u< a_i}     & \mbox{if $x_i<a_i$ } 
  \end{array}
\right.
$$
and $sgn$ is the signature function defined by
$$
sgn(r)= \left\{
  \begin{array}{cc}
\frac{r}{|r|}    & \mbox{for $r\neq 0$} \\
0    & \mbox{otherwise} 
  \end{array}
\right.
$$
\end{Lemma}
\begin{eproof}
Obviously we have
$$
f(x)-f(a)=\sum_{i=1}^{n}\left( f(x[i])-f(x[i-1]) \right).
$$
Put $g(t)=f(\seq{i-1}{x},t,a_{i+1},\cdots,a_n)$. Then
$g(a_i)=f(x[i-1])$ and $g(x_i)=f(x[i])$ and
$\sabun{g}(u)=\hensabun{f}{i}(x[i-1]+u\be_i)$. 
If $a_i\leq x_i$
$$
g(x_i)=g(a_i)+\sum_{a_i\leq u< x_i} \sabun{g}(u)\varepsilon.
$$
If $x_i<a_i$, then
$$
g(a_i)=g(x_i)+\sum_{x_i\leq u< a_i} \sabun{g}(u)\varepsilon,
$$
whence
$$
g(x_i)=g(a_i)+\sum_{u\in I_i(a,x)}sgn(x_i-a_i)\sabun{g}(u)\varepsilon.
$$
\end{eproof}

\begin{Proposition}
\label{prop:diff-quot-estimate}
If $f$ is a rational valued function on $X=\kukanN^{n}$ and
$$|\hensabun{f}{i}(x)|<M$$
for all $i\in [1..n]$ and $x\in X$. 
Then
$|f(x)-f(a)|\leq nMd(x,a)$
for all $x,a\in X$.
\end{Proposition}
\begin{eproof}
By Lemma~\ref{lemma:differene-quotient}, we have
$$|f(x)-f(a)|
\leq \varepsilon\sum_{i=1}^n\sum_{u\in I_i(a,x)}M=M\sum_{i=1}^n|x_i-a_i|\leq nMd(x,a).
$$
\end{eproof}

\subsection{Differentiability}
\begin{Definition}
\label{def:differentiabliblity-1st-multiple} 
A real function $F$ on $[0,1]^n$ is called \textit{differentiable} 
if it has a representation whose partial difference quotients are continuous, 
namely, $F$ has a representation
$(f,\kukanN^{n})$ such that the partial difference quotients
$\hensabun{f}{i}$ are continuous. 
By Corollary\ref{cor:taylor-type-exp-for-every-representation} below,
this definition does not depend on the choice of representations.
\end{Definition}

The real function on $[0,1]^n$ represented by
$(\hensabun{f}{i},[0,1]_{\varepsilon}^n)$ is called the $i$-th partial
derivative and is denoted by $\hendoukansuu{F}{i}$. Sometimes we
write it as $\partial_{i}F$ for brevity. 
These functions are independent not only 
of the choice of $f$ as is seen by the arguments in Lemma~\ref{lem:表現からの独立性}
 but also of the choice of $\varepsilon$ 
by Corollary~\ref{cor:taylor-type-exp-for-every-representation}.

\begin{Theorem}[Infinitesimal Taylor formula of first order]
\label{thm:taylor-1st-multivar}
Let $f$ be a continuous rational valued function on $\kukanN^{n}$ with continuous
partial difference quotients $\hensabun{f}{i}$ ($i\in [1..n]$). Then
\begin{equation}
\label{eq:910-taylor-1st-multivar-a}
f(x)\equiv_1 f(a)+\sum_{i=1}^n \hensabun{f}{i}(a)(x_i-a_i) \myif x\approx a.
\end{equation}
In particular the following also holds.
\begin{equation}
\label{eq:910-taylor-1st-multivar}
f(x)\approx_1 f(a)+\sum_{i=1}^n \hensabun{f}{i}(a)(x_i-a_i) \myif x\approx a
\end{equation}
\end{Theorem}
\begin{eproof}
For simplicity we consider the case $ a_i\leq x_i$.

By the continuity of $\hensabun{f}{i}$,
\begin{eqnarray*}
f(x[i])-f(x[i-1])&=&\sum_{a_i\leq u<x_i}\hensabun{f}{i}(x[i-1]+u\be_i)\Delta x
\\
&\equiv_1&
\hensabun{f}{i}(a)\sum_{a_i\leq u<x_i}\Delta x
=\hensabun{f}{i}(a)(x_i-a_i) \myif x\approx a.
\end{eqnarray*}
by Proposition~\ref{prop:landau-integral}.
\end{eproof}

\begin{Corollary}
\label{cor:taylor-type-exp-for-every-representation}
If $F$ is a differentiable function on $[0,1]^{n}$, 
then for any representation $(f,\kukanN^{n},\alpha)$, 
there are continuous functions $\seq{n}{g}$ such that
\begin{equation}
\label{eq:912-01aa}
f(x)\approx_1  f(a)+\sum_{i=1}^n g_i(a)(x_i-a_i)\myif x\approx a.
\end{equation}
and $(g_i,\kukanN^n,\alpha)$ represents $\partial_iF$.
\end{Corollary}

\begin{eproof}
Suppose $(h,[0,1]_{\delta})$ is a representation of $F$ with continuous
partial difference quotients so that
\begin{equation}
\label{eq:912-10}
h(y)\approx_1 h(b)+\sum_{i=1}^n \hensabun{h}{i}(b)(y_i-b_i) \myif y\approx b.
\end{equation}
Let $\alpha:[0,1]^n_{\varepsilon}\rightarrow [0,1]^n_{\delta}$ be the restriction of
$\kappa^n_{\delta}$ and put
$g_i(x)=\hensabun{h}{i}(\alpha(x))$
for $x\in [0,1]^n_\varepsilon$. Since $\alpha$ is a 
quasi-identity, \eqref{eq:912-10} implies \eqref{eq:912-01aa} 
by Proposition~\ref{prop:invariance-of-landau-by-normal-equivalence}.

Since $g_i\simeq \hensabun{h}{i}$, the continuous function $g_i$ represents $\partial_iF$.
\end{eproof}

\begin{Proposition}[Characterization of diffentiability]
\label{prop:characterization-of-diff}
Suppose a real function $F$ on $[0,1]^{n}$  
has a representation $(f,\kukanN^n)$ 
with continuous rational valued functions $\seq{n}{g}$ such that
\begin{equation}
\label{eq:912-01}
f(x)\approx_1  f(a)+\sum_{i=1}^n g_i(a)(x_i-a_i)\myif x\approx a.
\end{equation}
Then $F$ is differentiable.
\end{Proposition}

\begin{eproof}
By Proposition~\ref{prop:approx-imply-equiv}, 
there is a huge number $L$ such that $L\varepsilon\approx0$ 
and for $a\in [0,1]^n_{L\varepsilon}$
\begin{equation}
\label{eq:912-02}
f(x)\equiv_1 f(a)+\sum_{i=1}^n g_i(a)(x_i-a_i)\myif x\approx a,
\end{equation}
on $[0,1]_{L\varepsilon}^n$. Then 
$$\hensabun{f}{i}(a)=\frac{f(a+L\varepsilon\be_i)-f(a)}{L\varepsilon} = g_i(a)$$
for $a\in [0,1]_{L\varepsilon}^n$ hence
the partial difference quotients of $f$ restricted on $[0,1]_{L\varepsilon}^n$ 
are continuous.
\end{eproof}

\subsection{Chain Rule}
Let $F:[0,1]^{n_1}\rightarrow [0,1]^{n_2}$ be a morphism. Then
$F=(\seq{n_1}{F})$ with real functions $F_i$ on $[0,1]^{n_1}$.  
We call $F$ is \textit{a differentiable morphism} if each $F_i$
is differentiable. 
Let $G$ be a real function on $[0,1]^{n_2}$. The composition
$G\circ F$ is a real function on $[0,1]^{n_1}$.
\begin{Theorem}
If $F$ and $G$ are differentiable, then the composition $G\circ F$ is
differentiable and for $i\in [1..n]$
\begin{equation}
\label{eq:912-100}
\hendoukansuu{(G\circ F)}{i}=\sum_{j=1}^{n_2}\left(
\hendoukansuu{G}{j}\circ F
\right)\hendoukansuu{F_j}{i}.
\end{equation}

\end{Theorem}
\begin{eproof}
Let $(g,\kukanN^{n_2})$ be a representation of $G$ with
continuous difference quotients. 
By \ref{thm:taylor-1st-multivar}, for $b\in [0,1]^{n_2}_\varepsilon$
\begin{equation}
\label{eq:911-15}
g(y)\equiv_1 g(b)+\sum_{i=1}^{n_2} \hensabun{g}{i}(b)(y_i-b_i) \myif y\approx b
\end{equation}
Let $F$ be represented by 
$f:[0,1]_{\varepsilon}^{n_1}\rightarrow [0,1]_{\rational}^{n_2}$.
Then the composition $G\circ F$ is represented by 
$g\circ \kappa_{\varepsilon}^n\circ f$

\begin{center}
$\xymatrix{ 
   [0,1]_{\varepsilon}^{n_1} 
      \ar[r]^f
   & [0,1]_{\rational}^{n_2} 
      \ar[d]^{\kappa_{\varepsilon}}
   & \rational \\
   & \kukanN^{n_2} \ar[ur]^{g} 
   &
}$

\end{center}

By Corollary~\ref{cor:taylor-type-exp-for-every-representation},
we have for $i\in [1..n_1]$ and $a\in \kukane{\varepsilon}^n$
\begin{equation}
\label{eq:912-20}
f_i(x)-\left(f_i(a)+\sum_{j=1}^{n_1}\hensabun{f_i}{j}(a)(x_j-a_j)\right)
\approx_1 0\myif x\approx a.
\end{equation}
Since $\kappa_{\varepsilon}$ is a quasi-identity, this implies
\begin{equation}
\label{eq:912-20-2}
\kappa_{\varepsilon}(f_i(x))
-\left( \kappa_{\varepsilon}(f_i(a))+\sum_{j=1}^{n_1}(\hensabun{f_i}{j}(a))(x_j-a_j)\right)
\approx_10\myif x\approx a,
\end{equation}
by Proposition~\ref{prop:equivalence-higher-infinitesimal}.
Put $\alpha=(\kappa_{\varepsilon},\cdots,\kappa_{\varepsilon}):[0,1]_{\rational}^n\rightarrow[0,1]_{\varepsilon}$. 
If $x\approx a$ then 
$\alpha(f(x))-\alpha(f(a))\approx 0$, 
whence by substituting $y=\alpha(f(x))$,$y_i=\kappa_{\varepsilon}(f_i(x))$  and $b=\alpha(f(a))$ in \eqref{eq:911-15} 
we obtain
\begin{equation}
\label{eq:911-30}
g(\alpha(f(x)))\equiv_1 g(\alpha(f(a)))
+\sum_{i=1}^{n_2} \hensabun{g}{i}(\alpha(f(a)))(\kappa_{\varepsilon}(f_i(x))-\kappa_{\varepsilon}(a))) 
\myif x\approx a
\end{equation}
Hence by \eqref{eq:912-20-2}, we obtain
for $a\in [0,1]_{\varepsilon}^{n_{1}}$, 
\begin{eqnarray*}
g(\alpha(f(x)))
&\approx_1& 
g(\alpha(f(a))) +\sum_{i=1}^{n_2} \hensabun{g}{i}(\alpha(f(a)))
\left(\sum_{j=1}^{n_1}(\hensabun{f_i}{j}(a))(x_j-a_j)\right)
\\
&=& 
g(\alpha(F(a))) +\sum_{j=1}^{n_1}g_{j}(a)(x_j-a_j)
\end{eqnarray*}
with
$$g_{j}(a)=\sum_{i=1}^{n_2}\hensabun{g}{i}(\alpha(F(a)))(\hensabun{f_i}{j}(a)).
$$
Hence by Proposition~\ref{prop:characterization-of-diff}, $G\circ F$ is 
differentiable. Moreover, since $\hensabun{g}{i}(\alpha(f(a))$ represents 
$\hendoukansuu{G}{i}\circ F$, we have \eqref{eq:912-100}.
\end{eproof}

\subsection{Implicit Function Theorem}

The following implicit function theorem is in essense
the inverse function theorem \ref{thm:inverse-funct-theor} of one variable
with parameters and is proved by similar arguments. 

\begin{Theorem}[Implicit Function Theorem]
\label{thm:ImpFn}
\newcommand{\tempx}{x'}
Let $F$ be a differentiable function on $[-1,1]^n$ such that
$$F(0)=0 \quad, \hendoukansuu{F}{n}(0)\neq 0.$$
Then there is a differentiable function $G$ on
$[-c,c]^{n}$ with some $c\succ0$ such that 
for $y\in [-c,c]$ and $x'=(\seq{n-1}{x}) \in [-c,c]^{n-1}$,
$$
F(\tempx ,G(\tempx,y))= y,
$$
and if $x \in [-c,c]^{n}$ satisfies $|F(\tempx,[x_n])|\leq c$, 
then 
$$
G(\tempx ,F(\tempx,x_{n}))=x_{n}.
$$
The partial deriatives are given by 

$$
\partial_i{G}(\tempx,y) = \left\{
  \begin{array}{cc}
-\frac{\partial_iF(\tempx,G(\tempx,y))}%
{\partial_n{F}(\tempx,G(\tempx,y))}
    & \mbox{for $i\in [1..n-1]$} \\
-\frac{1}{\partial_n{F}(\tempx,G(\tempx,y))}
    & \mbox{for $i=n$} 
  \end{array}
\right.
$$
\end{Theorem}
This follows from the following lemma.
\begin{Lemma}
\label{lemma:ImpFn}
Let $f$ be a rational valued function on $[-1,1]_{\varepsilon}^n$ with continuous partial difference
quotients and suppose
$$f(\bfzero)=0 \quad, \hensabun{f}{n}(\bfzero)\not\approx 0.$$
Then there is a continuous rational valued function $g$ on
$[-c,c]_\varepsilon^{n}$ with some $c\succ0$ such that for $y\in [-c,c]_\varepsilon$ and
$\bx' \in [-c,c]_\varepsilon^{n-1}$
$$
f(\bx',g(\bx',y))\approx y.
$$
Moreover if 
$(\bx',x_n)\in [-c,c]_{\varepsilon}^{n}$ 
satisfies 
$f(\bx',x_n)\in [-c,c]_{\rational}$, 
then
$$
g(\bx',\kappa_{\varepsilon}(f(\bx',x_{n})))\approx x_{n}.
$$

The difference quotients of $g$ at $(\bx',y)\in [-c,c]_{\varepsilon}^n$ is given by
$$
\hensabun{g}{i}(\bx',y)\approx \left\{
  \begin{array}{cc}
-\frac{\hensabun{f}{i}(\bx',g(\bx',y))}{\hensabun{f}{n}(\bx',g(\bx',y))}    &
 \mbox{for }i\in [1..n-1]\\
  \frac{1}{\hensabun{f}{n}(\bx',g(\bx',y))}  & \mbox{for $i=n$} 
  \end{array}
\right.
$$
and are continuous.
\end{Lemma}
\begin{eproof}
We may assume $\hensabun{f}{n}(\bfzero)\succ 0$. Then
there are $\alpha,c_1\succ 0$ such that 
$\hensabun{f}{n}(x)\geq \alpha$ 
if $x\in [-c_1,c_1]^n_{\varepsilon}$.
By Proposition~\ref{prop:bibun-hyouka},
$$
f(\bfzero,-c_1)\prec -\frac12\alpha c_1, \quad \frac12\alpha c_1\prec f(\bfzero,c_1),
$$
whence there is $0\prec c_2\leq  c_1$ such that if $\bx'\in [-c_2,c_2]^{n-1}_{N}$
then
$$
f(\bx',-c_1)\prec -\frac12\alpha c_1, \quad\frac12\alpha c_1\prec f(\bx',c_1).
$$
Put 
$c=\min\seti{c_2,\frac{\alpha c_1}4}$,
then $(\bx',y)\in [-c,c]^{n}$ implies
$$
f(\bx',-c_1)\prec y \prec f(\bx',c_1),
$$
since $f(\bx',c_1)\succ \frac12\alpha c_1\geq 2c>c>y$ for the first equality.
Define
$$
g(\bx',y)=\max\setii{u\in [-c_1,c_1]_{\varepsilon}}{f(\bx',u)\leq y}.
$$
Then
$$ f(\bx',g(\bx',y))\leq y< f(\bx',g(\bx',y)+\Delta x),$$
which implies
\begin{equation}
\label{eq:914a}
f(\bx',g(\bx',y))\approx y.
\end{equation} 

By Corollary~\ref{cor:inverse-continuity-of-function-with-positive-sabun},  
\begin{equation}
\label{eq:914b}
\mbox{if $(\bx',u),(\bx',v)\in [-c_1,c_1]_{\varepsilon}^{n}$ satisfies $f(\bx',u)\approx f(\bx',v) \mbox{ then } u\approx v$.}
\end{equation}
If $(\bx',u)\in [-c,c]_{\varepsilon}^n$ satisfies
$f(\bx',u)\in [-c,c]$, then (\ref{eq:914a}) with $y=f(\bx',u)$ implies
$$f(\bx',g(\bx',\kappa_{\varepsilon}(f(\bx',u)))) 
\approx \kappa_{\varepsilon}(f(\bx',u))
\approx f(\bx',u)
,$$
whence by (\ref{eq:914b})
$$u\approx g(\bx',\kappa_{\varepsilon}(f(\bx',u))),$$
since the values of $g$ is in $[-c_1,c_1]_{\varepsilon}$.

If $(\bx'_i,y_i)\in [-c,c]_{\varepsilon}^n$ \itwo satisfies 
$\bx'_1\approx \bx'_{2},y_1\approx y_2$, then 
$$f(\bx_1',g(\bx_1',y_1))
\approx y_1\approx y_2
\approx f(\bx_2',g(\bx_2',y_2))
\approx f(\bx_1',g(\bx_2',y_2))
,$$
whence, by (\ref{eq:914b}) again,
$g(\bx_1',y_1)\approx g(\bx_2',y_2)$. Thus $g$ is continuous.

Suppose $f(\bx',x_n)=y $ and $f(\bu',u_n)=w $, which imply respectively
$x_n\approx g(\bx',y)$ and $u_n\approx g(\bu',w)$.
By Theorem~\ref{thm:taylor-1st-multivar},
\begin{eqnarray*}
f(\bx',x_n)
&\equiv_1&
 f(\bu',u_n)+\sum_{i=1}^{n-1}\hensabun{f}{i}(\bu',u_n)(x_i-u_i)
+\hensabun{f}{n}(\bu',u_n)(x_n-u_n)\\
&& \myif (\bx',x_n)\approx (\bu',u_n).
\end{eqnarray*}
Hence, from Lemma~\ref{lemma:lipshitz-condition-for-differentiable-function} below
and Proposition~\ref{prop:equivalence-higher-infinitesimal} it follows
$$
y
\approx_{1}
w
+\sum_{i=1}^{n-1}\hensabun{f}{i}(\bu',u_n)(x_i-u_i)
+\hensabun{f}{n}(\bu',u_n)(g(\bx',y)-g(\bu',w)) 
\myif (\bx',y)\approx(\bu',w).
$$
Since $\hensabun{f}{n}(\bu',u_n)\not\approx0$, 
by Proposition~\ref{prop:equivalence-higher-infinitesimal}
we have 
$$
g(\bx',y)\approx_1 g(u',w)
-
\sum_{i=1}^{n-1}\frac{\hensabun{f}{i}(\bu',g(\bu',w))}{\hensabun{f}{n}(\bu',g(\bu',w))}(x_i-u_i)
+\frac{1}{\hensabun{f}{n}(\bu',g(\bu',w))}(y-w).
$$
Since 
$$
h_i(u',w):= \left\{
  \begin{array}{cc}
-\frac{\hensabun{f}{i}(u',g(u',w))}{\hensabun{f}{n}(u',g(u',w))} 
& \mbox{for $i\leq n-1$} \\
\frac{1}{\hensabun{f}{n}(u',g(u',w))}  & \mbox{for $i=n$} 
  \end{array}
\right.
$$
are continuous, Proposition~\ref{prop:characterization-of-diff}
implies $g$ represents a differentiable function whose $i$-th partial
derivative is represented by $h_i$.
\end{eproof}

The following was used in the proof of Lemma~\ref{lemma:ImpFn}.
\begin{Lemma}
\label{lemma:lipshitz-condition-for-differentiable-function-multiple-var}
Let $f$ be a rational valued function on $[0,1]^n_{\varepsilon}$
with continuous difference quotient. Define a function
$F:[0,1]_{\varepsilon}^n\rightarrow [0,1]_{\rational}^n$
by $$F(x)=(x_1,\cdots,x_{n-1},f(x)).$$
If $\hensabun{f}{n}(a)\not\approx 0$,
then for some rationals $K_1,K_2,c\succ 0$ 
$$
K_1d(x,a)\preceq d(F(x),F(a))\preceq K_2d(x,a)
$$
holds for if $d(x,a)<c$.
\end{Lemma}
\begin{eproof}
Put $M:=\max\setii{|\hensabun{f}{i}{x}|}{i\in [1..n], x\in [0,1]^n_{\varepsilon}}$.
Then by Proposition~\ref{prop:diff-quot-estimate}
$$ |f(x)-f(a)|\leq M\sum_{i\in [1..n]}|x_i-a_i|\leq nMd(x,a)$$
for $x,a\in [0,1]^n_{\varepsilon}$.
Put $x'=(\seq{n-1}{x})$ and $a'=(\seq{n-1}{a})$. Then
$$d(F(x),F(a))=\max\seti{d(x',a'),|f(x)-f(a)|}.$$ 
Since $d(x',a')\leq d(x,a)$, if we put $K_2=\max\seti{1,nM}$,
then
$$
d(F(x),F(a))\leq K_2d(x,a).
$$

Put $b=\hensabun{f}{n}(a)$. We may assume $b>0$. 
Define 
$$b_2:=\max_{1\leq i\leq n-1}\left|\frac{\hensabun{f}{i}(a)}{\hensabun{f}{n}(a)}\right|$$
Then we can choose $c_1\succ 0$ and $c_2\succ 0$ such that 
$d(x,a)<c_1$ implies 
$\hensabun{f}{n}(x)>\frac{b}2$
and
$|\hensabun{f}{i}(x)|<2b_2\hensabun{f}{n}(x)$ for $i\leq n-1$.

Put $L=\min\seti{2,\frac{b}{4},\frac{b}{8b_2M}}$ and suppose $d(x',a')<L|x_n-a_n|$.
Then by Lemma~\ref{lemma:differene-quotient},
\begin{eqnarray*}
|f(x)-f(a)| 
&>  &
|\sum_{u\in I_n(a,x)}\hensabun{f}{n}(x[n-1]+u\be_n)\varepsilon
|
-
\sum_{i=1}^{n-1}\sum_{u\in I_i(a,x)}|\hensabun{f}{i}(x[i-1]+u\be_i)|\varepsilon
 \\
&\geq&
\frac{b}{2}|x_n-a_n|-2b_2\sum_{i=1}^{n-1}\sum_{u\in I_i(a,x)}|\hensabun{f}{n}(x[i-1]+u\be_i)|\varepsilon
 \\
&\geq&
\frac{b}{2}|x_n-a_n|-2b_2M\sum_{i=1}^{n-1}|x_i-a_i|
 \\
&\geq&
(\frac{b}{2}-2b_2ML)|x_n-a_n|\geq \frac{b}{4}|x_n-a_n|
\end{eqnarray*}
Since
$$d(x,a)=\max\seti{d(x',a'),|x_n-a_n|}
\leq \max\seti{L|x_n-a_n|,|x_n-a_n|}=L|x_n-a_n|$$
we conclude that $d(x',a')<L|x_n-a_n|$ implies
$$|f(x)-f(a)| > \frac{b}{4L}d(x,a).$$
Hence
$$d(F(x),F(a))=\max\seti{d(x',a'),|f(x)-f(a)|}
>\max\seti{d(x',a'),\frac{b}{4L}d(x,a)}=\frac{b}{4L}d(x,a),
$$
since $\frac{b}{4L}\geq 1$.

Suppose $d(x',a')\geq L|x_n-a_n|$. 
Then 
$$
d(x,a)=\max\seti{d(x',a'),|x_n-a_n|}
\leq
\max\seti{d(x',a'),\frac1Ld(x',a')}=L_2d(x',a'),
$$
where $L_2:=\max\seti{1,\frac1L}$.
Hence
$$
d(F(x),F(a))\geq d(x',a')\geq \frac1{L_2}d(x,a).
$$
So if we put $K_1=\min\seti{\frac{b}{4L},\frac1{L_2}}$, then
$d(x,a)<c_1$ implies
$$K_1d(x,a)\leq d(F(x),F(a)).$$
\end{eproof}

\subsection{Inverse Mapping Theorem}
Let $F=(\seq{m}{F})$ be a differentiable morphism from $[-1,1]^n$ to $\virtualline^m$.
The matrix 
$$dF(p):=\left(\;\partial_j{f_i}(p)\;\right)_{i\in [1..n],j\in [1..m]}$$
is called the Jacobian of $F$ at a point $p\in (-1,1)^n$. 

A differentiable map $F:[-1,1]^{n}\rightarrow [-1,1]^{n}$ with $F(0)=0$
is called \textit{a local diffeomorphism}
\index{local diffeomorphism@local diffeomorphism}
at the point $0$ if 
there is a $c\succ0$ and a differentiable morphism
$$G:[-c,c]^n \rightarrow [-1,1]^n$$ 
such that for $y\in [-c,c]^n$ 
$$F(G(y))=y$$
and 
$$G(F(x))=x,$$
for $x\in [-c,c]^{n}$ such that $F(x)\in [-c,c]^n$.
$G$ is called a \textit{local inverse} of $F$.

\begin{Theorem}
\label{th:913-1}
Let $F$ be a differentiable mapping $[-1,1]^n$ to $\virtualline^n$ such that
$F(0)=0$ with invertible Jacobian at $0$.
Then $F$ is a local diffeomorphism at $0$.
\end{Theorem}

Theorem~\ref{thm:ImpFn} implies the following special case.
\begin{Lemma}
\label{lem:914a}
Let $F_n$ be a differentiable function on $[-1,1]^n$ such that
$F_n(0)=0$, $\partial_n{F_n}(0)=1$ and 
$\partial_i{F_n}(0)=0$ for $i<n$. 
Then the differntial morphism 
$F:[-1,1]^n\rightarrow \virtualline^n$ defined by $F(x)=(\seq{n-1}{x},F_n(x))$
is a local diffeomorphism at $[0]$.
Moreover the Jacobian matrix of every local inverses of $F$ at 
$0$ is the identity matrix.
\end{Lemma}
\begin{eproofi}{ of Lemma~\ref{lem:914a}}
Let $(f,[-1,1]_{\epsilon}^n])$ represents $F_n$ so that
$\hensabun{f}{n}(0)\approx 1$  and $\hensabun{f}{i}(0)\approx 0$ for $i\in [1..n-1]$.
By Theorem~\ref{thm:ImpFn}, 
there is a differentiable function $g$ on $[-c,c]_{\varepsilon}^n$ with $c\succ0$ 
such that for $x \in [-c,c]_{\varepsilon}^n$
$$  f(x',g(x))\approx x_n  $$
where $x':=(\seq{n-1}{x})$ and if $f(x)\in [-c,c]$ then
$$  g(x',f(x))\approx x_n.$$
Put, for $x\in [-c,c]_{\varepsilon}^n$, $G(x):=(x',g(x))$. 
Then it is obvious that $G$ satisfies
the conditions of local inverse. 
By Theorem~\ref{thm:ImpFn},
$$
\hensabun{g}{i}(\bfzero)\approx 
\left\{
  \begin{array}{cc}
-\frac{\hensabun{f}{i}(\bfzero)}{\hensabun{f}{n}(\bfzero)}\approx 0
    & \mbox{for } i<n \\
-\frac{1}{\hensabun{f}{n}(\bfzero)}\approx 1
    & \mbox{for }  i=n
  \end{array}
\right.
$$
Hence the Jacobian matrix of $G$ at $\bfzero$ is the identity matrix.

\end{eproofi}

\begin{eproofi}{ of Theorem~\ref{th:913-1}}
By applying the inverse of $dF(\bfzero)$ to $F$, 
we may suppose $dF(\bfzero)$ is the identity matrix $I_n$ of size $n$.

Put $\Phi_n(x):=(\seq{n-1}{x},F_n(x))$. 
Applying Lemma~\ref{lem:914a} for $F_n$, we obtain 
differential morphism $\Gamma_n:[-c_n,c_n]^n\rightarrow \real^n$ such that
for $x\in [-c_n,c_n]^n$,
$$ \Phi_n(\Gamma_n(x))= x$$
and hence 
$$ F_n(\Gamma_n(x))= x_n, $$
and if $\Phi_n(x)\in [-c_n,c_n]^{n}$ then 
$$ \Gamma_n(\Phi_n(x))= x.$$

We define inductively local diffeomorphisms $\Gamma_{n-1},\Gamma_{n-1},\cdots,\Gamma_{1}$ 
satisfying, for each $j$, 
\begin{equation}
\label{eq:20120308-1}
d\Gamma_j(\bfzero)= I_n
\end{equation}
and 
\begin{equation}
\label{eq:914c}
F_{k}(\widehat{\Gamma}_{j}(x))= x_k  \mbox{ for } k\in [j..n],
\end{equation}
with $\widehat{\Gamma}_j=\Gamma_n\circ \Gamma_{n-1}\circ\cdots\circ \Gamma_{j}$. 
Then $\widehat{\Gamma}_1$ is the required local diffeomorphism.

Suppose we have constructed $\Gamma_n,\cdots,\Gamma_{i+1}$ satisfying the
conditions \eqref{eq:20120308-1} and \eqref{eq:914c} for $j\geq i+1$.  
Define 
$$\widetilde{F}_i(x):=F_i(\widehat{\Gamma}_{i+1}(x)).$$ 
Since $d\widehat{\Gamma}_{i+1}(\bfzero)= I_n$, we have
$\partial_j{\widetilde{F}_i}(\bfzero)=\partial_jF_i(\bfzero) =\delta_{ij}$.  Hence by Lemma~\ref{lem:914a},
there is a local diffeomorphism $\Gamma_i$ such that $d\Gamma_i(\bfzero)=I_n$, 
$\widetilde{F}_i(\Gamma_i(x))= x_i$ and, for $j\neq i$, the $j$-th component
of $\Gamma_i(x)$ is $x_j$. Then for $k>i$
$$
F_k(\widehat{\Gamma_i}(x))
=F_k(\widehat{\Gamma}_{i+1}(\Gamma_i(x)))
=(\Gamma_i(x))_k=x_k
$$
and 
$$
F_i(\widehat{\Gamma_i}(x))
=F_i(\widehat{\Gamma}_{i+1}(\Gamma_i(x)))
=\widetilde{F}_i(\Gamma_i(x))=x_i,
$$
whence (\ref{eq:914c}) holds for $j\geq i$.
\end{eproofi}

\subsection{Second Order Differentiability}

\begin{Definition}
  A real function $F$ on $[0,1]^n$ is \textit{differentiable up to
  second order}
\index{differentiable, second order}
 if it is differentiable and its partial derivatives are differentiable. 
The partial derivatives $\partial_j\partial_iF$ are called
\textit{second order partial derivatives}. It will be shown that it is symmetric
with respect to $i,j$, namely, $\partial_j\partial_iF=\partial_i\partial_jF$.
\end{Definition}

\newcommand{\Deltaf}[2]{\Frac{\Delta_{#1} \Delta_{#2} f}{\Delta x^2}}
\begin{Theorem}
  \label{th:multivariate-2nd-Taylor}
Suppose a function $F$ on $[0,1]^n$ is differentiable up to second order.
Let $f$,$g_i$ and $g_{ij}$ are  rational valued functions on $[0,1]_{\varepsilon}^n$
representing respectively $F$, $\partial_iF$ and $\partial_i\partial_jF$.
Then there is an infinitesimal $\delta\in \varepsilon\integer$ 
such that if $a\in [0,1]_{\delta}^{n}$ then
\begin{equation}
\footnotesize
\label{eq:multivariate-2nd-Taylor}
    f(x)\equiv_2 f(a)+\sum_{i=1}^ng_i(a)(x_i-a_i)
    +\sum_{1\leq j<i\leq n}g_{ji}(a)(x_i-a_i)(x_j-a_j)
    +\frac12\sum_{i}^{n}g_{ii}(a)(x_i-a_i)^2
\myif x\approx a.
\end{equation}
on $[0,1]_{\delta}^{n}$.
\end{Theorem}

\begin{eproof}
By Theorem \ref{thm:taylor-1st-multivar}, we have 
\begin{eqnarray}
\label{eq:915-1}
f(x) &\approx_1 & f(a)+\sum_{i=1}^ng_i(a)(x_i-a_i) \myif x\approx a
 \\
\label{eq:915-2}
g_i(x) &\approx_1 & g_{i}(a)+\sum_{j=1}^ng_{ji}(a)(x_j-a_j)  \myif x\approx a
\end{eqnarray}
Hence by Corollary~\ref{cor:restriction},
there is a number $L$ with $\varepsilon':=L\varepsilon\approx 0$ such that 
if $a\in X:=[0,1]_{\varepsilon'}$ then
\begin{eqnarray}
\label{eq:915-3}
f(x) &\equiv_1 & f(a)+\sum_{i=1}^ng_i(a)(x_i-a_i)
\myif x\approx a
 \\
\label{eq:915-4}
g_i(x) &\equiv_1 & g_i(a)+\sum_{j=1}^ng_{ji}(a)(x_j-a_j)
\myif x\approx a
\end{eqnarray}

Assume $x\approx a\in X$ and $a_i\leq x_i$ for all $i$. 
In the other cases, the assertion can be proved similarly.
Put $\Delta x=\delta$. Using the notation in Lemma~\ref{lemma:differene-quotient}, 
\eqref{eq:915-3} implies
  \begin{eqnarray*}
    f(x)-f(a) 
    &= &
    \sum_{i=1}^n(f(x[i])-f(x[i-1])).
    \\
    &\approx &
    \sum_{i=1}^n\left(\;
      \sum_{0\leq u< \frac{x_i-a_{i}}{\Delta x}} g_i(x[i-1]+u\Delta x\be_i)\Delta x
      \;\right)
    \\
    &=&
    \sum_{i=1}^ng_{i}(x)(x_i-a_i)
    +\sum_{i=1}^n\left(\;
      \sum_{0\leq u< \frac{x_i-a_{i}}{\Delta x}} A_i(u)\Delta x
      \;\right),
  \end{eqnarray*}
where 
  \begin{eqnarray*}
    A_i(u) &:= &
    g_i(x[i-1]+u\Delta \be_i)-g_i(x)
    \\
    &=& 
    g_i(x[i-1]+u\Delta \be_i)-g_i(x[i-1])
    +\sum_{j=1}^{i-1}g_i(x[j])- g_i(x[j-1]).
  \end{eqnarray*}
By \eqref{eq:915-4} and Proposition~\ref{prop:landau-integral},  
\begin{eqnarray*}
A_i(u) &\approx& 
    \sum_{0\leq w<u}
    g_{ii}(x[i-1]+w\Delta x \be_i)\Delta x
    +
    \sum_{j=1}^{i-1}\sum_{0\leq  w< \frac{x_j-a_j}{\Delta x}}
    g_{ji}(x[j-1]+w\Delta x\be_j)\Delta x
    \\
    &\approx_2& 
    g_{ii}(x)\sum_{0\leq w<u}\Delta x 
    +\sum_{j=1}^{i-1}\left(g_{ji}(x)\sum_{0\leq  w< \frac{x_j-a_j}{\Delta x}}\Delta x\right)
\myif x\approx a
    \\
    &\approx & 
    g_{ii}(x)u\Delta x
    +\sum_{j=1}^{i-1}g_{ji}(x)(x_j-a_j)
  \end{eqnarray*}
Hence 
\begin{eqnarray*}
&&\sum_{i=1}^{n}\sum_{0\leq u< \frac{x_i-a_i}{\Delta x}} A_i(u)\Delta x
\\
&\approx_2 &
\sum_{i=1}^{n}\left(g_{ii}(x)\sum_{0\leq u< \frac{x_i-a_i}{\Delta x}}u(\Delta x)^2\right)
+
\sum_{i=1}^{n}\sum_{j=1}^{i-1}g_{ji}(x)(x_j-a_j)
\sum_{0\leq u< \frac{x_i-a_i}{\Delta x}}\Delta x
\\
&=&
\sum_{i=1}^{n}\left(g_{ii}(x)\frac12({x_i-a_i})({x_i-a_i}+\Delta x)\right)
+
\sum_{j<i}g_{ji}(x)(x_j-a_j)(x_i-a_i)
\\
&\approx_2&
\sum_{i=1}^{n}\frac12g_{ii}(x)(x_i-a_i)^{2}
+
\sum_{j<i}g_{ji}(x)(x_j-a_j)(x_i-a_i)
\end{eqnarray*}
Hence, we have (\ref{eq:multivariate-2nd-Taylor}) with $\equiv_2$ replaced with
$\approx_2$. 

By Proposition~\ref{prop:approx-imply-equiv}, there is an integer $L'>0$ with 
$\delta:=L'\varepsilon'\approx0$ such that \eqref{eq:multivariate-2nd-Taylor}
holds for $a\in [0,1]_{\delta}$.
\end{eproof}

Let $f$ be a rational valued function on $X=\kukanN^{n}$ with \accessible{} $n$.
The \textit{second order partial difference quotients} are defined by
$$
\hensabun{f}{ij}:=\hensabun{(\hensabun{f}{j})}{i}.
$$
More explicitly we have the following.
\begin{Lemma}
\label{lem:second-order-diff-symmetric}
\begin{equation}
\label{eq:915-6}
\hensabun{f}{ii}:=\frac{f(x+2\Delta x \be_i)-2f(x+\Delta x\be_i)+f(x)}{(\Delta x)^{2}},
\end{equation}
\begin{equation}
\label{eq:915-7}
\hensabun{f}{ij}:=\frac{f(x+\Delta x \be_i+\Delta x \be_j)
-f(x+\Delta x\be_i)
-f(x+\Delta x\be_j)
+f(x)}{(\Delta x)^{2}} \mbox{  if $i\neq j$}.
\end{equation}
In particular $\hensabun{f}{ij}=\hensabun{f}{ji}$.
\end{Lemma}

\begin{Proposition}
\label{prop:915a}
Suppose a rational valued function $f$ on $[0,1]_r^n$ with an infinitesimal positive rational $r$ 
satisfies \eqref{eq:multivariate-2nd-Taylor} with continuous $g_i$ and $g_{ij}$. Then
$$
\hensabun{f}{i}\approx g_i,
$$
and if $j<i$ 
$$
\hensabun{f}{ji}\approx g_{ji}.
$$
\end{Proposition}

\begin{eproof}
Put $\Delta x=r$. Substituting $x=a+\Delta x\be_i$ and 
$x=a+2\Delta x\be_i$ in \eqref{eq:multivariate-2nd-Taylor} 
we obtain
\begin{eqnarray}
\label{eq:915-11}
f(a+\Delta x\be_i) &\equiv_2& f(a)+ g_i(a)\Delta x+ \frac12 g_{ii}(a)(\Delta x)^2
\\
\label{eq:915-12}
f(a+2\Delta x\be_i) &\equiv_2& f(a)+ 2g_i(a)\Delta x+ 2 g_{ii}(a)(\Delta x)^2.
\end{eqnarray}
Hence \eqref{eq:915-11} implies $\hensabun{f}{i}(a)\approx g_i(a)$. 

By \eqref{eq:915-6},~\eqref{eq:915-11} and \eqref{eq:915-12}, we have
$$\hensabun{f}{ii}(a)\approx g_{ii}(a).$$ 

Putting $x=a+\Delta x\be_j+\Delta x\be_i$ we have if $j<i$
\begin{eqnarray*}
f(a+\Delta x\be_i+\Delta x\be_j) &\equiv_2& 
f(a)+ (g_i(a)+g_j(a))\Delta x+ \frac12 g_{ji}(a)(\Delta x)^2.
\end{eqnarray*}
Hence by \eqref{eq:915-7} if $j<i$, $\hensabun{f}{ji}(a)\approx g_{ji}(a)$.
\end{eproof}

By Theorem~\ref{th:multivariate-2nd-Taylor} and Proposition~\ref{prop:915a}
we have the following.
\begin{Corollary}
\label{cor:second-order-diff}
If a real function $F$ is differentiable up to second order, then
$F$ is represented by a rational valued function $f$ with continuous 
partial difference quotients up to second order and $\partial_i\partial_jF$
is represented by $\hensabun{f}{ij}$
\end{Corollary}

By Lemma~\ref{lem:second-order-diff-symmetric} and Corollary~\ref{cor:second-order-diff},
we have proved the following.
\begin{Theorem}
\label{thm:commutativity}
Let $F$ be a function on $[0,1]^n$ differentiable up to second order. Then
\begin{equation}
\label{eq:915-20}
\partial_i\partial_jF=\partial_j\partial_iF.
\end{equation}
\end{Theorem}

\begin{Corollary}
  \label{cor:multivariate-2nd-Taylor}
  Suppose a function $F$ on $[0,1]^n$ is differentiable up to second
  order.  Let $f$,$g_i$ and $g_{ij}$  are rational valued
  functions on $[0,1]_{\Omega}^n$ representing respectively $F$,
  $\partial_iF$ and $\partial_i\partial_jF$.  Then $g_{ij}\approx
  g_{ji}$ and there is an integer $L$ with $r=L\varepsilon\approx0$
  such that if $a\in [0,1]_{r}^{n}$ then
\begin{equation}
    \label{eq2:multivariate-2nd-Taylor}
    f(x)\equiv_2 f(a)+\sum_{i=1}^ng_i(a)(x_i-a_i)
    +\frac12\sum_{1\leq i,j\leq n}g_{ji}(a)(x_i-a_i)(x_j-a_j)
\myif x\approx a. 
  \end{equation} 
\end{Corollary}

\newpage

\newpage
\section{Measure}
\def\runningtitle{Measure}
\label{sec:measure}

In this section, we show how to obtain the basic tools of Lebesgue
integration in our framework.
We start with a set $X$ with positive
probability density $p$, which give measure $m(A)$ of subsets of
$A\subset X$. A condition is said to be true almost everywhere if
there are subsets with arbitrary small measure outside of which it
holds.  The integral of a rational valued function $f$ on $X$ is
defined by $E(f):=\sum_{x\in X}p(x)f(x)$. A function $f$ is called
$L^1$ function if $E(f)\approx E(f^{a})$ for all huge $a$ where
$f^{a}$ denotes the function obtained by modifying $f(x)$ to zero when
$f(x)>a$. Note that this concept get meaning since we have functions
with huge values. Then the $L^1$ functions form a complete metric
space with respect to the distance function $d_1(f,g)=E(|f-g|)$. A
concrete sequence of $L^1$ functions converges with respect to $d_1$
then a subsequence converges pointwise almost everywhere.

\newcommand{\aecong}{\stackrel{a.e.}{\approx}}
\newcommand{\unionP}{\bigcup_{i=1}^{\infty}P_i}
\subsection{Probability Density}
Let $X$ be a set and $p$ be a \textit{probability density function}, 
\index{probability density@probability density}
namely, a rational valued function on $X$ 
satisfying 
\begin{itemize}
\item $0\leq p(x)\leq 1$
\item $\sum_{x\in X}p(x)=1$.
\end{itemize}
For a subset $A\subset X$, define its \textit{measure} $m(A):=\sum_{a\in A}p(a)$.  
\index{measure@measure, subset}
\index{measure@measure, $m(A)$}
Then for
$A,B\subset X$, we have obviously
$$m(A\bigcup B)=m(A)+m(B)-m(A\bigcap B).$$
If $\seq{N}{A}$ is a huge sequence of mutually disjoint subsets, then we have
also  
$$m\left(\bigcup_{i=1}^N A_i\right)= \sum_{i=1}^{N}m(A_{i}).$$

\subsection{Null Semisets}
We call a subsemiset $P\sqsubset X$ \textit{null semiset}
\index{semiset, null}
\index{$P\aecong\emptyset$}
and write $P \aecong \emptyset$ if for each \accessible{}  $k$, there is a subset $A\subset X$
satisfying $P\sqsubset A$ and $m(A)<\frac1k$. 
Note that a subset$A$  is a null semiset
if and only if $m(A)\approx 0$.

Obviously intersection and union of two null semisets are null. 

We call a condition $Q$, not necessarily objective,
 holds almost everywhere (a.e. for short) if the
subsemiset defined by $\lnot Q$ is a null semiset.
\index{a.e.}
Obviously $Q$ holds a.e. if and only if for every \accessible{}  $k$
there is a subset $B\subset X$ such that $m(B)>1-\frac1k$ and
$Q(x)$ holds for all $x\in B$. 

A typical condition we encounter is $f(x)\approx 0$ for a rational valued
function $f$ on $X$.

\begin{Lemma}
\label{lem:a-e-Robinson}
Let $f$ be a rational valued function on $X$. Then $f(x)\approx 0$ a.e. if and 
only if there is a subset $A\subset X$ with $m(A)\approx 0$ such that
$f(x)\approx 0$ for all $x\not\in A$.
\end{Lemma}

\begin{eproof}
Suppose $f(x)\approx 0$ a.e.. For \accessible{}  $k$, there is a subset 
$A_k\subset X$ such that $m(A_k)<\frac1k$ and 
$f(x)\approx 0$ for $x\nin A_k$. Then the set of numbers
$$
\setii{k}{
\mbox{there is a subset $A\subset X$ such that 
$m(A)<\frac1k$ and $|f(x)|<\frac1k$ for all $x\nin A$}
}
$$
contains all \accessible{}  numbers and hence also a huge number $K$.
Hence there is a subset $A\subset X$ such that
$m(A)<\frac1K$ and 
$|f(x)|<\frac1K$ for $x\nin A$. 
Then $m(A)\approx 0$ and $f(x)\approx 0$ for $x\nin A$.

The converse is obvious.
\end{eproof}

The following properties hold obviously.
\begin{Proposition}
\label{prop:boolean-properties-null-parts}
\begin{itemize}
\item If $P$ is a null semiset and $Q\sqsubset P$, then $Q$ is also a null semiset.
\item If $P_i\aecong \emptyset$\itwo are null semisets, then their union is also a null semiset.
\end{itemize}
\end{Proposition}

The infinite union of null semisets is also a null semiset.
\begin{Theorem}
\label{th:infinite-conjunction}
Suppose $\seqinf{P}$ is a concrete  sequence of null semisets. Then their union
$$\bigcup_{i}P_i$$
is also a null semiset.
\end{Theorem}
\begin{eproof}
Let $k$ be an \accessible{} number. Since $P_i$ is null, for each \accessible{}  $i$, there is 
a set $A_i\sqsupset P_i$ with 
\begin{equation}
\label{eq:9-803}
m(A_i)<\frac{1}{k2^i}.
\end{equation}
By the over-spill axiom, the concrete sequence $\seqinf{A}$ can be
extended to a huge sequence $(\seq{N}{A})$ which satisfy ~\eqref{eq:9-803}.
Put 
$$A=\bigcup_{i=1}^{M}A_i.$$
Then obviously $A\sqsupset \bigcup_{i=1}^{\infty}P_i$ and 
$$m(A)\leq \sum_im(A_i)\leq \frac1k\left(1-\left(\frac12\right)^{M+1}\right)<\frac1k.$$
\end{eproof}

For subsemisets  $P_i\sqsubset X$ \itwo, we write 
$P_1\aecong P_2$ if $P_1\Delta P_2\aecong \emptyset$.
\begin{Lemma}
\label{lemma:equivalence-relation-aecong}
The relation $\aecong$ is an equivalence relation.
\end{Lemma}
\begin{eproof}
The reflexsivity and symmetricity are obvious.
For the transitivity, it suffices to show 
\begin{equation}
\label{eq:1107-10}
P_1\Delta P_3 \sqsubset (P_1\Delta P_2)\cup (P_2\Delta P_3).
\end{equation}
Suppose $x\in P_1\setminus P_3$. If $x\in P_2$ then $x\in P_2\setminus P_3\sqsubset P_2\Delta P_3$,
whereas if $x\nin P_2$ then $x\in P_1\setminus P_2\sqsubset P_1\Delta P_2$. Hence
$ P_1\setminus P_3$ is included in the right hand side of  (\ref{eq:1107-10}).
Similarly it can be shown that $ P_3\setminus P_1$ is included in the right hand side of (\ref{eq:1107-10}).
\end{eproof}

\subsection{Measurable Semisets} 
 A subsemiset $P\sqsubset X$ is called \textit{measurable}
\index{semiset, mesurable}
\index{$\overline{m}(P)$}
if there is
a subset $A\subset X$ with null subsemiset $A \Delta P$. We define
then $\overline{m}(P):=[m(A)]$, 
which is independent of the choice of $A$ 
by the following lemma~\ref{lem:817-1} 
hence uniquely defined as a real number.

Note that if $P$ is measurable, 
then $\overline{m}(P)= 0$ means that $P$ is null semiset.

\begin{Lemma}
\label{lem:817-1}
If $P\sqsubset X$ and $A_1,A_2\subset X$ satisfies
$A_i\Delta P\aecong \emptyset$ \itwo, then 
$m(A_1)\approx m(A_2)$.
\end{Lemma}
\begin{eproof} 
First we show that $m(A_1\Delta A_2)\approx0$, for which  
it suffices to show 
$$ A_1\Delta A_2\sqsubset (A_1\Delta P)\bigcup (A_2\Delta P),$$
since there are $B_i\subset X$ such that $A_i\Delta P\sqsubset B_i$ and $m(B_i)\approx 0$
\itwo. Suppose $x\in A_1\setminus A_2$. If $x \in P$ then 
$x\in P\Delta A_2$ and if $x\nin P$ then $x\in P\Delta A_1$, whence
\begin{equation}
\label{eq:817-1}
x\in (A_1\Delta P)\bigcup (A_2\Delta P).
\end{equation}
Similarly it can be shown that (\ref{eq:817-1}) holds for $x\in A_2\setminus A_{1}$.

Define $C=A_1\cap A_2$, $A'_i=A_i\setminus C$ \itwo.
Then
from $A_1\Delta A_2=A'_1\cup A'_2$ it follows $m(A'_i)\approx 0$\itwo. Hence
$$
m(A_1)=m(A_1')+m(C)\approx m(A_2')+m(C)=m(A_2).
$$
\end{eproof}

Measurability of $P$ can be rephrased as follows.
\begin{Proposition}
\label{prop:rephrase-of-measurability}
A subsemiset $P\sqsubset X$ is measurable if and only if 
for each \accessible{}  number $k$, there exist subsets
$A,B\subset X$ satisfying 
\begin{equation}
\label{eq:1107-1}
A\setminus B\sqsubset P\sqsubset A\bigcup B,
\end{equation}
and $m(B)<\frac1k$.  
\end{Proposition}
\begin{eproof}
Suppose $P$ is measurable. Then there is a subset $B\subset X$ with $m(B)<\frac1k$
satisfying $P\Delta A\sqsubset B$. Put 
\begin{equation}
\label{eq:1107-2}
C=P\cap A\sqsubset X,\quad P'=P\setminus C, \quad A'=A\setminus C.
\end{equation}
Then
$$A',P'\sqsubset A'\cup P'=P\Delta A \sqsubset B.$$
Hence 
$$ A\setminus B\subset A\setminus A'=C\sqsubset P
= P'\cup C \subset B \cup C\subset B\cup A.$$

Conversely suppose for each \accessible{}  $k$ there is a subset $A,B$ with $m(B)<\frac1k$
satisfying (\ref{eq:1107-1}). Then 
$$A\Delta P\subset B.$$
In fact, if we define $C,A',P'$ by (\ref{eq:1107-2}), then it suffices to show
$A',P'\sqsubset B$. Since $A'\cap P=\emptyset$, 
$$
A'\setminus B\sqsubset A\setminus B \sqsubset P
$$
implies $A'\setminus B=\emptyset$, namely, $A'\subset B$. On the other hand
$$P'=P\setminus A\sqsubset (A\cup B)\setminus A\subset B.$$
\end{eproof}

\begin{Proposition}
\label{prop:join-of-measurable-semisets}
If $k$ is \accessible{}  and $P_i\sqsubset X$ ($i\in
[1..k]$) are measurable, then their intersection and union
are measurable.
\end{Proposition}
\begin{eproof}
There are subsets $A_i$ with $P_i\Delta A_i\aecong \emptyset$ for $i\in [1..k]$.
It suffices to show the following.
\begin{equation}
\label{eq:817-3}
\left(\bigcap_{i=1}^kP_i\right)\Delta \left(\bigcap_{i=1}^kA_i\right) 
\sqsubset
\bigcup_{i=1}^k\left(P_i\Delta A_i\right)
\end{equation}
\begin{equation}
\label{eq:817-4}
\left(\bigcup_{i=1}^kP_i\right)\Delta \left(\bigcup_{i=1}^kA_i\right) 
\sqsubset
\bigcup_{i=1}^k\left(P_i \Delta A_i\right)
\end{equation}

 To show (\ref{eq:817-3}), let $x$ be an element of the left hand side.
If $x\in \bigcap_i P_i$ and $x\nin \bigcap_iA_i$, then there is a $j$
such that $x\nin A_j$, whence $x\in P_j\Delta A_j$ and $x$ belongs to
the right hand side, which holds also in the case 
$x\nin \bigcap_i P_i$ and $x\in \bigcap_iA_i$ by similar arguments.

To show (\ref{eq:817-4}), let $x$ be an element of the left hand side.
If $x\in \bigcup_i P_i$ and $x\nin \bigcup_iA_i$, then there is
a $j$ with $x\in P_j$ and $x\nin A_j$, whence $x\in P_j\Delta A_j$ 
and $x$ belongs to the right hand side, which holds also in the case
$x\nin \bigcup_{i}P_i$ and $x\in \bigcup_iA_i$.
\end{eproof}

\begin{Proposition}
\label{prop:disjoint-union-mesurable-set}
If $k$ is \accessible{}  and $P_i\sqsubset X$ ($i\in
[1..k]$) are measurable and \textit{mutually almost disjoint}
\index{mutually almost disjoint@mutually almost disjoint}
 in the sense that
$P_i\bigcap P_j$ is null for $i\neq j$, then 
$$ \overline{m}(\bigcup_i P_i)=\sum_i \overline{m}(P_i). $$
\end{Proposition}
\begin{eproof}
For $i\in [1..k]$, let $A_i\subset X$ be a set satisfying 
$P_{i}\Delta A_i\aecong \emptyset$. Define
$$ B_i=P_i\bigcap A_i, P_i'=P_i\setminus B_i, A'_i=A_i\setminus B_i. $$
Since 
$$P'_i\bigcup A'_i=P_i\Delta A_i\aecong\emptyset$$
$P'_i$,$A'_i$ are null semisets. Since 
$$B_i\bigcap B_j\sqsubset P_i\bigcap P_j \aecong \emptyset$$
if $i\neq j$, we have 
$$A_i \bigcap A_j \sqsubset (B_i\bigcap B_j)\bigcup (A'_i\bigcup A'_{j}) \aecong \emptyset
$$
Hence, if we put 
$$C_i:=A_i\setminus \bigcup_{j>i}(A_i\bigcap A_j),$$
then $C_i$'s are mutually disjoint and $\bigcup_i C_i=\bigcup_iA_i$. Since 
$$ m(C_i)\leq m(A_i)\leq m(C_i)+ \sum_{j>i}m(A_j\bigcap A_i)\approx m(C_i),$$
we have $m(A_{i})\approx m(C_i)$. Hence
$$m(\bigcup_{i}A_i)= m(\bigcup_{i}C_i)
=\sum_im(C_i)\approx \sum_im(A_i)\approx \sum_im(P_i).$$
Hence, by (\ref{eq:817-4}) and Proposition~\ref{prop:boolean-properties-null-parts},
$$
m(\bigcup_i P_i)\approx m(\bigcup_{i}A_i)\approx \sum_i m(P_i).
$$ 

\end{eproof}

\index{sigma additivity@sigma additivity}
\begin{Theorem} 
\label{th:sigma-additivity}
If $\seqinf{P}$ is a concrete sequence of measurable subsemisets
and mutually almost disjoint, then the subsemiset
$$\unionP$$
is also measurable and 
$$
\overline{m}(\unionP)=\sum_{i=1}^{\infty}\overline{m}(P_i).
$$
\end{Theorem}

\begin{eproof}
For \accessible{}  $i$, choose $A_i\subset X$ satisfying
$$P_i\Delta A_i \aecong \emptyset.$$
Sine $A_i$ is a set, we extend the concrete  sequence $(\seqinf{A})$ to
a huge sequence of subsets $(\seq{M_0}{A})$. 

Put $b_p:=m(\bigcup_{1\leq i\leq p}A_i)$ for $p\in [1..M_0]$. 
By Lemma~\ref{th:convergence-of-increasing-sequence}
the increasing sequence of
rationals $(\seq{M_0}{b})$ has an upper bound $1$, whence it converges and
there is an $M_1\leq M_0$ such that $\lim_p b_p\approx b_K$ for all huge $K\leq M_1$. 

If $i$ is \accessible{} , it follows from $P_i\aecong A_i$ and (\ref{eq:817-3}) 
$$(A_i\cap A_j)\; \Delta\; (P_i\cap P_j) \sqsubset (A_i\Delta P_i)\cup  (A_j\Delta P_j)\aecong \emptyset,$$
whence if $i\neq j$ then $A_i\cap A_j \aecong P_i\cap P_j \aecong \emptyset.$ 
By Lemma~\ref{lemma:equivalence-relation-aecong}, we conclude
$A_i\cap A_j\approx \emptyset$,
for \accessible{}  $i\neq j$. By Proposition~\ref{prop:disjoint-union-mesurable-set}, 
for \accessible{}  $k$,
$$
m(\bigcup_{i\in [1..k]}A_i)\approx \sum_{i\in [1..k]}m(A_i),
$$
whence for some huge $M_2\leq M_1$, for every huge $L\leq M_2$
$$
m(\bigcup_{i\in [1..L]}A_i)\approx \sum_{i\in [1..L]}m(A_i).
$$
Since 
\begin{equation}
\label{eq:1107-20}
(\bigcup_{i\in [1..p-1]}A_i)
\cap 
(\bigcup_{i\in [p..\ell]}A_i)
\subset
\bigcup_{i\in [1..p-1],j\in [p..\ell]}(A_i\cap A_j)\aecong \emptyset,
\end{equation}
holds for \accessible{}  $\ell$, there is a huge $M_3\leq M_2$ such that
\eqref{eq:1107-20} holds also for huge $\ell\leq M_3$.
Thus if $M\leq M_3$  is huge we have
\begin{eqnarray*}
m\left(\bigcup_{i\in [1..M]}A_i\right)
&=&
m\left(\bigcup_{i\in [1..p-1]}A_i\right)
+
m\left(\bigcup_{i\in [p..M]}A_i\right)
-
m\left(
(\bigcup_{i\in [1..p-1]}A_i)
\cap 
(\bigcup_{i\in [p..M]}A_i)
\right)\\
&\approx&
m\left(\bigcup_{i\in [1..p-1]}A_i\right)
+
m\left(\bigcup_{i\in [p..M]}A_i\right).
\end{eqnarray*}

Hence 
\begin{eqnarray}
\label{eq:2}
\lim_pm\left(\bigcup_{i\in [p..M]}A_i\right)
&=&m\left(\bigcup_{i\in [1..M]}A_i\right)-
\lim_pm\left(\bigcup_{i\in [1..p-1]}A_i\right)
\nonumber
\\
\label{eq:20120331-40}
&=&b_M-\lim_p(b_p)\approx 0.
\end{eqnarray}

Put $A=\bigcup_{i=1}^{M_3}A_i$. We show
$$ \left(\unionP\right) \Delta A \aecong \emptyset$$

Now, using \eqref{eq:817-4}, we have 
{\footnotesize
\begin{eqnarray}
\left(\;\unionP\;\right)\;
\Delta \;A & = & 
\left(
  \left(\bigcup_{i=1}^{\ell}P_i\right)
    \bigcup 
  \left(\;\bigcup_{i=\ell+1}^{\infty}P_i\;\right)
\right)
\Delta
\left(
  \left(\bigcup_{i=1}^{\ell}A_i\right)
    \bigcup 
  \left(\;\bigcup_{i=\ell+1}^{M_3}A_i\;\right)
\right)
\nonumber
 \\
&\sqsubset&
\left(\;
  \left(\bigcup_{i=1}^{\ell}P_i\right)
     \Delta 
  \left(\bigcup_{i=1}^{\ell}A_i\right)
\;\right)
\bigcup 
\left(\;
  \left(\;\bigcup_{i=\ell+1}^{\infty}P_i\;\right)
  \Delta
  \left(\;\bigcup_{i=\ell+1}^{M_3}A_i\;\right)
\;\right)
\nonumber
\\
&\sqsubset&
\label{eq:825}
\left(\;
  \bigcup_{i=1}^{\ell}P_i\Delta A_i
\;\right)
\bigcup 
  \left(\;\bigcup_{i=\ell+1}^{\infty}P_i\setminus A_i\;\right)
\bigcup 
  \left(\;\bigcup_{i=\ell+1}^{M_3}A_i\;\right)
\end{eqnarray}
}
Since for each $i$, 
$$P_i\setminus A\sqsubset P_i\setminus A_i\sqsubset P_i\Delta A_i\aecong \emptyset,$$
the first and the second term in ~\eqref{eq:825} are null semisets. 
The measure of the last term converges to zero 
when $\ell\rightarrow \infty$ by \eqref{eq:20120331-40}.
Hence the left hand side $\left(\unionP\right)\Delta A$ is a null semiset. 
Thus $\unionP$ is measurable. 
By Proposition~\ref{prop:disjoint-union-mesurable-set},
\begin{eqnarray*}
\overline{m}(\unionP)
&=&[m(A)]=[b_{M_3}]=\lim_{p}[b_p]=\lim_p [m(\bigcup_{i=1}^{p}A_i)] \\
&=& \lim_p \overline{m}(\bigcup_{i=1}^{p}P_i)
= \lim_p \sum_{i=1}^{p}\overline{m}(P_i)=\sum_{i=1}^{\infty}\overline{m}(P_i).
\end{eqnarray*}
\end{eproof}

\begin{Theorem}
\label{th:join-of-increasing-measurable}
If $(\seqinf{P})$
is an increasing concrete sequence of measurable subsemisets, then
$\unionP$ is also measurable and 
$$ \overline{m}(\unionP)=\lim_i \overline{m}(P_i). $$
\end{Theorem}
\begin{eproof}
Define $Q_1=P_1$ and $Q_i=P_i\setminus P_{i-1}$ for $i>1$. Since $(\seqinf{Q})$
is a concrete sequence of subsemisets which is mutually disjoint, 
Theorem~\ref{th:join-of-increasing-measurable} implies that 
$$\unionP=\bigcup_{i=1}^{\infty}Q_i$$
is measurable and 
$$\overline{m}(\unionP)=
\overline{m}(\bigcup_{i=1}^{\infty}Q_i)
=
\lim_{p}\sum_{i=1}^p\overline{m}(Q_i)
=\
\lim_p\overline{m}(\bigcup_{i=1}^{p}Q_i)
=
\lim_p\overline{m}(P_p).$$
\end{eproof}

\index{sigma additivity@sigma additivity}
\begin{Theorem}
If $(\seqinf{P})$ is a concrete sequence of measurable subsemisets, then 
$\unionP$ is measurable and 
$$\overline{m}(\unionP)
\leq 
\sum_{i=1}^{\infty} \overline{m}(P_i).$$
\end{Theorem}
\begin{eproof}
Put $Q_p=\bigcup_{i=1}^{p}P_{i}$. Then $(\seqinf{Q})$ is increasing and 
$$\unionP=\bigcup_{i=1}^{\infty}Q_i$$
is measurable and 
$$\overline{m}(\unionP)
=
\overline{m}(\bigcup_{i=1}^{\infty}Q_i)
=
\lim_p\overline{m}(Q_p)
=\lim_{p}\overline{m}(\bigcup_{i=1}^{p}P_i).$$
Put $P'_1=P_1$ and $P'_i=P_i\setminus \bigcup_{j\in [1..i-1]}P_j$ for $i\geq 2$.
Then 
$$
\overline{m}(\bigcup_{i=1}^{p}P_i)
=\overline{m}(\bigcup_{i=1}^{p}P'_i)
=\sum_{i=1}^{p}\overline{m}(P'_i)
\leq \sum_{i=1}^{p}\overline{m}(P_i).
$$
Taking the limit, we obtain
$$
\overline{m}(\unionP)
=\lim_{p}\overline{m}(\bigcup_{i=1}^{p}P_i)
\leq \lim_{p}\sum_{i=1}^{p}\overline{m}(P_i)
=\sum_{i=1}^{\infty}\overline{m}(P_i).
$$
\end{eproof}

\subsection{Integration}
Let $(X,p)$ be a set $X$ with a probability density.
For a huge number $M$, the collection of functions $F_M(X):=Fun(X,[-M,M]_{\frac1M})$ 
is a set with $(M(2M+1))^{|X|}$ elements. 

For $f\in F_M(X)$ define its integration by
$$  E(f):=\sum_{x\in X}{p(x)f(x)}$$
and for $f,g\in F_M$ define
  $$d_1(f,g):=E(|f-g|).$$
If $d_1(f,0)<\infty$, we say $f$ is {\em integrable}. 
The class of integrable elements of $F_M(X)$ is denoted by $F^{int}_M(X)$.

\begin{Proposition}
\begin{enumerate}
\def\labelenumi{(\theenumi)}
\item\label{item:212}If $f,f'\in F_M(X)$ satisfies $f\approx f'$ and $f$ is integrable, 
then $f'$ is also integrable and
$E(f)\approx E(f')$.
\item\label{item:214} If $f_i\approx f'_i$ \itwo then 
$$ d_1(f_1,f_2)\approx d_1(f_1',f_2').$$
\item If $f_1,f_2\in F^{int}_M(X)$, then $d_1(f_1,f_2)$ is finite. 
\end{enumerate}
\end{Proposition}

\begin{eproof}
Put $\varepsilon:=\max_{x\in X}|f(x)-f'(x)|\approx0$. Then
$$ E(|f-f'|)|\leq E(\varepsilon)=\varepsilon.$$
Hence $E(f)\simeq E(f')$ and if $E(f)<\infty$ then $E(f')<\infty$.

Since $|f_1-f_2|\approx |f_1'-f_2'|$, we have
$$d_1(f_1,f_2)\approx d_1(f_1',f_2').$$

If $f_1,f_2\in F^{int}_M(X)$, then
$d_1(f_1,f_2)=E(|f_1-f_2|)\leq E(|f_1|)+E(|f_2|) $ is finite.
\end{eproof}

\begin{Proposition}
$$E(f)=\sum_{\lambda}\lambda m(f^{-1}(\lambda))
$$
\end{Proposition}
\begin{eproof} Since $X=\coprod_{\lambda}f^{-1}(\lambda)$, 
we have
$$
E(f)=\sum_{\lambda}\sum_{f(x)=\lambda}f(x)p(x)
=\sum_{\lambda}\lambda \sum_{f(x)=\lambda}p(x)
=\sum_{\lambda}\lambda m(f^{-1}(\lambda)).
$$
\end{eproof}

\begin{Proposition}[Chebishev Inequatlity]
\label{prop:Chebishev} Suppose $f$ is integrable, 
$f\geq 0$ and $c$ is a positive rational.
Then
$$
m(\setii{x}{f(x)\geq c})\leq \Frac{E(f)}{c}.
$$
\end{Proposition}

\begin{eproof}
$$
E(f)\geq \sum_{\lambda\geq c}\lambda m(f^{-1}(\lambda))  
\geq c\sum_{\lambda\geq c}m(f^{-1}(\lambda))=c\;m(\setii{x}{f(x)\geq c}).
$$
\end{eproof}

\begin{Theorem}
\label{th:consequence-of-E(h)=0} If $f,g\in F^{int}_M(X)$ satisfies
$d_1(f,g)\approx 0$, then $f\approx g$ a.e..
\end{Theorem}
\begin{eproof}
It suffices to show that if $E(h)\approx0$ then $h\approx 0$ a.e..

Suppose $E(h)\approx 0$. For \accessible{}  $k$, 
$$m\left(\setii{x}{|h(x)|>\frac1k}\right)\leq k E(h)\approx 0.$$
Hence the set 
$$\setii{k}{m\left(\setii{x}{|h(x)|>\frac1k}\right)\leq \frac1k}  $$
includes all \accessible{}  numbers and hence also a huge number $K$.
Put $A:=\setii{x}{|h(x)|>\frac1K}$.
Then $m(A)\leq \frac1K \approx 0$. Moreover if $x\nin A$ then $|h(x)|\leq \frac1K\approx 0$. 
Hence $h\approx 0$ a.e..
\end{eproof}

\begin{Remark} The converse does not hold, namely,  
$E(f)\simeq E(g)$ does not hold always even if $f\approx g$ a.e.. 
For example, suppose $p$ is uniform distribution, namely, $p(x)=\frac1{|X|}$.
Define
$$
f(x):= \left\{
  \begin{array}{cc}
  |X|  & \mbox{for } x=x_0 \\
  0  & \mbox{otherwise} 
  \end{array}
\right.
$$
Then $f\approx0$ a.e. but $1=E(f)\not\approx E(0)=0$. 

If $f,g$ are $L^1$ functions defined in the next subsection, then
$d_1(f,g)\approx 0$ if and only if $f\approx g$ a.e..
\end{Remark}

\subsection{$L^1$-functions}

For $f\in F_M(X)$ and a rational $a\in \kukanii{M}{M}$, 
define 
$$
f^{a}(x):= \left\{
  \begin{array}{cc}
0    & \mbox{if $|f(x)|>a$} \\
f(x)    & \mbox{otherwise} 
  \end{array}
\right.
$$

\begin{Definition} 
An $f\in F_M(X)$ is called an $L^{1}$ function if $E(|f-f^{a}|)\approx 0$ for huge rationals $a$.
\end{Definition}
\begin{Lemma}
$f\in F_M(X)$ is an $L^1$ function if and only if for huge rationals $a$
$$
\sum_{|f(x)|>a}|f(x)|p(x)\approx 0,
$$
if and only if
$$ \sum_{\lambda>a}\lambda m(|f|^{-1}\lambda)\approx 0 $$
\end{Lemma}
\begin{eproof}
Obvious since 
$$
E(|f-f^a|)
=\sum_{|f(x)|>a}|f(x)|p(x)
=\sum_{\lambda>a}\sum_{|f(x)|=\lambda}\lambda p(x)
=\sum_{\lambda>a}\lambda m(|f|^{-1}\lambda).
$$
\end{eproof}

\begin{Proposition}
\label{prop:characterization-of-L1}
An $f\in F_M(X)$ is an $L^1$ function, if and only if 
the following conditions are satisfied.
\begin{enumerate}
\item\label{item:174} $f\in F_{M}^{int}$,
\item\label{item:175} $E(|f|\chi_{A})\approx 0$ if $m(A)\approx 0$,
where $\chi_A$ denotes the characteristic function of $A\subset X$.
\end{enumerate}
\end{Proposition}
\begin{eproof}
Suppose $f$ is an $L^1$ function. Then for every huge $K$,
 $$E(|f|)\leq E(|f-f^a|)+E(|f^a|)\approx E(|f^a|)\leq K.$$
Hence there is an \accessible{} $k$ with $E(|f|)<k$.

Suppose a subset $A\subset X$ satisfies $m(A)\approx 0$.
Then 
  $$
E(|f|\chi_A)\leq E(|f-f^a|\chi_A)+E(|f^{a}|\chi_A)
\leq E(|f-f^a|)+a\;m(A)\approx0.
$$

Conversely suppose the conditions~\ref{item:174},~\ref{item:175} are satisfied.
Suppose $E(|f-f^a|)\succ 0$ for some huge $a$.
Then the subset
$$A:=\setii{x}{f(x)\geq a}$$
satisfies, by Proposition~\ref{prop:Chebishev}, 
$$m(A)\leq \Frac{E(f)}{a}\approx 0.$$
On the other hand, since $f^a=0$ on $A$ and $f-f^{a}=0$ on $A^{c}$, we have 
$$ E(|f|\chi_A)=E(|f-f^a|\chi_A)=E(|f-f^a|)\succ 0.$$
This contradicts the latter assumption.
\end{eproof}

\begin{Corollary}
\begin{itemize}
\item If $f,g$ are $L^1$ functions, then $f+g$ is also an $L^1$ function.
\item If $f$ is an $L^1$ function and $g$ is finite, then $fg$ is also an $L^1$ function.
\item If $f$ is an $L^1$ function and $|g|\leq |f|$, then $g$ is also an $L^1$ function.
\end{itemize}
\end{Corollary}

\begin{eproof}
Suppose $f,g$ are $L^1$ functions. From $E(|f+g|)\leq E(|f|)+E(|g|)$, it follows 
$f+g$ is integrable. On the other hand if $m(A)\approx0$ then 
$E(|f+g|\chi_A)\leq E(|f|\chi_A)+E(|g|\chi_A)\approx0$. Hence $f+g$ is an $L^1$ function.

Suppose $f$ is an $L^1$ function and $g$ is finite. Let $k$ be an \accessible{} number such 
that $|g(x)|<k$ for all $x\in X$. Then
$
E(|fg|)\leq k E(|f|),
$
hence $fg$ is integrable. On the other hand if $m(A)\approx 0$ then 
$$ E(|fg|\chi_A)\leq k E(|f|\chi_A)\approx0.$$
Hence $fg$ is an $L^1$ function.

Suppose $f$ is an $L^1$ function and $|g|\leq |f|$. Then $E(|g|)\leq E(|f|)<\infty$.
If $m(A)\approx 0$, then 
$$E(|g|\chi_A)\leq E(|f|\chi_A)\approx 0.$$
\end{eproof}

\begin{Proposition}
If $f$ is an $L^1$  function and $d_1(f,g)\approx 0$, then $g$ is also an $L^1$ function.
\end{Proposition}

\begin{eproof} Suppose $f$ is an $L^1$ function and $d_1(f,g)\approx 0$. 
Then $g$ is integrable since $f$ is integrable. Moreover if
$m(A)\approx 0$ then 
$$E(|g|\chi_A)\leq E(|f-g|\chi_A)+ E(|f|\chi_A)\leq d_1(f,g)+E(|f|\chi_A)\approx 0.$$
\end{eproof}

The class of $L^1$ functions forms a subclass of the set $F_M(X),d_1)$, which defines
a continuum denoted by $L^1_M(X)$. 

\begin{Theorem}
\label{th:characterization-of-E(h)=0}
If $f,g\in L_{M}^{1}(X)$, then $d_1(f,g)\approx 0$ if $f\approx g$ a.e..
\end{Theorem}
\begin{eproof}
Suppose $h\in L^1_{M}(X)$ satisfies $h\approx 0$ a.e..
By Lemma~\ref{lem:a-e-Robinson}, there is a subset $A\subset X$ such that 
$m(A)\approx 0$ and $h(x)\approx 0$ if $x\nin A$. Then 
$$E(|h|)=E(|h|\chi_A)+E(|h|\chi_{A^{c}}) \approx E(|h|\chi_A)\approx 0.$$
\end{eproof}

\begin{Theorem}
\label{thm:l1-functions-pointwise-convergence}
If $(\seqinf{f})$ is a concrete sequence of $L^1$ functions and 
converges to an $L^1$ function $f$ pointwise a.e., then 
$$\lim_i E(|f_i-f|)\approx 0.$$
In particular
$$\lim_i E(f_i)\approx E(f).$$
\end{Theorem}

\begin{eproof}
Suppose $f_i$ converges a.e. to $f$ pointwise.
By Lemma~\ref{lem:207-807} below,
there is a subset $A\subset X$ such that
$m(A)\approx 0$ and if $x\nin A$ then $f_i(x)$ converges to $f(x)$.
Then by Lemma~\ref{lem:pointwise2uniform}, $f_i$ converges uniformaly 
to $f$ on $A^c$. Hence for each \accessible{} $k$, there is an \accessible{} $n$ such 
that if $i$ is \accessible{} with $i>n$ then, 
$|f-f_i|<\frac1k$ on $A^c$ and we have 
$$E(|f-f_i|)=E(|f-f_i|\chi_A)+E(|f-f_i|\chi_{A^c})
\leq \frac1k+ E(|f-f_i|\chi_A). $$
Since $|f-f_i|$ is an $L^1$ function, by Proposition~\ref{prop:characterization-of-L1} we have
$$E(|f-f_i|\chi_A)\approx 0.$$
Thus for \accessible{}  $k$, there is an \accessible{}  $n$ such that for every \accessible{} 
$i\geq n$, we have $E(|f-f_i|)\leq \frac2k$, namely, 
the concrete sequence $\seqinf{f}$ converges to $f$ with respect to the metric $d_1$.
\end{eproof}

\begin{Lemma}\label{lem:pointwise2uniform}
If a concrete sequence $\seqinf{f}$ of functions on a set $X$ converges 
pointwise to $f$ everywhere, then it converges uniformly to $f$.
\end{Lemma}
\begin{eproof}
Suppose $f_i$ converges to $f$ pointwise. For each \accessible{} $k$ and $x\in X$, 
there is an \accessible{} $\ell_x$ such that if $n>\ell_x$ then 
\begin{equation}
\label{eq:120612-1}
 |f_n(x)-f(x)|<\frac1k.
\end{equation}
If we put $\ell=\max\setii{\ell_{x}}{x\in X}$, then for $n>\ell$ 
conditions \eqref{eq:120612-1} holds for every $x\in X$. Namely the
sequence $f_i$ converges to $f$ uniformely. 
\end{eproof}

\begin{Lemma}
\label{item:207-807}\label{lem:207-807}
\label{lem:converge-exceptional-set}
If a concrete sequence $\seqinf{f}$ of functions on $X$
converges pointwise to $f$ almost everywhere,
then there is a subset $A\subset X$ such that $m(A)\approx 0$ and
the concrete sequence of rationals $f_i(x)$ converges to $f(x)$ for $x\nin A$.
\end{Lemma}

\begin{eproof}
Suppose a concrete sequence $f_i$ converges almost everywhere to $f$.

Then for each \accessible{}  $k$, there is a subset $A_k\subset X$ satisfying
\begin{enumerate}
\item\label{item:204} $m(A_k)<\frac1k$
\item\label{item:205} If $x\nin A_k$ then $(f_1(x),f_2(x),\cdots )$ converges to 
$f(x)$.
\end{enumerate}
We may suppose $A_k$ is decreasing since $A'_k=A_1\cap A_2\cap \cdots\cap A_k$ 
satisfies the same conditions.

Extend the concrete sequence $(\seqinf{f})$ to a sequence $(\seq{N}{f})$
in $F_M(X)$. If $x\nin A_k$, then there is a huge $M_x$ such that for
all huge $i\leq M_x$,
\begin{equation}
\label{eq:9-807}
  f_{i}(x)\approx f(x)  
\end{equation}
holds. Put $M_{k}=\min_{x\nin A_k}M_{x}$. Then for huge $i\leq M_k$,
~\eqref{eq:9-807} holds for all $x\nin A_k$. Since
\begin{equation}
\label{eq:10-807}
  |f_{i}(x)-f(x)|<\frac1k
\end{equation}
holds for all huge $i\leq M_k$, there is an \accessible{}  $m_k$ such that 
~\eqref{eq:10-807} holds for all $i$ with $m_k\leq i\leq M_k$. Note that we may take
$m_k>k$. 

Thus for every \accessible{}  $k$, there are an \accessible{} number $m_k$, a huge number $M_k$,
and subset $A_k\subset X$ satisfying the following conditions.
\begin{enumerate}
\item $k<m_k$,
\item $m(A_k)<\frac1k $,
\item the condition ~\eqref{eq:10-807} holds for all $i$ satisfying $m_k\leq i\leq M_k$,
\item if $i<j\leq k$ then $m_i\leq m_j<M_j\leq M_i$,$A_j\subset A_i$.
\end{enumerate}
Since these conditions on $k$  are objective, we have a huge $K$ such that 
there are numbers $m_K,M_K$ and a subset $A_K\subset X$ 
satisfying the above four conditions for $k=K$.

Then if $x\nin A_K$, for every huge $i\leq M_K$ and \accessible{}  $k$, 
we have  $i\leq M_K\leq M_k$ and $A_K\subset A_k$, whence
(\ref{eq:10-807}) holds whence (\ref{eq:9-807}) holds.  
Thus if $x\nin A_{K}$ then $\lim_if_i(x)\approx f(x)$. 
\end{eproof}

We note that if a concrete sequence of $L^1$ functions converges 
everywhere to a function $f$ then the following can be proved easily.

\begin{Proposition}
Let $(\seqinf{f})$ be a concrete sequence of $L^1$ functions 
converging to $f$ pointwise everywhere. Then $f$ is an $L^{1}$ function and
$$\lim_i E(|f-f_i|)\approx 0.$$
\end{Proposition}

\begin{eproof}
By Lemma~\ref{lem:pointwise2uniform}, $f_i$ converges uniformly to $f$. 
In particular, there is an $L^1$ function $g=f_i$ such that  $|f(x)-g(x)|<1$.
Hence 
$$|f(x)|\leq 1+ |g(x)|$$
for all $x$ and $f$ is an $L^1$ function. Since $f_i$ converges uniformly to $f$,
it is obvious that
$$\lim_i E(|f-f_i|)\approx 0.$$

\end{eproof}

\begin{Theorem}
\label{thm:l1-functions-convergence}
Suppose a concrete sequence $(\seqinf{f})$ of $L^1$ functions
converges to $g\in F_M(X)$ with respect to the distance $d_1$. 
Then $g$ is also an $L^{1}$ function and a subsequence of
$\seqinf{f}$ converges to $g$ pointwise a.e..
\end{Theorem}
\begin{eproof}
For each \accessible{}  $k$, there is an \accessible{} $n_k$ such that 
\begin{equation}
\label{eq:1111-1}
E(|f_i-g|)\leq \frac{1}{4^k} \mbox{ for \accessible{} }i\geq n_k. 
\end{equation}

The function $g$ is integrable since
$$E(|g|)\leq E(|g-f_{n_i}|)+E(|f_{n_i}|)\leq \frac1{4^i}+E(|f_{n_i}|).$$

Furthermore, if $C$ is any subset of $X$ with $m(C)\approx 0$, then 
since $f_{n_i}$ is an $L^1$ function , we have $E(|f_{n_i}|\chi_C)\approx0$,
and for every \accessible{}  $i$
$$E(|g|\chi_{C})\leq E(|g-f_{n_i}|\chi_{C})+E(|f_{n_i}|\chi_{C})
\approx E(|g-f_{n_i}|\chi_{C})
\leq E(|g-f_{n_i}|)\leq \frac1{4^i}.
$$
Hence 
$$E(|g|\chi_{C})\approx0.$$
By Proposition~\ref{prop:characterization-of-L1}, $g$ is also an $L^1$ function.

Put 
\begin{equation*}
 A_k:=\setii{x\in X}{|f_{n_k}(x)-g(x)|\geq \frac1{2^k}},
\end{equation*}
then by Proposition\ref{prop:Chebishev}
$$
m(A_k)\leq \frac{\frac{1}{4^k}}{\frac1{2^k}}
=\Frac1{2^k}.
$$
Extend the concrete sequence $\seqinf{A}$ to a huge sequence $\seq{K}{A}$ such that
for all $i\leq K$, $m(A_i)\leq \frac1{2^i}$.

Put $B_i=\bigcup_{i\leq j\leq K}A_j$.  
Then
\begin{equation}
 m(B_i)\leq \sum_{i\leq j\leq K}\Frac1{2^{j}}\approx \frac1{2^{i-1}}.
\label{1111-2}
\end{equation}

Put $g_i:=f_{n_i}$. 
Suppose $x\nin B_i$. Since for all concrete $j\geq i$, we have $x\nin A_j$, whence
$$
|g_j(x)-g(x)|<\frac1{2^{j}}\leq \frac1{2^i},
$$
which implies that $(g_1(x),g_2(x),\cdots)$ converges to $g(x)$.
Hence by (\ref{1111-2}), the subsequence 
$(f_{n_1},f_{n_2},\cdots)$ converges pointwise to $g$ almost everywhere.

\end{eproof}

\begin{Corollary}[Completeness of $L^1_M(X)$]
\label{cor:L1-is-complete}
The metric space $L^1(X)$ is complete, namely, every concrete Cauchy sequence of 
$L^{1}$ functions with respect to $d_1$ converges to an $L^1$ function.
\end{Corollary}
\begin{eproof}
Suppose  $\seqinf{f}$ is a concrete Cauchy sequence of 
$L^{1}$ functions with respect to $d_1$. Extend to a
huge sequence $(\seq{N}{f})$ in $F_M(X)$, which is
convergent by Proposition~\ref{prop:convergence}. Hence $(\seq{N}{f})$ 
converges to a $g\in F_M(X)$ with respect to $d_1$. Then
the concrete sequence $\seqinf{f}$ converges to $g$ by the same proposition.
By Theorem~\ref{thm:l1-functions-convergence}, $g$ is an $L^1$ function.
Hence $\seqinf{f}$ converges in $L^1_M(X)$.
\end{eproof}

\begin{Remark}
If $L^1_M(X)$ were a set, then completeness follows from Proposition~\ref{prop:convergence}
directly. Since $L^1_M(X)$ is a proper semiset which cannot be 
defined by objective conditions, 
the extention of a concrete sequence to huge sequence is not possible in $L^1_M(X)$.
\end{Remark}

\newpage
\section{Concluding Remarks}
\def\runningtitle{Concluding Remarks}
\subsection{Recapitulation}
We showed that basic mathematical concepts with infinitary aspects
such as real number, calculus, topology, measure can be developed by
replacing the infinite axiom by the ``sorites axiom'' giving
qualitative plurality of finiteness.  We gave the terminology
``standard'' a semantical meaning of accessibility in order to make the
``validity'' of basic axioms obvious. We hope by this strategy the
``over-technicality'' of the traditional axiomatic foundation of
nonstandard mathematics is reduced considerably.

\subsubsection{Vague Concepts and Semisets}
The crucial point in order to actualize directly
the qualitative plurality of finiteness is to use vague conditions such as
``accessibility'' side by side with the usual mathematical
conditions. The essential difference between these types of conditions
give rises to the so-called ``overspill phenomena'', which turn out 
to be one of the basic principles in nonstandard mathematics.

However, logical usage of vague concepts needs drastic change of
basic concepts and principles of mathematics, although it should be
stressed that the change is of such a kind as to reduce unnecessary
complication of current mathematics resulting from not discriminating
between theoretical possibility and actual possibility.  

The most radical change is the introduction of ``proper semisets'' which are
proper classes included in a huge finite set. A typical example of a
proper semiset is the collection of accessible natural numbers. We have seen
that semisets play the role of the infinite sets but in a more
appropriate way since inaccessible numbers are not separated from the
accessible numbers by virtue of the overspill principle.

On this account we introduced three types of collection, namely,
proper classes, semisets, and sets. Furthermore only finite sets are
entitled to be sets but finite sets are ramified to huge sets and
concrete ones which have actual enumeration. Parenthetically we note 
that we used the word ``collection'' only informally.

Furthermore we introduced two attributes of conditions, objective and
definite. A condition is definite if it can be stated without
unbounded quantifications and objective if it can be sated without
using accessibility. We restricted the separation axiom only to objective
conditions, which implied the overspill principle.

Another deviation from the current mathematics is the understanding of
functions. We require functions to have explicit objective definite
specification since the usual notion of function as mapping have no
definite semantical meaning for proper classes. We checked that this
restriction is void for functions defined on sets but that functions
defined on semisets are extended uniquely to a mapping on surrounding
sets, which is revealed to be another basic principle in addition to
the overspill principle.

All these changes might be appear unnecessary complication at first
sight but they reflect important aspects of our way of understanding
the world and as a result they enrich mathematics with more intuitive
ways of arguments hitherto considered as merely informal ones.

\subsubsection{Treatment as Naive theory}
We did not and do not intend to present our new framework as a formal
theory. The reasons are as follows.

Firstly the so called ``formal theory'' itself depends on the current
mathematics with the doctrine of ``the $\nat$''. For example, even
syntactic concept as the provability is ramified in our alternative
mathematics with multiple concepts of finiteness. Hence ``formal
theory'' itself is not reliable from the point of view of the alternative
approach presented here.

But the more decisive reason is that our intention is to present an
alternative approach to mathematics without technical artifacts so
that freshmen could follow 
without specialized training of specific area such as mathematical
logic. We intend to grow an alternative mathematics for ``doing
mathematics'' just as current mathematics are done mostly by 
naive set theory without exact knowledge of axiomatic set theory.

\subsubsection{Continua and Points}
A continuum is usually identified with the infinite set consisting of
its points and the topology is captured as an additional structure
given by metric or topological structure.  In contrast we captured a
continuum as a collection, often a huge finite set, endowed with
indistinguishability relation. Here distinguishability is understood
from practical point of view and discriminated from theoretical
distinguishability. As a result the indistinguishability is a vague
relation, and the identity of a point on a continuum, defined as the
collection of elements indistinguishable from a fixed element, has
persistent indefiniteness, which is embodied in the sorites paradox
that both $x_1\approx x_2\approx \cdots \approx x_N$ and
$x_1\not\approx x_N$ can hold if $N$ is not accessible. In addition to
the indefiniteness, a point of a continuum itself has a structure of a
continuum if the indistinguishability is properly sharpened, which
reflects the fractal nature of continuum. 

We did not however try to define general continua as primitive entity
but only defined the linear continuum $\real$ as 
the proper class of rational numbers with the indistinguishability
relation. Subcontinua of $\real$ such as the unit interval
$[0,1]$ can be represented by finite sets with indistinguishability. 

Since points have nontrivial extensions,
a morphism between continua cannot be determined as a correspondence
of points. Besides indistinguishable elements must correspond to
indistinguishable ones, whence morphisms are represented by continuous
maps in the usual sense. Discontinuity simply means ill-definedness. 

We may say that nonstandard mathematics, by embracing indefiniteness
via ``standardness'', has given alternative approach to continuum 
more appropriate not only than current mathematics 
but also than the invaluable intuitionistic mathematics.

We note in passing that there is constructive approaches to 
nonstandard mathematics, for example \cite{MR1336645}, based on
intuitionistic type theory \cite{martin1990mathematics}, 
which however appear to be rather too formal to be relevant 
to the above mentioned intention to develop
alternatives mathematics at the same naive level as the usual one, 
based on intuitively clear simple semantics.

\subsubsection{Idealization}
We postulated that a number less than an accessible number is also
accessible. However, for example, most of the numbers less than
$10^{10^{10}}$, accessible by the exponentiation, cannot be described
concretely in any fixed notational system.  Thus our accessibility is
too idealized to have something to do with actual accessibility from
the radical ultrafinistic point of view. 

Moreover the existence of the huge numbers inaccessible by any
concrete methods might seem similar idealization as in the
introduction of the infinite sets. However the character of
idealization is utterly different. The idealization of infinity
as infinite sets is to regard essentially indefinite objects
as definite ones whereas  idealization of infinite as
huge sets keeps the indefiniteness so that it has potentially vast
superiority over the infinite sets. For example it is intuitively more
acceptable and more importantly it is safer from contradiction.
This allows us not to consider the coherence problem so seriously.

On might think that the usage of vague concepts might be a new
potential source of incoherency. In this respect, the usual coherent
usage of the terminology ``standard'' which is vague, in the sense that
it has no extension, gives us psychological assurance of our treatment of
accessibility, since ours are in a sense a tiny portion of most
axiomatic systems of nonstandard mathematics currently used.

\subsubsection{Transfer Principle}
We did not mention ``the transfer principle'' usually considered as
the key point of nonstandard mathematics. However we found that its
importance comes only from the requirement for nonstandard mathematics
to be conservative extension of current mathematics and transfer
principle is not necessary in developing mathematics itself.

However we briefly show that the transfer principle for definite
objective conditions is a trivial
consequence of the concept of accessibility if it is 
interpreted as the possibility of explicit specification.

Suppose $P(x)$ is a definite objective condition with all the
parameters accessible\footnote{Note that if $\Omega$ is a huge number then the definite condition ``$x<\Omega$'' is satisfied by all accessible number but is not by
$\Omega$. 
} and is satisfied by all \accessible{} numbers.
If some inaccessible numbers does not satisfy it, then the 
minimal numbers which do not satisfy $P$ is accessible by definition,
which contradicts to the assumption.

However we did not restrict the meaning of accessibility in order only
to get the transfer principle, which is not necessary if we do not
insist on the conservativeness of the nonstandard mathematics.

\subsection{Future Direction}
We remark on some of the important aspects not touched here and some
of the future promising directions .

\subsubsection{Accessibility of Higher Order Objects} 
Many arguments of nonstandard mathematics are carried over to our
framework except for those dependent on the standard-part operation.
For example, the compactness is usually defined by the condition that
every element is near-standard, that is, indistinguishable from a
standard one. 
However we have more intuitive characterization of compactness, 
as a sort of ``pigeon principle'', 
namely a continuum is compact if every huge subset has at least two 
mutually indistinguishable elements. 

Since not only ``standard-part operation'' but also 
the concept standardness itself applied to higher order objects 
such as sets and functions 
seems to result in undesirable technicalities in usual nonstandard mathematics.

We guess that the counterpart of the ``standard part arguments'' in our
framework is given by a sort of constructivity as is seen in the 
following examples. 

Suppose an increasing family of
sets $X(n)$ parametrized by accessible numbers $n$ are constructed by
a method independent of the specificity of the number $n$.  By
extension principle, we can substitute a huge number $\Omega$ to
obtain a huge set $X(\Omega)$. Such huge sets are considered to be
constructed by the series $\seti{X(n)}$.

Two different constructions of a set can be considered as its 
different structures. For example the huge set $[1..2^{\Omega}]$
have two constructions, one by substituting $2^{\Omega}$ to $n$ in $[1..n]$ 
and the other by substituting $\Omega$ in $n$ of $[1..2^{n}]$. 
The former series is constructed by adding $n+1$ to $[1..n]$.
The latter is constructed first by regarding $[1..2^{n}]$ 
as the set of infinitesimal intervals of width $2^{-n}$ 
and the step from $[1..2^{n}]$ to $[1..2^{n+1}]$ 
is done by halving all the intervals.

The vectors in $\rational^{[1..{\Omega}]}$ are sequences of rationals 
$(\seq{\Omega}{a})$ of length $\Omega$. 
The meaning of their accessibility depends on to the way the huge set 
$[1..\Omega]$ is constructed. 
For example, if $[1..\Omega]$ is constructed by the series $\seti{[1..n]}$ as above, 
then a vector $(\seq{\Omega}{a})$ is accessible 
if its essential is captured by the subsequences $(\seq{n}{a})$. 
These vectors form the $L^1$ space of convergent sequences. 

If $[1..2^\Omega]$ is the collection of infinitesimal intervals 
obtained by the halving processes as above,
then the vector $(\seq{2^{L}}{a})$, 
considered as a function which is constant on the intervals of width $2^{-L}$, 
is accessible if its essential part is captured by the functions constant
on the intervals of width $2^{-n}$ with accessible $n$. 
These vectors form the space of measurable functions on $[0,1]$.

The above sort of ``accessible'' sequences might be considered 
as typical ones of those usually called near-standard.  We guess that standard-part
operation can be captured by incorporating the constructivity of huge
finite sets as above. In this way, it seems that we arrive at
a mathematics which have much in common with the constructive
nature of the intuitionistic approach.

\subsubsection{Relative Accessibility}

There are now various frameworks which relativize standardness such as 
RST (Relative Set Theory ) of P\'eraire\cite{peraire1992théorie},
EST (Enlargement Set Theory ) of D. Ballard \cite{ballard1994foundational}, 
relative arithmetic of S. Sandars \cite{sanders2010relative} to mention a few.

Similarly we can relativize accessibility as follows.  Define a
binary relation ``number $y$ is \textit{accessible from a number} $x$
or simply \textit{$x$-accessible}'' if there is some method of
reaching $y$ using $x$ and the numbers less than $x$.  A number $y$ is
called $x$-inaccessible and written $y\gg x$ if it is not accessible
from $x$. The $1$-accessible numbers are accessible numbers of
\S~\ref{subsub:accessibility}.

The axioms are relativized as follows. 
A rational number $x$ is called \textit{infinitesimal at the level} $x$, or
simply \textit{$x$-infinitesimal} if $|x|<\frac1k$ for all numbers $k$
accessible from $x$. Two rational numbers $y,z$ is called \textit{$x$-indistinguishable}
and written $y\approx_x z$ if $y-z$ is $x$-infinitesimal.

We postulate that for every number $x$, there are $x$-inaccessible
numbers.  A condition is called \textit{objective} if it is defined
without using the binary relation $x\gg y$. We postulate the
separation axiom for objective conditions so that the collection of
elements in a set satisfying an objective condition is a set. This
implies the general overspill principle to the effect that if an
objective condition is satisfied by all the $x$-accessible numbers
then it is satisfied also by an $x$-inaccessible number.

\subsubsection{Continua of Syntactic Objects}
One of the innovative aspects of our approach is the possibility of
using huge syntactic objects, such as huge words, huge terms. In
contrast to the "infinite words", every operations on finite words
carry over to them and it is expected the mathematical world of huge
syntactic objects has new phenomena with both aspects of finite and
infinite and give new insights into the mysteries of the "complex
systems" for whose understanding 
the dichotomy between finite and infinite is a severe barrier.

\begin{figure}[htbp]
\includegraphics[width=0.3\textwidth]{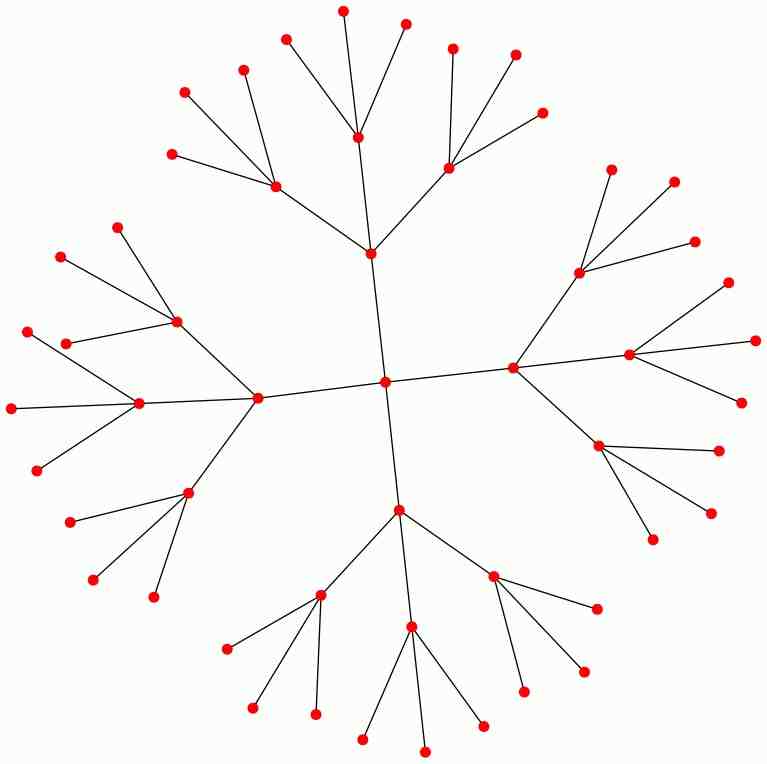}
\includegraphics[width=0.3\textwidth]{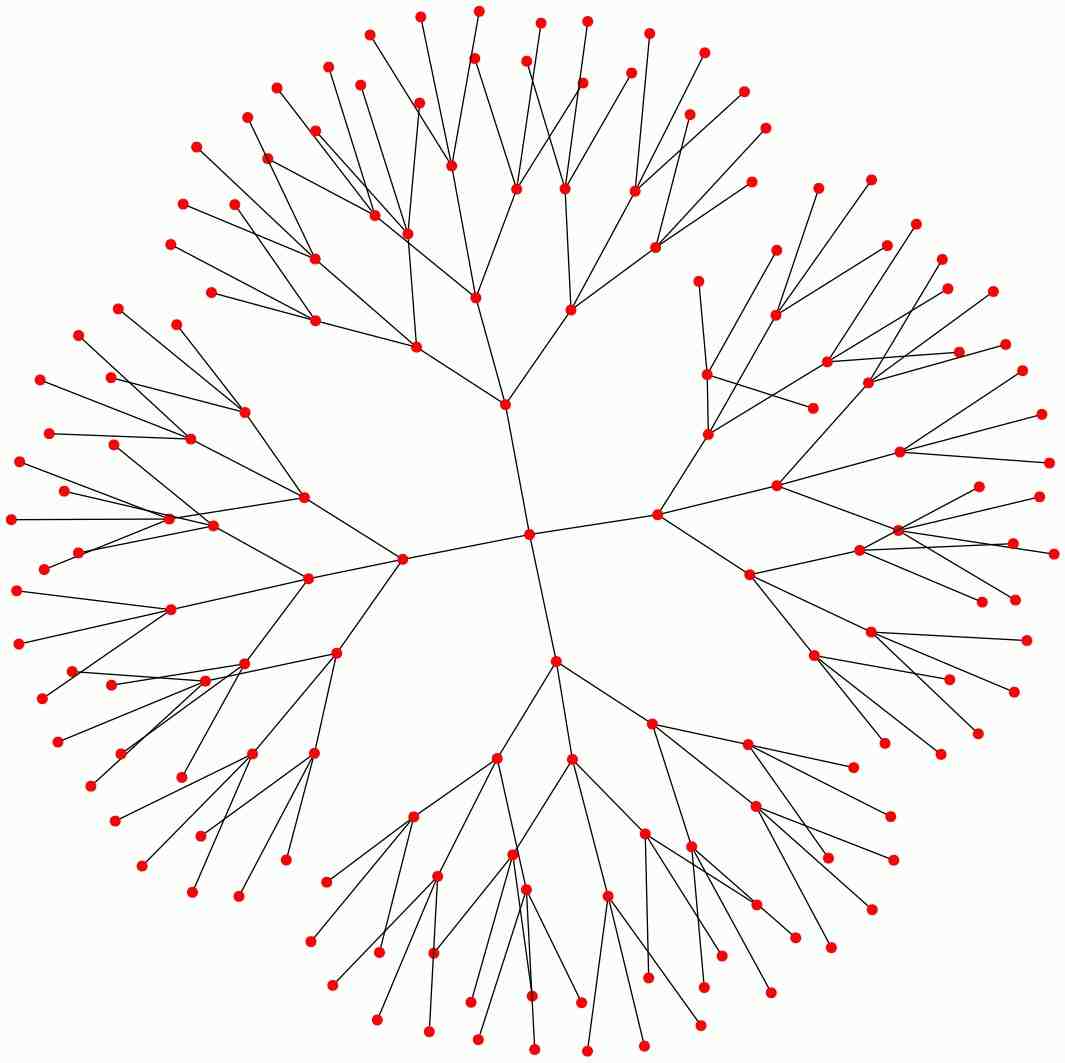}
\includegraphics[width=0.3\textwidth]{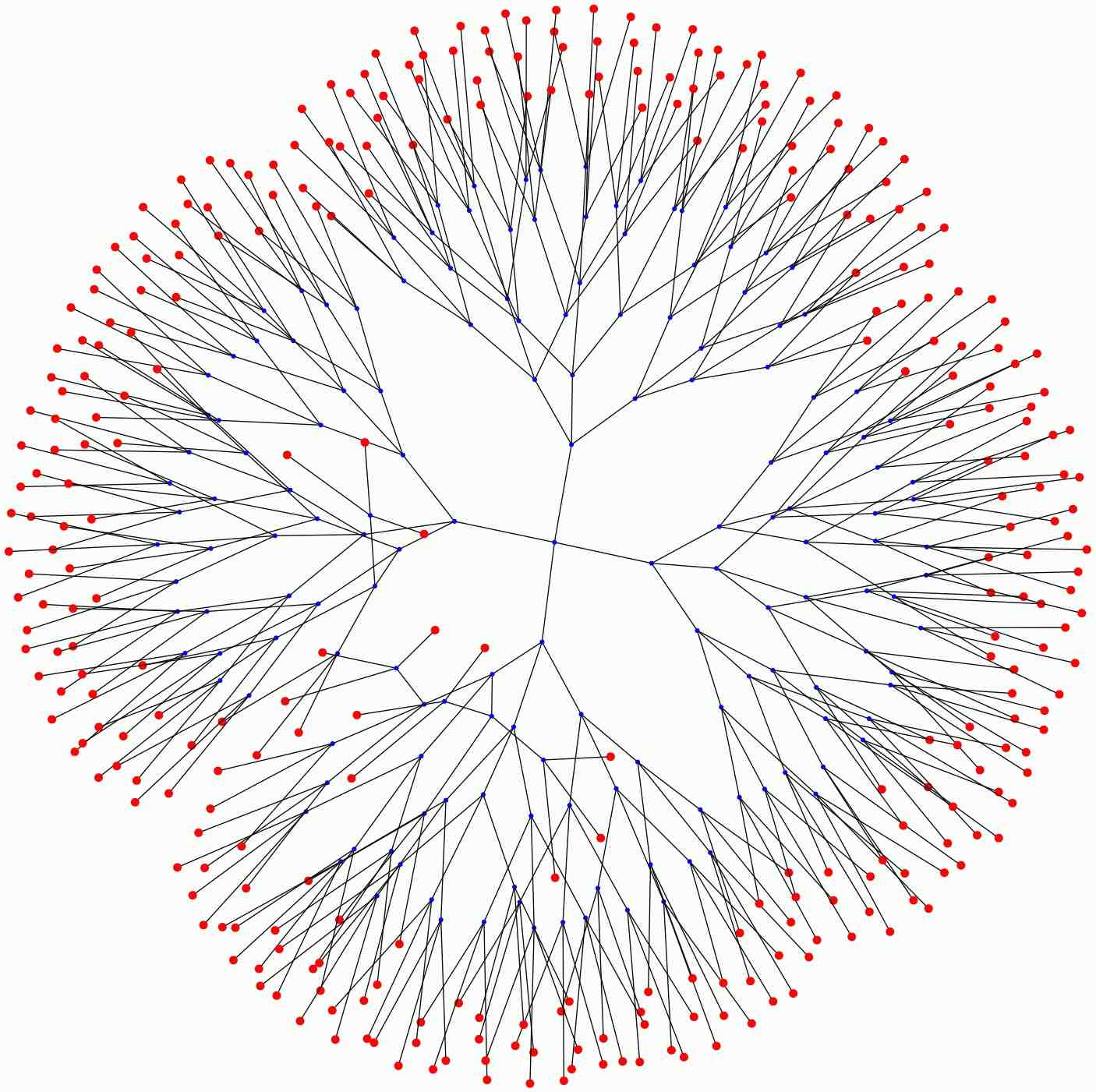}\\
\includegraphics[width=0.5\textwidth]{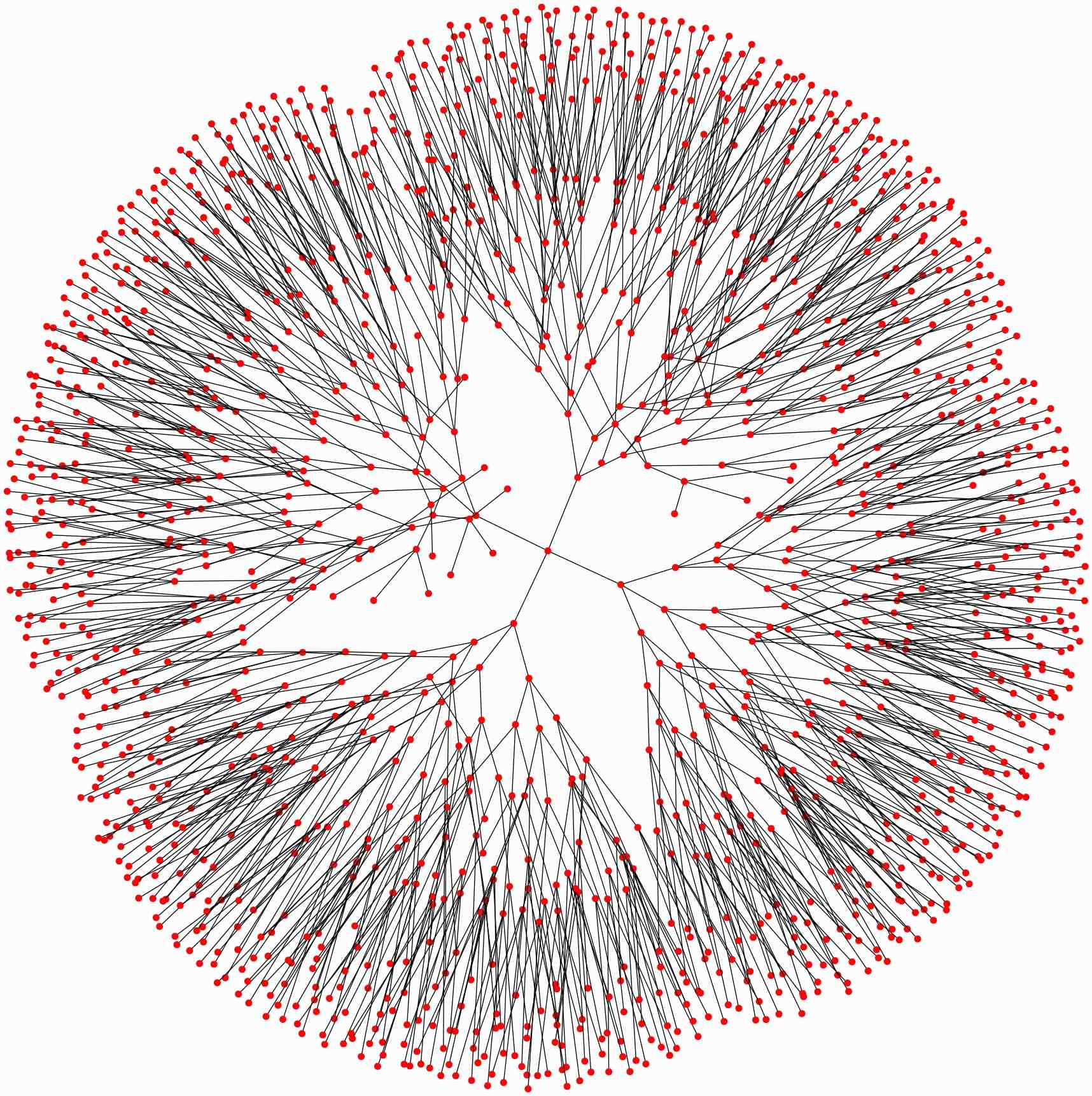}
\includegraphics[width=0.5\textwidth]{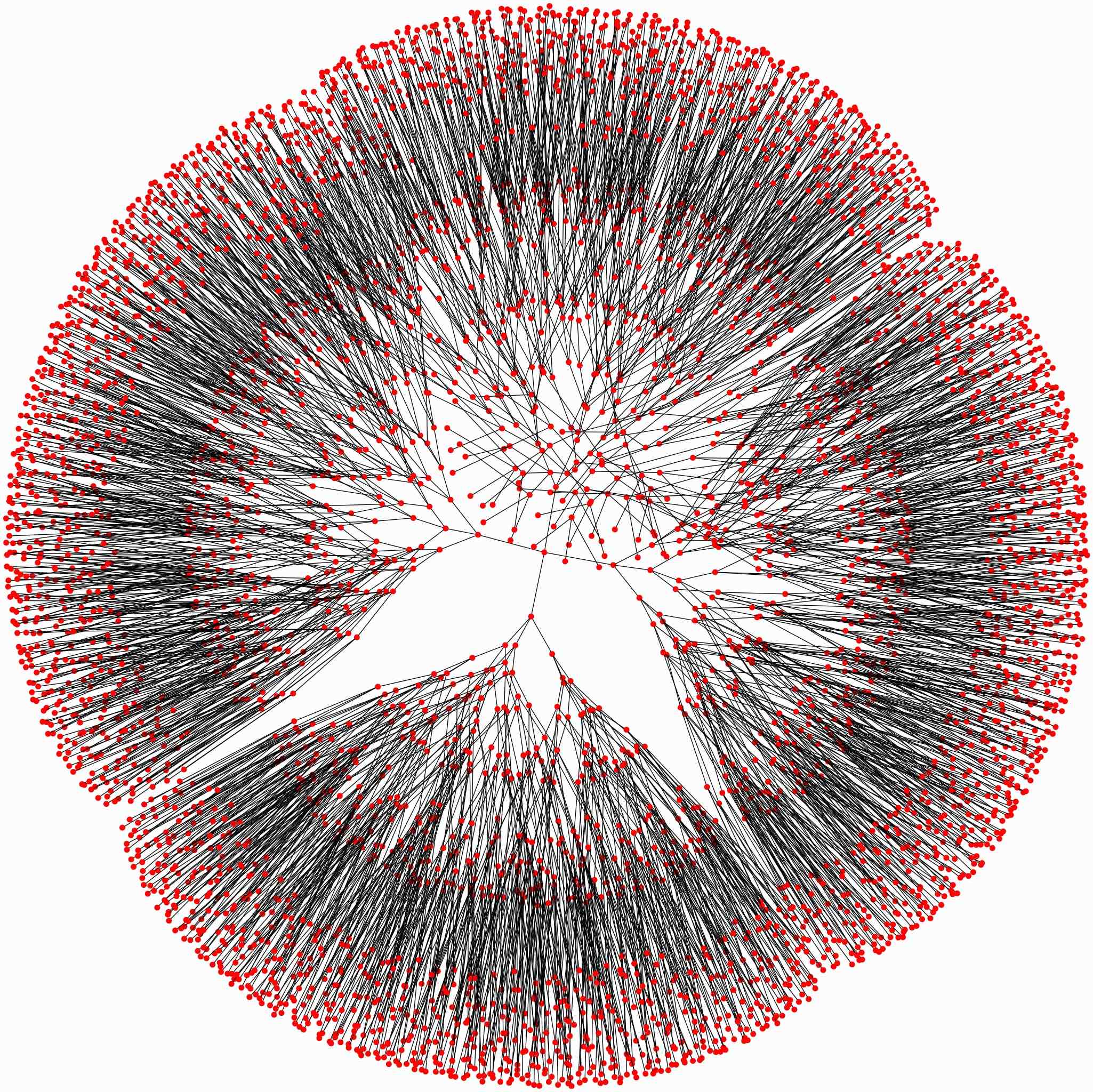}
\begin{center}
\caption{Cayley Graphs of Free Group on two generators. Parts
formed by words of length $\leq n$ ($n=3,4,5,6,7$)}
\label{fig:free-groups}       
\end{center}
\end{figure}

For example we consider it one of the main innovative points that 
continua can be directly constructed from huge syntactic objects. 
As a simplest example, we studied
in \S \ref{sec:binary-words} topological properties of the continua
formed by huge binary words with respect to a few distance
functions. Similarly the Cayley graphs of infinite groups 
define
directly complete metric spaces whose topological properties are
closely related to the algebraic properties of the groups.  The
investigation of these relations has been one of hot topics 
since 1980s as is exposed in 
\cite{MR1699320} and \cite{Epstein:1992:WPG:573874}.  
Fig~\ref{fig:free-groups} shows parts of the Cayley graph 
of the free group with two generators. 

Another example are associahedra. The terms of one binary symbol $b$
and one variable $x$ with $N$ occurrences of $b$ form a connected graph
when a term is connected to another if one is obtained from the other
directly by the associativity rule. The number of nodes is the $N$-th
Catalan number $c_N=\frac{1}{N+1}\binom{2N}{N}$. \par
\begin{figure}
\includegraphics[width=0.3\textwidth]{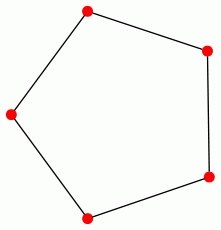}
\includegraphics[width=0.3\textwidth]{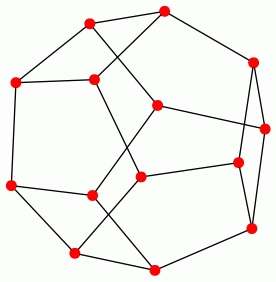}
\includegraphics[width=0.3\textwidth]{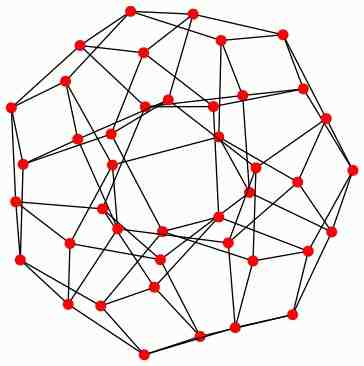}
\includegraphics[width=0.4\textwidth]{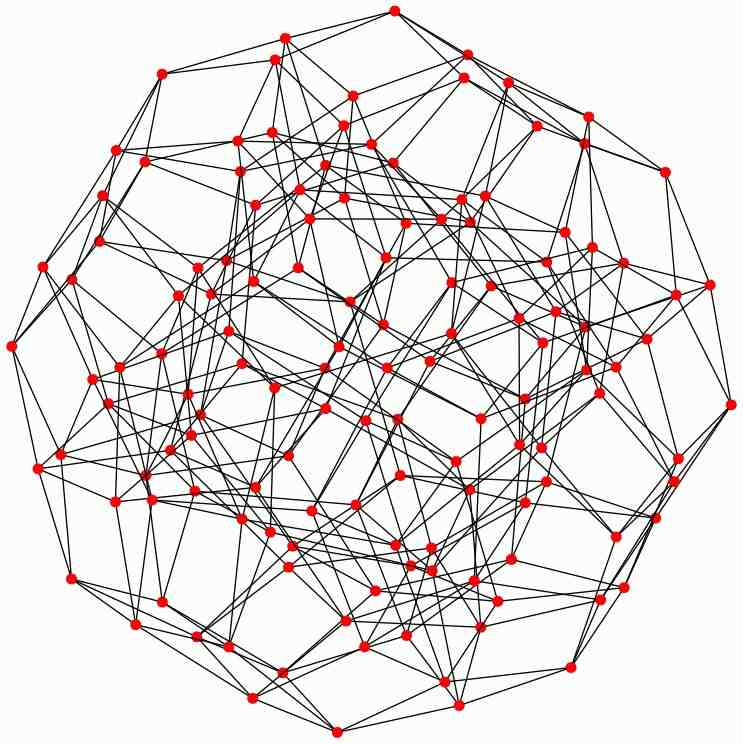}
\includegraphics[width=0.55\textwidth]{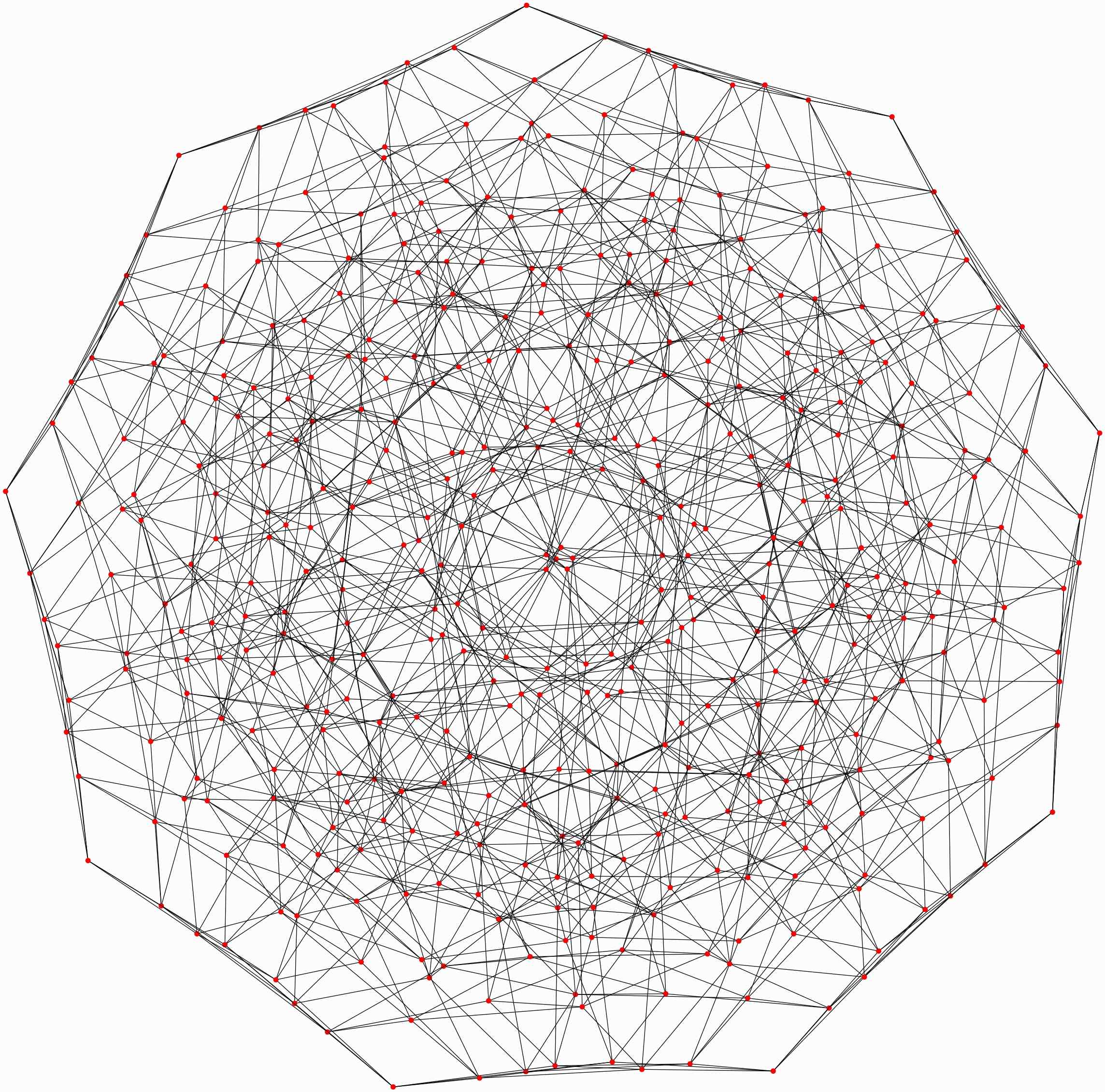}
\begin{center}
\caption{Associahedra (n=3,4,5,6,7)}
\label{fig:associahedron}       
\end{center}
\end{figure}

\begin{center}
\begin{tabular}[tb]{|c||*{10}{c|}}
N&3&4&5&6&7&8&9&10&11&12\\
\hline
$c_N$&5& 14& 42& 132& 429& 1430& 4862& 16796& 58786& 208012
\end{tabular}
\end{center}
\noindent These graphs are called associahedra for accessible $N$. See Fig~\ref{fig:associahedron}. We do not know yet much about the
topological property of this huge associahedron except that the
diameter is $N$ since the shortest path between the farthest pair
$t:=r_x^{N}x$ and $s:=l_x^{N}x$ is of length $N$, where $r_xw=b(w,x)$
and $l_xw=b(x,w)$.  There are many paths connecting $t$ and $s$ but
not so many compared with the hypercube treated in
\S~\ref{subsec:hypercube}.

The continua of huge syntactic objects seem to have potentiality to
become topics of productive investigation.

\subsection{Philosophical Implications}
\subsubsection{Rejection of Existence Absolutism}
\label{subsec:absolutism}
Although modern mathematics has acquired a sort of autonomy in human
intellectual activities and seems to continue flourishing without end,
it has a problem deeply embedded in its core which seems not only to
diminish its cultural value but also to endanger its very existence in
future. The problem is that modern mathematics allows no indefiniteness 
whatsoever and considers even infinity as a definite entity\footnote{
In \cite{vopenka1991philosophical}, P. Vop\v{e}nka emphasizes 
the close
connection between ``natural infinity'' and indefiniteness is
broken by the ``classical infinity'' of modern mathematics. He 
claims that the study of natural infinity via alternative set theory
is also a possible foundation for the science of indefiniteness.
}.
For example, although ``real numbers'' can not be determined by actual
calculations and have essential epistemological indefiniteness, they
are regarded as entities ontologically determined exactly\footnote{
  This determinism is connected with what H. Weyl calls the existence
  absolutism:
\begin{quote}
  "Before God", or "in itself", everything is determined into the last
  detail.  This existence absolutism is governed by a belief similar
  to the one that a process in the external world that we experience
  does not, in itself, carry any vagueness, even though our intuition
  can always only pick out spatial points and qualities in an
  approximate manner, and never delimit them with absolute sharpness.
  (``The Current Epistemological Situation in Mathematics'' \cite[p128]{bo:Mancosu98})
\end{quote}
}. The ``current state'' of the universe
is conceptually conceived as a definite element of
a huge mathematical space and the ideal of scientists is to find a law
which determines the future state from the current state.  In a sense,
the scientific determinism usually refered to by the name ``Laplace's
demon'' is implicitly intertwined to the basic way of thinking of
modern people. This determinism cannot be formulated without the
determinism of modern mathematics\footnote{ It may be objected that
  indeterminacy can be covered by probability theory, which however is
  based on the definiteness of the probability densities.}.

The world is the physical universe with complex phenomena, some of
which has the appearance of living entities and among them the human
beings appear to possess mind and free will, which is nothing but
associated phenomena of the physical activity of the brain. This naive
reductionism joined with the above scientific determinism seems to 
be deeply embedded in the quotidian view of the world. 

For example, there are many informal talks on the possibility of
``artificial reality'' and of immortality of a person in artificial
reality transferring his complete data into computer.  Even admitting
the possibility of describing the universe completely by mathematics,
such arguments lose sense if mathematics considers it non-sense to
think of an entity as determinate if it has potentially infinite details.

This radical determinism of modern mathematics seems not only to
narrow its influence on other disciplines such as biology and
sociology but also underestimate the potentiality of human being. 
Moreover, the determinism makes mathematics more difficult to be
learned than necessary. This difficulty is caused by 
the conflicting position treating indefinite entities as definite ones.  

The mathematical determinism results in the sharp dichotomy between
finite and infinite, which undervalue finiteness as trivial 
compared with ``unfathomable infinity''. This dichotomy mismatches the
naive insight that infinity is an aspect of huge finiteness recognized
under the limitation of our cognitive ability. Similarly the mathematical
determinism results in the sharp dichotomy between discrete and
continuous which also mismatches the naive insight that real continuum
arises from the congregation of finite but huge number of discrete objects
seen under the limitation of the granularity of our cognitive organs.

Most actual things are both finite and infinite depending on our
standpoints. In everyday life, incessant switching of viewpoints 
is indispensable to comprehend the world well. 
Similarly, by incorporating two incommensurable standpoints in mathematics, 
infinity and continuity emerge as the phenomena resulting 
from the interplay between two viewpoints. This is what is achieved by
incorporating sorites paradox in mathematics.

\subsubsection{Internal Measurement}
The indefiniteness thus introduced in mathematics reflects
appropriately the aspect so called
\textit{internal measurement} \cite{gunji04en,matsuno1995quantum} 
of our cognition of the world. Internal measurement refers to the 
stance  taking seriously into account the inevitable temporality and incompleteness of
the interaction between the observer and the observed behind the cognition. 
The ``objects'' are not fixed entities independent
of observation but are phenomena acquiring more clear features
through observation. 

In mathematics, axiomatics correspond to the way of observing
mathematical objects and each theorem may be regarded as an
observation. In contrast to the usual view that the objects are 
definite and immutable, we think that a new theorem deforms the
essence of the object. Such a view is seen in the predicative
reformulation of mathematics proposed by E. Nelson.
\begin{quote}
Let $C$ be an inductive formula. \dots
We can replace our concept of number (any $x$ ) by a more refined concept of number (any
$x$ such that $C^3[x]$). We can read $C^{3}[x]$ as ``$x$ is a number`` (leaving
open the possibility of formalizing an even more refined concept of number
at some time in future) \cite[p14]{nelson}. 
\end{quote}

Here is a new view that mathematical investigation change nature of
the objects by proving theorems. The induction inference is not
considered as an axiom but as the decision of adding new axioms which
refines the concept of numbers by a condition proven to be inductive.
This might be said to take account of internal measurement in
mathematical investigation and the indefiniteness of the totality of
the objects under study is crucial to support such view.

\clearpage
\addcontentsline{toc}{subsection}{Acknowledgement}
\paragraph{Acknowledgement}
I would like to express my deep gratitude to Shuichiro Tsunoda for
showing the crucial and ubiquitous roles played by indefiniteness in
mathematics, to Yukio-Pegio Gunji and Kouichiro Matuno for showing
essential aspects of the viewpoint of internal measurement convincing
me its importance in mathematics, to Taichi Haruna for carefully 
reading the manuscript and for pointing out many errors both mathematical 
and typographical, which improved the manuscript considerably.
I am indebted to Ichiro Tsuda, 
Kunihiko Kaneko and Takashi Ikegami whose stimulant
unprecedented activites on complex systems led me to notice the great
blind spot of modern mathemtaics supposedly arising from the peace and
easiness in the "Cantor's paradice" where infinities are treated as
handy definite entities thereby closing the way and making
indifference to understanding essential aspects of living being where
indefiniteness plays vital roles.  I am also indebted to Yoshitsugu
Oono, Gen Kuroki and Toshio Sunada for their various strong critisims
on the alternative view of mathematical science expressed in
\cite{tjst-1998} which emphasized the importance of indefiniteness in
mathematics based on the internal measurement point of view. 
Their criticism helped me to pursue the alternative way more concretely. 
I am thankful to many mathematicians for stimulative discussions,
especially to Shinsuke Shimogawa, Yohe Yamasaki, Akihiko Gyoja, 
Yoshifumi Takeda, Makoto Kikuchi, Shunsuke Yatabe. 
I am also greatfull for many researchers 
whose concern, sympathy and encouragement has been supporting
my research activity in this isolated direction, especially to 
Yoshihiro Fukumoto, Kazuyuki Tanaka, Yoshinori Shiozawa, Hideo
Mori, Kazufumi Nakajima, Isao Naruki, Etsuro Date, Noriaki Kawanaka,
Toshio Mikami. 
Finally I would like to express my deep gratitude to Shunichi Tanaka 
who turned my attention to "complex systems", to the late Kunihiko Kodaira 
who encouraged me to proceed to new domain of research, and 
to Koji Shiga who constantly
showed me, from my student days, the open-mindedness to mathematics
and invariable passionate concern with the mystery of infinity, which
formed in me the courage to pursue freed from the past without worry
whatever topics I judge important.

\clearpage
\def\oddpagehead{\runningtitle\hfill \thepage}
\def\evenpagehead{\oddpagehead}
\def\runningtitle{References}
\phantomsection
\addcontentsline{toc}{section}{References}
\bibliographystyle{amsalpha}

\begin{thebibliography}{CWF{\etalchar{+}}09}

\bibitem[AG06]{andreev2006theory}
P.V. Andreev and E.I. Gordon, \emph{A theory of hyperfinite sets}, Annals of
  Pure and Applied Logic \textbf{143} (2006), no.~1-3, 3--19.

\bibitem[Bal94]{ballard1994foundational}
D.~Ballard, \emph{Foundational aspects of "non" standard mathematics}, vol.
  176, Amer Mathematical Society, 1994.

\bibitem[Bec79]{beck79:_simpl_sets_and_found_of_analy}
Jon~M. Beck, \emph{Simplicial sets and the foundations of analysis},
  Proceedings of Conference on Sheaf Theory,Durham,England(July 1977), Lecture
  Notes in mathematics, vol. 753, Springer, Berlin, 1979, pp.~113--124.

\bibitem[Bec80]{beck80:_relat_between_algeb_and_analy}
\bysame, \emph{On the relationship between algebra and analysis}, Journal of
  Pure and Applied Algebra \textbf{19} (1980), 43--60.

\bibitem[BMW10]{brown2010step}
A.~Brown, M.A. McDonald, and K.~Weller, \emph{Step by step: Infinite iterative
  processes and actual infinity}, Research in collegiate mathematics education
  \textbf{7} (2010), 115.

\bibitem[Bor52]{borel-1952}
Emil Borel, \emph{Les nombres inaccessibles}, Gauthier-Villars, 1952.

\bibitem[CS95]{chuaqui1995free}
R.~Chuaqui and P.~Suppes, \emph{Free-variable axiomatic foundations of
  infinitesimal analysis: a fragment with finitary consistency proof}, The
  Journal of Symbolic Logic \textbf{60} (1995), no.~1, 122--159.

\bibitem[CT08]{MR2435141}
Peter~J. Cameron and Sam Tarzi, \emph{Limits of cubes}, Topology Appl.
  \textbf{155} (2008), no.~14, 1454--1461. \MR{2435141 (2010c:54040)}

\bibitem[CWF{\etalchar{+}}09]{chollet2009insight}
A.~Chollet, G.~Wallet, L.~Fuchs, G.~Largeteau-Skapin, and E.~Andres,
  \emph{Insight in discrete geometry and computational content of a discrete
  model of the continuum}, Pattern recognition \textbf{42} (2009), no.~10,
  2220--2228.

\bibitem[Die92]{diener1992application}
M.~Diener, \emph{Application du calcul de harthong-reeb aux routines
  graphiques}, Le Labyrinthe du Continu (1992), 424--435.

\bibitem[Dra85]{dragalin1985correctness}
A.~Dragalin, \emph{Correctness of inconsistent theories with notions of
  feasibility}, Computation theory (1985), 58--79.

\bibitem[Dum75]{dummett1975}
M.E. Dummett, \emph{Wang's paradox}, Synthese \textbf{30} (1975), 301--324.

\bibitem[EPC{\etalchar{+}}92]{Epstein:1992:WPG:573874}
David B.~A. Epstein, M.~S. Paterson, G.~W. Camon, D.~F. Holt, S.~V. Levy, and
  W.~P. Thurston, \emph{Word processing in groups}, A. K. Peters, Ltd., Natick,
  MA, USA, 1992.

\bibitem[Gro99]{MR1699320}
Misha Gromov, \emph{Metric structures for {R}iemannian and non-{R}iemannian
  spaces}, Progress in Mathematics, vol. 152, Birkh\"auser Boston Inc., Boston,
  MA, 1999, Based on the 1981 French original [ MR0682063 (85e:53051)], With
  appendices by M. Katz, P. Pansu and S. Semmes, Translated from the French by
  Sean Michael Bates. \MR{1699320 (2000d:53065)}

\bibitem[Gun04]{gunji04en}
Yukio-Pegio Gunji, \emph{原生計算と存在論的観測―生命と時間, protocomputing and
  ontological measurement, in japanese}, University of Tokyo Press, 2004, ISBN
  4130100971.

\bibitem[Har83]{harthong1983^^c3^^a9l^^c3^^a9ments}
J.~Harthong, \emph{{\'E}l{\'e}ments pour une th{\'e}orie du continu},
  Ast{\'e}risque \textbf{109} (1983), no.~110, 235--244.

\bibitem[HLO10a]{hrbacek2010analysis}
K.~Hrbacek, O.~Lessmann, and R.~O'Donovan, \emph{Analysis with ultrasmall
  numbers}, American Mathematical Monthly \textbf{117} (2010), no.~9, 801--816.

\bibitem[HLO10b]{Hrbaceck-ultrasmall-2010}
Karel Hrbacek, Olivier Lessmann, and Richard O'Donovan, \emph{Analysis with
  ultrasmall numbers}, Amer. Math. Monthly \textbf{117} (2010), no.~9,
  801--816. \MR{2760381 (2011j:26047)}

\bibitem[Isl80]{isles}
David Isles, \emph{Remarks on the notion of standard non-isomorphic natural
  number series}, Constructive Mathematics,Proceedings of the New Mexico State
  University Conference Held at Las Cruces, New Mexico, August 11--15,1980
  (F.Richman, ed.), Lecture Notes in Mathematics, vol. 873, Springer, 1980,
  pp.~111--134.

\bibitem[Lau92]{Laugwitz-omega-calculus}
Detlef Laugwitx, \emph{Leibniz' principle and omega calculus}, Le Labyrinthe du
  Continu. Paris: Springer France (1992), 144--154.

\bibitem[Lut87]{Lutz1987reveries}
Robert Lutz, \emph{R\^everies infinit\'esimales}, Gaz. Math. (1987), no.~34,
  79--87. \MR{918184 (89f:03065)}

\bibitem[Lut92]{lutz1992force}
R.~Lutz, \emph{La force des th{\'e}ories infinit{\'e}simales faibles}, le
  labyrinthe du continu, Springer France, Paris (1992), 414--423.

\bibitem[Mag07]{magidor2007strict}
O.~Magidor, \emph{Strict finitism refuted?}, Proceedings of the Aristotelian
  Society (Hardback), vol. 107, Wiley Online Library, 2007, pp.~403--411.

\bibitem[Man98]{bo:Mancosu98}
Paolo Mancosu (ed.), \emph{From {B}rouwer to {H}ilbert. {T}he debate on the
  foundations of mathematics in the 1920s}, Oxford University Press, Oxford and
  New York, 1998.

\bibitem[Mat95]{matsuno1995quantum}
K.~Matsuno, \emph{Quantum and biological computation}, BioSystems \textbf{35}
  (1995), no.~2-3, 209--212.

\bibitem[May00]{mayberry}
J.~P. Mayberry, \emph{The foundations of mathematics in the theory of sets},
  Encyclopedia of Mathematics and its Applications, no.~82, Cambridge
  University Press, Cambridge, 2000, ISBN 0-521-77034-3.

\bibitem[ML90]{martin1990mathematics}
P.~Martin-L{\"o}f, \emph{Mathematics of infinity}, COLOG-88, Springer, 1990,
  pp.~146--197.

\bibitem[Mon01]{monaghan2001young}
J.~Monaghan, \emph{Young peoples' ideas of infinity}, Educational Studies in
  Mathematics \textbf{48} (2001), no.~2, 239--257.

\bibitem[Myc81]{mycielski81:_analy_without_actua_infin}
Jan Mycielski, \emph{Analysis without actual infinity}, J. Symbolic Logic
  \textbf{46} (1981), 625--633.

\bibitem[Nel77]{nelson77:_inter_set_theor}
E.~Nelson, \emph{Internal set theory: a new approach to nonstandard analysis},
  Bulletin of the American Mathematical Society \textbf{83} (1977), 1165--1198.

\bibitem[Nel87a]{nelson}
\bysame, \emph{Predicative arithmetic}, Princeton University Press, 1987, ISBN
  0-691-08455-6.

\bibitem[Nel87b]{nelson-nonstandard}
\bysame, \emph{Radically elementary probability theory}, Annals of Mathematics
  Studies, no. 117, Princeton University Press, 1987, ISBN 0-691-08455-6.

\bibitem[Nel04]{nelson04:_bookr}
Edward Nelson, \emph{Bookreview:gnomes in the fog: The reception of brouwer's
  intuitionism in the 1920s, by dennis e. hesseling, science networks
  --historical studies, vol.28 birkhauer,base,2003,xxviii + 447 pp.,isbn
  3-7643-6536-6}, Bulletin of the AMS (2004).

\bibitem[Nel07]{nelson-2007-simplicity}
\bysame, \emph{The virtue of simplicity}, The strength of nonstandard analysis,
  SpringerWienNewYork, Vienna, 2007, pp.~27--32. \MR{2341412}

\bibitem[Pal95]{MR1336645}
Erik Palmgren, \emph{A constructive approach to nonstandard analysis}, Ann.
  Pure Appl. Logic \textbf{73} (1995), no.~3, 297--325. \MR{1336645
  (96c:03123)}

\bibitem[Par71]{r.71:_exist_and_feasib_in_arith}
R.~Parikh, \emph{Existence and feasibility in arithmetic}, Journal of Symbolic
  Logic \textbf{36} (1971), 494--508.

\bibitem[P{\'e}r92]{peraire1992th^^c3^^a9orie}
Y.~P{\'e}raire, \emph{Th{\'e}orie relative des ensembles internes}, Osaka J.
  Math \textbf{29} (1992), no.~2, 267--297.

\bibitem[P{\'e}r05]{peraire2005replacement}
\bysame, \emph{Le replacement du r{\'e}f{\'e}rent dans les pratiques de
  l’analyse issues de e. nelson et de g. reeb}, Philosophia Scienti{\ae}.
  Travaux d'histoire et de philosophie des sciences (2005), no.~CS 5, 257--273.

\bibitem[Ras73]{rashevskii1973dogma}
P.K. Rashevskii, \emph{On the dogma of the natural numbers}, Russian
  Mathematical Surveys \textbf{28} (1973), no.~4, 143--148.

\bibitem[Ree81]{reeb1981math}
G.~Reeb, \emph{Mathematique non standard (essai de vulgarisation)}, Bulletin
  APMEP \textbf{328} (1981), 259--273.

\bibitem[Rob66]{robinson1966}
A.~Robinson, \emph{Non-standard analysis}, North-Holland Publishing
  Company,Amsterdam, 1966.

\bibitem[RR96]{Reveilles-Richard1996}
Jean-Pierre Reveill^^c3^^a8s and Denis Richard, \emph{Back and forth between
  continuous and discrete for the working computer scientist}, Annals of
  Mathematics and Artificial Intelligence \textbf{16} (1996), no.~1, 89--152.

\bibitem[RS10]{DBLP:journals/corr/abs-1005-5685}
Ana Romero and Francis Sergeraert, \emph{Discrete vector fields and fundamental
  algebraic topology}, CoRR \textbf{abs/1005.5685} (2010).

\bibitem[San10]{sanders2010relative}
S.~Sanders, \emph{Relative arithmetic}, Mathematical Logic Quarterly
  \textbf{56} (2010), no.~6, 564--572.

\bibitem[Saz95]{sazonov-feasible}
Vladimir~Yu. Sazonov, \emph{On feasible numbers}, Logic and computational
  complexity (Leviant D, ed.), Lecture Notes in computer sicence, vol. 960,
  Springer, 1995, pp.~30--51.

\bibitem[SLSZ]{sochor-differential-calculus-AST}
Anton{\'\i}n Sochor, Alistair Lachlan, Marian Srebrny, and Andrzej Zarach,
  \emph{Differential calculus in the alternative set theory}, pp.~273--284,
  Springer Berlin / Heidelberg.

\bibitem[Tal80]{tall1980notion}
D.~Tall, \emph{The notion of infinite measuring number and its relevance in the
  intuition of infinity}, Educational Studies in Mathematics \textbf{11}
  (1980), no.~3, 271--284.

\bibitem[Tho92]{MR1413523}
Ren{\'e} Thom, \emph{L'ant\'eriorit\'e ontologique du continu sur le discret},
  Le labyrinthe du continu ({C}erisy-la-{S}alle, 1990), Springer, Paris, 1992,
  pp.~137--143. \MR{1413523}

\bibitem[Tra98]{tragesser98:_part_i}
Robert Tragesser, \emph{Part i:ultrafinitism,naturalism,vagueness},
  \url{http://www.cs.nyu.edu/pipermail/fom/1998-April/001825.html}, 4 1998.

\bibitem[Tsu98]{tjst-1998}
Toru Tsujishita, \emph{生命と複雑系, life and complex systems, in japanese},
  Science of Complex Systems and Modern Thought, pp.~75--225, Seidosha, 1998,
  ISBN 4-7917-9145-2.

\bibitem[TT01]{tall2001infinity}
D.~Tall and D.~Tirosh, \emph{Infinity--the never-ending struggle}, Educational
  studies in Mathematics \textbf{48} (2001), no.~2, 129--136.

\bibitem[vdDW84]{MR751150}
L.~van~den Dries and A.~J. Wilkie, \emph{Gromov's theorem on groups of
  polynomial growth and elementary logic}, J. Algebra \textbf{89} (1984),
  no.~2, 349--374. \MR{751150 (85k:20101)}

\bibitem[Ver98]{Vershik-Urysohn-1998-MR1691182}
A.~M. Vershik, \emph{The universal {U}ryson space, {G}romov's metric triples,
  and random metrics on the series of natural numbers}, Uspekhi Mat. Nauk
  \textbf{53} (1998), no.~5(323), 57--64. \MR{1691182 (2000b:53055)}

\bibitem[Vol70]{volpin}
A.S.~Essenin Volpin, \emph{The ultra-intuitionistic criticism and the
  anti-traditional programme for foundations of mathematcs}, Intuitionism and
  Proof Theory (J.~Myhill A.~Kino and R.E. Vesley, eds.),
  North-Holland,Amserdam, 1970, pp.~3--45.

\bibitem[Vop79]{vopenka}
Petr Vop$\check{\mbox{e}}$nka, \emph{Mathematics in the alternative set
  theory}, Teubner, 1979.

\bibitem[Vop91]{vopenka1991philosophical}
P.~Vop{\v{e}}nka, \emph{The philosophical foundations of alternative set
  theory}, International Journal Of General System \textbf{20} (1991), no.~1,
  115--126.

\bibitem[vS90]{MR1075018}
Walter~P. van Stigt, \emph{Brouwer's intuitionism}, Studies in the History and
  Philosophy of Mathematics, vol.~2, North-Holland Publishing Co., Amsterdam,
  1990. \MR{1075018 (92d:01054)}

\bibitem[Wey49]{weyl-philosophy}
Hermann Weyl, \emph{Philosophy of {M}athematics and {N}atural {S}cience.
  {R}evised and {A}ugmented {E}nglish {E}dition {B}ased on a {T}ranslation by
  {O}laf {H}elmer}, Princeton University Press, Princeton, N. J., 1949.
  \MR{0029851 (10,670c)}

\bibitem[Wey94]{Weyl-continuum}
\bysame, \emph{The continuum}, Dover Publications Inc., New York, 1994, A
  critical examination of the foundation of analysis, Translated from the
  German by Stephen Pollard and Thomas Bole, With a foreword by John Archibald
  Wheeler and an introduction by Pollard, Corrected reprint of the 1987
  translation [Thomas Jefferson Univ. Press, Kirksville, MO; MR1040831
  (91h:01105)]. \MR{1280464}

\bibitem[Yat09]{MR2500986}
Shunsuke Yatabe, \emph{Comprehension contradicts to the induction within \l
  ukasiewicz predicate logic}, Arch. Math. Logic \textbf{48} (2009), no.~3-4,
  265--268. \MR{2500986 (2010f:03025)}

\end{thebibliography}
\newcommand{\etalchar}[1]{$^{#1}$}
\providecommand{\bysame}{\leavevmode\hbox to3em{\hrulefill}\thinspace}
\providecommand{\MR}{\relax\ifhmode\unskip\space\fi MR }
\providecommand{\MRhref}[2]{%
  \href{http://www.ams.org/mathscinet-getitem?mr=#1}{#2}
}
\providecommand{\href}[2]{#2}

\clearpage
\pagestyle{empty}             %
\def\oddpagehead{}
\def\evenpagehead{\oddpagehead}
\addcontentsline{toc}{section}{Index}
\def\runningtitle{Index}
\printindex
\vfill
\noindent  Toru Tsujishita \\
Department of Mathematics, Ritsumeikan University, Shiga 525-8577, JAPAN \\
email: \texttt{tjst@se.ritsumei.ac.jp}

\end{document}